\documentclass[12pt,a4paper]{book}
\usepackage{latexsym} 
\usepackage{rawfonts} 
\usepackage{indentfirst} 
\title{SYSTEM ZEROS}
\author{Ye. M. Smagina}
\begin{document}
\maketitle
\newpage
\pagenumbering{roman}
\tableofcontents
\newpage
\pagenumbering{arabic}
\chapter*{Introduction}
\addcontentsline{toc}{chapter}{Introduction}
\markboth{Introduction}{Introduction}

     By the late 1950-s control methods based on the
state-space approach (i.e. optimal control, filtering and so on ) have been begun to develop and gave  excellent
results in control of complicated aerospace and industrial objects, which are described in the state-space by
multi-input and multi-output  systems.  In  view of the success  of the state-space approach this period
characterized by decreasing  the interest to the classic control design methods. Meanwhile  optimal control
revealed some disadvantages, which were only inherent to the state-space method but absent in the
frequency-response approach, for example, problems with  response analysis, difficulties with robustness and so
on.

     It is known that control problems in single-input/single-output
systems  are successfully solved by  classic frequency-response methods, which are based on notions of poles,
zeros and etc. Significant interest to the classical methods was appeared once again in the mid-1960. Many
researches attempted to extend the fundamental  concepts  of  the classic theory, such as a transfer function,
poles,  zeros, a frequency response  and etc. to linear multi-input/multi-output multivariable systems described
in  the state-space. For example, the well known method of modal control may be considered as an extension of
the classic method shifting poles.

     The main  difficulties  encountered in reaching this goal were the
generalization of the concept of a zero of a transfer function. Indeed, a classic transfer function of a
single-input/single-output system represents a rational function of the complex variable, which is a ratio of
two relatively prime polynomials. A zero of the classic transfer function is equal to a zero of a polynomial in
a numerator of the transfer function and coincides with a complex variable for which the numerator (and the
transfer function) vanishes.

A transfer function of a  multi-input/multi-output multivariable system  represents  a  matrix  with  elements
being rational functions  i.e. every element is an ratio of  two  relatively prime polynomials. In this case it
was very difficult to extend the classical zero definition to multivariable case. Only in 1970 H.H. Rosenbrock
introduced the notion of a zero of a multivariable  system, which  was  equivalent  to the classic one in the
physical meaning [R1]. Then this notion has been improved as Rosenbrock [R2], [R3] as others researchers [M3],
[W2], [M1], [M2], [D4], [A2], [P3], [K2]. As a result the main classic notions: minimal and nonminimal phase,
invertibility, the root-locus method, the integral feedback and etc. were extended to multivariable control
systems.

     The  first  review  devoted to the algebraic,  geometric  and  complex
variable properties of  poles and zeros of linear multivariable systems was published by MacFarlane and
Karcanias  in  1976  [M1]. The fullest survey devoted to definitions, classification and applications of system
zeros  was appeared in 1985 [S7]. The detailed review about system zeros was also published in [S2].

The present book is the first publication in English considered the modern problems of control theory and
analysis  connected with a concept of  system zeros. The previous book by Smagina [S9] had been written in
Russian  and it is inaccessible to  English speaking researchers.

The purpose of the offered  book is to systematize and consistently to state basic theoretical results connected
with properties  of  multivariable system zeros. Different zeros definitions and different types of  zeros are
studied. Basic numerical algorithms of  zeros computing and the zero assignment problem are also presented.
 All results  are illustrated by examples.

The book contains ten chapters. The first and second chapters are devoted to  different descriptions of a linear
multivariable dynamical system. They are   linear   differential equations   (state-space description) and
transfer function matrices. Few canonical forms having a companion matrix of dynamics are presented in the first
chapter. The second chapter is devoted to several basic  properties of transfer function matrices that related
with  controllability and observability notions. Also the Smith-McMillan canonical form of a transfer function
matrix and the Smith canonical form of its a numerator are studied.

Notions of transmission and invariant zeros are introduced in the third chapter. The physical interpretation of
these notions are explained. It is shown that  transmission and invariant  zeros are related to  complete
blocking some inputs that  proportional  to  $exp(zt)$ where  $z$  is a invariant (transmission) zero.

     In the fourth chapter the complete set  of   transmission
zeros  is  defined  via a transfer  function  matrix. Several methods of  transmission zeros calculation are
studied. These methods are based on the Smith-McMillan canonical form, transfer function matrix minors and
invariant polynomials of a numerator of the transfer function matrix. Also a new original method for
factorization of the transfer function matrix  is suggested.

Invariant  and system zeros are calculated via the system matrix in the fifth chapter. Notions of  decoupling
zeros are introduced. Also in this chapter we analyze relationships between zeros of different types.

     In the sixth chapter we study properties of zeros, i.e. it has been shown that  zeros are invariant under several
nonsingular transformations and the state and/or output feedback.

In the seventh chapter zeros of a controllable system are calculated via a special polynomial matrix (matrix
polynomial) formed by using the special canonical  representation  of a linear multivariable system. Proposed
method discovers relationships between zeros and the input-output  structure of a system. Several general
estimations of a number of zeros are obtained. Also it is presented a method of zero calculating via a matrix
pencil of the reduced order.

The computer-aided methods of  zeros computing and several  methods of zeros  assignment  are described  in the
eighth and ninth chapters.

     The applications of transmission zeros in the servomechanism problems and for maximally
achievable accuracy of an optimal system are included  in  the  tenth chapter.

\chapter[System description by differential equations]
        {System description by differential equations}

To control design we usually study a mathematical model obtained as a result of experiment or studying physical
laws. Depending on a way of obtaining the mathematical model can be represented as a set of differential
equations  and also through transfer functions. At first let us consider the description through differential
equations.

\section[State space representation]{State space representation}

Such representation is based on deduction of differential equations that describe dynamical behavior  of a
object  by studying physical laws. The  equations reveal internal correlation between all physical variables
that govern a work of the object. The set of these physical variables at any time $t$ is termed as a state of
the dynamical system and denoted by a vector $x(t)$. Individual physical variables and/or their linear
combinations are termed as state variables of the state vector $x(t)$ and denoted by $x(i),i=1,...,n$ where $n$
is a number of state variables, a dimension of the state-space.

   Let $u(t)$ is an $r$ dimensional vector-valued function of time
that  is called as an input of a dynamical system,  $y(t)$ is an $l$ dimension vector-valued function of time
that is called as an output of a dynamical system $(r,l\le n)$. The following  set of first order linear
vector-matrix differential equations presented in a vector-matrix form is named as a linear model of a dynamical
system in the state-space
$$  \dot{x}(t) = Ax(t) + Bu(t) \eqno(1.1) $$
$$  y(t) =Cx(t)              \eqno(1.2)$$
where $A,B,C$ are $n\times n$, $n\times r$  and $l\times n$ matrices respectively. If  elements of $A,B,C$ are
 functions of time then Eqns (1.1),(1.2) describe a time-depend linear dynamical model, otherwise if $A,B,C$
are constant matrices then (1.1),(1.2) is named as a time-invariant model.

  In some cases it is desirable to augment  equation  (1.2)  to
allow the output $y(t)$ to depend also  on  the  input  vector $u(t)$.  So, a general form of the linear
dynamical model is
$$ \dot{x}(t) = Ax(t) +Bu(t) $$
$$y(t) = C(t) +Du(t) \eqno(1.3)$$
where $D$ is an $r\times l$  matrix.

  In the following text we shall denote: $x=x(t), u=u(t),
y=y(t)$ and imply that vectors $x,u$, and $y$ are functions of time.

    The general  solution $x(t)$ of  the  linear  time-invariant
nonhomogeneous (forced) vector-matrix  differential  equation (1.1)
with initial state $x(t_o)=x_o$ is defined as [A1,W1]
$$  x(t) = e^{A(t-t_o)}x_o +
  \int_{t_o}^{t}{e^{A(t-\tau )}Bu(\tau)}d\tau \eqno(1.4)$$
where $e^{At}$   is  the  conventional  notation  of  the $n\times n$  matrix being termed as a matrix
exponential and defined by the formula
$$ e^{At} =I_n + At + \frac{A^2}{2!}t^2 + \cdots
\eqno(1.5)$$
Here $I_r$ is an $r\times r$  unity matrix.

    Let us recall  that the matrix $e^{At}$ is the state transition matrix
[W1] of the linear time-invariant homogeneous vector-matrix differential equation $\dot{x} = Ax$ with $t_o = 0$.
The matrix $e^{At}$ has the following properties
$$
  a) e^{At} e^{-At} =I_n , \qquad d) e^{A(t+t_o)} = e^{At} e^{At_o},
$$
$$
  b)e^{A(t)^{-1}} = e^{-At},\qquad  e) e^{A(t-t_o)} = e^{At} e^{-At_o},
$$
$$
  c)e^{I_nt} =I_n e^{t}, \qquad    f) \frac{d}{dt}e^{At}=Ae^{At}=
     e^{At}A  \eqno(1.6)$$

     The substitution of (1.4) into (1.2) gives the
output $y=Cx$ in the form

$$  y(t) =Ce^{A(t-t_o)}x_o +
 C\int_{t_o}^{t}{e^{A(t-\tau )}Bu(\tau)}d\tau
\eqno(1.7)$$

     Let the matrix $A$ has $n$ distinct eigenvalues $\lambda_1,\cdots,
\lambda_n$ with  corresponding linearly independent right eigenvectors $w_1,\cdots,w_n$ and dual left
eigenvectors $v_1,\cdots,v_n$. These vectors satisfy the relations [G1]
$$ Aw_i = \lambda_i w_i, \qquad v_i^TA = \lambda_iv_i^T,$$
 $$ v_i^Tw_j =\delta_{i,j} $$
where $\delta_{i,j} =1$  if $i=j$,  otherwise  zero. In this case the
matrix exponential can be presented as
$$
e^{At} =\sum_{i=1}^n{e^{\lambda_it}w_iv_i^T} \eqno(1.8)
$$
Substituting (1.8) into (1.7) enables to express
$y(t)$ as
$$
 y(t)= \sum_{i=1}^n\gamma_ie^{\lambda_i(t-t_o)}v_i^Tx_o +
\sum_{i=1}^n\gamma_i\int_{t_o}^{t}e^{\lambda_i(t-\tau )} \beta_i^Tu(\tau)d\tau \eqno(1.9) $$
 where column
vectors $\gamma_i$, $ i=1,2,...,n$ and row vectors $\beta_i^T$, $i=1,2,...,n$ are defined as follows
$$
        \gamma_i =Cw_i, \qquad  \beta_i^T = v_i^TB
\eqno(1.10) $$

     The  notions of controllability  and  observability  are
fundamental  ones  of   linear  dynamical   system (1.1),(1.2) [K1].

     {\it DEFINITION 1.1.} \rm [W1]: System (1.1),(1.2) is said to be
completely state controllable  or controllable if and only  if control $u(t)$ transferring any initial state
$x(t_o)$ at any time $t_o$ to any arbitrary final  state $x(t_1)$ at any  finite  time $t_1$ exists. Otherwise,
the system is said to  be uncontrollable.

     {\it DEFINITION 1.2.} \rm [W1]: System (1.1),(1.2)  is  said  to  be
completely state observable  or observable  if  and  only  if  the state $x(t)$ can be reconstructed over any
finite  time  interval $[t_o,t_1]$ from complete  knowledge  of  the  system  input $u(t)$ and output $y(t)$
over the time interval $[t_o,t_1]$ with $t_1>t_o \geq 0$.

     Let us introduce   algebraic  conditions  of  complete
controllability and observability,  which will be used late on.

     {\it THEOREM 1.1.}  \rm System (1.1), (1.2)  or,  equivalently,
the pair of matrices $(A,B)$ is controllable if and only if

$$  rank Y = rank \left[ \begin{array}{cccc} B,& AB,& \cdots, & A^{n-1}B
\end{array} \right] = n
\eqno(1.11)  $$ where an $n\times nr$  matrix $Y = [B,AB,\cdots,A^{n-1}B]$  is called by the controllability
matrix of the pair $(A,B)$.

     {\it THEOREM 1.2.}  \rm System (1.1), (1.2)  or,  equivalently,  the
pair of matrices $(A,C)$ is observable if and only if

$$
   rank Z = rank \left[ \begin{array}{c} C \\ AC \\ \vdots \\ A^{n-1}C
\end{array}  \right] = n
\eqno(1.12)  $$ where an $n\times nl$ matrix $Z^T = [C^T,A^TC^T,\cdots,(A^T)^{n-1}C^T]$  is called by the
observability matrix of the pair $(A,C)$.

     Proofs of these theorems may be found in [A1], [V1], [O1].

     Let us consider also the following simple algebraic conditions of
controllability and observability.

  {\it THEOREM 1.3. \rm   System (1.1),(1.2) is  controllable  if  and
only if
$$
     rank \left[ \begin{array}{cc} \lambda_iI_n-A,& B \end{array} \right]
 = n \eqno(1.13)  $$
where $\lambda_i$ is an eigenvalue of $A$, $i=1,...,n$.

     {\it THEOREM 1.4.}  \rm  System (1.1), (1.2)  is  observable  if  and
only if
$$     rank \left[ \begin{array}{c} \lambda_iI_n-A\\ C \end{array} \right]
 \eqno(1.14)  $$
where $\lambda_i$ is an eigenvalue of $A$, $i=1,...,n$.

     The proof is given in [R2].

     Dynamical behavior of a linear time-invariant
system may be described also via input-output variables by a set
of  differential equations of an order $p$
$$
  F_py^{(p)} + F_{p-1}y^{(p-1)} +\cdots + F_oy =
  B_ku^{(k)} + B_{k-1}u^{(k-1)} + \cdots + B_ou
\eqno(1.15)   $$ where $y^T=[y_1,\ldots ,y_r]$ is an $r$-dimensional vector of the output, $u^T=[u_1,\ldots
,u_r]$ is an $r$-dimensional vector of the  input $(r \ge 1)$, $F_i$ , $B_i$  are constant $r\times r$ matrices,
$rankF_p=r$.

     Let's note that
we can transfer from the input-output representation (1.15) to the state-space representation (1.1),(1.2) by an
linear combination of input and output variables [A1], [S23], [M5].

\section[Block companion canonical forms of time-invariant
             system]{Block companion canonical forms of \\
            time-invariant  system}
\subsection[Companion and block companion matrix]{Companion
            and block companion matrix}

At first we consider a sense of the term 'companion matrix'. Let us introduce a monic polynomial in s of an
order  $n$  with  real coefficients $a_1,a_2,\ldots ,a_n$
$$ \phi (s) = s^n + a_1s^{n-1} + \cdots + a_{n-1}s + a_n
\eqno(1.16) $$ and an  $n\times n$  matrix $P$  of the following
structure
$$ P\; = \;
\left[ \begin{array}{ccccc} 0 & 1 & 0 & \cdots & 0 \\ 0 & 0 & 1 & \cdots & 0 \\
\vdots & \vdots & \vdots & \ddots & \vdots \\  0 & 0 & 0 & \cdots & 1 \\
-a_n & -a_{n-1} & -a_{n-2} & \cdots & -a_1 \end{array} \right]
\eqno(1.17) $$
The matrix $P$ is known as the \it{companion
matrix} \rm of the polynomial $\phi(s)$ [L1]. Indeed, the
following equality is true
$$
 \phi(s) = det(sI_n -P) \eqno (1.18)
$$
Let us introduce a regular\footnote{A polynomial matrix $ L(s) = L_os^p + L_1s^{p-1} +\cdots +L_p$ is termed a
regular one if $ L_p $  is a nonsingular matrix.} $r\times r$ polynomial  matrix
$$
 \Phi (s) = I_rs^p + T_1s^{p-1} + \cdots + T_{p-1}s + T_p \eqno(1.19)
$$
whose elements are polynomials in $s$, matrices $T_1,T_2,\cdots ,T_p$ are $r\times r$ constant matrices  with
real  elements. Matrix $\Phi(s)$ is called [G2] as the monic matrix polynomial of degree $p$. We define an
$rp\times rp$ block matrix $P^*$  as follows:
$$
 P^* =
\left[ \begin{array}{ccccc} O & I_r & O & \cdots & O\\ O& O& I_r& \cdots & O\\
\vdots & \vdots & \vdots & \ddots & \vdots \\  O& O& O& \cdots & I_r\\
-T_p & -T_{p-1} & -T_{p-2} & \cdots & -T_1 \end{array} \right]
\eqno(1.20)
$$
The matrix $P^*$ is called the \it{block companion matrix} \rm [B1] of the matrix polynomial $\Phi(s)$. The
following assertion reveals a relation between $P^*$ and $\Phi(s)$.

{\it ASSERTION 1.1.}
$$
     det(sI_{rp}-P^*) = det\Phi(s)
\eqno(1.21)  $$

    {\it PROOF. }\rm Indeed,
$$ det(sI_{rp}-P^*) = det \left[ \begin{array}{ccccc}
 sI_r & -I_r & O & \cdots & O\\ O& sI_r & -I_r& \cdots & O\\
\vdots & \vdots & \vdots & \ddots & \vdots \\  O& O& O& \cdots & -I_r\\
 T_p &  T_{p-1} & T_{p-2} & \cdots & sI_r+T_1 \end{array} \right]
$$
The matrix $sI-P^*$ is partitioned into four blocks
$$
sI-P_{11} =
\left[ \begin{array}{ccccc} sI_r & -I_r & O & \cdots & O\\ O & sI_r& -I_r&
\cdots & O\\ \vdots & \vdots & \vdots & \ddots & \vdots \\  O& O& O& \cdots &
sI_r  \end{array} \right] , -P_{12}= \left[ \begin{array}{c} O\\O\\ \vdots \\
-I_r \end{array} \right] ,
$$
$$
-P_{21} = \left[ \begin{array}{ccccc} T_p & T_{p-1} & T_{p-2} & \cdots &
 T_2 \end{array} \right] , sI-P_{22} =  sI_r + T_1
$$
Assuming $s\ne 0$ and using formulas of Schur  [G1] we can reduce a determinant
of the block matrix \\ $ P^* = \left[ \begin{array}{cc} sI-P_{11} & -P_{12}\\
-P_{21} & sI-P_{22} \end{array} \right] $   to the form
$$
  det(sI_{rp}-P^*) = det(sI-P_{11})det(sI-P_{22} - P_{21}(sI-P_{11})^{-1}P_{12})
\eqno(1.22)
$$
It is easy to verify that
$$
(sI-P_{11})^{-1} = \left[ \begin{array}{cccc} s^{-1}I_r & s^{-2}I_r & \cdots &
s^{(p-1)}I_r \\ O & s^{-1}I_r & \cdots & s^{(p-2)}I_r \\ \vdots &  \vdots &
\ddots & \vdots \\ O & O & \cdots & s^{-1}I_r \end{array} \right]
$$
$$
  (sI-P_{11})^{-1}P_{12} = \left[ \begin{array}{c} s^{(1-p)}I_r \\ s^{(2-p)}I_r
\\ \vdots \\s^{-1}I_r \end{array} \right]
$$
Substituting these relationships and  blocks $P_{21}$, $sI-P_{22}$ into  (1.22) gives
$$
det(sI_{rp}-P^*)  = s^{r(p-1)}det(sI_r + T_1 + [T_p,T_{p-1},\ldots ,T_2]
\left[ \begin{array}{c} s^{(1-p)}I_r \\ s^{(2-p)}I_r \\ \vdots \\s^{-1}I_r
\end{array} \right]) = $$
$$
s^{r(p-1)}det(sI_r + T_1 + T_2s^{-1} +T_3s^{-2},\ldots ,T_ps^{1-p}]) =
det(I_rs^p + T_1s^{p-1} + \cdots + T_p) = det\Phi(s)$$
The assertion is proved.

     {\it REMARK 1.1. \rm
    It is  evident  that  the relationship  (1.21)  is  true  for $s=0$.
Actually, let us find $detP^* = det(sI - P^*)/_{s=o} =det(T_p)$.
The right-hand side of the  last  expression  coincides  with  the
right-hand side of (1.21) for $s=0$.

     Now we consider the  modification  of  the  block  companion  matrix
(1.20). Let   matrices $T_i$  ($i=1,2,\ldots,p$)  in (1.19) have the following structure
$$
  T_i =[O,\hat{T}_i], \qquad i=1,2,\ldots ,p
\eqno(1.23) $$
where $\hat{T}_i$ are  $r\times l_{p-i+1}$ non-zero matrices, integers $l_1, l_2,\ldots ,l_{p-1}$
satisfy the following inequality
$$
  l_1 \le l_2 \le \cdots \le l_{p-1} \le r,  \qquad l_p=r
\eqno(1.24) $$ We denote $\bar{n} = l_1 +l_2 + \cdots +l_p $ and
define an $\bar{n}\times \bar{n}$ matrix
$$  \hat{P} =
\left[ \begin{array}{ccccc} O & E_{1,2} & O & \cdots & O\\ O& O& E_{2,3}&
\cdots & O\\\vdots & \vdots & \vdots & \ddots & \vdots \\  O& O& O& \cdots &
E_{p-1,p}\\ -\hat{T}_p & -\hat{T}_{p-1} & -\hat{T}_{p-2} & \cdots & -\hat{T}_1
\end{array} \right]
\eqno(1.25)
$$
where  $l_i\times l_{i+1}$ blocks $E_{i,i+1}$ have forms
$$
           E_{i,i+1} = [O, I_{l_i}] ,\qquad i=1,2,\ldots p-1
$$
     We will call the matrix $\hat{P}$ by the \it{generalized block companion
matrix} \rm of the matrix polynomial $\Phi(s)$ having the matrix coefficients $T_i =[O,\hat{T}_i]$.

  {\it ASSERTION 1.2.} \rm  For $s\ne 0$ the following equality is true
 $$
     det(sI_{\bar{n}} -\hat{P}) =s^{\bar{n}-rp}det\Phi(s)
\eqno(1.26)$$

    {\it PROOF.}  \rm The matrix $sI-\hat{P}$ is partitioned into four blocks
$$
sI-P_{11} =
\left[ \begin{array}{ccccc} sI_{l_1} & -E_{1,2} & O & \cdots & O\\ O &
sI_{l_2}& -E_{2,3} & \cdots & O\\ \vdots & \vdots & \vdots & \ddots & \vdots \\
O& O& O& \cdots & sI_{l_p-1}  \end{array} \right] , -P_{12}=
\left[ \begin{array}{c} O\\O\\ \vdots \\ -E_{p-1,p} \end{array} \right] ,
\eqno(1.27)$$
$$
-P_{21} = \left[ \begin{array}{cccc} \hat{T}_p & \hat{T}_{p-1} & \cdots &
  \hat{T}_2 \end{array} \right] , \qquad sI-P_{22} =  sI_r + \hat{T}_1
\eqno(1.28)$$
We assume $s\ne 0$ and will use formula (1.22). At
first we calculate determinants $det(sI-P_{11})$ and
$(sI-P_{11})^{-1}P_{12}$. Using (1.24) we find
$$
    det(sI-P_{11}) = s^{l_1+l_2+\cdots +l_{p-1}} = s^{\bar{n}-r}
\eqno(1.29)$$
Then we determine
$$
       (sI-P_{11})^{-1} = s^{1-p}
\left[ \begin{array}{ccccc} s^{p-2}I_{l_1} & s^{p-3}E_{1,2} &
s^{p-4}E_{1,2}E_{2,3} & \cdots & s^0E_{1,2}E_{2,3}\cdots E_{p-2,p-1}\\ O &
s^{p-2}I_{l_2} & s^{p-3}E_{1,2} & \cdots & s^1E_{2,3}E_{3,4}\cdots E_{p-2,p-1}\\
\vdots & \vdots & \vdots & \ddots & \vdots \\  O& O& O& \cdots & s^{p-2}I_{l_p-1}  \end{array} \right]
\eqno(1.30)$$ and using the structure of $P_{12}$ present the product $(sI-P_{11})^{-1}P_{12}$ as
$$
(sI-P_{11})^{-1}P_{12} =  (sI-P_{11})^{-1} \left[ \begin{array}{c} O \\ O \\
\vdots \\ E_{p-1,p} \end{array} \right] = s^{1-p} \left[ \begin{array}{c}
s^0E_{1,2}E_{2,3}\cdots E_{p-1,p}\\ s^1E_{2,3}E_{3,4}\cdots E_{p-1,p}\\
\vdots \\ s^{p-2}E_{p-1,p} \end{array} \right] \eqno(1.31)$$
 Let
us calculate  terms  $E_{i,i+1}E_{i+1,i+2}\cdots E_{p-1,p}$, which are products of the appropriate matrices
$E_{ij}$. Substituting $E_{i,i+1} = [O, I_{l_i}]$ yields
$$
E_{i,i+1}E_{i+1,i+2}\cdots E_{p-1,p} = \underbrace{[O,
I_{l_i}]}_{l_{i+1}} \underbrace{[O, I_{l_{i+1}}]}_ {l_{i+2}}\cdots
\underbrace{[O, I_{l_{p-1}}]} _{l_p=r} = \\ \underbrace{[O,
I_{l_i}]}_{r}
$$
Then varying $i$ from 1 to $p-1$ we obtain
$$
\begin{array}{rcc} E_{1,2}E_{2,3}\cdots E_{p-1,p} & = & [O,I_{l_1}] \\
                 E_{2,3}E_{3,4}\cdots E_{p-1,p}& = &[O,I_{l_2}]  \\
                                           &   \vdots   &    \\
                  E_{p-2,p-1}E_{p-1,p} & = & [O,I_{l_{p-2}}] \\
                  E_{p-1,p} &= & [O,I_{l_{p-1}}] \end{array}
\eqno(1.32) $$
Substituting  (1.32) in (1.31) gives the following expression
$$
  (sI-P_{11})^{-1}P_{12} = \left[ \begin{array}{c} s^{1-p}[O,I_{l_1}] \\
s^{2-p}[O,I_{l_2}] \\ \vdots \\s^{-1}[O,I_{l_{p-1}}] \end{array} \right] \eqno(1.33)$$
 Then inserting the
right-hand sides of  (1.29),(1.33) into (1.22)  and using the blocks $P_{21}$ and $ sI-P_{22} $ (1.28) we obtain
$$
det(sI_{\bar{n}}-\hat{P})  = s^{\bar{n}-r}det(sI_r + \hat{T}_1 +
[\hat{T}_p,\hat{T}_{p-1},\ldots ,\hat{T}_2] \left[
\begin{array}{c} s^{1-p}[O,I_{l_1}] \\ s^{2-p}[O,I_{l_2}] \\
\vdots \\s^{-1}[O,I_{l_{p-1}}]
\end{array} \right]) =$$
$$ =\; s^{\bar{n}-r}det(sI_r + \hat{T}_1 +
 [O,\hat{T}_2]s^{-1} +[O,\hat{T}_3]s^{-2},\ldots ,[O,\hat{T}_p]s^{1-p})=
s^{\bar{n}-r}det((I_rs^p + T_1s^{p-1} + \cdots + T_p)s^{1-p}) =
$$
$$ = s^{\bar{n}-r}(s^{r(1-p)}det(I_rs^p + T_1s^{p-1} + \cdots + T_p)) =s^{\bar{n}-rp}
det\Phi(s)
$$
     The assertion is proved.

     Further  we  consider  several canonical  forms having the companion  (block  companion,  general
block companion) matrix of dynamics.

\subsection[Controllable (observable) companion canonical
   form of single-input (output) systems]{Controllable (observable)
   companion canonical form of single-input (output) systems}

     Let us consider controllable system (1.1),(1.2) with a scalar input $u$
$$ \dot{x} = Ax +bu $$
$$y(t) = Cx
\eqno (1.34)  $$ where   $b$ is a nonzero constant column vector.  We will find a linear nonsingular
transformation of  state variables
$$    z = Nx
\eqno(1.35) $$
with an nonsingular $n\times n$  matrix $N$ that transforms (1.34) to the controllable canonical
form [M4]
$$ \dot{z} = \hat{A}z +\hat{b}u $$
$$y = CN^{-1}z
\eqno (1.36)  $$
 where $\hat{A} $ is the companion matrix of the
characteristic polynomial $ \phi (s) = s^n + a_1s^{n-1} + \cdots +
a_{n-1}s + a_n$  of the matrix $A$ , i.e. $\hat{A}$ has the form
(1.17), $\hat{b}$ is the $n$ column vector
 $$ \hat{b} = \left[ \begin{array}{c} 0 \\ \vdots \\ 0 \\ 1 \end{array} \right]
\eqno (1.37)$$ The matrix $CN^{-1}$ has no  the special structure. For uniformity we will call (1.36) as the
\it{controllable companion canonical} \rm form.

 \bf {Determination of transformation matrix} \it{N}. \rm  Let us calculate
the controllability matrix of the pair$(\hat{A},\hat{b})$
$$ \hat{Y} = [\hat{b},\hat{A}\hat{b}, \cdots ,\hat{A}^{n-1}\hat{b}] =
\left[ \begin{array}{ccccc} 0 & 0 & 0 & \cdots & 1 \\ 0 & 0 & 0 & \cdots
 & -a_1 \\ \vdots & \vdots & \vdots & \ddots & \vdots \\ 0 & 1 & -a_1 & \cdots&
\cdots \\ 1 & -a_1 & -a_2+a_1^2 & \cdots & \cdots \end{array} \right]
\eqno(1.38) $$
Alternatively, substituting (1.35) into  the  first equation of
(1.34) we obtain
$$  \dot{z} =NAN^{-1}z + Nbu       $$
Thus, $\hat{A}$ and $\hat{b}$ are expressed via $A$, $N$ and $b$
as follows
$$   \hat{A} = NAN^{-1},  \;\;\;   \hat{b} = Nb  $$
Writing the controllability matrix of the pair $(NAN^{-1}, Nb)$ as
$$ \hat{Y} = [Nb,NAb, \cdots ,NA^{n-1}b] = NY
\eqno(1.40)$$ where $Y = [b,Ab, \cdots ,A^{n-1}b]$ is the $n\times n$ controllability matrix of the pair $(A,b)$
we can express   $N$ from (1.40) as
$$ N = \hat{Y}Y^{-1}
\eqno(1.41)$$ Since the matrix $\hat{Y}$ is the lower triangular
matrix then $rank\hat{Y}= n$ and $rank N= n$.
     In the literature it is usually used the matrix $N^{-1} = Y\hat{Y}^{-1}$ [A1], [M4]
having the following simple structure
$$ \hat{Y}^{-1} =
\left[ \begin{array}{ccccc} a_{n-1} & a_{n-2}& \cdots & a_1 & 1 \\
a_{n-2} & a_{n-3}& \cdots & 1 & 0 \\ \vdots & \vdots & \vdots & \ddots
& \vdots \\ a_1 & 1 &  \cdots & 0 & 0 \\ 1 & 0 & \cdots & 0 & 0
 \end{array} \right]
$$

     Thus, we show: if the pair $(A,b)$ is controllable then
the nonsingular transformation (1.35) with $N$ from  (1.41) always exists. This transformation reduces  system
(1.34) to the controllable companion canonical form (1.36). To calculate $N$ it is enough to know the
controllability  matrix of the pair $(A,b)$ and the characteristic polynomial of the dynamics matrix $A$.

     Let us discover a structure of  canonical system  (1.36).
Denoting  variables of the vector $z$ by $z_i$, $i=1,2,...,n$ we can rewrite the first equation in (1.34) as
follows
$$ \begin{array}{ccl} \dot{z}_1 & = & z_2 \\ \dot{z}_2 & = & z_3 \\
& \vdots & \\  \dot{z}_{n-1} & = &z_n \\ \dot{z}_n & = & -a_nz_1 - a_{n-1}z_2 \cdots -a_1z_n + u  \end{array}
\eqno(1.42)$$ It is evident from   (1.42)   that   the each   state   variable $z_i$, $i=1,2,...,n-1$ is the
integral of the following state variable $z_{i+1}$ and $z_{n}$ is  the  integral  of  control $u$ and signals
$a_iz_j$ ($i=n,n-1,...,1;j=1,2,...,n$).
%

If  $l=1, y=z_1$ then we can directly pass
 from  the state  space  representation  (1.42)  to the
 input-output representation (1.15)  with $r=1,\; p=n,\;  k=1$ and
 $ F_i =1,\;i=0,1,\ldots,n;\; B_1 =1$. Indeed, let us denote
$$ y = z_1,\;\; \dot{y} = z_2,\;\; y^{(2)} = z_3,\;\;  \cdots,\;  y^{(n-1)} = z_n   $$
Since $y^{(n)} = \dot{z}_n$ then substituting $y^{(i)}$ ($i=1,2,...,n$) in the last equation of (1.42) gives a
linear differential equation
$$ y^{(n)} + a_1y^{(n-1)} + \cdots +a_ny = u
\eqno(1.43) $$
     The dual result can be obtained  for  a  observable
system. If the system (1.34) has a scalar output $y$, i.e. $C$ is
an $n$ row vector, and the  pair  $(A,C)$  is  observable  then
(1.34) can be transformed  into the observable (companion)
canonical form
$$ \dot{z} = \tilde{A}z +\tilde{B}u $$
$$   y = \tilde{c}z
\eqno (1.44)  $$ where $ \tilde{A} =\bar{A}^T$ , $\tilde{c} =  [0
\; 0 \;\cdots\; 0 \; 1]$,  matrix $ \tilde{B} $  has  no  special
structure.

\subsection[Controllable (observable) block companion
 canonical form of multi-input(output) systems]{Controllable
 (observable) block companion\\
 canonical form of multi-input(output) systems}

  \bf{Asseo's  form } \rm [A4]. Let us consider the controllable  system
(1.1),(1.2) with an $r$  input  vector $u$ ( $r > 1$ ). We will  find a nonsingular transformation of  state
variables  (1.35) which reduces the system to the canonical form having the block companion matrix of dynamics
(see 1.20). This canonical form have been first obtained by Asseo [A4] in 1968.

     Let us propose that $rank B=r$. We define the controllability index of the  pair
$ (A,B)$  as  the smallest integer $\nu (\nu \le n)$ such that
$$  rank \left[ \begin{array}{cccc} B,& AB,& \cdots, & A^{n-1}B
\end{array} \right] =  rank \left[ \begin{array}{cccc} B,& AB,& \cdots, &
A^{\nu-1}B \end{array} \right] \eqno(1.45)  $$ and consider a system with $n = r\nu $ , i.e. $r$ is the divisor
of $n$. Only such type a system is reduced to the canonical form with the block companion matrix of dynamics
(1.20). This canonical form (Asseo's form) is the particular case of Yokoyama's canonical form where $r\nu>n$.

     Let  the  transformation $ z = Nz $ reduces system (1.1),
(1.2) to a canonical form
$$  \dot{z} = A^*z + B^*u  $$
$$  y =C^*z
\eqno(1.46)$$
where $A^*=NAN^{-1}$ is  block companion matrix (1.20): $ A^* = P^* $ with  $p = \nu $ and
$$  B^* = NB = \left[ \begin{array}{c} O \\ I_r  \end{array}  \right]
\eqno(1.47) $$ The matrix $ C^* = CN^{-1} $  has no special structure. We  will   call (1.46) as the
\it{controllable block companion  canonical} \rm form a the multi-input system.

 \bf{ Determination of transformation matrix \it N. \rm  The  matrix $N$ is
partitioned as
$$ N = \left[ \begin{array}{c}  N_\nu \\ N_{\nu-1} \\ \vdots \\ N_1
\end{array}  \right]
 \eqno(1.48) $$
where $N_i$ are $r\times n$ submatrices. Since $A^* =NAN^{-1}$
then substituting  (1.20) for $P^*=A^*$, $p = \nu$  and (1.48) for
$N$ into the equality  $A^*N =NA$ gives the following matrix
equation
$$
\left[ \begin{array}{ccccc} O & I_r & O & \cdots & O\\ O& O& I_r& \cdots & O\\
\vdots & \vdots & \vdots & \ddots & \vdots \\  O& O& O& \cdots & I_r\\
-T_p & -T_{p-1} & -T_{p-2} & \cdots & -T_1 \end{array} \right]
\left[ \begin{array}{c} N_\nu \\ N_{\nu-1} \\ \vdots \\ N_2 \\ N_1
\end{array} \right]
=\left[ \begin{array}{c} N_\nu \\ N_{\nu-1} \\ \vdots \\ N_2 \\ N_1
\end{array} \right]  A   $$
from which blocks $N_i$ ($i=1,2,\ldots ,\nu $)  are defined as follows
$$
\begin{array}{ccccc} N_{\nu-1} & = & N_{\nu}A & &  \\
N_{\nu-2} & = & N_{\nu-1}A & = & N_{\nu}A^2 \\  & \vdots & & & \\
N_1 & = & N_2A & = & N_{\nu}A^{\nu -1} \end{array} \eqno(1.49) $$
Thus, any  block $N_i$ ($i=1,2,\ldots ,\nu$)
is  defined via the block $N_\nu$. To determine $N_\nu$ we shall use the approach of Sec.1.2.2. Since $A^*
=NAN^{-1}$, $B^* =NB$ then we can express  blocks of the  controllability matrix of the pair $(A^*,B^*)$ via
matrices  $A,B,C$ as follows
$$ Y^* = [B^*,A^*B^*,\ldots , (A^*)^{n-1}B^*] = [NB,NAB, \ldots ,NA^{n-1}B] =
N[B,AB,\ldots ,A^{n-1}B] $$
From this equality we obtain
$$ (A^*)^{i-1}B^* = NA^{i-1}B,  \;\;\; i=1,2, \ldots ,n
\eqno(1.50)$$
  Blocks  $A^kB$  $(k> \nu -1)$  are
linearly  dependent on $A^kB$  $(k\le \nu -1)$ because the pair $(A,B)$ has the controllability index $\nu$ and
satisfies the condition (1.45).

     Let us consider the $ n\times n$  matrix $\bar{Y}^* = [B^*,A^*B^*,\ldots ,
(A^*)^{\nu-1}B^*]$. Using (1.50) we have
$$
\bar{Y}^* = N[B,AB,\ldots ,A^{\nu-1}B] \eqno (1.51)$$ Calculating products $(A^*)^iB^* $ with $A^*$ and $B^*$
from (1.20) and (1.47) we reveal the structure of the matrix $\bar{Y}^*$
$$ \bar{Y}^* =
 \left[ \begin{array}{cccc} B^*, & A^*B^*, & \ldots ,& (A^*)^{\nu-1}B^*
\end{array} \right] = \left[ \begin{array}{ccccc} O & O& \cdots & O & I_r\\
O & O &\cdots  & I_r & -T_1 \\ \vdots & \vdots & \ddots & \vdots & \vdots \\
0 & I_r & \cdots & X & X \\ I_r & -T_1 & \cdots & X & X
\end{array} \right]
\eqno(1.52) $$ where $X$ are  some unspecified  submatrices.  The matrix $\bar{Y}^*$ is nonsingular one because
it has unity blocks on the diagonal, i.e. $rank\bar{Y}^* = n$. Substituting the right-hand side of (1.52) into
the left-hand side of (1.51) we obtain the equality
$$
\left[ \begin{array}{ccccc} O & O& \cdots & O & I_r\\
O & O &\cdots  & I_r & -T_1 \\ \vdots & \vdots & \ddots & \vdots & \vdots \\
0 & I_r & \cdots & X & X \\ I_r & -T_1 & \cdots & X & X
\end{array} \right]  =
\left[ \begin{array}{c} N_\nu \\ N_{\nu-1} \\ \vdots \\ N_2 \\ N_1
\end{array} \right]  \left[ \begin{array}{cccc} B,& AB,& \cdots, &
A^{\nu-1}B \end{array} \right]
\eqno(1.53) $$
from which it
follows the equation for  $N_\nu$
$$  N_\nu
\left[ \begin{array}{cccc} B,& AB,& \cdots, & A^{\nu-1}B
\end{array} \right] = \left[ \begin{array}{ccccc} O,& O,& \cdots,
& O, & I_r  \end{array} \right] \eqno(1.54) $$
Thus
$$
 N_\nu = \left[ \begin{array}{ccccc} O,& O,& \cdots, & O, & I_r \end{array}
\right] \left[ \begin{array}{cccc} B,& AB,& \cdots, & A^{\nu-1}B
\end{array} \right] ^{-1} \eqno(1.55)$$ Others blocks
$N_i$, $i=1,2,\ldots,\nu-1$   are  calculated  by formulas (1.49).
     It should be noted that obtained blocks $N_i$,$i=1,2,\ldots,\nu-1$
satisfy   relation  (1.53).  Indeed,   the following  equalities
take place from  (1.49)  and (1.54)
$$
N_{\nu-1} \left[ \begin{array}{cccc} B,& AB,& \cdots, & A^{\nu-1}B
\end{array} \right] \;= \; N_\nu A \left[ \begin{array}{cccc}
B,& AB,& \cdots, & A^{\nu-1}B \end{array} \right] \;=
$$
$$
=\; N_\nu \left[ \begin{array}{cccc} AB,& A^2B,& \cdots, & A^\nu B \end{array}
 \right] \;= \;\left[ \begin{array}{cccccc} O,& O,& \cdots, & O, & I_r,& X
\end{array} \right] ,
$$
$$  N_{\nu-2} \left[ \begin{array}{cccc} B,& AB,& \cdots, & A^{\nu-1}B
\end{array} \right] \;=
\;= N_\nu A^2 \left[ \begin{array}{cccc} B,& AB,& \cdots, & A^{\nu-1}B
\end{array} \right] \;=
$$
$$
=\; N_\nu \left[ \begin{array}{cccc} A^2B,& A^3B,& \cdots, & A^{\nu+1}B
\end{array} \right] \;=\; \left[ \begin{array}{ccccccc}
 O,& O,& \cdots, & O, & I_r,& X,& X  \end{array} \right]
$$
and so on.

The matrix $N$ is nonsingular one. It follows from  nonsingularity of the  matrix  in the left-hand  side  of
(1.53)  and  the matrix $[B, AB, \cdots,  A^{\nu-1}B]$.

     Thus, we show: if the pair matrix $(A,B)$ is
controllable with the controllability index  $\nu =  n/r$ and $rankB = r$  then the nonsingular  transformation
$z = Nx $ with $N$ from (1.48), (1.49), (1.54) always  exists. This transformation reduces  system  (1.1),(1.2)
to the canonical form (1.46).

     The analogous dual result can be obtained  for  an  observable
system. Let $rankC =l$ and the pair $(C,A)$ is  observable with the  observability index $\alpha $, which is a
smallest integer such as
$$
 rank [C^T,A^T,C^T,..., (A^T)^{n-1}C^T] \;=
\;rank [C^T,A^T,C^T,..., (A^T)^{\alpha -1}C^T] = n
$$
Let $\alpha = n/l$.  Then system (1.1), (1.2) can be transformed
into  the observable block companion canonical form
$$  \dot{z} = \tilde{A}z + \tilde{B}u  $$
$$    y =\tilde{C}z  $$
where $\tilde{A}=(P^*)^T$ ,$\tilde{C}=[ O,O,\ldots ,O,I_l]$ and
matrix $\tilde{B}$ has  no special structure.

     Let us consider  the  structure  of  the canonical  system  (1.46).
We introduce  subvectors $\bar{z}_i$, $i=1,2,\ldots,\nu$
$$
\bar{z}_1 = \left[ \begin{array}{c} z_1 \\ z_2 \\ \vdots \\ z_r \end{array}
\right], \bar{z}_2 = \left[ \begin{array}{c} z_{r+1} \\ z_{r+2} \\ \vdots \\
z_{2r} \end{array} \right], \cdots
$$
where $z_i$, $i=1,2,\ldots ,n$  are  components of the vector $z$. Using these notions and the block structure
of $A^*$ we  rewrite the first equation in (1.46) as follows
$$ \begin{array}{lcl} \dot{\bar{z}}_1 & = &\bar{z}_2 \\
\dot{\bar{z}}_2 & = & \bar{z}_3 \\
& \vdots & \\  \dot{\bar{z}}_{\nu -1} & = &\bar{z}_\nu \\
\dot{\bar{z}}_\nu & = & -T_\nu \bar{z}_1 - T_{\nu -1}\bar{z}_2 \cdots -T_1\bar{z}_\nu + I_ru  \end{array}
\eqno(1.56)$$
In (1.56) the each group of state variables $\bar{z}_i$ ($i=1,2,\ldots,
 \nu -1$) is the integral of the next group $\bar{z}_{i+1}$  and $\bar{z}_\nu $
is the integral of the control vector $u$ and vectors $T_i\bar{z}_j$ ($i=\nu, \nu -1, \ldots ,1$;
$j=1,2,\ldots,\nu$). The general structure of (1.56) coincides with the structure of (1.42) with $n = \nu$, $z_i
=\bar{z}_i$, $ a_i = T_i$,  $CN^{-1} = [ C_1,C_2,\ldots,C_\nu ]$  where $C_i$ are $l\times r$ submatrices.

Let us show that for $l =r$  we  can  pass from the state-space  representation  (1.56)  to the  input-output
representation  (1.15). In  fact, defining  the output  vector for (1.56) as $\bar{y} =\bar{z}_1$ and using
(1.56) we obtain
$$ \bar{y} = \bar{z}_1,\qquad \dot{\bar{y}} = \bar{z}_2 , \;\;\cdots,\;\;
 \bar{y}^{(\nu-1)} = \bar{z}_\nu   $$
So far as
$$ \bar{y}^{(\nu )} =\bar{z}_\nu^{(1)} = -T_\nu \bar{z}_1 - T_{\nu-1}\bar{z}_2
- \cdots - T_1\bar{z}_\nu + I_ru$$
then replacing $\bar{z}_i$ by  $\bar{y}^{(i)}$ in the last expression  we obtain
$$ \bar{y}^{(\nu )} = -T_\nu \bar{y} - T_{\nu-1}\bar{y}^{(1)} -
T_{\nu-2}\bar{y}^{(2)} - \cdots - T_1\bar{y}^{(\nu-1)} + I_ru \eqno(1.57)$$ The vector differential equation
(1.57) coincides with (1.15) when $ p = \nu $, $F_p = I_r$, $B_o = I_r$, $B_i = O$ ($i=1,2,\ldots,k$).

 \bf {Yokoyama's form } \rm [Y1], [Y2]. Let us consider the general case of
system (1.1), (1.2) with $r$ input vector $u$ $(r > 1)$, $rankB = r$ and the controllability index $\nu \ne
n/r$, i.e. $r$ does not the divisor of $n$: $n < r\nu $ .  Using  the nonsingular transformation of  state
variables (1.35)  and  input variables
$$                v = M^{-1}u   \eqno(1.58)   $$
where $M$ is an $r\times r$ permutation\footnote{A permutation matrix has a single unity element in each row
(column) and zeros otherwise.} matrix we  can  reduce   system (1.1),  (1,2)  to  the canonical  form  with the
general  block companion matrix of dynamics (1.25). This canonical form have been worked out by Yokoyama in 1972
[Y1].

     For a pair of matrices $A$ and  $B$  with the controllability
index $\nu$ we define the integers $l_1,l_2, \ldots, l_\nu$ by the
rule
$$ l_1 = rank [B, AB, \cdots,  A^{\nu-1}B] -  rank[B, AB, \cdots, A^{\nu-2}B] $$
$$ l_2 =  rank[B, AB, \cdots,  A^{\nu-2}B] -  rank[B, AB, \cdots, A^{\nu-3}B] $$
$$   \cdots   \eqno(1.59)$$
$$ l_{\nu -1} = rank [B,AB] - rankB  $$
$$l_\nu = rankB = r$$
From (1.59) it follows that
$$  l_1 \le l_2 \le \cdots \le l_\nu                \eqno(1.60)$$
Let us determine the sum  of $l_i$ ,$i=1,2, \ldots, \nu$ by adding the left-hand and the right-hand sides of
(1.59). We obtain the relation
$$ l_1 + l_2 + \cdots l_\nu = rank[B,AB, \ldots , A^{\nu -1}B] = n $$

     Now we use  transformation  (1.35), (1.58) to reduce    system
(1.1), (1.2) to  Yokoyama's canonical form
$$ \dot{z} = Fz +Gv $$
$$ y = CN^{-1}z
\eqno (1.61)  $$ where the matrix $F = NAN^{-1}$ is the general block companion  matrix of the structure (1.25)
with   $p =\nu $, $\bar{n} = n$, $ -\hat{T}_p = F_{\nu1}$, $ -\hat{T}_{p-1} = F_{\nu2}$, $\ldots$, $ -\hat{T}_1
= F_{\nu\nu}$
$$
 F =
\left[ \begin{array}{ccccc} O & E_{1,2} & O & \cdots & O\\ O& O& E_{2,3}& \cdots & O\\\vdots & \vdots & \vdots &
\ddots & \vdots \\  O& O& O& \cdots & E_{\nu-1,\nu}\\ F_{\nu1} & F_{\nu2} & F_{\nu3} & \cdots & F_{\nu\nu}
\end{array} \right]
\eqno(1.62)
$$
and  blocks $E_{i,i+1}$ of the form
$$ E_{i,i+1} = [O, I_{l_i}] ,\qquad i=1,2,\ldots, \nu -1
\eqno(1.63)$$
 In (1.62)  blocks $F_{\nu i}$ have the  increasing numeration for convenience. In
(1.61) the matrix $G = NBM$  has the form
$$  G = \left[ \begin{array}{c} O \\ G_\nu  \end{array}  \right]
\eqno(1.64) $$ where an $r\times r$ block $G_\nu$ is a lower triangular  matrix  with   unity diagonal elements
$$  G_\nu = \left[ \begin{array}{cccc} I^1 & O & \cdots & O \\
X & I^2 & \cdots & O \\ \vdots & \vdots & \ddots & \vdots \\
X & X & \cdots & I^{\nu}  \end{array}  \right]
 \eqno(1.65)$$
 Here $I^i$, $i=1,2, \ldots, \nu$  are  unity matrices of
orders $l_{\nu -i+1} -l_{\nu -i}$, $l_o =0$, $X$  are some unspecified  submatrices. Let us note that the matrix
$CN^{-1}$ has no special structure. We will  call  (1.61) as the \it{controllable generalized block companion
canonical} \rm  form or Yokoyama's form.

\bf{Determination of transformation matrix} \it N \rm [S3]. At first we construct the $n\times r\nu $ matrix
containing the first $\nu $ blocks $A^iB$, $i=0,1, \ldots,\nu -1 $ of the controllability matrix $B,AB,
\ldots,A^{n-1}B$  satisfying the equality
 $$  rank[B,AB, \ldots,A^{\nu-1}B]   = n  $$
Let us find a permutation matrix $M$ rearranging columns  of $B$ such that the matrix $$[BM,ABM,
\ldots,A^{\nu-1}BM]$$  has linearly independent columns in last columns of $A^iBM$ ,$i=0,1, \ldots, \nu -1 $.
From (1.59) it  follows  that the blocks $A^iBM$  maintain $l_{\nu -i}$ linearly  independent columns.

    The matrix N is partitioned as follows
$$ N = \left[ \begin{array}{c}  N_\nu \\ N_{\nu-1} \\ \vdots \\ N_1
\end{array}  \right]  $$
where  $l_{i+1}\times n$  blocks $N_{\nu -i}$, $i=0,1, \ldots, \nu -1 $ have  the following structure
$$  N_{\nu -i} = \left[ \begin{array}{c} P_{\nu -i} \\
\tilde{P}_{\nu -i} \end{array}  \right]
 \eqno(1.66) $$
 In (1.66)
$P_{\nu -i}$ are $(l_{i+1} - l_i)\times n$ submatrices ($i=0,1, \ldots,\nu -1$), $l_o = 0$. Using the equality $
FN = NA $ and the structure of $F$ (1.62)  we present the blocks $N_\nu, N_{\nu -1}, \ldots N_1$ as
$$ E_{1,2}N_{\nu -1} = N_\nu A, \qquad E_{2,3}N_{\nu -2} = N_{\nu -1}A,
\qquad \ldots \qquad E_{\nu -1,\nu}N_1 = N_2 A
$$
and by (1.63) obtain $\tilde{P}_{\nu -i}$
$$\tilde{P}_{\nu -1}  = N_\nu A, \qquad \tilde{P}_{\nu -2} = N_{\nu -1}A,
\qquad \ldots \qquad \tilde{P}_1 = N_2 A
$$
Substituting  last expressions into (1.66)  we find the structure of  blocks $N_i$, $i = \nu, \nu -1, \ldots, 1$
$$   N_\nu =P_\nu,$$
$$  N_{\nu-i} = \left[ \begin{array}{c} P_{\nu -i} \\  \dotfill \\
N_{\nu -i+1}A \end{array}  \right] =  \left[ \begin{array}{c} P_{\nu -i} \\
\dotfill \\ P_{\nu -i+1}A \\ P_{\nu -i+2}A^2 \\ \vdots \\ P_\nu A^i
\end{array}  \right] \begin{array}{cl} \} & l_{i+1} - l_i \\  & \\
\} & l_i - l_{i-1} \\ \} & l_{i-1} - l_{i-2} \\ \vdots & \\  \} & l_1
\end{array}
\eqno(1.67) $$

     Thus, the blocks $N_{\nu -1}$ of the matrix $N$ are  defined  via  the  blocks
$P_\nu,P_{\nu -1}, \ldots, P_1$. To find these blocks we use controllability matrices of pairs $(A,BM)$ and
$(F,G)$ denoted as $Y_{FG}$ and $Y_{A,BM}$ respectively. Since $F = NAN^{-1}$, $G = NBM$ then $Y_{FG}$ is
expressed via $Y_{A,BM}$ as follows
$$ Y_{FG} = [NBM, NABM, \ldots, NA^{n-1}BM] = N[BM, AMB, \ldots, A^{n-1}BM] =
N Y_{A,BM} $$
Let us denote by
$$ \tilde{Y} =  N[BM, AMB, \ldots, A^{\nu -1}BM]
\eqno (1.68) $$
the $n\times r\nu $ submatrix of the matrix $Y_{FG}$. On the other hand the matrix $\tilde{Y}$
can be constructed from matrices $F$ and $G$ as follows
$$  \tilde{Y} = [ G, FG, \ldots, F^{\nu -1}G] = [ \tilde{Y}_1, \tilde{Y}_2,
\ldots, \tilde{Y}_\nu]
\eqno(1.69) $$
where
$$ \tilde{Y}_1 = G,  \qquad  \tilde{Y}_i = F\tilde{Y}_{i-1}, \qquad
i =2,3, \ldots, \nu $$
Using  formulas (1.62)  and  (1.65)  we   can  find  $n\times r$ matrices $\tilde{Y}_i$
$$ \tilde{Y}_1 = \left[ \begin{array}{c}  O \\ \Theta_1 \end{array}  \right]
$$
$$ \tilde{Y}_i =\left[ \begin{array}{cc} O & O \\ \dotfill & \dotfill \\
X & \Theta_i \\ \dotfill & \dotfill \\ X & X \end{array}  \right]
\begin{array}{cl} \} & l_1 + l_2 + \cdots +l_{\nu -i} \\ & \\
\} & l_{\nu -i+1} \\  &  \\ & \end{array} \eqno (1.70)$$
 where  square matrices $\Theta_i$ of the order $l_{\nu
-i+1}$ are lower triangular matrices
$$  \Theta_1 =G_\nu, \qquad \Theta_i
 = \left[ \begin{array}{cccc} I^i & O & \cdots & O \\
X & I^{i+1} & \cdots & O \\ \vdots & \vdots & \ddots & \vdots \\
X & X & \cdots & I^{\nu}  \end{array}  \right] ,\qquad i=2,3,\ldots ,\nu \eqno(1.71)$$ In (1.71) $I^i$  are
unity matrices of the order $l_{\nu -i+1}-l_{\nu -i}$, $l_o = 0$, $X$ are some matrices. From (1.68) and (1.69)
we obtain the relation
$$ \left[ \begin{array}{c}  N_\nu \\ N_{\nu-1} \\ \vdots \\ N_1
\end{array} \right] [BM, AMB, \ldots, A^{\nu- 1}BM] \qquad=\qquad
[ \tilde{Y}_1, \tilde{Y}_2, \ldots, \tilde{Y}_\nu] \eqno(1.72) $$ Let us denote  $l_{\nu -i+1}$ last (linearly
independent) block columns of the submatrix $ A^{\nu -1}BM$ by $ V_i$, $i=2,\ldots,\nu $
$$
  V_i = A^{i- 1}BM \left[ \begin{array}{c} O \\ I_{l_{\nu -i+1}}
\end{array} \right] , i = 2, \ldots, \nu ,\qquad V_1 = BM
\eqno (1.73) $$
The matrix $V = [ V_1, V_2, \ldots, V_\nu ] $ of the  size  $n\times \sum_{i=1}^\nu l_i =
n\times n$ is  the nonsingular  square  matrix. Using  (1.72),(1.70) we  find structure of the product $NV$
$$ \left[ \begin{array}{c}  N_\nu \\ N_{\nu-1} \\ \vdots \\ N_1
\end{array} \right] [ V_1, V_2, \ldots, V_\nu ] \: =
\:  \left[ \begin{array}{cccllc} O & O & \cdots & O & O & \Theta_\nu \\
O & O & \cdots  & O & \Theta_{\nu -1} & X \\ O & O & \cdots & \Theta_{\nu -2} &
 X &  X \\ \vdots & \vdots & \cdots & \vdots & \vdots & \vdots \\
O & \Theta_2 & \cdots & X & X & X \\
\underbrace{\Theta_1}_{l_{\nu}} & \underbrace{X}_{l_{\nu-1}} &
\cdots & X & \underbrace{X}_{l_{2}} & \underbrace{X}_{l_{1}}
\end{array} \right] \begin{array}{cl} \} & l_1 \\  \} & l_2 \\ \} & l_3 \\
\vdots & \\ \vdots & \\ \vdots & \\
\end{array}
\eqno(1.74) $$
where $X$ are some unspecified submatrices. Taking
into account that $P_i$ are  $(l_{\nu -i+1}-l_{\nu -i}) \times n$
upper blocks of $N_i$  and  matrices $ \Theta_i$ have the
structure (1.71) as well as using the equality $NV_1 = G$ we can
rewrite  relation (1.74) in terms of blocks $P_i$, $i=1,2,\ldots,
\nu$
$$ \left[ \begin{array}{c}  P_\nu \\ P_{\nu -1} \\ \vdots \\ P_1
\end{array} \right] [ V_1, V_2, \ldots, V_\nu ] = \:
\: \left[ \begin{array}{lcccclcclccc} O & O & \vdots &\cdots & \vdots &
O & O & \vdots & O & O & \vdots & I^\nu \\   O & O & \vdots & \cdots & \vdots &
O & O & \vdots & I^{\nu -1} & O & \vdots & X \\   O & O & \vdots & \cdots &
\vdots &  I^{\nu -2} & O & \vdots & X &  X & \vdots & X \\
\vdots & \vdots & \vdots & \cdots & \vdots & \vdots & \vdots & \vdots &
\vdots &\vdots & \vdots & \vdots \\
I^1 & O & \vdots & \cdots & \vdots & X & X & \vdots & X & X &
\vdots & X
\end{array} \right] \begin{array}{cl} \} & l_1 \\  \} & l_2-l_1 \\ \} & l_3-l_2 \\
\vdots & \\ \} & l_{\nu}-l_{\nu-1}\\
\end{array}
\eqno(1.75) $$ where   matrix block columns have sizes $n\times l_\nu, n\times l_{\nu -1}, \ldots, n\times l_1$
respectively, blocks  $X$ are  some unspecified  submatrices. Equation (1.75) may be used for calculating
blocks $P_i$. Then matrices $N_i$ are obtained by relation (1.67).

     Let us demonstrate that these $N_i$ satisfy (1.74). At first we
evaluate the block $N_\nu V$
$$ N_\nu V\: = \: P_\nu V\: = \: P_\nu[V_1,V_2, \ldots, V_\nu] \: = \:
 [ O,O,\ldots, O,I^\nu]
\eqno(1.76)  $$ Then we find the block $N_{\nu -1} V$
$$ N_{\nu -1}V \: = \: \left[ \begin{array}{c} P_{\nu -1}V \\ P_\nu AV
\end{array} \right] \: = \:\left[ \begin{array}{c} P_{\nu -1}
 [V_1,V_2, \ldots, V_\nu] \\ P_\nu [AV_1,AV_2, \ldots, AV_\nu]
\end{array} \right]
 $$
Using (1.75) we obtain
$$ P_{\nu -1}[V_1,V_2, \ldots, V_\nu] \: = \: [ O,O,\ldots, I^{\nu -1}, O, X]
\eqno (1.77)$$ and find  blocks $P_\nu AV_i$, $i=\nu,\nu
-1,\ldots,1$ of the matrix $P_\nu AV = P_\nu [AV_1,AV_2, \ldots,
AV_\nu]$.
 From (1.73)  it follows
 $  V_\nu = AV_{\nu-1} \left[ \begin{array}{c} O \\ I_{l_1}
\end{array} \right] $. Thus
$$
P_\nu AV_{\nu -1} = P_\nu (AV_{\nu-1} \left[ \begin{array}{c} I_{l_2 -l_1} \\ O
\end{array} \right],AV_{\nu-1} \left[ \begin{array}{c} O \\ I_{l_1}
\end{array} \right]) \: = \:
 P_\nu (AV_{\nu-1} \left[ \begin{array}{c} I_{l_2 -l_1} \\ O \end{array} \right]
, V_\nu) \: = \: [ X, I^\nu]
$$
Other products $P_\nu AV_i$ , $i=1,2,\ldots,\nu -2$ are equaled to
zeros because  columns of matrices $AV_1$,\ldots, $AV_{\nu -1}$
are linearly dependent on blocks $V_1,V_2,\ldots, V_{\nu -1}$ for
which   equality (1.76) is true. We result in
$$
  P_\nu AV \: = \: [ O,O,\ldots, O, X, I^\nu, X] $$
 Uniting (1.77) with the last expression we find
$$ N_{\nu -1}V \: =\: \left[ \begin{array}{cccllc} O & O & \cdots &
I^{\nu -1} & O & X \\ O & O & \cdots & X & I^\nu & X  \end{array} \right]
\: = \: \left[ \begin{array}{ccclc} O, & O, & \cdots, &  \Theta_{\nu -1}, X]
\end{array} \right]$$
Now it is evident that the right-hand  side  of  the  last  formula coincides with  the  second  block row of
the  matrix  in   the right-hand side of (1.74). And so on.

Then we need to show that the matrix $G =NBM$ coincides  with (1.64). Calculating  $N_\nu BM,N_{\nu -1}BM,
\ldots, N_1BM$ and using $BM = V_1$ we obtain from   (1.74) that
$$  N_iBM = O, \;\; i=1,\nu, \ldots, 2, \qquad N_1BM = \Theta_1 = G_\nu $$

 {\it REMARK 1.2.} \rm  For $l_1 = l_2 = \cdots = l_{\nu} = r$ (Asseo's  form)
 we  have $\nu = n/r, \;\; V_1 = B,\; \; V_2 = AB,  \ldots, V_\nu = A^{\nu -1}B,
 \;\;
 l_i = r,\;\; l_{i+1} -l_i = 0\;\; (i=1,2, \ldots, \nu -1), \;\; M = I_r,\; \; N_\nu =P_\nu$.
Therefore, equation (1.75) may be rewritten as
$$  P_\nu
\left[ \begin{array}{cccc} B,& AB,& \cdots, & A^{\nu -1}B \end{array} \right]
\: = \: \left[ \begin{array}{ccccc} O,& O,& \cdots, & O, & I_r
\end{array} \right]   $$
The matrix $N$ has the following simple structure
$$ N = \left[ \begin{array}{c}  N_\nu \\ N_{\nu -1} \\ \vdots \\ N_1
\end{array}  \right] \; = \; \left[ \begin{array}{c}  P_\nu \\ P_\nu A \\
\vdots \\ P_\nu A^{\nu -1} \end{array}  \right]
$$
Let us note that the last  formula  coincides  with  (1.54), (1.49) respectively. So, Asseo's form is the
particular case of Yokoyama's form.

{\it REMARK 1.3.}  \rm If $rankC = l$ and  the  pair  $(A,C)$ is observable with  the  observability  index
$\alpha < n/l$ then system  (1.1), (1.2) can be transformed into the observable generalized  block companion
canonical form
$$  \dot{z} = \tilde{A}z + \tilde{B}u  $$
$$    y =\tilde{C}z  $$
where $\tilde{A}$ and   $\tilde{C}$ are
$$  \tilde{A}=F^T , \qquad \tilde{C}=[ O,O,\ldots ,O,\tilde{G}_\alpha^T]
$$

     Let us show that the structure of  canonical  system  (1.61)
resembles with (1.42) or (1.56). We combine components $z_i$, $i=1,\ldots,n$  of the vector $z$ into subvectors
$\tilde{z}_1,\tilde{z}_2, \ldots, \tilde{z}_\nu $    by the  rule
$$
\tilde{z}_1 = \left[ \begin{array}{c} z_1 \\ z_2 \\ \vdots \\ z_{l_1}
 \end{array} \right], \tilde{z}_2 = \left[ \begin{array}{c} z_{l_1+1} \\
z_{l_1+2} \\ \vdots \\ z_{l_1 +l_2} \end{array} \right], \cdots
$$
Using the block structure of $F$  we can rewrite the first equation in (1.61) as
$$ \begin{array}{lcl} \dot{\tilde{z}}_1 & = & [O,I_{l_1}]\tilde{z}_2 \\
\dot{\tilde{z}}_2 & = & [O,I_{l_2}]\tilde{z}_3 \\ & \vdots & \\
\dot{\tilde{z}}_{\nu-1} & = & [O,I_{l_{\nu -1}}]\tilde{z}_\nu \\
\dot{\tilde{z}}_\nu & = & F_{\nu 1}\tilde{z}_1 + F_{\nu
2}\tilde{z}_2 + \cdots + F_{\nu \nu}\tilde{z}_\nu + G_\nu v
\end{array}
\eqno(1.78)$$ In (1.78) the each subvector $\tilde{z}_i$ is the integral  of   last components of the subvector
$\tilde{z}_{i+1}$ and $\tilde{z}_\nu $ is the integral of the control vector $G_\nu v$ and   vectors $F_{\nu
i}\tilde{z}_i$, $i=1,2,\ldots,\nu $.

Now let us show that  for $l = r$  it  is  possible  to  transfer from the state-space   representation   (1.78)
to the input-output representation (1.15), which is set of  linear differential equations of  the  order $\nu$.
For this purpose we introduce the $r$ vector $y$ containing  $l_i-l_{i-1}$ subvectors $y_i$ ($i=1,\ldots, \nu
$), $l_o = 0$
$$
y = \left[ \begin{array}{c} y_\nu \\ y_{\nu -1} \\ \vdots \\ y_1
\end{array} \right]
\eqno(1.79) $$
  Denoting
$$
 \tilde{z}_1 = y_1, \: \tilde{z}_2 = \left[ \begin{array}{c} \dot{y}_2 \\
\dot{y}_1 \end {array} \right ], \:   \tilde{z}_3 = \left[ \begin{array}{c}
 y^{(2)}_3 \\ y^{(2)}_2 \\ y^{(2)}_1  \end{array} \right ], \:
\cdots, \:  \tilde{z}_\nu =  \left[ \begin{array}{c}  y^{(\nu -1)}_\nu \\
 y^{(\nu -1)}_{\nu -1} \\ \vdots \\ y^{(\nu -1)}_1 \end{array} \right ]
$$
and using (1.79) we  express
$$
\tilde{z}_1  =  [O,I_{l_1}]y, \;\; \tilde{z}_2  =
[O,I_{l_2}]\dot{y}, \;\; \cdots , \;\; \tilde{z}_{\nu-1} =
[O,I_{l_{\nu -1}}]y^{(\nu -2)}, \;\; \tilde{z}_\nu  = y^{(\nu -1)}
\eqno(1.80)$$
Since
$$ y^{(\nu)} = \dot{\tilde{z}_\nu} =  F_{\nu 1}\tilde{z}_1 +
F_{\nu 2}\tilde{z}_2 + \cdots + F_{\nu \nu}\tilde{z}_\nu + G_\nu v  $$ then replacing $\tilde{z}_i$ by $y^{(i)}$
in the last equation
$$ y^{(\nu)} =  F_{\nu 1}[O,I_{l_1}]y + F_{\nu 2}[O,I_{l_2}]\dot{y} + \cdots
+ F_{\nu \nu}y^{(\nu -1)} + G_\nu v  $$
 and performing  multiplications we get
$$ y^{(\nu)} = [O, F_{\nu 1}]y + [O,F_{\nu 2}]\dot{y} + \cdots
+ [O,F_{\nu \nu}]y^{(\nu -1)} + G_\nu v \eqno(1.81) $$ This vector differential equation coincides with (1.15)
when $ p = \nu,\;\; F_p = I_r,\;\; F_i = -[O,F_{\nu i}] \;\;(i=1,\ldots, \nu -1),\;\; B_i = O (i\ne 0),\;\; B_o
= G_\nu. $

 Let us consider several examples.

{\it EXAMPLE 1.1.} \rm

 We need to find the controllable canonical form of the system
with $n=3$ and $r=1$
$$ \dot{x} = \left[ \begin{array}{ccc} 2 & 1 & 0 \\ 0 & 1 & 1 \\
1 & 0 & 0 \end{array} \right ]x + \left[ \begin{array}{c} 1 \\ 0 \\ 0
\end{array} \right ]u
\eqno(1.82)$$
 Since $detY  = det [b,Ab,A^2b] = - 1 \ne 0$ then the system is controllable. We calculate a
characteristic polynomial of $A$:  $ \phi (s) = det(sI - A) = s^3 + a_1s^2 + a_2s + a_3 = s^3 - 3s^2 + 2s - 1$
and find $a_3 = -1$, $a_2 = 2$, $ a_1 = -3$. Now we construct matrices $Y = [b,Ab,A^2b]$ and $\hat{Y}^{-1} =
\left[
\begin{array}{ccc} a_2 & a_1 & 1 \\ a_1 & 1 & 0 \\ 1 & 0 & 0 \end{array}
\right ]$  and calculate
$$ Y = \left[ \begin{array}{ccc} 1 & 2 & 4 \\ 0 & 0 & 1 \\ 0 & 1 & 2 \end{array}
\right ] , \qquad \hat{Y}^{-1} = \left[ \begin{array}{rrr} 2 & -3 & 1 \\
-3 & 1 & 0\\ 1 & 0 & 0 \end{array} \right ]
$$
By formula (1.41) we find
$$ N^{-1} = \left[ \begin{array}{rrr} 0 & -1 & 1 \\ 1 & 0 & 0 \\ -1 & 1 & 0
\end{array} \right ] $$
 Since
$$  N = \left[ \begin{array}{ccc} 0 & 1 & 0 \\ 0 & 1 & 1\\ 1 & 1 & 1
\end{array} \right ]
$$
then using formula (1.39) we obtain
$$ \hat{A} = NAN^{-1} = \left[ \begin{array}{rrr} 0 & 1 & 0 \\ 0 & 0 & 1\\
1 & -2 & 3 \end{array} \right ], \hat{b} = \left[ \begin{array}{c}
0 \\ 0 \\ 1 \end{array} \right ]
\eqno(1.83) $$
It is evident that the matrix $\hat{A}$ in  (1.83)  is  the  companion
matrix of the polynomial  $\phi(s) =  s^3 - 3s^2 + 2s - 1$.

\it{EXAMPLE 1.2.} \rm

Let us consider the following system with $n = 4,\;\; r = 2$
$$ \dot{x} = \left[ \begin{array}{cccc} 2 & 1 & 0 & 1 \\ 1 & 0 & 1 & 1 \\
 1 & 1 & 0 & 0 \\ 0 & 0 & 1 & 0 \end{array} \right ]x +
\left[ \begin{array}{cc} 1 & 0 \\ 0 & 0 \\ 0 & 0 \\ 0 & 1
\end{array} \right ]u
\eqno (1.84) $$
Since $det[B,AB] = det \left[ \begin{array}{cccc} 1 & 0 & 2 & 1 \\
0 & 0 & 1 & 1 \\ 0 & 0 & 1 & 0 \\ 0 & 1 & 0 & 0 \end{array} \right ] = -1$
then  the  system  is controllable with the controllability index $\nu = 2
 = n/r$ and can  be  transformed  into  the  controllable  block
companion canonical form (Asseo's  form). Let us find the related transformation matrix $N$. It has the
following structure
$$ N = \left[ \begin{array}{cc} N_2 \\ N_1 \end{array} \right ] =
\left[ \begin{array}{cc} N_2 \\ N_2A \end{array} \right ] \eqno(1.85)$$ where the $2\times4$ submatrix $N_2$ is
calculated from equation (1.54)
$$ N_2[B,AB] = [O,I_2]  $$
Since
$$ [B,AB]^{-1} \: = \: \left[ \begin{array}{rrrr} 1 & -1 & -1 & 0 \\ 0 & 0 & 0 & 1 \\
0 & 0 & 1 & 0 \\ 0 & 1 & -1 & 0 \end{array} \right ]
$$
then we can find
$$ N_2 = [O,I_2][B,AB]^{-1} = \left[ \begin{array}{rrrr} 0 & 0 & 1 & 0 \\
 0 & 1 & -1 &  0  \end{array} \right ] $$
and using (1.85) calculate
$$ N = \left[ \begin{array}{rrrr}   0 &  0 &  1 &  0 \\   0 &  1 & -1 &  0 \\
 1 &  1 &  0 &  0 \\ 0 & -1 &  1 &  1  \end{array} \right ]
\eqno(1.86) $$
As
$$ N^{-1} = \left[ \begin{array}{rrrr} -1 & -1 &  1 &  0 \\ 1 &  1 &  0 &  0 \\
 1 &  0 &  0 &  0 \\ 0 &  1 &  0 &  1   \end{array} \right ]
\eqno(1.87) $$ then we  find
$$ A^* = NAN^{-1} =  \left[ \begin{array}{rrrr} 0 &  0 &  1 &  0 \\
 0 &  0 &  0  & 1 \\ -1 &  0 &  3 &  2 \\  1 &  0 &  0 & -1
\end{array} \right ] , \qquad B^* = NB = \left[ \begin{array}{rr} 0 & 0 \\
0 & 0 \\ 1 & 0 \\ 0 & 1 \end{array} \right ]
\eqno(1.88) $$
Matrix $A^*$ in (1.88) is  the  block  companion  matrix  for  the
matrix polynomial  $\Phi(s) = I_2s^2 + T_1s + T_2$ with
$$ T_1 =\left[ \begin{array}{rr}  -3 & -2 \\ 1 &  0 \end{array} \right ],
\qquad  T_2 = \left[ \begin{array}{rr} 1 & 0 \\ -1 &   0 \end{array} \right ].
$$

{\it EXAMPLE 1.3.} \rm

 Let us find   Yokoyama's  canonical  form  for  the controllable
system with  $n = 4,\;\; r = 2$ [S3]
$$ \dot{x} = \left[ \begin{array}{cccc} 2 & 1 & 0 & 0 \\ 0 & 1 & 0 & 1 \\
 0 & 2 & 0 & 0 \\ 1 & 1 & 0 & 0 \end{array} \right ]x +
\left[ \begin{array}{cc} 1 & 0 \\ 0 & 0 \\ 0 & 0 \\ 0 & 1
\end{array} \right ]u
\eqno (1.89) $$ At first we build the controllability matrix
$Y_{AB} = [B,AB,A^2B,A^3B] $. As $rank[B,AB,A^2B] = 4$  then  $\nu
= 3$. Using  formulas  (1.59) we calculate $l_1 = l_2 = 1, l_3
=2$. Since
$$ rank [B,AB,A^2B] \: = \: rank [b_1,b_2,Ab_2,A^2b_2] $$
where $b_1$, $b_2$ are columns of the matrix $B$ then $M$ is the
unity matrix
$$ M =  \left[ \begin{array}{cc}  1 & 0 \\ 0 &  1 \end{array} \right ]
\eqno(1.90)$$ Let us find the matrix $V = [V_1, V_2, V_3]$. Using (1.73) we find
$$ V_1 = [b_1,b_2] =  \left[ \begin{array}{cc}  1 & 0 \\ 0 & 0 \\
0 & 0 \\ 0 & 1  \end{array} \right ], \: V_2 = AV_1 \left[ \begin{array}{c}
0 \\ 1  \end{array} \right ]  = \left[ \begin{array}{c}  0 \\ 1 \\ 0  \\ 0
\end{array} \right ], \: V_3 = AV_2 = \left[ \begin{array}{c} 1 \\ 1 \\ 2 \\ 1
\end{array} \right ]
$$
and by  formulas (1.67) discover the structure of the matrix $N$
$$  N \: = \: \left[ \begin{array}{c} N_3 \\  \dotfill \\  N_2 \\ \dotfill \\
N_1 \end{array}  \right ] \: = \: \left[ \begin{array}{c} P_3 \\  \dotfill \\
 P_3A \\  \dotfill \\  P_1 \\ P_3A^2 \end{array} \right ]
\eqno(1.91) $$ where  $N_3 = P_3$ is an $1\times 4$ submatrix ($l_1 = 1$). Here a submatrix $P_2$ does't exist
because $l_2 - l_1 = 0$ and $N_2 = P_3A$. Submatrices $P_1$ and $P_3$ are satisfied the following equation
$$ \left[ \begin{array}{c}  P_3 \\  P_1 \end{array} \right]
[ V_1, V_2, V_3 ]  \: = \: \left[ \begin{array}{lcclcl} 0 & 0 & \vdots & 0 &
\vdots & 1 \\  1 & 0 & \vdots & x_1 & \vdots & x_2 \end{array} \right]
\begin{array}{cl} \} & l_1=1 \\  \} & l_3-l_2 = 1\end{array}
\eqno(1.92) $$ that follows from formula (1.75) for the concrete
$l_1 =1$, $l_2 = 1$, $l_3 = 2$. In (1.92) $x_1$ and $x_2$  are any
numbers. Assigning  $x_1 = 0$, $x_2 = 1$  we calculate from (1.92)
$$ P_3 =  [ 0 \: \:  0 \: \: 0.5 \: \: 0 ]  $$
$$ P_1 =  [ 1 \: \: 0 \: \:  0  \: \:  0]   $$
and find $N$  from (1.91)
$$ N = \left[ \begin{array}{cccc}   0 &  0 &  0.5 &  0 \\   0 &  1 & 0 &  0 \\
 1 &  0 &  0 &  0 \\ 0 & 1 &  0 &  1  \end{array} \right ]
\eqno(1.93) $$
 Thus
$$ F = NAN^{-1} =  \left[ \begin{array}{rrrr} 0 &  1 &  0 &  0 \\
 0 &  0 &  0  & 1 \\ 0 &  1 &  2 &  0 \\  0 &  1 &  1 & 1
\end{array} \right ] , \qquad G = NB = \left[ \begin{array}{rr} 0 & 0 \\
0 & 0 \\ 1 & 0 \\ 0 & 1 \end{array} \right ]
 \eqno(1.94) $$
 It is
evident that the structure of matrices $F$ and $G$ corresponds
formulas  (1.62)-(1.65).  Indeed, in (1.62)
$$
  E_{1,2} = 1, \: E_{2,3} = [0 \: \: 1], \: F_{31}  = \left[ \begin{array}{c} 0 \\
0 \end{array} \right ],\: F_{32}  = \left[ \begin{array}{c} 1 \\ 1
\end{array} \right ],\: F_{33}  = \left[ \begin{array}{cc} 2 & 0 \\ 1 & 1
\end{array} \right ]
$$
and in (1.64), (1.65) $I^1 = 1$, $I^2$  does not  exist,  $I^3 = 1$ .  The matrix $F$ in (1.94) is the general
block companion  matrix  for  the matrix polynomial
$$   \Phi(s) = I_2s^3 + T_1s^2 + T_2s + T_3   $$
with
$$ T_1 = -F_{33} = \left[ \begin{array}{rr} -2 & 0 \\ -1 & -1 \end{array}
 \right ], \: T_2 = -[O, \: F_{32}] = \left[ \begin{array}{rr} 0 & -1 \\
0 & -1 \end{array} \right ],\:T_3 = -[O, \: F_{31}] = \left[ \begin{array}{rr}
0 & 0 \\ 0 & 0 \end{array} \right ]
$$

For testing we find $ det \Phi(s) \:=\: det \left[
\begin{array}{cc} s^3 - 2s^2 & -s \\ -s^2 & s^3 -s^2 -s
\end{array} \right ] \:= \: s^3(s^3 - 3s^2 + s + 1)$ \qquad and
\qquad $det (sI_4 - F) = s^4 - 3s^3 + s^2 + s$. It is evident that $s^{-2}det\Phi(s) \:= \:det(sI_4 - F)$. The
last equality corresponds to  Assertion 1.2 (formula (1.26)).

\chapter[System description by transfer function matrix]{System
 description by transfer function matrix}

\section[The Laplace transform]{The Laplace transform}

     Let's consider a scalar function $f(t)$ of a real variable  $t$
such that the function $f(t)e^{-st}$ where $s$ is a complex variable has a convergent integral
$$
   \bar{f}(s) = \int_0^\infty f(t)e^{-st}dt
\eqno(2.1)$$ This integral is known as a direct one-sided Laplace
transform of a time-dependent function or a Laplace integral. It
is calculate, by definition, as follows
$$
  \int_0^\infty f(t)e^{-st}dt = \begin{array}{c} \\ lim \\ \scriptstyle
T \rightarrow \infty ,\epsilon \rightarrow  0 \end{array} \displaystyle
 \int_\epsilon^T f(t)e^{-st}dt $$

If a limit exists  then  the  Laplace  integral  is  a  convergent integral. These questions are studied detail
in any textbooks, for example in [B3].

     A function $\bar{f}(s)$ of a  complex  variable $s$ is  called
the  Laplace transform of $f(t)$  and  denoted as $L[f(t)] =
\bar{f}(s)$. Let us write the  main properties of the Laplace
transform which will be useful in  the present study.
     Let $f(t)$, $f_i(t)$, $i=1,2$ are scalar functions of time and $a$,
 $b$ are constant variables. We have

1. $$  L [af_1(t) + bf_2(t)] = a\bar{f}_1(s) +  b\bar{f}_2(s) $$
     where $ L[f_1(t)] = \bar{f}_1(s)$ , $ L[f_2(t)] = \bar{f}_2(s)$

2. $$  L[f^{(i)}(t)] = s^i\bar{f}(s) - f(+0)s^{i-1} -\cdots- f^{(i-1)}(+0),
       \:  f(+0) =  \begin{array}{c} \\ lim \\ \scriptstyle t \rightarrow +0
         \end{array} \displaystyle f(t) $$

3. $$  L[ \int f(t)dt] = \bar{f}(s)/s + (\int f(t)dt/_{t=+0})/s $$

4. $$  L [f(t-a)] = e^{-as}\bar{f}(s), \;\; {\rm for} \;\; a > 0 $$

5. $$  L[\vec{f}(t)] = \vec{\bar{f}(s)}  \eqno(2.2)   $$
where $\vec{f}(t)$, $\vec{\bar{f}(s)} $  are $n$-vectors.

Thus, the  Laplace transform makes possible to replace a differential equation  in $f(t)$   by an algebraic
equation in $\bar{f}(s)$. Solving the algebraic equation  we can find $\bar{f}(s)$.  For obtaining $f(t)$ we
should   use the inverse Laplace transform ( $L^{-1}$ - transform). For more details   see, for example, [B3].

\section[Transformation from state-space to frequency domain
 representation. Transfer function matrix]{Transformation
 from state-space to \\frequency domain representation. \\
 Transfer function matrix}

     We study equation  (1.1) with  $t_o = 0$.   Taking the Laplace
transforms of both sides of (1.1), (1.2) and using properties
(2.2) gives
$$ s\bar{x}(s) - x(0) = A\bar{x}(s) + B\bar{u}(s)  \eqno(2.3) $$

$$ \bar{y}(s) = C \bar{x}(s)  \eqno(2.4)$$
where vectors $\bar{x}(s)$, $\bar{u}(s)$, $\bar{y}(s)$  are the
Laplace transforms of the vectors $x(t)$, $u(t)$, $y(t)$
respectively. Assuming $ s \ne \lambda_i$ ( $\lambda_i$    are
eigenvalues of $A$) we express $\bar{x}(s)$ in equation (2.3) as
follows
$$
 \bar{x}(s) = (sI_n - A)^{-1}\{ B\bar{u}(s) + x(0) \}$$
The last relation is true for all $s\neq \lambda_i$,
$i=1,2,\ldots,n$. Substituting $\bar{x}(s)$ in (2.4) we  get  the
expression for $ \bar{y}(s)$
$$
\bar{y}(s) = C(sI_n - A)^{-1} B\bar{u}(s) + C(sI_n - A)^{-1}x(0) \eqno(2.5)$$ where the first term depends on
the  input vector and the second one  depends on the initial state vector. Taking the inverse Laplace transform
of (2.5) we get formula (1.7) where the inverse Laplace transform of the second term is $Ce^{At}x(0)$ and the
first one is $C\int_{t_o}^{t}{e^{A(t-\tau )}Bu(\tau)d\tau}$.

     When $x(0)$ is equal to zero then
$$
\bar{y}(s) = C(sI_n - A)^{-1} B\bar{u}(s)
\eqno(2.6) $$
The matrix
$$
    G(s) = C(sI_n - A)^{-1} B  \eqno(2.7)  $$
is called as a transfer function matrix.
     Similarly, for  system (1.3) we can get
$$
    G(s) = C(sI_n - A)^{-1} B  + D  \eqno(2.8)  $$
     The elements $ g_{ij}$  , $i=1,\ldots l;\;\; j=1,\ldots,r$ of $G(s)$
are  rational functions of $s$. Each element $g_{ij}$ is the transfer function from  $j$-th  component  of  the
output  to $i$ -th component of the input. As $(sI - A)^{-1} = adj(sI - A)/ det(sI - A)$ then a numerator degree
of each elements $g_{ij}$ of $G(s)$  (2.7)  is strictly less  than  a denominator degree. Then the following
condition is true
$$  \begin{array}{c} \\ lim \\ \scriptstyle t \rightarrow \infty
      \end{array} \displaystyle G(s) = O $$
Such  $G(s)$ is  known  as  a  strictly  proper  transfer  function matrix\footnote{ Further the abbreviation
TFM will be used}.
     If  $D \ne O$ (2.8) then $\begin{array}{c} \\ lim \\ \scriptstyle t
 \rightarrow \infty \end{array} \displaystyle G(s) = D \ne O $.
This TFM is known as a proper transfer function matrix. It possess several  (or  single)  elements having equal
degrees of  a numerator and a denominator.

     Let us consider an element $ g_{ij}$  of the strictly proper TFM (2.7).

     \it{DEFINITION 2.1.} \rm A complex $s_i$ is called a pole of  $G(s)$
if several (or single) elements of $G(s_i)$  are equal to  $ \infty$.

     Zeros of a  least common  denominator of $ g_{ij}$ form  a
subset of the complete set of the TFM poles. The complete  set  of
poles coincides with zeros of a polynomial being the least common
denominator of all nonzero minors of all orders of $G(s)$ [M2].

     For example, let's determine poles of the following TFM
$$
     G(s) =  \left[ \begin{array}{cc} \frac{s-1}{(s+2)(s+3)}& 0 \\
0 & \frac{1}{s+2} \end{array} \right]     $$
Zeros of the least
common denominator of $ g_{ij}$ $i=1,2;j=1,2$ are $s_1 = -2, s_2 =
-3$. They form the  subset  of  the complete set of poles: $s_1 =
-2, s_2 = - 2, s_3 = - 3$.

     Now we consider a definition of  system poles.

     \it{DEFINITION 2.2.} \rm  A complex $s$ which is a some zero  of the
polynomial $det(sI_n -A)$ is called as a system pole.

     The complete set of system poles coincides with eigenvalues of the matrix  $A$.

     If  all  elements $ g_{ij}$ of $G(s)$ have  relatively   prime
numerators and denominators  then  the set of TFM poles coincides
with the set of  system poles.

\section[Physical interpretation of transfer function matrix]{Physical
 interpretation of transfer function matrix}

\subsection[Impulse response matrix]{Impulse response matrix}

     Let   system  (1.1),(1.2)  has  been  completely  at  rest
$(x(0) = 0)$  when a delta-function impulse  $\delta(t)u_o$ is
applied  where
$$
   \delta(t) = \left \{ \begin{array}{ccc} 0 & , & t\ne 0 \\ \infty & , & t = 0
\end{array} \right. , \qquad   \int_{-\epsilon}^{+\epsilon} \delta(t)dt = 1,
\qquad \epsilon > 0 $$ and $u_o$ is a constant  $r$-vector having only unit  element  and zeros otherwise. Since
$$
  L[\delta(t)] = \int_{0}^{\infty} \delta(t)e^{-st}dt = 1   $$
then the Laplace transform of the output with  $x(0)=0$ is equal
to
$$
\bar{y}(s) = C(sI_n - A)^{-1} Bu_o
\eqno(2.9) $$
 Let us find
$y(t)$. For this purpose we consider  general solution  (1.7) of differential equations (1.1),(1.2) with
$x(t_o)|_{t_o =o}= 0$
$$
  y(t) =  C\int_{0}^{t}{e^{A(t-\tau )}Bu(\tau)d\tau}  $$
Setting  $u(\tau) = \delta(\tau)u_o $ in the last equation we
obtain
$$
\begin{array}{c} y(t) =  C\int_{0}^{t}{e^{A(t-\tau )}B\delta(\tau)u_od\tau}
=
  C(\int_{o}^{t}{e^{A(t-\tau )}\delta(\tau)d\tau)}Bu_o = \\ \\
  C\int_{0}^{t}{e^{At}\delta(t-\tau)d\tau}Bu_o = Ce^{At}Bu_o
  \end{array}$$
The matrix
$$ G(t) = Ce^{At}B  \eqno(2.10)     $$
is called as the impulse response matrix of a system [W1].

Using the Laplace transform of $e^{At} $: $L[e^{At}] = (sI_n - A)^{-1}$ we find
$$  L[G(t)] = C(sI_n - A)^{-1}B = G(s) $$
So, TFM $G(s)$  is the Laplace transform  of  a  impulse response
matrix.

\subsection[Frequency response matrix]{Frequency response matrix}

     It is known that  exponential functions $e^{st}$ with a
complex parameter $s$   describe  oscillatory  signals  of  all
frequencies with a constant or exponential amplitude. Indeed, if
  $s = j\omega$  is an imaginary variable then $e^{j\omega t} =
cos\omega t +jsin\omega t$   and   $cos\omega t = 0.5(e^{j\omega t} + e^{-j\omega t})$. So, we  have an
oscillatory  function with the frequency $\omega $ . If $s$ is a complex variable: $ s = \bar{s} + j\omega $ (
$\bar{s}$ - real variable ) then $e^{st} = e^{\bar{s}t} (cos\omega t +jsin\omega )$. Thus we have an oscillatory
signal with the exponential increasing or decreasing amplitude and the frequency $\omega $. Applying an
exponential input signal we can reveal a relationship between TFM  and  the transient response of a system.

     Suppose we use an exponential input having the following
form
$$
   u(t) = \left \{ \begin{array}{ccc} 0 & , & t\le 0 \\ u_o e^{st} & ,
& t > 0 \end{array} \right.
\eqno(2.11) $$
The function (2.11) can
be rewritten as follows
$$     u(t) = u_o e^{st}1(t)    \eqno(2.12) $$
where $s$  is a complex variable, $u_o$  is a constant  $r$ vector, $1(t)$ denotes a unit  step function of time
$$
   1(t) = \left \{ \begin{array}{ccc} 0 & , & t\le 0 \\ 1 & , & t > 0
\end{array}   \right. $$
Let us write a general solution of (1.1), (1.2) with the input (2.12). Assuming that $s$  does not coincide with
any eigenvalue of $A$   we obtain
$$ y(t) =  Ce^{At}x(0) + C\int_{0}^{t}{e^{A(t-\tau )}Bu_o e^{st}d\tau}
    =  Ce^{At}x(0) + Ce^{At}\int_{0}^{t}{e^{-A\tau}Bu_o
e^{s\tau}d\tau} $$
 We consider the second item. Since
$$ Bu_oe^{s\tau} = Bu_oI_n e^{s\tau} =  I_n e^{s\tau}Bu_o
= e^{I_ns\tau}Bu_o $$ then we can write
$$ y(t) =  Ce^{At}x(0) +Ce^{At}\int_{0}^{t}{e^{(sI_n-A)\tau}Bu_od\tau}
\eqno(2.13)  $$
The vector $Bu_o$ does not depend in $\tau$, therefore, it should be taken out from the integral
$$ y(t) =  Ce^{At}x(0) +Ce^{At} \{ \int_{0}^{t}{e^{(sI_n-A)\tau}d\tau}
\} Bu_o $$
 Integrating we have
$$ \int_{0}^{t}{e^{(sI_n-A)\tau}d\tau} = \{
e^{(sI_n-A)\tau}\arrowvert_0^t \} (sI_n - A)^{-1} = (e^{(sI_n - A)t} - I_n)(sI_n - A)^{-1}    $$ Substituting
the right-hand side of the last relation in (2.13) we obtain
$$
y(t) =
 Ce^{At}x(0) +Ce^{At}(e^{(sI_n - A)t} - I_n)(sI_n - A)^{-1}Bu_o =
$$
$$
Ce^{At}x(0) + Ce^{st}(sI_n - A)^{-1}Bu_o - Ce^{At}(sI_n - A)^{-1}Bu_o =
$$
$$
Ce^{At} \{ x(0) - (sI_n - A)^{-1}Bu_o \} + C(sI_n - A)^{-1}Bu_o e^{st}
$$
In this expression the second term is equal to $G(s)u(t)$. The first term is determined due  system response in
the initial time. If  the system is asymptotically stable ( $Re \lambda_i(A) < 0$ )  and $Re s > Re \lambda_i$,
$i=1,\ldots,n$ then we have for large values of  $t$
$$
     y(t) \cong G(s)u(t),\qquad t \gg 0   \eqno(2.14)  $$
Thus, the  transfer  function  matrix  $G(s)$ describes asymptotic behavior of a system in response to
exponential inputs of the complex frequency $s$.

     Let us consider  the   oscillatory  input
$$ u(t) = u_o e^{j\omega t}1(t), \qquad t \geq  0
\eqno(2.15) $$ where  a real value $\omega$ is the frequency of the oscillation. Substituting (2.15) in (2.14)
yields
$$  y(t) \cong G(j\omega)u(t)   $$
where $ G(j\omega ) = C(j\omega I_n - A)^{-1}B$ is called as the
frequency response matrix.

     So, it has been shown that  TFM with  $s = j\omega$  coincides
with the frequency response matrix.

\section[Properties of transfer function matrix]{Properties
  of transfer function matrix}

\subsection[Transformation of state, input and output vectors]
 {Transformation of state, input and output vectors}

     Let's carry  out  a  nonsingular  transformation  of  the state
vector $x$
$$   \hat{x} =Nx    \eqno(2.16) $$
where $\hat{x}$  is a new state  vector, $N$ is  a transformation
$n\times n$ matrix.
     Expressing $x = N^{-1}\hat{x}$  and  substituting  into  (1.1),
(1.2) we obtain a new system
$$  \dot{\hat{x}} = \hat{A}\hat{x} +\hat{B}u $$
$$  y = \hat{C}\hat{x}
   \eqno (2.17)  $$
where  $\hat{A} = NAN^{-1}$, $\hat{B} = NB$, $\hat{C} =CN^{-1}$.
     System (2.17) has the transfer function matrix
$$  \hat{G}(s) = \hat{C}(sI - \hat{A})^{-1}\hat{B}
 \eqno(2.18)  $$
Substituting  $\hat{A} = NAN^{-1}$, $\hat{B} = NB$, $\hat{C} =CN^{-1}$
in (2.18) yields
$$  \hat{G}(s) = CN^{-1}(sI_n - NAN^{-1})NB = C(sI_n - A)^{-1}B = G(s) $$

 \it{PROPERTY 2.1.} \rm A transfer function matrix is  invariant  under
a nonsingular state transformation.  The  physical  meaning  of
input and output vectors is preserved.

     Let's carry out a  nonsingular transformation  of  the  input  and
output vectors
$$ \hat{u} = Mu, \qquad \hat{y} = Ty   \eqno(2.19)   $$
where $M$ and $T$ are nonsingular matrices  of  dimensions $r\times r$ and $l\times l$ respectively.
Substituting (2.19) into (1.1), (1.2) yields
$$  \dot{x} = Ax + \hat{B}\hat{u} $$
$$  \hat{y} = \hat{C}x $$
where $\hat{B} = BM^{-1}$, $\hat{C} = TC$. The transfer function matrix of this system is
$$ \hat{G}(s) = TC(sI_n - A)^{-1}BM^{-1} = TG(s)M^{-1}
\eqno(2.20) $$
The following property follows from (2.20)

\it{PROPERTY 2.2.} \rm  A transfer function  matrix  does  not invariant under nonsingular input and output
transformations. The physical meaning of  input and output vectors is not preserved.

     Let $M$ and $T$ are  permutation matrices. Multiplying  by permutation matrices permutes   columns and rows
of $G(s)$ and, in fact, changes a numeration of input and output variables.

\it{PROPERTY 2.3.} \rm  A transfer function matrix is  invariant under the  permutation  transformation of input
and/or   output. This transformation rearranges  columns and/or rows.  The physical meaning  of  input and
output vectors is preserved.

     From  Property  2.1  follows that TFM  describes  only  the
external  (input-output) behavior of  a system and does not depend on a  choice of the state vector. In what
follows we define the   relationship   between     TFM   and    controllability/observability characteristics of
a system.

\subsection[Incomplete controllable and/or observable system]
 {Incomplete controllable and/or observable system}

     Let uncontrollable and observable system (1.1), (1.2)
has the controllability matrix $Y_{AB}$ with
$$  rank Y_{AB} = rank[B,AB,\ldots,A^{n-1}] = m < n
\eqno (2.21) $$
 A subspace \bf{N} \rm is the controllability one
if every state $x \in \bf{N} $  \rm can be reached from  the
initial  state  along  a controllable state  trajectory  during a
finite interval of  the time. The subspace \bf{N} \rm   has the
dimension coinciding  with $rank Y_{AB}$. For case (2.21) the
dimension of \bf{N} \rm is equal to $m$.

     Let vectors  $e_1,e_2,\ldots, e_m$  are the basis of the
controllability  subspace \bf{N}. \rm  We define   $n-m$ linearly independent   vectors $e_{m+1},e_{m+2},\ldots,
e_n$ which form the orthogonal complement of the controllability subspace basis. All vectors  $e_1,e_2,\ldots,
e_n$  form the basis of the state-space. We consider the nonsingular $n\times n$ matrix
$$     N = [N_1, N_2] $$
with  $n\times m$  and  $n\times (n-m)$  submatrices
$$  N_1 = [e_1,e_2,\ldots, e_m],  \qquad
N_2 = [e_{m+1},e_{m+2},\ldots, e_n]  $$
 and introduce the transform state vector $\hat{x}$
$$   \hat{x} =N^{-1}x     $$
Since
$$ x = N\hat{x} = [N_1, N_2] \left [ \begin{array}{c} \hat{x}_1 \\
\hat{x}_2 \end{array} \right ] \;=\; N_1\hat{x}_1 + N_2\hat{x}_2
     \eqno (2.22)  $$
then  substituting (2.22) in (1.1) yields [K5]
$$  \dot{\hat{x}} = \hat{A}\hat{x} +\hat{B}u
\eqno (2.23)  $$
where
$$ \hat{A} =  \left [ \begin{array}{cc} \hat{A}_{11} & \hat{A}_{12} \\
O & \hat{A}_{22} \end{array} \right ], \qquad  \hat{B} = \left [
\begin{array}{c}  \hat{B}_1 \\ O \end{array} \right ], \qquad
\hat{x} = \left [ \begin{array}{c} \hat{x}_1 \\\hat{x}_2
\end{array} \right ]   $$

Let us rewrite  (2.23) as  two equations
$$ \dot{\hat{x}_1} = \hat{A}_{11}\hat{x}_1 + \hat{A}_{12}\hat{x}_2 + \hat{B}_1u$$
$$ \dot{\hat{x}_2} = \hat{A}_{22}\hat{x}_2
\eqno (2.24)  $$
 In (2.24) the first subsystem with  the $m\times m$ dynamics matrix  $\hat{A}_{11}$ is completely
controllable. This follows from analysis  of  the controllability matrix $Y_{\hat{A}\hat{B}}$. Indeed, since
$rank N = n$ then the following rank equalities  take place
$$ rank[\hat{B}_1,\hat{A}_{11}\hat{B}_1, \ldots, \hat{A}_{11}^{m-1}\hat{B}_1] =
\left [ \begin{array}{cccc} \hat{B}_1 & \hat{A}_{11}\hat{B}_1 & \cdots &
\hat{A}_{11}^{n-1}\hat{B}_1 \\ O & O & \cdots & O \end{array} \right ] =
$$
$$
=\;rank[\hat{B},\hat{A}\hat{B}, \ldots, \hat{A}^{n-1}\hat{B}] =
rank \{ N^{-1}[B,AB,\ldots,A^{n-1}B] \} = rank[B,AB,\ldots,A^{n-1}B] = m $$

     Eigenvalues of the matrix $\hat{A}_{11}$ are refereed to  controllable
 poles of system (1.1),(1.2). Eigenvalues of the
matrix $\hat{A}_{22}$  are called as uncontrollable poles (input
decoupling poles)  of  system (1.1), (1.2).

\it {ASSERTION 2.1.}  \rm  Eigenvalues $\lambda_i^*$  of the
matrix $A$
 for which the equality $rank[\lambda_i^*I_n - A, B] < n$ is satisfied
coincide with   eigenvalues of $\hat{A}_{22}$ ( uncontrollable
poles of (1.1), (1.2)).

\it{PROOF}. \rm It follows from the following rank equalities
$$ rank[\lambda_i^*I_n - A, B] = rank \left [ \begin{array}{ccc} \lambda_i^*I_m -
\hat{A}_{11} & -\hat{A}_{12} & \hat{B}_1 \\  O & \lambda_i^*I_{n-m} -
 \hat{A}_{22} & O \end{array} \right ] = $$
$$ =\;rank \left [ \begin{array}{ccc} \lambda_i^*I_m -  \hat{A}_{11} & \hat{B}_1 &
-\hat{A}_{12} \\ O & O & \lambda_i^*I_{n-m} - \hat{A}_{22} \end{array} \right ]
=  rank[ \lambda_i^*I_m -  \hat{A}_{11} , \hat{B}_1] +  rank[ \lambda_i^*I_{n-m}
 -  \hat{A}_{22}] $$
Since the pair of matrices $(\hat{A}_{11},\hat{B}_1)$  is completely controllable then $ rank[ \lambda_i^*I_m -
\hat{A}_{11} , \hat{B}_1] = m$. Hence,  the  rank of the  matrix $[\lambda_i^*I_n - A, B]$  is reduced if and
only if $\lambda_i^*$ are eigenvalues of $\hat{A}_{22}$.

\it{ASSERTION 2.2.} \rm   The number  of uncontrollable poles  is
equal to the rank deficient of controllability matrix $Y_{AB}$
(2.21).

\it{PROOF} \rm. Let consider the controllability matrix of system (2.23)
 $Y_{\hat{A}\hat{B}} = \\
rank[\hat{B},\hat{A}\hat{B},
 \ldots, \hat{A}^{n-1}\hat{B}] $.
Using the structure of matrices $\hat{A}$ and  $\hat{B}$ we obtain
$$  Y_{\hat{A}\hat{B}} \;=\; \left [ \begin{array}{cccccc} \hat{B}_1 &
\hat{A}_{11}\hat{B}_1 & \cdots & \hat{A}_{11}^{m-1}\hat{B}_1 & \cdots &
\hat{A}_{11}^{n-1}\hat{B}_1 \\ O & O & \cdots & O & \cdots & O \end{array}
 \right ]  \begin{array}{cl} \} & m \\ \} & n-m \end{array}
$$
Since  the  pair $(\hat{A}_{11},\hat{B}_1)$  is completely
controllable then
$$rank[\hat{B}_1,\hat{A}_{11}\hat{B}_1, \ldots, \hat{A}_{11}^{n-1}\hat{B}_1]
= rank[\hat{B}_1,\hat{A}_{11}\hat{B}_1, \ldots,
\hat{A}_{11}^{m-1}\hat{B}_1] = m $$ Hence $rank Y_{\hat{A}\hat{B}}
= m < n$ and the rank deficient of $Y_{\hat{A}\hat{B}}$  coincides
with the number of uncontrollable poles. The last one is equal to
$n - m$  . Then from relations
$$ Y_{\hat{A}\hat{B}} = N^{-1} Y_{AB},\qquad rankN = n $$
we have
$$  rankY_{\hat{A}\hat{B}} = rank(N^{-1} Y_{AB}) = rank Y_{AB} $$

The last equality completes the proof.

Let's find the output vector of the transformed system (2.23). As $x = N \hat{x}$   then  $ y = Cx =
CN^{-1}\hat{x}$. We denote  $\hat{C} = CN$ and  rewrite the vector $y$ as follows
$$  y =  \hat{C}\hat{x} = CN\hat{x} = CN_1\hat{x}_1 + CN_2\hat{x}_2 =
\hat{C}_1\hat{x}_1 + \hat{C}_2\hat{x}_2   \eqno (2.25)$$

The following subsystems $S_1$ and $S_2$
$$ S_1 : \qquad \dot{\hat{x}_1} = \hat{A}_{11}\hat{x}_1 + \hat{A}_{12}\hat{x}_2
+ \hat{B}_1u , \qquad y_1 = \hat{C}_1\hat{x}_1 $$
$$  S_2 : \qquad \qquad \qquad \dot{\hat{x}_2} = \hat{A}_{22}\hat{x}_2,
    \qquad  y_2 = \hat{C}_2\hat{x}_2 $$
have properties: subsystem $S_1$ is completely controllable and
observable, subsystem $S_2$ is uncontrollable.
%

Now we find the transfer function matrix of  system (2.23), (2.25)
$$  \hat{G}(s) = \hat{C}(sI - \hat{A})^{-1}\hat{B} =
CN \left [ \begin{array}{cc} sI -\hat{A}_{11} & - \hat{A}_{12} \\
O & sI -\hat{A}_{22}  \end{array} \right ] \left [ \begin{array}{c} \hat{B}_1\\
 O \end{array} \right ] $$
 Since for $detA \ne 0, \; detC \ne 0$
$$  \left [ \begin{array}{cc}  A & -B \\ O & C \end{array} \right ]^{-1} =
\left [ \begin{array}{cc}  A^{-1} & A^{-1}BC^{-1} \\ O & C^{-1} \end{array}
 \right ] \qquad \qquad    $$
then
$$  \hat{G}(s) = CN \left [ \begin{array}{cc} (sI_m -\hat{A}_{11})^{-1} &
 (sI_m -\hat{A}_{11})^{-1} \hat{A}_{12}(sI_{n-m} -\hat{A}_{22})^{-1} \\
O & (sI_{n-m} -\hat{A}_{22})^{-1}  \end{array} \right ] \left [ \begin{array}{c}
 \hat{B}_1\\ O \end{array} \right ] = $$
$$ =\;[\hat{C}_1,\hat{C}_2] \left [ \begin{array}{c}
(sI_m -\hat{A}_{11})^{-1}\hat{B}_1 \\ O \end{array} \right ] =
\hat{C}_1 (sI_m -\hat{A}_{11})^{-1}\hat{B}_1  $$

So,  transfer function matrix $\hat{G}(s)$  coincides with one for the completely controllable subsystem $S_1$.
Since from Property 2.1 $\hat{G}(s) = G(s)$ then we have the following assertion.

\it{ASSERTION 2.3.} \rm   TFM  of  an  uncontrollable  and completely observable  system  coincides  with  the
TFM of the completely controllable and observable subsystem. Poles of TFM  coincide with controllable poles.

     The  analogous  result  can  be  obtained  for  a  completely
controllable and unobservable system.  Such system is transformed
into the following canonical form [K5]
$$  \dot{\hat{x}} = \left [ \begin{array}{cc} \hat{A}_{11} & O \\
\hat{A}_{21} & \hat{A}_{22} \end{array} \right ] \hat{x} + \left [
\begin{array}{c}  \hat{B}_1 \\ \hat{B}_2 \end{array} \right ] u, $$
$$ y = \left[ \begin{array}{cc}  \hat{C}_1 & O \end{array} \right ] \hat{x}
\eqno (2.26)  $$
System (2.26) is decomposed into two subsystems
$$ S_1 : \qquad \dot{\hat{x}_1} = \hat{A}_{11}\hat{x}_1 +  \hat{B}_1u ,
\qquad y_1 = \hat{C}_1\hat{x}_1 $$
$$  S_2 : \qquad \dot{\hat{x}_2} = \hat{A}_{21}\hat{x}_2 +
\hat{A}_{22}\hat{x}_2 +  \hat{B}_2u $$ with  completely controllable and observable subsystem $S_1$ and
unobservable subsystem $S_2$.
     Eigenvalues of the matrix $\hat{A}_{11}$ are named as  observable
poles.  Eigenvalues  of   matrix  $\hat{A}_{22}$ are named as
unobservable poles (output decoupling poles).
It can be shown that  TFM of controllable  and  unobservable system (1.1), (1.2)  is
$$ G(s) = \hat{C}_1 (sI_m -\hat{A}_{11})^{-1}\hat{B}_1     $$
So, we have the following result.

\it{ ASSERTION  2.4.} \rm   TFM  of  an  unobservable  and completely controllable system coincides  with TFM
of the completely controllable and observable subsystem. TFM poles coincide with observable poles.

     Similar to Assertions 2.1, 2.2 we can obtain

\it{ASSERTION  2.5.} \rm   Eigenvalues $\lambda_i^*$ of the matrix $A$ for which the equality
$rank[\lambda_i^*I_n - A^T, C^T] < n$ is satisfied coincide with   eigenvalues of $\hat{A}_{22}$ (unobservable
poles of (1.1),(1.2)).

\it{ASSERTION  2.6.} \rm  A number of unobservable poles  is  equal to a rank deficient of the unobservability
matrix.

Let's consider  uncontrollable  and  unobservable  system  (1.1), (1.2). Using a nonsingular transformation of
the  state vector  this system can be reduce to the following block form [K1]
$$ \left [ \begin{array}{c} \dot{\hat{x}_1} \\ \dot{\hat{x}_2} \\
\dot{\hat{x}_3} \\ \dot{\hat{x}_4} \end{array} \right ] =
\left [ \begin{array}{cccc} \hat{A}_{11} &  \hat{A}_{12} & \hat{A}_{13} &
\hat{A}_{14}  \\ O & \hat{A}_{22} &  \hat{A}_{23} & \hat{A}_{24} \\
O & O &   \hat{A}_{33} & \hat{A}_{34} \\  O & O & O & \hat{A}_{44}  \end{array}
\right ] \left [ \begin{array}{c} \hat{x}_1 \\ \hat{x}_2 \\ \hat{x}_3 \\
\hat{x}_4 \end{array} \right ]  + \left [  \begin{array}{c}  \hat{B}_1 \\
 \hat{B}_2  \\ O  \\ O  \end{array} \right ] u, $$
$$ y = \left[ \begin{array}{cccc}  O & \hat{C}_2  & O & \hat{C}_4  \end{array}
 \right ] \hat{x} , \qquad \hat{x}^T = [ \hat{x}^T_1, \hat{x}^T_2,\hat{x}^T_3,
\hat{x}^T_4] \eqno (2.27)  $$
We may rewrite (2.27) as four
connected subsystems
$$ \begin{array}{cccccl}  S_1 &  : &   & \dot{\hat{x}_1} &  = &
\hat{A}_{11}\hat{x}_1 + \hat{A}_{12}\hat{x}_2  + \hat{A}_{13}\hat{x}_3 +
\hat{A}_{14}\hat{x}_4 + \hat{B}_1u \\
 S_2  & : &  & \dot{\hat{x}_2} &  =  &  \hat{A}_{22}\hat{x}_2  +
\hat{A}_{23}\hat{x}_3  + \hat{A}_{24}\hat{x}_4 +  \hat{B}_2u, \qquad
 y_1 = \hat{C}_2\hat{x}_2 \\
 S_3  & : &  & \dot{\hat{x}_3} & = & \hat{A}_{33}\hat{x}_3 + \hat{A}_{34}
\hat{x}_4  \\
S_4  & : &  & \dot{\hat{x}_4} & = & \hat{A}_{44}\hat{x}_4,\qquad
 y_2 = \hat{C}_4\hat{x}_4 \end{array}
\eqno(2.28) $$
 where $S_1$  is  completely  controllable  but
unobservable, $S_2$   is completely  controllable  and observable,
$S_3$   is  completely uncontrollable and unobservable, $S_4$ is
completely uncontrollable but observable.
Eigenvalues of the matrix $\hat{A}_{22}$ are  simultaneously controllable and observable poles. Eigenvalues of
matrix $\hat{A}_{11}$  are unobservable poles (output decoupling poles). Eigenvalues of   matrix $\hat{A}_{44}$
are uncontrollable poles (input decoupling poles). Eigenvalues of  matrix $\hat{A}_{33}$ are simultaneously
uncontrollable and unobservable poles (input - output decoupling poles).

     Let's find  TFM of system (2.27). Using the block  structure
we can determine
$$ \hat{G}(s) = \hat{C}_2 (sI -\hat{A}_{22})^{-1}\hat{B}_2
 \eqno(2.29) $$
Since $G(s) = \hat{G}(s)$ then we obtain the assertion.

\it{ASSERTION  2.7.} \rm   TFM of an incompletely  controllable and observable system coincides with TFM of the
completely controllable and observable subsystem $S_2$. Poles of TFM coincide with  simultaneously  controllable
and observable poles.

\it{CONCLUSION } \rm

\qquad 1. Poles of a completely controllable and observable system
(1.1), (1.2) coincide with  eigenvalues of the dynamics matrix
$A$.

\qquad 2. Poles of an incompletely controllable and/or observable system are  eigenvalues of the dynamics matrix
$A$ without decoupling poles.

\qquad 3. TFM of an incompletely controllable and observable system coincides with TFM of a completely
controllable  and observable subsystem.

\section[Canonical forms of transfer function matrix]{Canonical
    forms of transfer function \\ matrix}

  \subsection[Numerator of transfer function matrix]{Numerator
  of transfer function matrix}

     Let single-input/single-output  system (1.1),(1.2)
has a strictly proper scalar rational transfer function
$$ g(s) = \frac{\psi(s)}{\phi(s)} $$
with $\psi(s)$ and $\phi(s)$ relatively prime \footnote{ The polynomials
$\psi(s)$ and $\phi(s)$ relatively prime if they have not any common
 multipliers being a polynomial in $s$.} polynomials in a  complex
variable $s$ having real coefficients. Orders  of $\psi(s)$ and $\phi(s)$ are $m$ and $n$   respectively $(m <
n)$. Let's present $g(s)$ as
$$  g(s) = \psi(s)[\phi(s)]^{-1} = [\phi(s)]^{-1} \psi(s)$$
We may say that the transfer function $g(s)$ is factorizated as the product of the  polynomial $\psi(s)$  and
the inverse  of the other polynomial $\phi(s)$. The polynomial $\psi(s)$ is known as a numerator of transfer
function $g(s)$.

     We try to get the similar factorization of the strictly  proper
transfer function matrix $G(s)$  with elements are strictly proper rational functions in a complex variable $s$
with   real coefficients. We need to factorizate $G(s)$ into a product  $G(s) = C(s)P(s)^{-1}$ where $C(s)$ and
$P(s)$ are relatively prime polynomial matrices in $s$.\footnote{ A polynomial matrix has polynomials in complex
variable $s$ as elements. We consider only polynomials with real coefficients.}

     At first we introduce some  definitions. In  the  product  $P(s) = Q(s)R(s)$ a  matrix $R(s)$   is
called a right divisor of the matrix $P(s)$  and a matrix $P(s)$ is called a left multiple of the matrix $R(s)$

\it{DEFINITION 2.3.} \rm A square polynomial  matrix  $D(s)$ is called as a common right divisor of  matrices
$C(s)$ and $P(s)$ if and only if
$$ C(s) = C_1(s)D(s), \qquad  P(s) = P_1(s)D(s) \eqno(2.30)  $$
where $C_1(s)$, $P_1(s)$ are some polynomial matrices.

\it{DEFINITION  2.4.} \rm  A square polynomial  matrix $D(s)$ is called as  a greatest common right divisor of
matrices $C(s)$  and $P(s)$
 if  and  only if

a) the matrix $D(s)$ is the common right divisor of matrices
$C(s)$ and $P(s)$,

b) the matrix $D(s)$  is the left multiple of every common right
divisor  of matrices $C(s)$ and  $P(s)$.

     Similarly we can define the greatest common left  divisor  of
matrices $C(s)$ and $P(s)$.

     Let us consider the important  particular  case.

\it{DEFINITION  2.5.} \rm  A square polynomial matrix $U(s)$  is
called as unimodular matrix if and only if $detU(s)  $  is  an
nonzero scalar that  independent in $s$.

     An inverse of the unimodular matrix  is also a polynomial
matrix.

\it{DEFINITION  2.6.} \rm   Two polynomial matrices  $C(s)$  and
$P(s)$ are called as relatively right(left) prime ones if a
greatest common right (left) divisor of $C(s)$ and  $P(s)$ is a
unimodular matrix.

\it{THEOREM  2.1.} \rm  [D2] Any  proper $l \times r$ rational
function matrix (having  rational functions as elements)   always
may be  (nonuniquely) represented as  the  product
$$ G(s) = C(s)P(s)^{-1}   \eqno (2.31)   $$
where  $C(s)$  and $P(s)$   are  relatively prime  polynomial
matrices of dimensions $l\times r$ and $r\times r$ respectively.

     The representation (2.31) is known as a factorization  of a
 transfer function matrix $G(s)$ [W1].

     Similarly  the matrix $G(s)$  may  be  factorizated   into
a product of relatively left prime polynomial matrices $N(s)$ and $Q(s)$  of dimensions $l\times l$  and
$l\times r$ respectively
$$ G(s) = N(s)^{-1}Q(s)   \eqno  (2.32)  $$

It is significant that although matrices $C(s)$ and $P(s)$ in  (2.31) are relatively prime but polynomials
$detC(s)$ and $det P(s)$ are  not  relatively  prime for $l = r $.   For  example,   this result is observed for
the following matrices [D2]
$$ C(s) =  \left [ \begin{array}{cc}  s-1 & 0 \\ 0 & s-2 \end{array} \right ],
\qquad  P(s) = \left [ \begin{array}{cc}  s-2 & 0 \\ 0 & s-1 \end{array}
\right ]  $$

\it{DEFINITION 2.7.} \rm    If   polynomial matrices $C(s)$ and $P(s)$ in (2.31) are relatively right prime then
an $l\times r$ polynomial matrix $C(s)$  is called as a  \it{numerator} \rm of TFM $G(s)$ [W2].

     Similarly if   polynomial   matrices $N(s)$ and $Q(s)$
in (2.32) are relatively left prime then an $l\times r$  polynomial matrix $Q(s)$ is called as a  numerator of
TFM .

     In the following we will show that all numerators of any
TFM  can be transform into  the  uniquely  standard  canonical form.  This canonical form is known as Smith's
form.

    \subsection[Smith form for numerator]{Smith form of numerator}

     For a polynomial matrix with real coefficients  we  introduce
notions of  \it{elementary} \rm  row  (column)  operations [R1], [W1]:

1. interchanging any two rows(columns),

2. multiplication any row (column) by a nonzero real scalar,

3. Adding to any row (column) a product of any other row (column)
by any polynomial or real scalar.

     We need to note that an unimodular matrix is obtained from  the  identity
matrix $I$ by  a  finite  number  of  elementary   row   and column operations  on $I$. Therefore, a determinant
of  an unimodular matrix is a nonzero scalar.

     It  follows from the definition of an unimodular matrix    that  any
sequence of  elementary row (column) operations on a polynomial matrix is equivalent to the premultiplication
(postmultiplication) this matrix by appropriate an unimodular matrix $U_L(s)(U_R(s))$. Such operations we will
call as \it{equivalent} \rm  operations.

\it{DEFINITION 2.8.} \rm   Two polynomial matrices $P(s)$ and
$Q(s)$ will be called  as \it{equivalent} \rm  polynomial matrices
if and only if the first one can be obtained from the second one
by a sequence of  equivalent operations.

     Equivalent polynomial matrices $P(s)$  and $Q(s)$ satisfy the
following relation
$$  P(s) = U_L(s)Q(s)U_R(s)    \eqno(2.33)  $$
where $U_L(s)$ and $U_R(s)$ are unimodular matrices.

     Since an unimodular matrix is nonsingular one then   it
follows  from  (2.33) that equivalent operations do not change the
rank of a polynomial matrix, i.e. $rankP(s) = rankQ(s)$.

     Let's consider reducing an $l\times r$ polynomial matrix
$P(s)$ of the rank $m \le min(r,l)$ to  the  Smith  form  [M1] (or normal  form). We denote polynomial elements
of the matrix $P(s)$ by $p_{ij}(s)$.  Let $p_{jd}(s)$  and  $p_{hk}(s)$ are  two nonzero elements.  We need to
show that if neither  of these elements is a divisor of the other,  then carrying out  only  equivalent
operations we can obtain  a  new  matrix $P_1(s)$, which  contains a  nonzero element of  lower  degree than
either $p_{jd}(s)$ or $p_{hk}(s)$.

We will analyze the three  cases:

     1. If $p_{jd}(s)$  and  $p_{hk}(s)$  are in the same column $( d = k )$,
and $\rho(p_{jd}(s)) \le \rho(p_{hk}(s))$ where $\rho$ is  a degree of a polynomial element, then subtracting
$g(s)$  times the $j$-th row of $P(s)$ from the  $h$-th row we obtain the following relation
$$ p_{hk}(s) =g(s)p_{jk}(s) + r(s) $$
where $g(s)$  is a nonzero polynomial and $r(s)$ is a polynomial with $\rho(r(s)) < \rho(p_{jk}(s))$  or  $r(s)
= 0$. That is $r(s)$ is the lowest degree polynomial remainder after division of the  polynomial $p_{hk}(s)$ by
the polynomial $g(s)$.

If we assume that
\newline
a)  $p_{jk}(s)$ and  $p_{hk}(s)$ are relatively prime
\newline
b) $\rho(p_{jk}(s)) \le \rho(p_{hk}(s))$ \newline then $r(s)$ must be nonzero polynomial and $\rho(r(s)) <
\rho(p_{hk}(s))$, $\rho(r(s)) < \rho(p_{jk}(s))$. Therefore, using only equivalent operations we can decrease
the degree of a element $p_{hk}(s)$ changing this element by $r(s)$,  which is the remainder from the division
$p_{hk}(s)$ by $p_{jk}(s)$.

     2. If $p_{jd}(s) $  and $p_{hk}(s)$  are in the same row $(h=j)$ then
the similar procedure  may be applied where the role of a row/column is interchanged.

     3. If $k \neq d$, $h \neq j$   then the same procedure can  be  applied
to both a row and a column by comparing  $p_{jd}(s) $  and $p_{hk}(s)$ with a element $p_{hd}(s)$.

     Thus we has shown

\it{ASSERTION 2.8. } \rm If $p_{ij}(s)$   is the least degree
element of $P(s)$ and  it does not divide every element  of $P(s)$
then equivalent operations, as just consider, will lead to a new
matrix $P_1(s)$ containing  elements of lower degrees.

Carrying out this procedure many times we can obtain matrices  $P_i(s),i=2,3,\ldots $ containing elements of
more lower degrees. This process be finished after a finite number of steps since a degree of a  polynomial is a
finite positive integer.

     Let's suppose that the process is finished by  a matrix $\bar{P}(s)$. We can
permute  rows and columns of $\bar{P}(s)$ until the element $\bar{p}_{11}(s)$ becomes nonzero and of  a least
degree. Let emphasize that $\bar{p}_{11}(s)$ must divide every element of $\bar{P}(s)$. This important property
of $\bar{p}_{11}(s)$ follows from constructing  matrices  $\bar{P}_1(s)$,
 $\bar{P}_2(s)$,\ldots,  $\bar{P}(s)$. Indeed, let some
$\bar{p}_{st}(s)$ does not  divided by $\bar{p}_{11}(s)$   without a remainder. Then we can represent
$\bar{p}_{st}(s)$ as  $\bar{p}_{st}(s) =
 \bar{g}(s)\bar{p}_{11}(s)+ \bar{r}(s) $  where $\bar{r}(s)$ is a polynomial
with $\rho(\bar{r}(s)) \le \rho(\bar{p}_{st}(s))$. Therefore, we obtain the contradiction with the assumption
that the matrix $\bar{P}(s)$ has the element $\bar{p}_{11}(s)$  of the least degree.

     Taking into account this property of the element $\bar{p}_{11}(s)$ and using
elementary row (column) operations of the third type can reduce the matrix $\bar{P}(s)$ to the following matrix
$\bar{P}_1(s)$
$$     \bar{P}_1(s) = \left[ \begin{array}{ccc} \bar{p}_{11}(s) &  \vdots &
 0  \cdots  0 \\ \dotfill & \dotfill &  \dotfill \\
0 & \vdots &  \\  \vdots  & \vdots & X(s)  \\ 0 & \vdots &
\end{array}  \right]
\eqno (2.34)$$ where $X(s)$ is a some $(l-1)\times (r-1)$ polynomial submatrix.

     Repeating the whole operation with the smaller matrix $X(s)$
without changing the first row or column of $\bar{P}_1(s)$ we get the following matrix  $\bar{P}_2(s)$
$$
    \bar{P}_2(s) = \left[ \begin{array}{cccc} \bar{p}_{11}(s) & 0 & \vdots &
 0  \cdots  0 \\ 0 & \bar{p}_{22}(s) & \vdots & 0 \cdots 0 \\
\dotfill & \dotfill &  \dotfill & \dotfill \\ 0 & 0 & \vdots & \\
\vdots  & \vdots & \vdots & X(s)  \\ 0 & 0 & \vdots &
\end{array}  \right]
\eqno (2.35)$$
 Continuing this  process   we  reduce finally the matrix $P(s)$ to the form
$$
    S(s) = \left[ \begin{array}{cccccc} s_1(s) & 0 & \cdots &  0 & \cdots & 0 \\
0 & s_2(s) & \cdots & 0 &  \cdots &  0 \\ \vdots & \vdots & & \vdots & & \vdots \\
0 & 0 & \cdots & s_m(s) & \cdots & 0 \\ \vdots & \vdots & & \vdots & & \vdots \\
0  & 0 & \cdots & 0 & \cdots & 0 \end{array}  \right] \eqno
(2.36)$$
 where $m$  is a normal rank\footnote{ The normal rank
(or rank) of a polynomial matrix is the order of the largest minor, which is not identically zero [B2]}  of the
$l\times r$ polynomial matrix  $P(s)$ and  every  $s_i(s)$  divides $s_j(s)$ without a remainder $(i<j)$.
Divisibility $s_j(s)$ by $s_i(s)$ follows from the construction of the matrix  $S(s)$  because the element
$\bar{p}_{11}(s)$  divides all elements of $X(s)$ and so on.

     Since the matrix $S(s)$ is resulted from the matrix $P(s)$ by a
sequence of equivalent operations which could be  realized  by unimodular matrices $U_L(s)$ and $U_R(s)$ (i.e.
$S(s) = U_L(s)P(s)U_R(s)$) then
$$   rank S(s) = rank P(s) = m   $$
The matrix $S(s)$ (2.36) is known as the Smith canonical form for a polynomial matrix or the Smith form [L2].

\subsection[Invariant polynomials]{Invariant polynomials}

     Now we show that only equivalent polynomial matrices  may  be
reduced to identical Smith forms. Let  a  polynomial matrix $P(s)$ has  a normal  rank $m \le min(r,l)$. A
greatest common divisor of all  $j$-th order minors $(j=1,2,\ldots,m)$ of $P(s)$ we denote  by $d_j(s)$. Since
any $j$-th order minor $(j>2)$ would be expressed as a linear combination of ($j-1$)-th order minors then
$d_{j-1}(s)$  is  the divisor of $d_j(s)$. If we denote $d_o(s) = 1$ then every $d_j(s)$  is divided by
$d_{j-1}(s),\;\; j=1,\ldots,m $   in the sequence $d_o(s), d_1(s) ,\ldots,d_m(s)$. We define
$$   \epsilon_m(s) = \frac{d_m(s)}{d_{m-1}(s)} , \qquad \epsilon_{m-1}(s) =
\frac{d_{m-1}(s)}{d_{m-2}(s)} , \qquad \cdots ,\qquad
\epsilon_1(s) = \frac{d_1(s)}{d_o(s)} \eqno(2.37) $$ Polynomials
$\epsilon_1(s) ,\ldots,\epsilon_m(s)$ are called as \it{
invariant} \rm polynomials of $P(s)$.  It  can  be  shown these
polynomials are  invariants  under  equivalent operations.

\it{ASSERTION 2.9.} \rm  Two equivalent polynomial matrices $P(s)$
and $Q(s)$  have equal invariant polynomials.

\it{PROOF}. \rm If two matrices $P(s)$ and $Q(s)$ are equivalent then two unimodular matrices $U_L(s)$ and
$U_R(s)$ exist such that $ Q(s) = U_L(s)P(s)U_R(s)$.  Unimodular   matrices $U_L(s)$ and $U_R(s)$ have nonzero
scalar determinants  then $rank Q(s) = rankP(s)$. Let's denote greatest common divisors of all $j$-th order
minors of matrices $P(s)$ and $Q(s)$ by $d_k(s)$ and $\delta_k(s)$ respectively.  From the equality $ Q(s) =
U_L(s)P(s)U_R(s)$  we obtain  that every $k$-th order minor $(1 \le k \le m )$ of the matrix $Q(s)$ should be
expressed by the formula Caushy-Binet [G1] as a linear combination of  $k$-th order minors of $P(s)$. Hence
$\delta_k(s)$ is divided by $d_k(s)$. Vise versa: from the equality $ P(s) = U_L(s)^{-1}Q(s)U_R(s)^{-1}$  it
follows  divisibility of  $d_k(s)$ by $\delta_k(s)$. That is why
$$    d_k(s) =  \delta_k(s) , \qquad k=1,\ldots,n $$
Consequently $P(s)$  and $Q(s)$  have  equal invariant polynomials. This completes the proof.

     Let us calculate  invariant polynomials  of the matrix   $S(s)$
(2.36). Since
$$ d_1(s) = s_1(s), \qquad  d_2(s) = s_1(s)s_2(s), \qquad \ldots, \qquad
d_m(s) = s_1(s)s_2(s) \cdots s_m(s)  $$
then  we have from (2.37)
$$  \epsilon_m(s) = s_m(s), \qquad \epsilon_{m-1}(s) = s_{m-1}(s), \qquad
\ldots, \qquad \epsilon_1(s) = s_1(s)      \eqno(2.38) $$ Hence,
diagonal elements of the Smith form coincide with invariant
polynomials of $S(s)$.

     As the matrix  $S(s)$ is obtained from the  matrix  $P(s)$  by  the
sequence of equivalent operations then these matrices are equivalent and have similar invariant polynomials. The
last follows from Assertion 2.9. We have obtained the following assertion.

\it{ASSERTION 2.10.} \rm  Any polynomial matrix $P(s)$ is reduced
to Smith form (2.36) with diagonal elements  $s_i(s)$ that are
invariant polynomials of $P(s)$ .

       Since any two equivalent  polynomial  matrices  have  equal
invariant polynomials then the following corollary is true.

\it{COROLLARY  2.1. } \rm Any two polynomial matrices have the
unique Smith form.

     Let us consider two any  numerators $P(s)$ and  $Q(s)$  of a proper
transfer function matrix $G(s)$. It is evident that  $P(s)$ and $Q(s)$ are $l\times r$ polynomial matrices. As
it has been shown in [W2]  if $P(s)$ and  $Q(s)$ are two numerators of a rational function matrix $G(s)$ then
they are equivalent. This means that polynomial matrices $P(s)$ and $Q(s)$ satisfy  the relation (2.33) and by
Corollary 2.1 they have the unique Smith form and equal invariant polynomials.

     So, we conclude that all numerators of  TFM
$G(s)$  have equal invariant polynomials and may be reduced to the unique Smith canonical form.

         \it{ EXAMPLE 2.1. } \rm

 Let's calculate the  Smith  form of the following matrix
$$      P(s) = \left[ \begin{array}{ccc} s & 0 & 0 \\ 0 & s & s+1 \\
s & s-1 & 0  \end{array}  \right] \eqno (2.39)$$
 Here  $ r = l =3$
and $m = rankP(s) = 3$.

    At first  we  construct  the matrix $P_1(s)$ (2.24)
by subtracting the second  row of the matrix (2.39) from the third one
$$      P_1(s) = \left[ \begin{array}{ccc} s & 0 & 0 \\  0 & s & s+1 \\
s & -1 & -s-1  \end{array}  \right]  $$
Interchanging  rows and
columns in $P_1(s)$ we can reduce it to the following form with
$\bar{p}_{11}(s) = -1 $
$$  \bar{P}_1(s) = \left[ \begin{array}{ccc} -1 & s & -s-1 \\ s & 0 & s+1\\
 0 &  s &  0  \end{array}  \right]  $$
Using the second and third types elementary operations  we  obtain the matrix $\bar{P}_1(s)$  in  the form
(2.34)
$$     \bar{P}_1(s) = \left[ \begin{array}{ccc} -1 &  \vdots &
 0 \qquad 0 \\ \dotfill & &  \dotfill \\  0 & \vdots & X(s) \\
 0 & \vdots & \end{array}  \right]  \eqno (2.40)$$
with
$$ X(s) = \left[ \begin{array}{cc} s^2 & -s^2 +1 \\ s & 0 \end{array} \right] $$
 Adding the first column to second one in $X(s)$ and
interchanging  columns obtained we result in
$$ \bar{X}_1(s) =  \left[ \begin{array}{cc} 1 & s^2 \\ s & s\end{array} \right] $$
Using the third and second type elementary operations we transform $\bar{X}_1(s)$ into the form (2.34)
$$ \bar{X}_1(s) =  \left[ \begin{array}{cc} 1 & 0 \\ s & s(s^2-1)
\end{array} \right] $$
So, we have reduced the matrix $\bar{P}_1(s)$ to the following one
$$ \bar{P}_2(s) =  \left[ \begin{array}{ccc} 1 & 0 & 0 \\ 0 & 1 & 0 \\
0 & 0 &  s(s^2-1) \end{array} \right] \eqno(2.41)  $$
 This matrix
has  the Smith form (2.36)  with $s_1(s) =1$, $s_2(s) =1$, $s_3(s) =s(s^2-1)$.

     For checking  we calculate  invariant  polynomials  of the
matrix (2.39). It has three nonzero first order minors: $s$, $\;s-1$, $\;s+1$ with the great common divisor
$d_1(s) =1$,  four nonzero second order minors: $s^2$,  $\;s(s+1)$,  $\;s(s-1)$,  $\;(s+1)(s-1)$ with the great
common divisor $d_2(s) =1$  and  one third  order nonzero minor   $d_3(s) =s(s-1)(s+1)$. Using (2.37) we
determine invariant polynomials
$$   \epsilon_3(s) = \frac{d_3(s)}{d_2(s)} = s(s^2-1), \qquad
\epsilon_2(s) = \frac{d_2(s)}{d_1(s)}= 1, \qquad \epsilon_1(s) =
\frac{d_1(s)}{d_o(s)} = 1  $$

     So,  invariant polynomials of
$P(s)$   coincide with  diagonal elements of the form (2.41). This result is adjusted with Assertion 2.10.

\subsection[Smith-McMillan form of transfer function
   matrix]{Smith-McMillan form of transfer function matrix}

     Now we demonstrate that the canonical form of a rational function matrix
(the Smith-McMillan canonical form)  can  be  obtained by using  the Smith  canonical  form.  We consider any $l
\times r$ rational function matrix $W(s)$ having a rank $m \le min(r,l)$. Let a polynomial $\phi(s)$ is a monic
least common denominator of all elements of $W(s)$. We form  the  matrix $T(s) = \phi(s)W(s)$ that is a
polynomial $l\times r$ matrix. Let $S_T(s)$ is the Smith form of the matrix $T(s)$, i.e.
$$  S_T(s) = U_L(s)T(s)U_R(s)   \eqno(2.42)   $$
where $U_L(s)$ and $U_R(s)$  are  unimodular $l\times l$ and
$r\times r$ matrices respectively  and $S_T(s)$ has the following
structure
$$     S_T(s) = \left[ \begin{array}{ccc} diag \{ s_{T1}(s),s_{T2}(s),\ldots,
s_{Tm}(s) \} & \vdots & \; O \; \\ \dotfill & &  \dotfill \\  O &
\vdots & \; O \; \end{array}  \right] \eqno (2.43)$$
 where
$s_{Ti}(s)$  are  invariant polynomials of $T(s)$ . Since $T(s) = \phi(s)W(s)$  then dividing both left-hand and
right-hand sides of (2.42) by $\phi(s)$  ($s \neq s_i$, $s_i$ is a zero of $\phi(s)$)  yields
$$     \frac{S_T(s)}{\phi(s)} =  U_L(s)W(s)U_R(s)
 \eqno(2.44)  $$
To discover  a structure of $ S_T(s)/ \phi(s) $ we  use (2.43)
$$   \frac{S_T(s)}{\phi(s)}   = \left[
\begin{array}{ccc} diag
 \{  \frac{s_{T1}(s)}{\phi(s)},\frac{s_{T2}(s)}{\phi(s)}, \ldots,
 \frac{s_{Tm}(s)}{\phi(s)} \} & \vdots & \; O \; \\ \dotfill & &  \dotfill \\
O & \vdots & \; O \; \end{array}  \right]
 \eqno (2.45) $$
 Since  elements $s_{Ti}(s)$, $i=1,2,\ldots,m$ are  invariant
polynomials of $T(s)$  then they satisfy the following condition: $s_{Ti}(s)$ is divided by $s_{T,i-1}(s)$
without a remainder. Therefore,  elements $s_{Ti}(s) /\phi(s) $ satisfy the same requirement. Carrying  out all
possible cancellations in $ s_{Ti}(s)/ \phi(s) $ we result in a ratio of two monic polynomials $ \epsilon_i(s)/
\psi_i(s) $, which have to satisfy  conditions:
 $$    \frac{\epsilon_i(s)}{\psi_i(s)} \qquad  \rm{must \; \; divide}  \qquad
   \frac{\epsilon_{i+1}(s)}{\psi_{i+1}(s)}, \qquad \it{i=1,\ldots,m-1}
\eqno (2.46) $$
 It follows from (2.46) that
 $$    \epsilon_i(s) \qquad \rm{ must \;\; divide}  \qquad \epsilon_{i+1}(s)
\qquad \it{i=1,\ldots,m-1} \eqno (2.47) $$

 $$    \psi_{i+1}(s) \qquad \rm{must \; \;divide}  \qquad \psi_i(s),
 \qquad \it{i=1,\ldots,m-1} ,\qquad  \psi_1(s) = \phi(s)
\eqno (2.48) $$
 Denoting  $M(s) = S_T(s)/ \phi(s)$ we obtain from
(2.45) and the last relation
$$   M(s) = \left[ \begin{array}{ccc} diag
 \{  \frac{\epsilon_1(s)}{\psi_1(s)},\frac{\epsilon_2(s)}{\phi_2(s)}, \ldots,
 \frac{\epsilon_m(s)}{\phi_m(s)} \} & \vdots & \; O \; \\ \dotfill & &
\dotfill \\ O & \vdots & \; O \; \end{array}  \right]
 \eqno (2.49) $$
This matrix is known as the Smith-McMillan canonical  form  of  a
rational function matrix $W(s)$. Using (2.44) we can write
$$     M(s) =  U_L(s)W(s)U_R(s)
 \eqno(2.50)  $$
Thus, a rational function matrix $W(s)$  is reduced to the Smith-McMillan canonical form by the sequence of
elementary operations.

                 \it{EXAMPLE   2.2} \rm [S14]

 Let a $4\times 3$ transfer function matrix has the form
$$ W(s) =  \left[ \begin{array}{ccc} \frac{1}{s(s+1)^2} & \frac{s^2+2s-1}
{s(s+1)^2}  &  \frac{s+2}{s-2} \\ \\ 0 & \frac{s+2}{(s+1)^2} & 0 \\
\\ 0 & 0 &  \frac{3(s+2)}{s+1} \\\\  \frac{s+3}{s(s+1)^2} &
\frac{2s^2+3s+3}
 {s(s+1)^2}  &  \frac{s+2}{s-2}  \end{array} \right]
\eqno(2.51)  $$
 We calculate  the  least common  denominator of
all elements of $W(s): \phi(s) = s(s + 1)^2$  and form the matrix
$T(s) = \phi(s)W(s)$
$$ T(s) =  \left[ \begin{array}{ccc} 1  &  s^2 +2s-1 &
(s+2)s(s+1)\\ \\
0 &  s(s+2)  &   0 \\  \\ 0  &  0  &  3s(s+1)(s+2) \\  \\ s+3 &
2s^2 +3s-3 & s(s+1)(s+2) \end{array} \right]  $$

 Let's  find the
Smith form of $T(s)$. By  (2.36) we need to know invariant  polynomials $\epsilon_i(s) = s_{Ti}(s)$ of $T(s)$.
For this  we calculate a greatest common divisor $d_i(s)$ of the  $i$ order minors ($i=1,2,3$)  of the matrix
$T(s)$: $d_o  = 1$, $d_1(s)  = 1$, $d_2(s) = s(s + 2)$, $d_3(s) = s^2(s + 1)(s + 2)^2$ and using (2.37)  find
$$
   s_{T1}(s) = \frac{d_1(s)}{d_0(s)} = 1, \qquad  s_{T2}(s) =
  \frac{d_2(s)}{d_1(s)} = s(s+2), \qquad    s_{T3}(s) = \frac{d_3(s)}{d_2(s)}
   = s(s+1)(s+2) $$
Thus, the Smith form $S_T(s)$  of  the matrix  $T(s)$ is
$$ S_T(s) =  \left[ \begin{array}{ccc} 1  &  0 & 0 \\ \\
0 &  s(s+2)  &   0 \\ \\  0  &  0  &  s(s+1)(s+2) \\ \\ 0 & 0 & 0
 \end{array} \right]  $$
 Calculating  $ \epsilon_i(s)/ \psi_i(s) =  s_{Ti}(s)/ \phi(s) $
 with $\phi(s)=s(s+1)^2$
$$ \frac{\epsilon_1(s)}{\psi_1(s)} = \frac{ s_{T1}(s)}{\phi(s)} = \frac{1}
{s(s+1)^2}, \qquad  \frac{\epsilon_2(s)}{\psi_2(s)} = \frac{s_{T2}(s)}{\phi(s)}
 = \frac{s+2}{(s+1)^2}, \qquad  \frac{\epsilon_3(s)}{\psi_3(s)}  =
\frac{s_{T3}(s)}{\phi(s)} = \frac{s+2}{(s+1)}  $$
 and using the
formula (2.49) we  obtain  the  matrix  $M(s)$ being  the
Smith-McMillan form of  $W(s)$
$$ M(s) =  \left[ \begin{array}{ccc} \frac{1}{s(s+1)^2} & 0  & 0
\\\\
0 & \frac{s+2}{(s+1)^2} & 0 \\ \\ 0 & 0 &  \frac{s+2}{s+1} \\ \\ 0
& 0 & 0  \end{array} \right] \eqno(2.52)  $$

    Now we reveal the relationship between the  matrix $M(s)$  and the Smith
canonical form of any $l\times r$ numerator of a transfer function matrix of a rank $m$.
     For simplicity we let $ m = min(r,l)$ and   denote
$$ E(s) = diag(\epsilon_1(s), \epsilon_2(s),\ldots,
\epsilon_m(s))$$
$$ \Psi(s) = diag(\psi_1(s),\psi_2(s),\ldots, \psi_m(s))
\eqno (2.53) $$
 Using (2.53) we can present the matrix $M(s)$ (2.49)
as
$$ M(s) \; = \;  E(s)\Psi(s)^{-1} \; = \; \Psi(s)^{-1} E(s)
\eqno(2.54)$$
 Thus, we obtain that the Smith-McMillan form of a rational function matrix $M(s)$
is factorized into the product of the $l\times r$ polynomial   matrix  $E(s)$ and an inverse of the $r\times r$
polynomial    matrix $\Psi(s)$. Matrices  $E(s)$ and $\Psi(s)$ are the relatively right prime ones and $E(s)$ is
a numerator of the rational function matrix $M(s)$. The matrix $E(s)$ is in the Smith form that is unique one
for all numerators of TFM (2.51) by Corollary 2.1.

                  \it{CONCLUSION} \rm

     Polynomials $\epsilon_i(s)$ of  the  Smith-McMillan
  can be calculated by the factorization of TFM  into
a product of relatively prime polynomial matrices and following calculating the Smith form of the numerator.

\chapter[Notions of transmission and invariant zeros]{Notions
   of transmission and invariant zeros}

  \section[Classic definition of zeros]{Classic definition of zeros}

    Let's consider  system  (1.1), (1.2)  with  a single input and
single output $(r=l=1)$. The transfer function of this system is defined from  formula (2.7)  as
$$  g(s) = c(sI_n - A)^{-1}b
     \eqno(3.1)   $$
where $c$ is an $n$-dimensional row vector, $b$ is an $n$-dimensional column vector. We propose that system
(1.1),(1.2) is completely controllable and observable. Then its a transfer function $g(s)$ is the ratio of two
relatively prime polynomials
$$   g(s) = \frac{\psi(s)}{\phi(s)}
 \eqno(3.2)    $$
where  $\phi(s) = det(sI_n - A)$  is a polynomial of order $n$ with zeros $\lambda_1, \lambda_2, \ldots,
\lambda_n$  being poles of $g(s)$ and $\psi(s) = c (adj(sI_n -A))b$  is a polynomial of order $m <n$. Zeros
 $z_1, z_2,\ldots,z_m$  of $\psi(s)$  are called as zeros of  the
scalar transfer function (TF) $g(s)$. Since $g(s)$ is  an irreducible rational function then  $\lambda_i \neq
z_j$ $(i=1,2,\ldots,n; \;j=1,2,\ldots,m)$.

      Now we study the physical interpretation of zeros of TF.  Let
a scalar exponential  signal  is  applied  at the input of system (1.1),(1.2)
$$         u(t) = u_o\exp(jwt)1(t), \qquad  t>0
 \eqno(3.3)  $$
where  $jw$ is a complex frequency, $w$ is a real value, $u_o \neq 0$ is a scalar constant value, 1(t) is the
unit step function\footnote{ $1(t)= \left\{ \begin{array}{rl} 0,& \mbox{if} \;\;t\leq 0\\1 ,& \mbox{if}\;\; t >
0 \end{array} \right. $}. According to formulas (2.14),(3.2)  we obtain the following output response for
$x(t_o) = 0$
$$ y(t) \cong g(jwt)u(t) = \frac{\psi(jw)}{\phi(jw)}u_o\exp(jwt)\eqno  (3.4) $$
It is  follows from (3.4)
that $\psi(jw) = 0$
 $(\phi(jw) \neq 0)$  if the complex frequency $jw$ coincides with a
certain zero $z_i$ of $g(s)$. In this case the output response be trivial (identically zero).

      So, we  conclude:  In the classic
 single input/output controllable  and observable system a transmission zero is defined as a value of
a complex  frequency $ s =z_i$ at which  the transmission of  the exponential signal $\exp(z_it)$ is blocked. It
is evident that the transmission zero coincides is any zero of the numerator of the transfer function $g(s)$.

     Similarly we may  define
a transmission zero of a  linear multivariable system (having a vector input/output) as a complex frequency at
which  the transmission of a signal is 'blocked'. In the following section we will define  transmission zeros
via the transfer function matrix $G(s)$.

 \section[Definition of transmission zero via transfer function matrix]
            {Definition of transmission zero via transfer function matrix}

     Lets's consider completely controllable  and  observable  system
(1.1), (1.2) with $r$ inputs and $l$ outputs $(r,l> 1)$. The transfer function matrix $G(s)$ of this  system is
a rational function matrix  with elements being rational irreducible scalar functions. We propose that $l \ge r$
and  $rank G(s) = min(r,l) = r$. Let a exponential signal
$$         u(t) = u_o\exp(jwt)1(t), \qquad  t>0
\eqno(3.5)  $$
 is applied at the input of (1.1),(1.2). In (3.5) $jw$ is a complex frequency, $w$ is a  real
value,  $u_o \neq 0$ is a nonzero constant $r$  vector. According to formula (2.14) we can write  the following
output steady response for $x(t_o) = O $
$$  y(t) \cong G(jwt)u(t) = G(jw)u_o\exp(jwt)\eqno(3.6)  $$
 We assume that  $jw \neq \lambda_i$ where
 $\lambda_i$ ($i=1,2,\ldots,n$)  are poles of the system, which
coincide with zeros of the polynomial $\phi(s) = det(sI_n - A)$
for the completely controllable and observable system.

     If  the rank of $G(s)$ is  locally  reduced  at $s = jw$ then
 a nonzero constant $r$ vector $v$ exists, which is a nontrivial solution of the linear system
$$            G(jw)v = O  \eqno (3.7)   $$
Setting $u_o =v$ in (3.5) and taking into account (3.7) we have  $y(t) = G(jw)v\exp(jwt)= O$. This fact means
that there exists the exponential input  vector (3.5) of a corresponding complex frequency $jw$ such that the
output steady response is identically zero. The value of a complex frequency  $s =jw$ that locally reduces the
rank of $G(s)$ is the transmission zero.

\it{DEFINITION 3.1.} \rm  A complex frequency $s =z$ at which a rank of the transfer function matrix is locally
reduced
$$    rankG(s)/_{s=z} < min(r,l)
\eqno(3.8) $$
is called as a  transmission zero.

     Thus, if  a complex
frequency $s$ coincides with a transmission zero $z$ of a multi-input/multi-output system then there exists some
nonzero proportional $exp(zt)$  input vector such that its propagating through the system is blocked.

     Let's consider the difference between the multi-input/multi-output
case and the classic one. The  scalar  transfer  function $g(s)$ vanishes at $s=z$ in the
single-input/single-output case, then  any proportional $exp(zt)$ input signal does not transmit through the
system ('blocked'). In the multi-input/multi-output case the transfer function matrix does not become the zero
matrix at $s=z$  but its the rank is locally reduced. This fact means that there exists a nonzero proportional
$exp(zt)$ input vector signal,  which may propagate through the system.

 \it{REMARK 3.1.} \rm    The inequality (3.8) for $l \ge r$ is the necessary and sufficient
condition  for  'blocking'  the   transmission of a exponential signal. It is only sufficient condition for $r >
l$. Namely,  if the condition (3.8) is satisfied for $s = z$ then  there always exists a proportional to
$exp(zt)$ input signal such that $y(t) = O, t \gg 0$. The inverse proposition ( if $y(t) = O$ for a proportional
$exp(zt)$ input signal then  $rankG(s) < l$) is not true. Indeed, consider the following $2\times 3$ transfer
function matrix
$$ G(s) = \frac{1}{(s+1)(s+2)} \left [ \begin{array}{ccc}  s & 0 & s\\
0 & s-0.5 & 0  \end{array} \right ]  $$
 of the normal rank $2$. The steady response $y(t)$
to the input signal  $u(t) = [1, 0, -1]^T \exp^{jwt}1(t) (w \neq 0)$ is defined by formula (3.6)
$$  y(t) = \frac{1}{(jw+1)(jw+2)} \left [ \begin{array}{ccc}  jw & 0 & jw \\
0 & jw-0.5 & 0  \end{array} \right ] \left [ \begin{array}{r} 1 \\0 \\ -1
\end{array} \right ]\exp(jwt) = \left [ \begin{array}{c} 0 \\0 \end{array}
 \right ] $$
Hence the steady response is equal to zero although the rank of
$G(s)/_{s=jw}$ is not reduced.

     The following  simple  example illustrates blocking an
oscillatory input signal of the frequency coincided with a transmission zero.

             {\it EXAMPLE 3.1.} \rm

Let completely controllable and observable  system (1.1),(1.2)
with  $n = 3$, $l = r = 2$ has the transfer function matrix
$$ G(s) = \frac{1}{(s+1)(s+2)(s+3)} \left [ \begin{array}{cc}  s^2 + 0.25 & 0 \\
0 & 1  \end{array} \right ]
 \eqno(3.9)$$
One can  see that the rank of $G(s)$ is reduced at $s = \pm 0.5j$. Let's inject in the system the following
oscillatory signal of frequency 0.5
$$ u(t) = \left [ \begin{array}{c} v_1\\ v_2  \end{array} \right ] 2cos(0.5t),
\qquad v_1 \neq v_2 \neq 0 $$
This signal can be represented as
the sum of  two complex exponential signals
$$ u(t) = \left [ \begin{array}{c} v_1\\ v_2  \end{array} \right ] 2cos(0.5t)
= \left [ \begin{array}{c} v_1\\ v_2  \end{array} \right ]
(\exp(0.5jt)+ \exp(0.5jt)) = u_1(t) +u_2(t) $$
 Using  formula
(3.6) we find the output steady response to zero initial conditions
$$ y(t) = y_1(t) + y_2(t) = G(0.5j)u_1(t) + G(-0.5j)u_2(t) =$$
$$=\frac{1}{(0.5j+1)(0.5j+2)(0.5j+3)} \left [ \begin{array}{cc} 0 &
0 \\ 0 & 1 \end{array} \right ] \left [ \begin{array}{c} v_1\\v_2
\end{array} \right ]
 \exp(0.5jt) + $$
$$ +\frac{1}{(-0.5j+1)(-0.5j+2)(-0.5j+3)} \left [ \begin{array}{cc}
 0 & 0 \\ 0 & 1 \end{array} \right ] \left [ \begin{array}{c} v_1\\v_2
\end{array} \right ] \exp(-0.5jt) $$
It is easy to verify  that $y(t) \equiv O$  for $v_1 = q \neq 0, v_2 = 0$ where $q$ is any constant value. This
fact means that the steady output response to the nonzero input signal
$$  u(t) =  \left [ \begin{array}{c} 2q \\ 0 \end{array} \right ] 2cos(0.5t) $$
 is the identically zero.

       \section[Transmission zero and system response]
                  {Transmission zero and system response}

     Now we consider the output response of system (1.1),(1.2) with
the nonzero initial state  $x(t_o) \neq O $. We  assume that matrix $A$ has $n$  distinct eigenvalues
$\lambda_1,\ldots,\lambda_n$. Let's apply  the   following the exponential type input
$$         u(t) = u_o\exp(\alpha^*t)1(t)
\eqno(3.10) $$ where $u_o \neq O$ is a constant  $r$  vector, $\alpha^* \neq \lambda_i$, $i=1,\ldots,n$  is a
complex number. Substituting (3.10) in  formula (1.9) (see Sec 1.1) and assuming $t_o = 0$  yields
$$
 y(t)= \sum_{i=1}^n\gamma_i\exp(\lambda_it)v_i^Tx_o +
\sum_{i=1}^n\gamma_i\int_{0}^{t}\exp(\lambda_i(t-\tau )) \beta_i^Tu_o\exp(\alpha^*\tau)d\tau $$
  where $\gamma_i
=Cw_i, \;\; \beta_i^T = v_i^TB, \;\; x_o = x(t_o)=0$.

 Taking out from the integral the terms that
independent on $\tau$
$$
 y(t)= \sum_{i=1}^n\gamma_i\exp(\lambda_it)v_i^Tx_o +
\sum_{i=1}^n\gamma_i\beta_i^Tu_o\exp(\lambda_it)\int_{0}^{t}
\exp((\alpha^*- \lambda_i)\tau)d\tau      $$
 and integrating we obtain
$$ y(t)= \sum_{i=1}^n\gamma_i\exp(\lambda_it)v_i^Tx_o +
\sum_{i=1}^n\frac{\gamma_i\beta_i^Tu_o\exp(\lambda_it)}{\alpha^* -
\lambda_i} (\exp((\alpha^* - \lambda_i)t) - 1) = $$
$$= \sum_{i=1}^n\gamma_i\exp(\lambda_it)v_i^Tx_o +
    \sum_{i=1}^n\frac{\gamma_i\beta_i^Tu_o}{\alpha^* - \lambda_i}
    (\exp(\alpha^*t) - \exp(\lambda_it)) = $$
$$ = \sum_{i=1}^n\gamma_i\exp(\lambda_it)v_i^Tx_o -
\sum_{i=1}^n\frac{\gamma_i\beta_i^Tu_o}{\alpha^* -
\lambda_i}\exp(\lambda_it) +
\sum_{i=1}^n\frac{\gamma_i\beta_i^Tu_o}{\alpha^* - \lambda_i}
\exp(\alpha^*t)
 =  y_o(t) + y_1(t) + y_2(t)
\eqno(3.11) $$
 The analysis of this expression shows that the
response  $y(t)$ is the sum of the following components: the response to the initial conditions $y_o(t)$, the
term  associated with the free motion $y_1(t)$ and the forced response $y_2(t)$. If all $\lambda_i$ have
negative real part then $y(t) \to y_2(t)$ as $t \to \infty $. For small time $y(t)$ depends on $y_o(t)$,
$y_1(t)$ and $y_2(t)$.

     We always may choose  vector $x_o$ such that the components
$y_o(t)$ and $y_1(t)$ are mutual eliminated. Indeed, from the equation
$$  y_o(t) + y_1(t) = \sum_{i=1}^n\gamma_i\exp(\lambda_it)v_i^Tx_o -
\sum_{i=1}^n\gamma_i\frac{\exp(\lambda_it)v_i^T}{\alpha^* -
\lambda_i}Bu_o = O $$ we can write
$$        v_i^Tx_o = \frac{v_i^T}{\alpha^* - \lambda_i}Bu_o  $$
Changing  $i$  from  1  to  $n$  and denoting  $V = [v_1, v_2,
\ldots, v_n]$ we obtain the system of  linear equations in $x_o$
$$     V^Tx_o = diag(\alpha^* - \lambda_1, \alpha^* - \lambda_2, \ldots,
  \alpha^* - \lambda_n)^{-1}V^TBu_o $$
from which $x_o$ is calculated as
$$ x_o = (V^T)^{-1}diag(\alpha^* - \lambda_1, \alpha^* - \lambda_2, \ldots,
 \alpha^* - \lambda_n)^{-1}V^TBu_o $$
 Since $(V^T)^{-1} = W$ where the matrix $W$ consists of  right
eigenvectors $w_i$ ($W = [w_1, w_2, \ldots, w_n]$) then  the right-hand side of the last expression becomes
$$ W diag(\alpha^* - \lambda_1, \alpha^* - \lambda_2, \ldots,
\alpha^* - \lambda_n)^{-1}V^TBu_o \;\;=\;\; (\alpha^*I_n - A)^{-1}Bu_o \eqno (3.12) $$
 and
$$         x_o \; = \; (\alpha^*I_n - A)^{-1}Bu_o \eqno(3.13) $$
Substituting this $x_o$ in (3.11) we obtain the output response
containing only the forced response
$$ y(t) = y_2(t) = \sum_{i=1}^n\frac{\gamma_i\beta_i^T}
    {\alpha^* - \lambda_i} u_o\exp(\alpha^*t)
\eqno(3.14a) $$
Using   notions  (1.10) ( $\gamma_i =Cw_i, \;
\beta_i^T = v_i^TB$) and the relation (3.12) we can rewrite
(3.14a) as
$$ y_2(t) = C[w_1, w_2, \ldots, w_n]diag(\alpha^* - \lambda_1,
   \alpha^* - \lambda_2, \ldots, \alpha^* - \lambda_n)^{-1}[v_1,v_2,\ldots,
   v_n]^TBu_o \exp(\alpha^*t) = $$
$$ = CW diag(\alpha^* - \lambda_1,\alpha^* -
\lambda_2, \ldots, \alpha^* - \lambda_n)^{-1}V^TBu_o \exp(\alpha^*t) = C(\alpha^*I_n -
A)^{-1}Bu_o\exp(\alpha^*t) $$
 Since $C(\alpha^*I_n - A)^{-1}B$ is the transfer function matrix of  system
(1.1),(1.2) at $s = \alpha^*$  then
$$y_2(t) = G(\alpha^*)u_o\exp(\alpha^*t)
\eqno(3.14b)  $$

If the complex frequency $\alpha^*$ coincides with  a  transmission zero then the rank of $G(\alpha^*)$  is
reduced and  there  exists a  nonzero vector $u_o$  such as  $y_2(t) \equiv O$.

     So, it has been shown: if $\alpha^*$ is a transmission  zero  then
there exists a nonzero vector $u_o$ and a nonzero  initial  state condition $x_o(\alpha^*I_n - A)^{-1}Bu_o$ such
as the  output response to the input (3.10) is identically zero, i.e. $y(t) = y_o(t) + y_1(t) + y_2(t) \equiv
O$.

               { \it CONCLUSION}  \rm

     If a proportional $ exp(\alpha^*t)$ signal is applied to an input
of a completely controllable  and  observable  system  where
$\alpha^*$ is  a transmission zero, $\alpha^* \neq \lambda_i$ (
 distinct eigenvalues of $A$) then

     a. there exists an initial state conditions $\bar{x}_o$ such that the
transmission of this signal through the system is blocked: $y(t) \equiv O$,

     b.  for  $x(t_0) \neq \bar{x}_o  \neq O$ the output response $y(t)$
  consists of the sum  $y_0(t) + y_1(t)$;  the transmission of forced
response $y_2(t)$ is blocked,

     c. for  $x(t_0) \neq \bar{x}_o = O$ the output response $y(t)$ contains
only the free motion term $y_1(t)$; the transmission of forced
response $y_2(t)$ is blocked.

     For the illustration we consider the following example.

               {\it EXAMPLE 3.2.} \rm

Let's completely  controllable  and  observable system (1.1),(1.2)
with $ n = 3$, $r = l  = 2$ has the following matrices $A$, $B$,
$C$
$$ A = \left[ \begin{array}{rrr} 0 & 1 & 0 \\ 0 & 0 & 1 \\-6 & -11 & -6
 \end{array} \right ],  \qquad B = \left[ \begin{array}{rr} -1 & 0 \\
0 & 0 \\ 0 & -1 \end{array} \right ], \qquad C = \left[
\begin{array}{rrr} 0 & -1 & 1 \\ -1 & -1 & 0 \end{array} \right
]$$
 Poles of this system coincide with  eigenvalues of $A$:
  $\lambda_1 = -1$,
$\lambda_2 = -2$, $\lambda_3 = -3$. At first we find eigenvectors
$w_i$ and  $v_i^T$, $ i =1,2,3$. Since $A$ has distinct
eigenvalues then the matrix $W$ is the Vandermonde matrix of the
structure
$$
   W = [w_1, w_2, w_3] = \left[ \begin{array}{rrr} 1 & 1 & 1 \\
\lambda_1 & \lambda_2 & \lambda_3  \\ \lambda_1^2 & \lambda_2^2 & \lambda_3^2
\end{array} \right ] = \left[ \begin{array}{rrr} 1 & 1 & 1 \\ -1 & -2 & -3 \\
1 & 4 & 9 \end{array} \right ]   $$
 Determining the matrix
$$ V^T \; = \; W^{-1} \; = \; \left[ \begin{array}{r} v_1^T \\ v_2^T \\ v_3^T
 \end{array} \right ] \; = \; \left[ \begin{array}{rrr} 3 & 2.5 & 0.5 \\
-3  & -4 & -1 \\  1 &  1.5 & 0.5 \end{array} \right ]   $$
 and vectors $\gamma_i =Cw_i, \;  \beta_i^T = v_i^TB$,
$i=1,2,3$
$$ [\gamma_1, \gamma_2, \gamma_3] \; = \; C[w_1, w_2, w_3] = \left[
\begin{array}{rrr} 0 & -1 & 1 \\ -1 & -1 & 0 \end{array} \right ]
\left[ \begin{array}{rrr} 1 & 1 & 1 \\ -1 & -2 & -3 \\ 1 & 4 & 9 \end{array}
 \right ] \; = \; \left[ \begin{array}{rrr} 2 & 6 & 12 \\ 0 & 1 & 2 \end{array}
 \right ], $$
$$ \left[ \begin{array}{r} \beta_1^T \\ \beta_2^T \\ \beta_3^T \end{array} \right ] \;
 = \;\left[ \begin{array}{r} v_1^T \\ v_2^T \\ v_3^T \end{array} \right ] B \;
= \; \left[ \begin{array}{rrr} 3 & 2.5 & 0.5 \\ -3  & -4 & -1 \\
 1 &  1.5 & 0.5 \end{array} \right ] \left[ \begin{array}{rr} -1 & 0 \\
0 & 0 \\ 0 & -1 \end{array} \right ]  \; = \; \left[
\begin{array}{rr} -3 & -0.5 \\ 3 & 1 \\ -1 & -0.5 \end{array}
\right ]  $$
 and using formulas (3.14a),(3.14b) we  calculate
the transfer function matrix $G(s)$ at $s = 1$
$$  G(1) \; = \; \sum_{i=1}^3\frac{\gamma_i\beta_i^T}{1- \lambda_i} \; = \;
[\gamma_1, \gamma_2, \gamma_3] diag(1-\lambda_1,  1-\lambda_2, 1-\lambda_3)^{-1}
\left[ \begin{array}{r} \beta_1^T \\ \beta_2^T \\ \beta_3^T \end{array} \right ]
= $$
$$ = \; \left [ \begin{array}{rrr} 2 & 6 & 12 \\ 0 & 1 & 2 \end{array} \right ]
   \left[ \begin{array}{rrr} 0.5 & 0 & 0 \\ 0 & 0.33 & 0 \\ 0 & 0 & 0.25
\end{array} \right ] \left[ \begin{array}{rr} -3 & -0.5 \\ 3 & 1 \\
-1 & -0.5 \end{array} \right ] \; = \; \left[ \begin{array}{rr} 0 & 0 \\
\frac{1}{2} & \frac{1}{12} \end{array} \right ]$$
Since the rank
of the matrix $G(s)$ is reduced  at $s=1$ then this frequency $s$
is the transmission zero.

Then we find the forced response  $y_2(t)$ to the input signal
$$ u(t) =  \left [ \begin{array}{c} v \\ -6v \end{array} \right ]
\exp(1t),\;\;\;\; v \neq O \eqno(3.15) $$ as follows
$$ y_2(t) \; = \; G(1)u(t) \; = \;\left[ \begin{array}{rr} 0 & 0 \\
\frac{1}{2} & \frac{1}{12} \end{array} \right ] \left [
\begin{array}{r} v \\ -6v \end{array} \right ] \exp(1t)    \;
\equiv \; 0 $$
 and evaluate the free motion component
$$ y_1(t) \; = \; \left [ \begin{array}{rrr} 2 & 6 & 12 \\ 0 & 1 & 2 \end{array}
 \right ] \left[ \begin{array}{ccc} 0.5\exp(-t) & 0 & 0 \\ 0 & 0.33\exp(-2t) &
 0 \\ 0 & 0 & 0.25\exp(-3t) \end{array} \right ] \left[ \begin{array}{rr}
 -3 & -0.5 \\ 3 & 1 \\ -1 & -0.5 \end{array} \right ] \left[ \begin{array}{r}
  v \\ -6v \end{array} \right ] \; = $$
$$= \; \left[ \begin{array}{r} -6\exp(-2t)+ 6\exp(-3t)\\ -6\exp(-2t)+
\exp(-3t) \end{array} \right ]v $$ We can see that $y_1(t) \to O$
as $t \to \infty$.
  Hence, if the initial state condition is zero
$(x_0 = O)$ then the output response  $y(t) = y_1(t) + y(_2(t) \to O$ as $t \to \infty$. We  have been observed
the interesting phenomenon: the growing signal (3.15) is applied to the system input but the output response is
vanished. This phenomenon is stipulated by coincidence of the input signal frequency with the transmission zero.
If  $x_o \neq O$  eliminates free motion component $y_1(t)$ then the output response remains zero even for a
small $t$.

    \section[Definition of invariant zero  by state-space representation]
               {Definition of invariant zero by state-space representation}

     Formerly we have been shown that a transmission zero is defined  via
the transfer function matrix $G(s)$ of the completely controllable and observable system.

     Let's consider an incomplete controllable
and/or observable system  described in the state-space by linear differential equations (1.1),(1.2). We will
study a response to the input signal
$$         u(t) = u_o\exp(\alpha^*t)1(t)
\eqno(3.16) $$
 and demonstrate that if $\alpha^*$ is an invariant
zero then the output response to (3.16) may be zero
$$      y(t)  \; = \; O, \qquad t > 0
\eqno(3.17) $$
 for the following state motion
$$         x(t) = x_o\exp(\alpha^*t)1(t)
\eqno(3.18) $$
 namely, we show that the  invariant  zero  coincides
with a frequency $s = \alpha^*$  at which the transmission of the exponential signal
 $exp(\alpha^*t)$ through the system  is blocked.

At first we show that an invariant  zero associates with reducing a rank of the $(n+l)\times (n+r)$ system
matrix [R1]
 $$  P(s) \; = \; \left[ \begin{array}{cc} sI_n-A & -B \\ C & O
\end{array} \right ]
\eqno(3.19) $$
     Taking  the  Laplace  transform of
(1.1), (1.2) and expressing $\bar{x}(s)$ via  $\bar{u}(s)$ we get
$$
\bar{y}(s) =  C(sI_n - A)^{-1}x_o + C(sI_n - A)^{-1} B\bar{u}(s)
$$
For a proper system the following equalities are follows from $y(t)=O$
$$      y(t_o) = O    \qquad  \rm{or} \it \qquad   Cx_o = O
\eqno(3.20)$$
$$     \bar{y}(s) = O
\eqno (3.21) $$
  Substituting  $\bar{y}(s)$ in (3.21) we obtain the relation
$$
  C(sI_n - A)^{-1}\{ x_o + B\bar{u}(s)) \}   = O
 \eqno(3.22) $$
which using
$$ \bar{u}(s)  = \frac{u_o}{s - \alpha^*}
 \eqno (3.23) $$
can be rewritten for $s \neq \alpha^* $ (3.22) as follows
$$
  C(sI_n - A)^{-1}\{(s - \alpha^*)x_o + Bu_o \} = O
  \eqno(3.24) $$
Applying the obviously identity
$$   (s - \alpha^*)I_n \;=\; (sI_n - A) - (\alpha^*I_n - A)
\eqno (3.25) $$ to (3.24) we obtain  series of equalities
$$ C(sI_n - A)^{-1}\{((sI_n - A) - (\alpha^*I_n - A))x_o + Bu_o \} \;= $$
$$ =\; C(sI_n - A)^{-1}( - (\alpha^*I_n - A)x_o + Bu_o) + C(sI_n - A)^{-1}
     (sI_n - A)x_o \; =$$
$$ =\;C(sI_n - A)^{-1}( - (\alpha^*I_n - A)x_o + Bu_o) + Cx_o \; = \;O  $$
that can be written for $Cx_o=O$ (see (3.20)) as
$$ C(sI_n - A)^{-1}((\alpha^*I_n - A)x_o - Bu_o) \; = \;O
\eqno (3.26) $$
Since the multiplier $(\alpha^*I_n - A)x_o - Bu_o$ is independent
 on $s$ then the following  condition
$$   (\alpha^*I_n - A)x_o - Bu_o  \; = \; O   \eqno(3.27) $$
is the necessary one to  fulfilling (3.26) for any $s$. Uniting (3.20) and (3.27) we obtain the equality
 $$  \left[ \begin{array}{cc} sI_n-A & -B \\ C & O \end{array} \right ]
     \left[ \begin{array}{c} x_o \\ u_o \end{array} \right ] \; = \; O
\eqno(3.28) $$
 that is a necessary  condition for existence of
$\alpha^*$, $x_o$, $u_o$ such that   $y(t) = O$.

     Let's consider relation (3.28) as a linear matrix  equation  in the
vector $[x_o^T, u_o^T]$. It has a nontrivial solution  if a rank of the matrix
 $$  P(\alpha^*) \; = \; \left[ \begin{array}{cc} \alpha^*I_n-A & -B \\ C & O
\end{array} \right ]  $$
is smaller then $min(n+r,n+l)$  (it  is the  necessary and sufficient condition for $l \ge r$ and only the
sufficient one  for $l < r$). Therefore, for $l \ge r$  the column rank reduction of the matrix$P(\alpha ^* )$
ensures existence of nonzero vectors $x_o$ and/or $u_o$ such that the transmission of the signal (3.16) through
the system is blocked: $y(t) = O$. The corresponding nonzero initial state $x_o$ is defined for $s \neq \alpha^*
$  from (3.27) as
$$ x_o \; = \; (\alpha^*I_n - A)^{-1}Bu_o
\eqno (3.29) $$
 Then from  (2.3)  and  (3.23)  we calculate the
state response $\bar{x}$
$$ \bar{x}(s) \; = \; (sI_n - A)^{-1}(x_o + \frac{Bu_o}{s-\alpha^*}) \; =
\; \frac{1}{s-\alpha^*}(sI_n - A)^{-1}((s-\alpha^*)x_o + Bu_o) $$ that can be presented by using  identity
(3.25) as
$$ \bar{x}(s) \; = \; \frac{1}{s-\alpha^*} \{ (sI_n - A)^{-1}(
    (-\alpha^*I_n - A)x_o + Bu_o)  + x_o \} $$
Then  from (3.27) we find  obtain $\bar{x}$
$$ \bar{x}(s) \; = \; \frac{x_o}{s-\alpha^*}
\eqno (3.30)  $$
 that is transformed by the inverse Laplace into the form
$$  x(t) = x_o\exp(\alpha^*t), \; t > 0     $$
Thus, the state vector is the nonzero exponential vector of the  frequency coincided with the input signal
frequency.

     The condition (3.28) is also the sufficient one for existence of
$\alpha^*$, $x_o$ and $u_o$ such that $y(t) \equiv O$ because all steps  of the proof can be reversed [M1]. So,
it has been stated:   to block the transmission of a proportional $exp(\alpha^*t)$ signal  it is necessary and
sufficient that a rank of the matrix $P(s)$ is locally reduced at $s = \alpha^*$.

{\it DEFINITION 3.2.} \rm  A  complex  frequency  $s = \alpha^*$
at which the column rank of $P(s)$ is locally reduced
$$   rankP(s)/_{s=\alpha*} < min(n+r, n+l)  \eqno (3.31) $$
is called as an  invariant zero [M1].

                    {\it EXAMPLE 3.3.} \rm

Let's calculate an   invariant  zero  of  the following system
with $n = 2$, $r =l = 1$
$$ \dot{x} = \left[ \begin{array}{cc} 2 & 0 \\ 1 & 1 \end{array} \right ]x
+ \left[ \begin{array}{c} 1 \\ 0 \end{array} \right ]u, \; y =
\left[ \begin{array}{cc} 1 & 1 \end{array} \right ]x
 \eqno(3.32) $$
This system is controllable  and  observable and has  two poles $\lambda_1 = 1$, $ \lambda_2 = 2 $. Constructing
the $3\times 3$ system matrix (3.19)
$$  P(s) \; = \; \left[
\begin{array}{cc} sI_2-A & -B \\ C & O
\end{array} \right ] \; = \; \left[ \begin{array}{ccc} s-2 & 0 & -1 \\
-1 & s-1 & 0 \\ 1 & 1 & 0  \end{array} \right ] $$
 we  reveal
that  the column (and row) rank of $P(s)$ is locally reduced from $3$ to $2$ at  $s=0$. Hence, $\alpha^* = 0$ is
the invariant zero. From
 equation (3.28) with $\alpha^* = 0$
$$
 \left[ \begin{array}{rrr} -2 & 0 & -1 \\ -1 & -1  &  0 \\ 1 &  1 &  0
\end{array} \right ] \left[\begin{array}{c} x_{10} \\ x_{20} \\ u_o \end{array}
\right ]  = O $$
 we find  the vector $x_0^T =[x_{10}, x_{20},
u_o]$ with  $x_o^T = [1, 1]$,  $u_0 = -2$. So, if the nonzero signal $u(t) = -2\exp^01(t)=2\cdot 1(t)$ is
applied in the input of system (3.32) then the output response is zero for  initial conditions $x_1(0) = 1$,
$x_2(0) = -1$.

{\it REMARK 3.1.} \rm The vector $x_0$ calculating from (3.27) coincides with  (3.13).  Therefore, condition
(3.28) generalizes conditions  (3.8)  and  (3.13)  to   an incomplete controllable and/or  observable system
with the matrix $A$ of a general structure and without the restriction: $\alpha^* \neq \lambda_i$.

{\it ASSERTION  3.1.} \rm  If $\alpha^* \neq \lambda_i(A)$ is a transmission zero then it is  an invariant zero,
the converse is not true.

{\it  PROOF}. \rm $\;$ Let $\alpha^*$ is  a  transmission  zero that satisfies  condition (3.8). If  system
(1.1), (1.2)  is completely controllable   and observable then its TFM is $G(s) = C(sI_n - A)^{-1}B$, otherwise
TFM of (1.1),(1.2) coincides with TFM of the completely controllable  and observable subsystem. Transmission
properties of such system depend on  transmission properties of a completely controllable and observable
subsystem. Therefore,  we will propose that system (1.1),(1.2) is completely controllable and observable and
$G(s) = C(sI_n - A)^{-1}B$.

     Since we study invariant zeros then we must consider the case $l \ge r$.
Let $rankG(s) = min(r,l) = r$. If $\alpha^*$ is a transmission zero then by Definition 3.1 $rankG(\alpha^*) < r$
and hence all minors $det(C^{i_1,i_2,\ldots,i_r}(\alpha^*I_n - A)^{-1}B)$ of the matrix $G(\alpha^*)$ are equal
to zeros
$$
det(C^{i_1,i_2,\ldots,i_r}(\alpha^*I_n - A)^{-1}B) \; = \; 0,
      \qquad  i_k \in \{1,2,\ldots,l \}, \; k=1,2,\ldots,r
\eqno(3.33) $$
 Here  $C^{i_1,i_2,\ldots,i_r}$ denotes a  $r\times
n$ matrix formed from $C$ by deleting all rows except rows $i_1,
i_2,\ldots,i_r$.

Then let us calculate  $(n+r)$ order minors $\tilde{P}$ of the $(n+l)\times (n+r)$  matrix $ P(\alpha^*) \; = \;
\left[
\begin{array}{cc} \alpha^*I_n-A & -B \\ C & O \end{array} \right ]
$
$$det\tilde{P} \; = \; det\left[ \begin{array}{cc}\alpha^*I_n-A & -B \\
C^{i_1,i_2,\ldots,i_r} & O \end{array} \right ] \; = \;
det(\alpha^*I_n - A) det(C^{i_1,i_2,\ldots,i_r}(\alpha^*I_n -
A)^{-1}B) $$
 It follows from (3.33)  that all $det\tilde{P}= 0$. Hence
  $rankP(\alpha^*) <n+r$ and according  Definition 3.2 $\alpha^*$ is an invariant zero.

     The following example illustrates Assertion 3.1.

                  {\it EXAMPLE 3.4.} \rm

 Let us calculate transmission and  invariant  zeros  of the following
controllable and unobservable system with $ n=3,r=l=1$

$$ \dot{x} = \left[ \begin{array}{rrr} 1 & 4 & 0 \\ 0 & -1 & 0 \\
0 & 2 & -3 \end{array} \right ]x + \left[ \begin{array}{c} 0 \\ -1 \\ -1
\end{array} \right ]u, \qquad y = \left[ \begin{array}{ccc} -1 & -1 & 0
\end{array} \right ]x
 \eqno(3.34) $$
 At first we define TFM of this system
$$ G(s) = [ -1 \; -1\; 0]\left[ \begin{array}{ccc} s-1 & -4 & 0 \\0 & s+1 & 0\\
0 & -2 & s+3 \end{array} \right ]^{-1} \left[ \begin{array}{c} 0 \\
-1
\\ -1
\end{array} \right ] = \frac{s+3}{(s-1)(s+1)} $$
The system has two poles $(1, -1)$  and one  transmission  zero $ (-3)$. We need to note that  the dynamics
matrix of (3.34) has three eigenvalues  $1, -1, -3$  but eigenvalue $\lambda_1 =-3$ coincides with the
unobservable pole and it is cancelled in the transfer function matrix $G(s)$.
 To determine invariant zeros we
construct the $4\times4$ system matrix  $P(s)$  (3.19)
$$
P(s) = \left[ \begin{array}{cccc} s-1 & -4 & 0 & 0 \\ 0  & s+1 & 0 & 1 \\
0 & -2 & s+3 & 1\\ -1  &   -1 &   0 & 0  \end{array}\right ] $$
 and observe that its the column rank is locally
reduced at $s = -3$  from $4$  to  $2$.  Therefore, $\alpha^*_1 = -3$, $\alpha^*_2 = -3$ are  invariant zeros.
Moreover,  the first  invariant zero is the transmission zero simultaneously but the second invariant zero is
not the transmission zero and it does not appeared in $G(s)$. This result corresponds to Assertion 3.1 because
the transmission zero of $G(s)$ is simultaneously  the invariant zero but the converse does not held.

     Let's find the initial state $x_0^T =[x_{10}, x_{20},x_{30}]$ and
$u_o$, which assure equality (3.28). From the equation
$P(s|_{s=3})\left[
\begin{array}{c} x_o \\ u_o \end{array} \right ] = O $ or
$$
 \left[ \begin{array}{rrrc} -4 & -4 & 0 & 0  \\ 0 & -2  &  0 & 1 \\
  0 & -2 &  0 & 1 \\ -1 & -1 & 0 & 0 \end{array} \right ]
 \left[\begin{array}{c} x_{10} \\ x_{20} \\ x_{20} \\ u_o \end{array} \right ]
  = O $$
we  calculate $x_o^T =[ 1, -1, \beta ], \;u_o = -2$  where $\beta$
is  any real number.

                    {\it CONCLUSIONS}

     In present section we demonstrate that zeros associate with
transmitting  an exponential signal. Namely,

     1. A transmission zero $\alpha^*$ is defined from the condition (3.8).
This  zero  associates  with the transmission-blocking  [M1] properties of the system. If the condition (3.8) is
carried  out  then the   steady  forced  output  response  to an exponential input of
 frequency $\alpha^*$   is blocked.

     2. An invariant zero $\alpha^*$ is defined from  the condition (3.28).
This  zero  associates  with the zero-output [M1] behavior of the system. If the condition (3.28) is carried out
then  an initial state  $x_o$ and a vector $u_o$   exist such that the whole output response  to an exponential
input of frequency $\alpha^*$ is blocked. The  output of the system at frequency $\alpha^*$ is identically equal
to zero.

\chapter[Determination of transmission zeros
         via TFM]
        {Determination of transmission zeros \\
         via TFM}

\section[Calculation of poles and zeros via Smith-McMillan form]
        {Calculation of poles and zeros via Smith-McMillan form}

We will seek a complete set of transmission zeros as a set of a complex $z_i$  at which the transmission of
steady exponential signals is absent. Let's consider the Laplace transform $\bar{y}(s)$ of the output $y(t)$ of
system (1.1), (1.2)  with $x(t_o) = O$
$$      \bar{y}(s) = G(s)\bar{u}(s)  \eqno(4.1)  $$
where $\bar{u}(s)$  is the Laplace transform of  $u(t)$, $G(s)$ is a matrix of a rank $\rho = min(r,l)$ (the
normal rank).

     At first we consider the case  $l \ge r$. Denoting
 the Smith-McMillan canonical form of the matrix $G(s)$ by $M(s)$
 we obtain from (2.50) $G(s) =
U_L^{-1}(s)M(s) U_R^{-1}(s)$. Therefore, relation (4.1) can be rewritten  as follows
$$   \bar{y}(s) = U_L^{-1}(s)M(s) U_R^{-1}(s)\bar{u}(s) =
U_L^{-1}(s) \left [ \begin{array}{c} M^*(s) \\ O \end{array}
 \right ] U_R^{-1}(s)\bar{u}(s)
\eqno(4.2) $$
 where  $U_L^{-1}(s), U_R^{-1}(s)$ are unimodular
matrices of dimensions $l\times l$ and  $r\times r$ respectively, an $\rho\times \rho$ matrix $M^*(s)$ has the
following form
$$ M^*(s) = diag \left \{ \begin{array}{cccc} \frac{\epsilon_1(s)}{\psi_1(s)},&
\frac{\epsilon_2(s)}{\psi_2(s)},& \ldots,& \frac{\epsilon_{\rho}(s)}
{\psi_{\rho}(s)} \end{array}  \right \}
 \eqno (4.3) $$
We denote  zeros of polynomials $\epsilon_i(s)$ ($i=1,2,\ldots,\rho$) taken  all together by  $z_j$,
$j=1,2,\ldots,\eta $.  As $U_L^{-1}(s)$ and $U_R^{-1}(s)$ are unimodular matrices of a full rank for any $s$
then
$$    rankG(s) \; = \; rankM^*(s)  \eqno(4.4)  $$
It is evident that the rank of $M^*(s)$ is reduced below the normal rank ($\rho $) if and only if a complex
variable $s$ coincides with some of $z_j$, $j=1,2,\ldots,\eta $.  It follows from (4.4)  that the rank of $G(s)$
is also reduced below the normal rank $\rho $ if and only if   $s=z_j$, $j=1,2,\ldots, \eta$ . Therefore, there
exists a nonzero vector $u^*(s)$ such that $G(z_i)u^*(s) = O$. The latest leads to blocking the transmission of
a proportional $exp(z_jt)$ steady signal  at $s=z_j$, $j=1,2,\ldots,\eta$ .

     Then let $l < r$. The Laplace transform of  $y(t)$ is
$$   \bar{y}(s) = U_L^{-1}(s)M(s) U_R^{-1}(s)\bar{u}(s) =
U_L^{-1}(s) \left [ \begin{array}{cc} M^*(s), & O \end{array}
 \right ] U_R^{-1}(s)\bar{u}(s) $$
Thus  equality  (4.4)  is fulfilled  and if the row rank of the matrix $G(s)$ is reduced at $s=z_j$
$(j=1,2,\ldots,\eta )$  then there exists a nonzero vector $u^*$ such that $\bar{y}(z_j) = O$. Hence, for $l<r$
 the rank reducing is only the sufficient  condition to 'block' transmission of an exponential signal at $s=z_j$.

     So, we have been shown: if a complex  frequency $s$  coincides
with  any  zero of the invariant polynomial $\epsilon_i(s)$,
$i=1,2,\ldots, \rho$ then we can find a nonzero vector
$\bar{u}(s)$ for which $\bar{y}(s) = O$. All zeros of invariant
polynomials $\epsilon_i(s)$ ( $i=1,2,\ldots, \rho$)  form a set of
frequencies at which the transmission of steady  exponential
signals may be absent.

\it{DEFINITION 4.1.} \rm [M1] Zeros of polynomials $\epsilon_i(s)$, $j=1,2,\ldots,\rho$,  taken all together,
form the set of  transmission zeros.\footnote{In [M1] these zeros are termed by zeros of TFM.}

\it{DEFINITION 4.2.} \rm [M1] Zeros of polynomials $\psi_i(s)$, $j=1,2,\ldots,\rho$, taken all together, form
the set of poles of the transfer function matrix.

\it{DEFINITION 4.3.} \rm  A polynomial $z(s)$ having transmission zeros as zeros is called as a zero polynomial
of $G(s)$.

     For illustration we consider the transfer function matrix
 (2.51) (see example 2.2). Using the Smith-McMillan
form of $G(s)$ we can calculate two transmission zeros $z_1 = -2$,
$z_2 = -2$. Substituting $s=-2$ into (2.51)  yeilds
$$  G(s)/_{s=-2} \; = \; \left [ \begin{array}{ccc}   -0.5  & 0.5 & 0 \\
0 & 0 & 0 \\  0 & 0 & 0 \\  -0.5 & 0.5 & 0 \end{array} \right ] $$
 Thus, the rank of $G(s)$ is
 reduced from 3 to 1 at $ = -2$. This fact confirms  presence the transmission zero of double multiplicity at
 $z = -2$. The zero polynomial of $G(s)$
is $z(s) =(s+2)^2$. From (2.52) we can find poles of $G(s)$. They are  $ -1,-1,-1,0,-1,-1.$

\section[Transmission zero calculation via minors
         of TFM]
        {Transmission zero calculation via minors
         of TFM}

     Applying  the Smith-McMillan canonical form for calculating
 transmission zeros is rather uncomfortable, especially, for
manual operations.  We consider  the  alternative  method used minors of the matrix $G(s)$ [M3]. It is the
direct method and may be applied  for a system with  a few number of inputs and outputs.

     Frequencies $z_i$ are
transmission zeros if a rank of the matrix $G(z_i)$ is locally reduced ( it is the necessary and sufficient
condition for $l \ge r$ and only the sufficient one for $l < r$). Let the matrix $G(s)$ has the normal rank
$\rho =
 min(r,l)$. We consider all non identically zero minors
 $G(s)^{i_1, i_2,\ldots,i_{\rho}}_{j_1, j_2,\ldots,j_{\rho}}$
 of order $\rho$
of the matrix $G(s) = C(sI -A)^{-1}B$, which are formed from $G(s)$ by  deleting all rows except rows $i_1,
i_2,\ldots,i_{\rho}$ and all columns except columns $j_1, j_2,\ldots,j_{\rho}$ .  It is evident that
$$ G(s)^{i_1, i_2,\ldots,i_{\rho}}_{j_1, j_2,\ldots,j_{\rho}} \; = \;
    det \{ C^{i_1, i_2,\ldots,i_{\rho}} (sI_n - A)^{-1}
     B_{j_1, j_2,\ldots,j_{\rho}}   \} \eqno (4.5) $$
where an $\rho \times n$  matrix $C^{i_1, i_2,\ldots,i_{\rho}}$ is formed from $C$ by deleting all rows except
rows $i_1, i_2,\ldots,i_{\rho}$  and  an $n\times \rho$ matrix $B_{j_1, j_2,\ldots,j_{\rho}}$ is formed from $B$
by the deleting all columns except columns $j_1, j_2,\ldots,j_{\rho}$.

     Let a polynomial $p(s)$ of the degree $k \le n$ is the least
common denominator of these minors. We  add on numerators of minors $G(s)^{i_1, i_2,\ldots,i_{\rho}}_{j_1,
j_2,\ldots,j_{\rho}}$ such a way that they will have the polynomial $p(s)$ as  common denominator.  Resulting
minors become
$$ Z^{i_1, i_2,\ldots,i_{\rho}}_{j_1, j_2,\ldots,j_{\rho}}\frac{1}{p(s)} $$
 where polynomials $Z^{i_1, i_2,\ldots,i_{\rho}}_{j_1,
j_2,\ldots,j_{\rho}}$ are numerators of  new minors. It is evident that the normal rank of $G(s)$ is locally
reduced at $s=z_i$ if all $Z^{i_1, i_2,\ldots,i_{\rho}}_ {j_1, j_2,\ldots,j_{\rho}}$ become equal to zeros at
$s=z_i$. Thus, these numerators must have the  divisor $(s-z_i)$ and we obtain the following definition of
transmission zeros  [M3].

\it{ DEFINITION  4.4.} \rm Transmission zeros are zeros of the
 polynomial $z(s)$ that is a greatest common  divisor  of
numerators $Z^{i_1, i_2,\ldots,i_{\rho}}_ {j_1, j_2,\ldots,j_{\rho}}$  of all non identically zero minors of
$G(s)$ of the order $\rho = min(r,l)$, which are constructed so that these numerators have polynomial $p(s)$ as
the common denominator.

                     \it{ EXAMPLE 4.1.} \rm

Let us calculate  transmission zeros of
 transfer function matrix (2.51). At first we find four minors of
order 3 of the form  $\;$ $ G^{i_1,i_2,i_3}_{1,2,3}, \;
i_1,i_2,i_3 \in \{1,2,3,4\}$
$$    G^{1,2,3}_{1,2,3}\; = \; \frac{s(s+2)^2}{s(s+1)^5}, \qquad
      G^{2,3,4}_{1,2,3}\; = \; \frac{3(s+3)(s+2)^2}{s(s+1)^5}, $$
$$    G^{1,2,4}_{1,2,3}\; = \; \frac{(s+2)}{s(s+1)^5} -
       \frac{(s+3)(s+2)^2}{s(s+1)^5} \; = \; -\frac{(s+2)^2}{s(s+1)^5},$$
$$    G^{1,3,4}_{1,2,3}\; = \; \frac{3(s+3)(s+2)(s^2+2s-1)}{s^2(s+1)^5} -
    \frac{3(s+2)(2s^2+3s-3)}{s^2(s+1)^5}\; = \; \frac{3(s+2)^2}{s(s+1)^5}
$$
The least common denominator of these minors is $    p(s) \; = \; s(s + 1)^5.  $ Adding on  numerators  of above
minors such that they  have polynomial $p(s)$  as the common denominator we obtain   $ Z^{i_1,i_2,i_3}_{1,2,3},
\; i_1,i_2,i_3 \in \{1,2,3,4\}$
$$  Z^{1,2,3}_{1,2,3}\; = \; 3(s+2)^2, \qquad Z^{2,3,4}_{1,2,3}\; =
       \; 3(s+3)(s+2)^2, $$
$$    Z^{1,2,4}_{1,2,3}\; = \; -(s+2)^3, \qquad Z^{1,3,4}_{1,2,3}\; =
       \; s(s+1)(s+2)^2 $$
The greatest common divisor of these numerators is  $(s+2)^2$.
Hence, transmission zeros coincide with  zeros of the polynomial
$z(s) = (s+2)^2$. Similarly result has been obtained above by
using Smith-McMillan form (2.52).

\it{REMARK 4.1.} \rm  If   system  (1.1), (1.2)  is controllable and observable and $r=l$  then $G(s)$  has the
only minor of order $r$, which is equal to $det(C(sI - A)^{-1}B)$. This minor can be represented in the form
$$ det(C(sI - A)^{-1}B) \; = \; \frac{\psi(s)}{det(sI - A)}, \qquad
     \psi(s) \;= \; \frac{det(Cadj(sI - A)B)}{det(sI - A)^{r-1}} $$
where the polynomial $\psi(s)$  is the  numerator of the minor
$det(C(sI - A)^{-1}B)$. Here transmission zeros coincide with
zeros of the polynomial $\psi(s)$.

\section[Calculation of transmission zeros via
         numerator of TFM]
        {Calculation of transmission zeros via
         numerator of TFM}

     To  calculate  transmission zeros we may also use  the factorization
(2.31) or (2.32) of the  transfer function matrix  $G(s)$.  As it has been shown in Sec.2.5.4  polynomials
$\epsilon_i(s)$ $(i=1,2,\ldots, \rho )$ of the Smith-McMillan form of $G(s)$ coincide with invariant polynomials
of the Smith form of any numerator of $G(s)$  (which are polynomial  matrices).  So, we can formulate the
equivalent definition of  transmission zeros [W2].

\it{DEFINITION 4.5.} \rm  Transmission zeros are equal to zeros of polynomials $\epsilon_i(s)$
$(i=1,2,\ldots,\eta )$  of the Smith form of  any  numerator of $G(s)$, taken all together.

     From this definition it follows  the following
procedure to compute transmission zeros:

     (i) factorize $G(s)$ into the product (2.31) or (2.32),

    (ii) find the Smith form (or invariant polynomials) of any
        numerator of the transfer function matrix.

     In the rest of the section we consider  a  simple
method of factorization of  TFM based on the block companion canonical form (Asseo's canonical form).

\subsection[Factorization of transfer function matrix\\
            by using Asseo's canonical form]
           {Factorization of transfer function matrix\\
            by using Asseo's canonical form }

    Using the nonsingular transformation
    $$ \hat{x} = Nx
\eqno(4.6) $$
 where the structure of the $n\times n$ matrix $N$ is defined by formulas
(1.48),(1.49),(1.55) we reduce completely controllable system (1.1), (1.2) to the block companion canonical form
(1.46) with the dynamic matrix $A^*$ (1.20) (where $p=\nu$) and the input matrix $B^*$ (1.47). It follows from
Property 2.1 (see Sec.2.4.1) that
$$ \hat{G}(s) = G(s)
\eqno(4.7) $$
 where $\hat{G}(s)$ is the  transfer  function matrix of  the  canonical system and $G(s)$  is TFM
of (1.1), (1.2).

      To calculate the
matrix $\hat{G} = C^*(sI_n - A^*)^{-1}B^*$  where  $C^* = CN^{-1}$ we use the structure of $A^*$ and $B^*$. At
first we find $(sI_n - A^*)^{-1}$. To this point we partition the  matrix $sI_n - A^*$  as
$$
    sI_n - A^* = \left[ \begin{array}{ccc} sI_{n-r} - A^*_{11} & \vdots &
 -A^*_{12} \\ \dotfill & \vdots &  \dotfill \\ A^*_{21}, A^*_{22},\ldots,
A^*_{2,\nu-1}& \vdots &  sI_r + A^*_{2,\nu}\end{array}  \right]
\eqno (4.8) $$
 where $A^*_{2i}$ are  square $r\times r$
submatrices  ($i=1,\ldots,\nu$) \footnote{For convenience we rename lower blocks of $A^*$ in (1.20) as
$T_p=A^*_{21}, T_{p-1}=A^*_{22},\ldots, T_1=A^*_{2\nu}$.}, matrices $sI_{n-r} - A^*_{11}$ and $A^*_{12}$ of the
dimensions $(n-r)\times (n-r)$ and $(n-r)\times r$ are
$$
sI_{n-r}-A^*_{11} =
\left[ \begin{array}{ccccc} sI_r & -I_r & O & \cdots & O\\ O & sI_r& -I_r&
\cdots & O\\ \vdots & \vdots & \vdots & \ddots & \vdots \\  O& O& O& \cdots &
sI_r  \end{array} \right] , \qquad A^*_{12}= \left[ \begin{array}{c} O\\O\\
\vdots \\ -I_r \end{array} \right] \eqno(4.9) $$
Assuming that
$sI_{n-r} - A^*_{11}$ is the nonsingular matrix $(s \neq 0)$ and
using the formula from [G1] we  calculate
$$
    (sI_n - A^*)^{-1} = \left[ \begin{array}{ccc} X & \vdots &
(sI_{n-r} - A^*_{11})^{-1}A^*_{12}T^{-1} \\ \dotfill & \vdots &  \dotfill \\
X & \vdots &  T^{-1} \end{array}  \right] \eqno (4.10) $$
 where $T
= sI_r + A^*_{2\nu} + [A^*_{21}, A^*_{22},\ldots, A^*_{2,\nu-1}] (sI_{n-r} - A^*_{11})^{-1}A^*_{12}$ and $X$ are
some submatrices. Using the left multiplication of both sides of Eqn.(4.10) by $C^* = CN^{-1}$ and the right
multiplication of that by $B^* = \left[
\begin{array}{c} O \\I_r \end{array}
 \right] $  we obtain
$$  \hat{G}(s) \;= \;C^*(sI_n - A^*)^{-1}B^* \;= $$
$$  =\; C^*\left[ \begin{array}{c}
    (sI_{n-r} - A^*_{11})^{-1}A^*_{12} \\ \dotfill \\ I_r \end{array}  \right]
    \{ sI_r + A^*_{2\nu} + [A^*_{21}, A^*_{22},\ldots, A^*_{2,\nu-1}]
     (sI_{n-r} - A^*_{11})^{-1}A^*_{12} \}^{-1}
   \eqno(4.11)  $$
To calculate $(sI_{n-r} - A^*_{11})^{-1}A^*_{12}$ we find at first
$$
(sI-A^*_{11})^{-1} = \left[ \begin{array}{cccc} s^{-1}I_r &
s^{-2}I_r & \cdots & s^{(1 -\nu )}I_r \\ O & s^{-1}I_r & \cdots &
s^{(2 -\nu)}I_r \\ \vdots &  \vdots & \ddots & \vdots \\ O & O &
\cdots & s^{-1}I_r \end{array} \right] \eqno (4.12) $$
 Thus using (4.9)  we get
$$
  (sI-A^*_{11})^{-1}A^*_{12} = (sI-A^*_{11})^{-1} \left[ \begin{array}{c} O \\
 O \\ \vdots \\ I_r\end{array} \right]  \; = \; \left[ \begin{array}{c}
 s^{(1-\nu)}I_r \\ s^{(2-\nu)}I_r \\ \vdots \\s^{-1}I_r \end{array} \right]
\eqno(4.13)$$
Substituting  (4.13)  in  (4.11) and taking out
$s^{1-\nu}$ gives
$$   \hat{G}(s) = C^* \left[ \begin{array}{c} s^oI_r \\ s^1I_r \\
\vdots \\ s^{\nu-2}I_r \\ s^{\nu-1}I_r \end{array} \right] s^{1-\nu}
\{ sI_r + A^*_{2\nu} + [A^*_{21}, A^*_{22},\ldots, A^*_{2,\nu-1}] s^{1-\nu}
 \left[ \begin{array}{c}  s^oI_r \\ s^1I_r \\ \vdots \\ s^{\nu-2}I_r \\
 s^{\nu-1}I_r \end{array} \right] \}^{-1}
\eqno(4.14)  $$
 Since we have proposed that $s \neq 0$ then the term  $s^{1-\nu}$ may be canceled.
  Multiplying matrices in the right-hand side of
(4.14)  and  partitioning the matrix $C^*$ as $C^* = [C_1,
C_2,\ldots,C_{\nu}]$ where $C_i$ are $ l\times r$ submatrices we
get the following expression for $\hat{G}(s)$
$$ \hat{G}(s)\; = \; (C_1 + C_2s + \cdots + C_{\nu-1}s^{\nu-2} +
C_{\nu}s^{\nu-1})(A^*_{21} + A^*_{22}s + \cdots + A^*_{2,\nu-1}s^{\nu-2} +
A^*_{2,\nu}s^{\nu-1} + I_rs^{\nu})^{-1}
   \eqno(4.15) $$
Denoting
$$  C(s) = C_1 + C_2s + \cdots + C_{\nu-1}s^{\nu-2} +C_{\nu}s^{\nu-1}
\eqno(4.16)  $$
$$       A_2^*(s) = A^*_{21} + A^*_{22}s + \cdots + A^*_{2,\nu-1}s^{\nu-2} +
            A^*_{2,\nu}s^{\nu-1} + I_rs^{\nu}   \eqno(4.17) $$
we present the transfer function matrix $\hat{G}$ as
$$        \bar{G}(s) \; = \; C(s)A^*_2(s)^{-1}   $$
Since $\hat{G}(s) = G(s)$  (see Eqn.(4.7)) then we obtain
$$
            G(s) \; = \;  C(s)A^*_2(s)^{-1}    \eqno(4.18) $$
Thus, it has been shown that a transfer  function  matrix  of the
proper  controllable  system  (1.1),(1.2)  is  factorizated   into
the product of the $l\times r$ matrix polynomial $C(s)$ of the
degree $\nu -1$ and the inverse  of  the $r\times r$  matrix
polynomial $A^*_2(s)$ of the degree $\nu$.

     Using similar way for the observable block companion  canonical
form we can factorizate a transfer function matrix of completely
observable  system (1.1),(1.2) as
$$     G(s) \; = \; N(s)^{-1}Q(s)
   \eqno(4.19)     $$
where $Q(s)$ is an $l\times r$  matrix polynomial of a degree $\alpha -1$  and $N(s)$ is an $l\times l$  matrix
polynomial  of a degree $\alpha$ where $\alpha$ is the  observability index of the pair $(A,C)$.

\it{REMARK} \rm 4.2. The factorization (4.18) takes place for $s=0$. Indeed, partitioning the matrix $A^*$ as
$$     A^* \; = \; \left[ \begin{array}{cc} O & I_{n-r} \\ -A^*_{21} &
 -\tilde{A}^*_{22} \end{array}  \right]
\eqno (4.20) $$
 where $-A^*_{21}$, $-\tilde{A}^*_{22}$  are
$r\times r$ and $r\times (n-r)$ submatrices  respectively and assuming that  $detA^*_{21} \neq 0$  we find
$$
     (A^*)^{-1} = \left[ \begin{array}{ccc}  - (A^*_{21})^{-1}\tilde{A}^*_{22} &
\vdots &  -(A^*_{12})^{-1} \\ \dotfill & \vdots &  \dotfill \\
 I_{n-r} & \vdots & O \end{array}  \right]
\eqno (4.21) $$
 Thus, TFM $\hat{G}(s)$  at  $s=0$ is $
\hat{G}(0) = C^*(-A^*)^{-1}B^* \;=\; C^*(-A^*)^{-1} \left[
\begin{array}{c} O \\I_r \end{array} \right]\; =
 C^*\left[ \begin{array}{c} (A^*_{21})^{-1} \\ O \end{array}
\right]  $.  Partitioning the matrix  $C^* = [C_1^*, C^*_2]$ we obtain
$$  \hat{G}(0)\; = \; C^*_1(A^*_{21})^{-1}
 \eqno(4.22)  $$
The right-hand side of (4.22) coincides with (4.18) at  $s=0$.

                 \it{EXAMPLE 4.2.}  \rm

     For illustration of the method we consider system (1.84) with
     the
output
$$ y \; = \; \left[ \begin{array}{cccc}   1 & 0 & 0 & 0\\ 0 & 1 & 1 & 0
\end{array} \right ] \eqno  (4.23)  $$
The transfer function matrix of this system is
$$ G(s) \; = \; C(sI - A)^{-1}B \; = \; \left[ \begin{array}{cccc}  1 & 0 & 0 & 0\\
 0 & 1 & 1 & 0 \end{array} \right ] \left[ \begin{array}{crrr}
s-2 & -1 & 0 & -1\\ -1 & s & -1 & -1 \\ -1 & -1 & s & 0 \\ 0 & 0 & -1 & s
\end{array} \right ]^{-1} \left[ \begin{array}{cc} 1 & 0 \\ 0 & 0 \\ 0 & 0 \\
0 & 1 \end{array} \right ] \; = $$
$$ = \; \frac{1}{s(s^3-2s^2-2s-1)} \left[ \begin{array}{cc} s^3-s-1 & s^2+s-1 \\
 2s^2+2s+1 & s^2+s+1 \end{array} \right ]
\eqno(4.24)$$
 Calculating  the matrix $C^* = CN^{-1}$  with
$N^{-1}$ from (1.87) gives
$$ C^* \; =  \; \left[ \begin{array}{cccc}  1 & 0 & 0 & 0 \\
 0 & 1 & 1 & 0 \end{array} \right ] \left[ \begin{array}{rrcc}
-1 & -1 & 1 & 0 \\ 1 & 1 & 0 & 0 \\ 1 & 0 & 0 & 0 \\ 0 & 1 & 0 & 1
\end{array} \right ] \; = \;  \left[ \begin{array}{rrcc} -1 & -1 & 1 & 0 \\
2 & 1 & 0 & 0 \end{array} \right ]  $$
 Thus we can find blocks
$C_1$, $C_2$
$$   C_1 = \left[ \begin{array}{rr} -1 & -1 \\ 2 & 1 \end{array} \right ],
 \qquad    C_2 =  \left[ \begin{array}{cc} 1 &  0 \\ 0 & 0 \end{array} \right ]
\eqno (4.25)  $$
 Using expression (1.88) we obtain  $2\times 2$ submatrices
$$   A^*_{21} = - \left[ \begin{array}{rc} -1 & 0 \\ 1 & 0 \end{array} \right ]
\; = \; \left[ \begin{array}{rc} 1 & 0 \\ -1 & 0 \end{array} \right ] $$
$$  A^*_{22} = - \left[ \begin{array}{cr} 3 &  2 \\ 0 & -1 \end{array} \right ]
\; = \; \left[ \begin{array}{rr} -3 & -2 \\ 0 & -1 \end{array}
\right ] \eqno (4.26)  $$
 and get the following matrix
polynomials $C(s)$ and $A^*_2(s)$ in the factorization (4.18)
$$   C(s) = C_1 + C_2s =\left[ \begin{array}{rr} -1 & -1 \\ 2 & 1 \end{array}
   \right ] +  \left[ \begin{array}{cc} 1 &  0 \\ 0 & 0 \end{array} \right ] s
   \; = \; \left[ \begin{array}{cc} s-1 &  -1 \\ 2 & 1 \end{array} \right ] ,$$
$$ A^*_2(s) = A^*_{21} + A^*_{22}s +I_2s^2 = \left[ \begin{array}{rc} 1 & 0 \\
-1 & 0 \end{array} \right ]  + \left[ \begin{array}{rr} -3 & -2 \\ 0 & -1
\end{array} \right ]s + \left[ \begin{array}{cc} 1 & 0 \\ 0 & 1 \end{array}
 \right ]s^2 \; = \;  \left[ \begin{array}{cc} 1-3s+s^2 &  -2s \\ -1 & s^2+s
 \end{array} \right ]
\eqno (4.27)  $$
 The matrix  $G(s)$ becomes
$$
     G(s) \; = \;  C(s)A^*_2(s)^{-1} \; = \; \left[ \begin{array}{cc}
           s-1 &  -1 \\ 2 & 1 \end{array} \right ] \left[ \begin{array}{cc}
          1-3s+s^2 &  -2s \\ -1 & s^2+s  \end{array} \right ]^{-1}
\eqno (4.28)  $$
      For checking we compute directly the product (4.28). Since
$$  A^*_2(s)^{-1} \; = \;  \frac{1}{s(s^3-2s^2-2s-1)} \left[ \begin{array}{cc}
                        s^2+s & 2s \\ 1 & 1-3s+s^2 \end{array} \right ]   $$
then
$$ C(s)A^*_2(s)^{-1} \; = \; \frac{1}{s(s^3-2s^2-2s-1)}\left[ \begin{array}{cc}
    s-1 & -1 \\ 2 & 1 \end{array} \right ] \left[ \begin{array}{cc} s^2 +s & 2s \\
    1 & 1-3s+s^2 \end{array} \right ]  \; = $$
$$   = \;  \frac{1}{s(s^3-2s^2-2s-1)} \left[ \begin{array}{cc} s^3-s-1 &
  s^2+s-1 \\ 2s^2+2s+1 & s^2+s+1
 \end{array} \right ]  $$

  \subsection[Calculation of numerator]{Calculation of numerator}

     Now we find conditions, which ensure that  polynomial matrices
 $C(s)$ and $A^*_2(s)$ in the factorization (4.18)  are  relatively right
prime. Such $C(s)$ is a numerator of  TFM  $G(s)$.

\it{THEOREM 4.1.} \rm  Let $\nu = n/r$ is the controllability index of $(A,B)$. If the pair of matrices $(A,B)$
is completely controllable and the pair of matrices  $(A,C)$ is  completely observable  then  matrices $C(s)$
and $A^*_2(s)$  are relatively  right prime.

\it{PROOF. } \rm  If the pair $(A,B)$ is controllable with the
controllability index $\nu = n/r$ then  system (1.1),(1.2) has the
controllable  block companion canonical form (Asseo's form) and
  matrices $C(s)$ and $A^*_2(s)$ in the factorization (4.18)
have forms (4.16), (4.17) respectively. Let's build the $l\nu
r\times \nu r$ matrix
$$   R(C^*,A^*) \; = \; \left[ \begin{array}{c} C^* \\ C^*A^* \\ \vdots \\
 C^*(A^*)^{\nu r-1} \end{array} \right ]
\eqno (4.29) $$
  where $C^* = CN^{-1} = [C_1,C_2,\ldots,C_{\nu}]$
and  the $\nu r\times \nu r$ matrix $A^*=NAN^{-1}$ has the form (1.20) with $p = \nu$, $T_p = A^*_{21}$,
$T_{p-1} = A^*_{22}$ , $\ldots$, $T_1 = A^*_{2\nu}$. Matrix polynomials $C(s) = C_1 + C_2s + \cdots +
C_{\nu-1}s^{\nu-2} +C_{\nu}s^{\nu-1}$ and $A_2^*(s) = A^*_{21} + A^*_{22}s + \cdots + A^*_{2,\nu-1}s^{\nu-2} +
A^*_{2,\nu}s^{\nu-1} + I_rs^{\nu}$ will be relatively right prime if and only if  the matrix (4.29) has the full
rank that equals to $\nu r$ [20].

     To calculate the rank of $R(C^*,A^*)$ we substitute
$A^* = NAN^{-1}$, $C^* = CN^{-1}$ into the right-hand side of
(4.29) and write series of the equalities
$$   R(C^*,A^*) \; = \; \left[ \begin{array}{c} CN^{-1} \\CN^{-1}NAN^{-1} \\
 CN^{-1}(NAN^{-1})^2 \\ \vdots \\   CN^{-1}(NAN^{-1})^{\nu r-1} \end{array}
  \right ] \;=\; \left[ \begin{array}{c} CN^{-1} \\CAN^{-1} \\
  CA^2N^{-1} \\ \vdots \\   CA^{\nu r-1}N^{-1} \end{array}  \right ] \;=\;
  \left[ \begin{array}{c} C\\CA\\ CA^2 \\ \vdots \\CA^{\nu r-1}\end{array}
   \right ] N^{-1} $$
Since $\nu r = n $ then the matrix $Z_{AC}=[C^T, C^TA^T,\ldots,
C^T(A^{\nu r-1})^T]$ =
 $[C^T, C^TA^T,\ldots, C^T(A^{n-1})^T]$  is the observability
matrix of the pair $(A,C)$. Thus,  $rank R(C^*,A^*)=rank(Z_{AC}N^{-1})=n$. Since $rank (N^{-1}) = n$ and $rank
(Z_{AC}) = n$ (the  pair $(A,C)$ is completely observable) then we obtain
$$  rank R(C^*,A^*) \; = \; n $$
     Therefore,  matrices $C(s)$   and  $A^*_2(s)$  are  relatively
right prime. The theorem is proved.

     Some important corollaries are follows from this theorem.

\it{COROLLARY 4.1.} \rm  If the pair $(A,B)$ is completely
controllable with  $\nu r = n$ and the pair $(A,C)$ is completely
observable then $C(s)$ is a numerator of the transfer function
matrix  $G(s)$.

\it{COROLLARY 4.2.} \rm   Transmission  zeros  of  the completely controllable and observable system (1.1),(1.2)
with $\nu r = n $ are equal to zeros of  invariant polynomials $\epsilon_i(s)$ $(i=1, 2,\ldots,\rho)$ of the
Smith form of the matrix polynomial $C(s) = C_1 + C_2s + \cdots + C_{\nu-1}s^{\nu-2} +C_{\nu}s^{\nu-1}$.

\it{COROLLARY 4.3.} \rm    Transmission  zeros of the completely
controllable and observable system (1.1),(1.2) with $\nu r =n$ and
$r=l$ are equaled to zeros of the polynomial $det C(s)$.

                 {\it EXAMPLE 4.3.}

 We calculate  transmission  zeros  of  the  system from Example
 4.2. This system is completely controllable
 and  observable. Moreover, this system
has the equal number of inputs and outputs. Thus, the  polynomial
matrix $C(s)$ in the factorization (4.28) is  the numerator  of
the transfer function matrix $G(s)$ and, according   Corollary
(4.3),
 transmission zeros are equal to zeros of the polynomial
$$ Z(s) \; = \; detC(s) \; = \; det \left [ \begin{array}{cr} s-1 & -1 \\
   2 & 1  \end{array} \right ] \; = \; s+1  $$
Hence, the system  has  the unique transmission zero that equal to -1.

     For verification we calculate the rank of $G(s)$ (4.24) at  $s =  -1$
$$   rankG(s)/_{s=-1} \; = \;  rank \left [ \begin{array}{rr} -0.5 & -0.5 \\
     0.5 &  0.5 \end{array} \right ] \; = \; 1  $$
As the rank of the transfer function matrix  $G(s)$ is  locally
reduced from 2 to 1 at $s=-1$ then $s = -1$ is the transmission
zero.

\chapter[Zero definition via system matrix]
        {Zero definition via system matrix}

   In present chapter we will based on the definition of zeros via
$(n+l )\times (n+r)$  Rosenbrock's system matrix
 $$  P(s) \; = \; \left[ \begin{array}{cc} sI_n-A & -B \\ C & O
\end{array} \right ]
\eqno(5.1) $$
 of the  normal rank  $n + min(r,l)$.

\it{DEFINITION 5.1.} \rm A complex frequency $s=z$ at which the
normal rank of the matrix $P(s)$ is reduced
$$   rankP(s)/_{s=z} < n + min (r,l)
  \eqno (5.2) $$
is named as a system zero of system (1.1), (1.2).

    First  definitions of system zeros have been introduced by prof. Rosenbrock
[R2], [R3]. The definition of system zeros as zeros of invariant polynomials of the Smith form of $P(s)$ was
introduced in  1973 [R2]. More recently these zeros were refereed to as invariant zeros. The complete set of
system zeros in terms of special formed minors of $P(s)$  have been introduced in 1974 [R3]. Later we consider
these  notions more detail. We will consider also other types of zeros defined in terms of the matrix $P(s)$.

\section[Complete set of invariant zeros]
           {Complete set of invariant zeros}

     In Section 3.4 we already have presented the definition
of a invariant zero for a system with $l \ge r $ where $r$ and $l$ are numbers of inputs and outputs
respectively. The invariant zero has been defined as a complex frequency that reduces a rank of the matrix
$P(s)$.
     The complete set of invariant zeros consists of the complete set
of complex frequencies $z_i$, $i=1,2, \ldots$  for which the rank inequality (5.2) is fulfilled. To find these
frequencies we will seek the Smith  form $S(s)$ of  the  matrix $P(s)$.  As it was shown in Section 2.5.2, a
polynomial $(n+l)\times (n+r)$ matrix with the normal rank  $n + min (r,l)$ has  the matrix $S(s)$ of the
following structure
$$   S(s) \; = \; U_L(s)P(s)U_R(s) \; = \; \left \{ \begin{array}{cc} \left [
\begin{array}{c} diag(s_1(s),s_2(s),\ldots,s_{n+r}(s)) \\
\dotfill\\ O \end{array} \right ], & l\ge r \\  \\ \protect [diag(s_1(s),s_2(s), \ldots,s_{n+l}(s)),O], & l\le r
\end{array} \right. \eqno(5.3) $$
 where $U_L(s)$ and $U_R(s)$  are unimodular  matrices  of dimensions
$(n+l)\times(n+l)$ and $(n+r)\times (n+r)$ respectively, $s_i(s)$
 are  invariant polynomials of $P(s)$. We present $P(s)$ as
$$   P(s) \; = \; U_L(s)^{-1}S(s)U_R(s)^{-1} \; = \; \left \{ \begin{array}{cc}
 U_L(s)^{-1}\left [ \begin{array}{c} diag(s_1(s),s_2(s),\ldots,s_{n+r}(s)) \\
\dotfill\\ O \end{array} \right ]U_R(s)^{-1}, & l\ge r \\  \\
\protect U_L(s)^{-1}[diag(s_1(s),s_2(s),
\ldots,s_{n+l}(s)),O]U_R(s)^{-1}, & l\le r \end{array} \right.
\eqno(5.4) $$
 Since matrices $U_L(s)$ and $U_R(s)$  are the
unimodular ones with constant  determinants then  inverse matrices $U_L(s)^{-1}$ and $U_R(s)^{-1}$ have similar
properties and a complex $s$ is an invariant zero if and only if  it is a zero of any polynomials $s_i(s)$,
$i=1,2\ldots,n+min(r,l)$. If $l\ge r$ the column rank of $P(s)$ is reduced, if $l\le r$ the row rank of $P(s)$
is reduced.

      Thus,  we  can  define  the  complete  set  of
invariant zeros as follows

\it{DEFINITION 5.1. } \rm Zeros of all invariant polynomials
$s_i(s)$, $i=1,\ldots,n+min(r,l)$, taken all together, form the
complete set of  invariant zeros.

\it{REMARK 5.1.} \rm  Davison and Wang in 1974 [D4] were
 defined an invariant zero via the inequality (5.2) and the
'complete set' of invariant zeros of a completely controllable and
observable system as zeros  of  the  highest  order  invariant
polynomial ( i.e. the polynomial $s_{n+\sigma}(s)$, $\sigma =
min(r,l)$).
 These zeros were named as 'transmission zeros' [D4] . It is evident that the
invariant zeros of Davison and Wang form a  subset of the
Rosenbrock's ones.

\it{REMARK 5.2.} \rm  It follows from  Definition 5.1 and Eqn.(2.37) that the complete set of invariant zeros
 coincides with  zeros of the monic largest common
divisor $\psi_I(s)$ of all   $n + min(r,l)$ order minors (non identically zero) of the matrix $P(s)$  of the
normal rank $n + min(r,l)$. The polynomial $\psi_I(s)$ is
$$  \psi_I(s)\; = \; s_1(s)s_2(s)\cdots,s_{n+\sigma}(s),\qquad \sigma = min(r,l)
\eqno(5.5) $$
     This remark may be used for the manual calculating invariant
zeros.

                      \it{EXAMPLE 5.1.} \rm

We consider the  system of the form
$$ \dot{x} = \left[ \begin{array}{rrr} 1 & 0 & 0 \\ 0 & -1 & 0 \\
0 & 0 & -3 \end{array} \right ]x + \left[ \begin{array}{c} 0 \\ -1 \\ -1
\end{array} \right ]u, \qquad y = \left[ \begin{array}{crc} 1 & -1 & 0 \\
0 & 2 & 0 \end{array} \right ]x
 \eqno(5.6) $$
To find the complete set of invariant zeros we construct the
matrix
$$ P(s) \; = \; \left [ \begin{array}{cccc}   s-1 & 0 & 0 & 0 \\
 0 & s+1 & 0 & 1 \\ 0 & 0 & s+3 & 1 \\  1 & -1 & 0 & 0 \\ 0 & 2& 0 & 0
 \end{array} \right ] \eqno(5.7) $$
and determine minors $P_i$ of the order  $n + min(r,l) = 3 + 1 = 4$ by deleting the row $i$ ($i=1,2,3,4,5$). As
a result we get
$$ P_1 = -2(s+3),\;\; P_2 = 0, \;\;P_3 = 0,\;\; P_4 =-2(s-1)(s+3),\;\;
 P_5 =(s-1)(s+3) $$
The monic largest common divisor of non identically zero minors $P_i$ is equal to  $s+3$. Therefore,  $\psi_I(s)
= s+3$ and the system (5.6) has the only invariant zero  $z = -3$.

\section[Complete set of system zeros]
        {Complete set of system zeros}

     Analysis of the matrix $[sI_n - A,B]$  with $A$  and $B$
from Example 5.1 reveals that  the system (5.6) is uncontrollable at   $s = 1$  because
$$      rank [sI_n - A,B]/_{s=1} \; = \; rank \left [ \begin{array}{cccc}
         0 & 0 & 0 & 0 \\  0 & 2 & 0 & 1 \\ 0 & 0 & 4 & 1 \end{array} \right ]
         \; = \; 2 < 3  \eqno(5.8) $$
Therefore, the  signal that proportional to $exp(1t)$  does not appear in the input of system (5.6) or in the
output of the dual system. Thus, the set of  invariant zeros do not include all frequencies for which signal
transmitting through the system is 'blocked'. The complex variable $s = 1$ is  called as 'decoupling zero'.
These zeros  will study later in Section 5.3.

     Now we study system zeros that form  a complete set of frequencies,
which are not propagated through a system. This definition of system zeros is based on  minors  of  $P(s)$
having the special form. Let's consider all $n+k$ order minors of the matrix $P(s)$ constructing  by deleting
all rows of $P(s)$ expect rows $1,2,\ldots,n,n+i_1,\ldots, n+ i_k$ and all columns  expect columns
$1,2,\ldots,n, n+j_1,\ldots,n+j_k$ where $i_k\in \{1,2,\ldots,r\}, \; j_k \in\{1,2,\ldots,l\}$. We denote these
minors as
$$  P(s)^{1,2,\ldots,n,n+i_1,\ldots, n+ i_k}_{1,2,\ldots,n, n+j_1,\ldots,n+j_k}
\eqno (5.9) $$
 where the integer $k$ varies from  1  to  $min(r,l)$. Let  $\delta \; (0 \le \delta \le
min(r,l))$ is a maximal value $k$ such that at least the only minor (5.9) of the order $\rho = n+\delta$  does
not identically zero. We denote this minor by $P_{\rho}(s)$. Let's suppose  that we obtain a few such minors and
the polynomial $\psi(s)$ is the greatest common divisor of these minors (if we obtain the only  minor then
$\psi(s) = P_{\rho}(s)$).

\it{DEFINITION  5.2.} \rm  The complete set of  system zeros
coincides with zeros of the polynomial $\psi(s)$ that is  the
greatest common divisor of non identically zero minors (5.9) of
the maximal order $\rho$.

\it{EXAMPLE  5.2.} \rm  To calculate  system zeros of the system (5.6) we find  two minors of the structure
(5.9) of the matrix (5.7): the first one with $i_1 = 1,j_1 = 1$ and  the second one with $i_1 = 2, j_1 = 1$

$$    P(s)^{1,2,3,4}_{1,2,3,4} = (s-1)(s+3), \qquad
      P(s)^{1,2,3,5}_{1,2,3,4} = -2(s-1)(s+3)    \eqno(5.10)  $$
The greatest common divisor of these minors  is  $\psi(s) = (s-1)(s+3)$. Therefore, the system (5.6) has two
system zeros: $ 1, -3$.

     It follows from  Examples 5.1 and 5.2  that invariant zeros
are a subset of system zeros. We prove this property in the general case

\it{ASSERTION  5.1.} \rm  The set of invariant zeros  is  a subset
of  system zeros.

\it{ PROOF}. $\;$ \rm  Without loss   of generality we  may assume that $l \ge r$. Then the normal rank of
matrix $P(s)$ is equal to $n + min(r,l) = n + r$. Consider all minors (5.9) of $P(s)$  of the order $n+r$  with
$i_k = r$
$$  P(s)^{1,2,\ldots,n,n+i_1,\ldots, n+ i_k}_{1,2,\ldots,n, n+j_1,\ldots,n+j_k},
        \qquad i_m \in\{1,2,\ldots,l\}, m = 1,3,\ldots,r $$
We denote these minors by  $Q_i(s), i=1,2,\ldots,n_c$  where the number of  minors ($n_c$) is  calculated as
follows
$$     n_c = \frac{l(l-1)\cdots (l-r+1)}{12\cdots r}    $$
As the normal rank of $P(s)$  is equal to  $n+r$  then there exists at least the only minor that is non
identically zero. Let $\psi_c(s)$ is the greatest common divisor of  nonzero minors $Q_i(s)$. By Definition 5.2
zeros of $\psi_c(s)$ form the set of
 system zeros.

     Then we consider all possible minors of  the  matrix $P(s)$
 of  order  $n+r$,
which  are  constructed from  $P(s)$ by  deleting superfluously
$l-r$ rows from  $n+l$   ones. We denote these minors by
$\beta_i(s), i=1,2,\ldots,n_I$ where $n_I$ is the  number of these
minors
$$     n_I = \frac{(n+l)(n+l-1)\cdots (l-r+1)}{12\cdots (n+r)}    $$
     It is clear that  $n_I \ge n_c$ and the set of $Q_i(s)$ is  a
subset of the set of $\beta_i(s)$. Let $\psi_I(s)$ is the greatest
common divisor of nonzero $\beta_i(s)$.  By Remark 5.2  zeros of
the polynomial $\psi_I(s)$ form the set of  invariant zeros.

     Since the set of  $Q_i(s)$  is the subset of
$\beta_i(s)$  then we can write the equality
$$       \psi_c(s) =  \psi_I(s)\psi_k(s)  $$
where $\psi_k(s)$ is non  identically  zero  polynomial.
Therefore, the degree of $\psi_I(s)$  is not greater than that of
 $\psi_c(s)$. The assertion is proved.

     We can illustrate this result using Examples 5.1, 5.2.
 Indeed system (5.6) has  $n_c = 2$
minors $Q_i(s)$ (see Example 5.1):
$$  Q_1(s) = (s-1)(s+3), \qquad  Q_2(s) = -2(s-1)(s+3)   $$
and $n_I = 5$ minors $\beta_i(s)$ (see Example 5.2):
$$ \beta_1 = -2(s+3), \; \beta_1(s) = \beta_1(s) = 0, \;\beta_4 = -2(s-1)(s+3),
   \; \beta_5 =(s-1)(s+3)    $$
The monic greatest common divisor of the  minors $Q_i (s)$ is
$\psi_c(s) = (s-1)(s+3)$. The monic greatest common divisor of the
nonzero minors $\beta_i$ is $\psi_I(s) = s+3$. It is evident that
the set of invariant zeros $\{-3\}$ is the  subset of the set of
system zeros $\{1,-3\}$.

\it{ ASSERTION 5.2.} \rm  If $l=r$  and $\rho = n + l = n +r$ then
sets of invariant zeros and system zeros  coincide.

     The proof follows from the  structure of  minors $Q_i(s)$
and $\beta_i(s)$. These minors are equal to the only non
identically zero minor of order $\rho = n + r = n + l$ .

     Let's find the minor of form (5.9) when  $r = l$ .
Using formula for the determinant  of a block matrix [G1] we get
$$ detP(s)\; =\; det \left[ \begin{array}{cc} sI_n-A & -B \\ C & O
   \end{array} \right ] \; = \; det(sI_n -A)det(C(sI_n -A)^{-1}B) \; = \;
    det(sI_n -A)detG(s)   $$
Thus, the complete set of  system zeros of a  system with equal number of inputs/outputs  coincides with  zeros
of the polynomial
$$ \psi(s) \;=\; det(sI_n -A)detG(s)   \eqno(5.11) $$
Such definition of system zeros was introduced in [K5].

\section[Decoupling zeros]
        {Decoupling zeros}

     In the general case  sets of invariant zeros and system zeros
 are
distinguished  by    presence  of   decoupling  zeros.  We observe this fact in Examples 5.1 and 5.2: the
complete set of system zeros contains the zero $z = 1$ that does not the invariant zero. This zero coincides
with the frequency at which the system is uncontrollable and the proportional $exp(1t)$ signal  does not appear
in the input. Such a zero is the decoupling zero. Now we study these  zeros.

     Let's study the structure of the matrix  $P(s)$  (5.1). If
system (1.1),(1.2) is unobservable or/and uncontrollable then
 there exists a complex variable $ s=z$  at which the normal rank of
the block column $ \left [ \begin{array}{c} sI_n-A \\C \end{array} \right ] $
 or/and the block row  $[sI_n-A,B]$ is locally reduced. This
complex variable is named as a decoupling zero. These  zeros have been introduced by  Rosenbrock  in  1970 [R1].
They associate with complex frequencies (modes, eigenvalues of A), which are decoupled from the input/output.

     \it{OUTPUT DECOUPLING ZEROS}. \rm They  appear  when  several  free  modal
(exponential type) motions of the system state  $x(t)$  are decoupled from the  output. Let's consider this
situation in detail for the matrix $A$ having $n$ distinct eigenvalues $\lambda_1,\ldots,\lambda_n$. Without
loss of generality we can assume that the forced response is absent ( $u(t) = 0$). Then using  (1.8)  we expand
the solution $x(t)$ of  linear time-invariant differential equation (1.1) with $u(t) = O,\; t_o =0,\; x(t_o) =
x_o \neq O $ as follows
$$        x(t) = \sum_{i=1}^n{w_i\exp(\lambda_it)\xi_i}, \qquad
          \xi_i = v_i^Tx_o     $$
 where $w_i$, $v_i^T$  are  right  and  left  eigenvectors of
the matrix $A$, $\xi_i$  is a nonzero scalar\footnote{ We consider nontrivial case when all  $\xi_i = v_i^Tx_o
\neq 0$ and the free modal motion of $x(t)$ has all modes.}.  We can express the output of the  system as
$$  y(t) \; = \; Cx(t) \; = \; C\sum_{i=1}^n{w_i\exp(\lambda_it)\xi_i} \; =
    \; C\sum_{i=1}^n{x_i(t)}  $$
where  $x_i(t)=w_iexp(\lambda_i t)\xi_i$ are modal components of the state $x(t)$. If certain  modal component
$x_l(t) = x_l = w_l\exp(\lambda_lt)\xi_l, \; l \in \{1,\ldots,n\}$ is decoupled with the output then the
following condition is satisfied
$$ Cx_l\; = \; C w_l\exp(\lambda_lt)\xi_l \; = O
\eqno (5.12)$$
 As $ \exp(\lambda_lt) \neq 0 $, $\xi_l \neq 0$ then
it follows from (5.12)  that $Cw_l = O $. Adding the last
expression to the following: $Aw_l =w_l\lambda_l$  or
$(\lambda_lI_n - A)w_l = O$ we obtain
$$ \left [ \begin{array}{c} \lambda_lI_n-A \\C \end{array} \right ]w_l
 \; = \; O  \eqno(5.13) $$
Considering (5.13) as an linear homogeneous equation in $w_l$ we recall that a nontrivial solution of (5.13)
exists if a column rank of the $(n+l)\times n$ block matrix in (5.13) is locally reduced below $n$. The
appropriate value of the complex frequency $(\lambda_l)$ is called as the output decoupling zero.

\it{DEFINITION 5.3.} \rm  Output decoupling zeros are formed by
the set of complex variables  $s$ at which the normal column rank
of the matrix
$$ P_o(s) \; = \; \left [ \begin{array}{c} sI_n-A \\C \end{array} \right ]
\eqno (5.14)$$
is reduced.

     The output decoupling zeros are calculated as   zeros  of
 invariant polynomials of $P_o(s)$.

\it{INPUT DECOUPLING ZEROS}. \rm They  appear  when  certain
 free modal (exponential type) motions of the state $x(t)$  are
decoupled from the  input. Considering the dual system we may show that there exist an eigenvalue $\lambda_l$
and a left eigenvector $v_l$ of the matrix $A$ such that  the following equalities take place
$$    v_l^T(\lambda_lI_n - A) = O, \qquad v_l^TB = O   $$
 Uniting these  equalities yields
$$    v_l^T[\lambda_lI_n - A, B]  = O  \eqno  (5.15) $$
Considering (5.15) as an equation in the vector $v_l^T$ we conclude that this equation has a nontrivial solution
in $v_l^T$ if the row rank of the $n\times(n+r)$ matrix in (5.15) is locally reduced below $n$. The appropriate
value of  the complex frequency $\lambda_l$ is called as an input decoupling zero.

\it{DEFINITION  5.4.} \rm Input decoupling zeros are formed by the
set of  complex variables $s$ at which the normal row  rank of the
matrix
$$         P_i(s) = [sI_n - A,B]
\eqno(5.16) $$
is reduced.

     The input decoupling zeros are calculated as   zeros  of
invariant polynomials of $P_i(s)$.

     In Section 2.4.2  we have introduced notions of uncontrollable and
unobservable  poles, which coincide with  eigenvalues of $A$ reducing the normal rank of  matrices  $[sI_n -
A,B]$ and $[sI_n - A^T,C^T]$.  Now we show  that  these poles are equal to decoupling zeros.

\it{THEOREM  5.1.} \rm  Output  decoupling zeros and input
decoupling zeros of system (1.1), (1.2) coincide with unobservable
and  uncontrollable poles of this system respectively.

\it{PROOF}. \rm Let $z$ is an output decoupling zero. Then the column rank of $P_o(s)$ at $s=z$ is locally
reduced below  $n$. We need to show that the complex variable $z$  coincides with an eigenvalue of the matrix
$A$. Indeed, if rank $P_o(s)/_{s=z} < n$ then there exists a nontrivial vector $f$ such as $P_o(z)f = O$. Using
the structure of $P_o(s)$ (5.14) we  write the following equations
$$         (zI_n - A)f = O      \eqno(5.17)   $$
$$         Cf = O     \eqno(5.18)  $$
It follows from (5.17)  that $z$ is the eigenvalue of $A$  and $f$
is the corresponding eigenvector. Thus we immediately obtain from
Assertion 2.5  that $z$ is the unobservable pole of (1.1), (1.2).

     A similar way may be used for the second part of the theorem.

\it{REMARK 5.1.} \rm  Uncontrollable (unobservable) poles
 are sometimes refereed as decoupling poles.

\it{ NUMBER OF DECOUPLING ZEROS}. \rm The following relations can be get from Assertions 2.2, 2.6 and Theorem
5.1:

      1. The number of input decoupling zeros is equal to the rank
deficient of the controllability matrix $Y_{AB} = [B,AB, \ldots, A^{n-1}B]$.

      2. The number of output decoupling zeros is equal to the rank
deficient of the observability matrix $Z_{CA}^T = [C^T,A^TC^T, \ldots,
(A^T)^{n-1}C^T]$.

\it{INPUT-OUTPUT DECOUPLING ZEROS}. \rm They appear when there
exist $\lambda_l$, $w_l$ and $v^T_l$ such as two equalities
(5.13), (5.15) are held simultaneously. Such $\lambda_l$ is named
as an input-output decoupling zero.

    \it{EXAMPLE 5.3.} \rm

We consider the system (5.6). If $s = 1 $ then the rank of the matrix
$$  P_i(s) \; = \; [sI_n - A,B] \; = \; \left [ \begin{array}{cccr} s-1 & 0 & 0 & 0 \\
     0 &  s+1 &  0 & -1 \\ 0 & 0 & s+3 & -1 \end{array} \right ] $$
is  reduced below $n = 3$
$$               rankP_i(s)/_{s=1} = 2 < 3   $$
 Hence, $z = 1$ is the input decoupling zero.
     Let's find a number of decoupling zeros of this system.
The rank deficient of the controllable matrix
$$    Y_{AB} = [B,AB,A^2B] = \left [ \begin{array}{rrr}  0 & 0 & 0 \\
      -1 &  1 & -1 \\ -1 &  3 & -9 \end{array} \right ] $$
 is  equal  to  1.
So, the system (5.6) has the only input decoupling zero.

     Then we find the rank deficient of the observability matrix
$$   Z_{CA}^T = [C^T,A^TC^T,(A^T)^2C^T] = \left [ \begin{array}{rrrrrr}
     1 &  0 &  1 &  0 &  0 &  0 \\ -1 &  2 &  1 & -2 & -1 &  4 \\
      0 &  0 &  0 &  0 &  0 &  0 \end{array} \right ] $$
It is  equal  to  1  then the system  (5.6) has the only output decoupling zero. To find this zero we construct
the matrix $P_o(s)$
$$ P_o(s) = \left [ \begin{array}{c} sI_n-A \\C \end{array} \right ] \; = \;
\left [ \begin{array}{ccc} s-1 &  0 &  0\\ 0 & s+1 &  0 \\ 0 & 0 & s+3 \\
 1 & -1 &  0 \\ 0 & 2 & 0  \end{array} \right ] $$
and discover that the column rank of $P_o(s)$  is reduced below $n= 3$  at $s=-3$ . Hence,  $z = -3$ is the
output decoupling zero.
     Thus,  the system  (5.6)  has  the   input
decoupling zero $z = 1$  and the output decoupling zero $z = -3$.

\section[Relationship between different zeros]
         {Relationship between different zeros}

     At first we introduce the following notations:

$$  \begin{array}{lcl}   \{ n \} &  - & \rm{a\; set\; of\; system\; zeros}, \\
    \{ i \} & - & \rm{a\; set\; of\; invariant \;zeros}, \\
      \{ p \}  &  - &  \rm{a\; set\; of\; transmission \;zeros},\\
     \{ i.d. \} & - & \rm{a\; set \;of \;input\; decoupling\; zeros}, \\
     \{ o.d. \} & - & \rm{a\; set\; of\; output\; decoupling\; zeros},  \\
      \{ i.o.d. \} & - & \rm{a\; set\; of\; input-output\; decoupling\; zeros}.
\end{array} \eqno(5.19)   $$

\subsection[Transmission and invariant zeros]
           {Transmission and invariant zeros}

     It has been shown that the invariant zeros are associated  with
 reducing a column or row rank of the matrix  $P(s)$.

     Let $l \ge r$. The normal rank of the $(n+l)\times(n+r)$ matrix $P(s)$
is not changed after the right multiplication of $P(s)$ by the nonsingular unimodular $(n+r)\times (n+r)$ matrix
 $$  L_1(s) \; = \; \left[ \begin{array}{cc} I_n & (sI_n-A)^{-1}B \\ O& I_r
\end{array} \right ]
\eqno(5.20) $$
 As the determinant of $L_1(s)$ does not depend on
$s$ then the following rank equalities are satisfied
$$  rankP(s) \; = \; rank\{P(s)L_1(s)\} \; = \; rank \left[ \begin{array}{cc}
sI_n-A & O \\ C & C(sI_n-A)^{-1}B \end{array} \right ] \; = $$
$$= rank \left[ \begin{array}{cc} sI_n-A & O \\ C & G(s) \end{array}
\right ] \eqno (5.21) $$
 Hence, a column rank of $P(s)$ is
depended  on the rank of $G(s)$: if the rank of $G(s)$ is reduced then the column rank of $P(s)$ is also
reduced.

     A similar way may be used for  $l\ge r$.  The normal rank of
$P(s)$ is not changed  after the left  multiplication of  $P(s)$  by the unimodular $(n+l)\times(n+l)$ matrix
 $$  L_2(s) \; = \; \left[ \begin{array}{cc} I_n & O \\ -C(sI_n-A)^{-1} & I_l
\end{array} \right ]
\eqno(5.22) $$
Thus
$$  rankP(s) \; = \; rank\{L_2(s)P(s)\} \; = \; rank \left[ \begin{array}{cc}
 sI_n-A & -B \\ O & G(s) \end{array} \right ]
\eqno (5.23) $$
     and  if the rank of $G(s)$ is reduced then the row rank  of $P(s)$
is also reduced.

     We result in that the set of transmission zeros (defined
via $G(s)$ ) is the subset of the set of invariant zeros (defined via $P(s)$). Using notations (5.19) we summary
this result as the inclusion
$$                  \{ p \} \subseteq \{ i \}    \eqno  (5.24) $$
The similarly result follows from Assertion 3.1.

\subsection[Invariant, transmission and decoupling zeros]
           {Invariant, transmission and decoupling zeros}
     Let   system (1.1), (1.2) with  $l>r$  possesses  invariant
zeros. If a complex variable $s=z$  coincides  with  an  invariant zero then the column rank of $P(s)$ is
reduced. Then there exists a nontrivial solution of equation (3.28) with respect to the vector $[x_o^T,u_o^T]$.
Equation (3.28) can be rewritten for $u_o = O$ as
$$        P(z)\left [ \begin{array}{c} x_o \\ O \end{array} \right ] = O $$
or in the equivalent form
$$    \left [ \begin{array}{c} zI_n - A \\ C \end{array} \right ]x_o = O $$
  This expression corresponds  to the condition of unobservability and
the  complex  variable $z$  coincides  with  an  unobservable pole that is equal to an  output decoupling zero.
Hence, $z$ is as well the invariant zero as the output decoupling zero. We conclude that
    if a system has more outputs  than inputs then several
invariant zeros may be output decoupling zeros simultaneously.

    The dual situation may take place for a system with  $l<r$,
when several zeros are as well invariant  zeros as  input
decoupling zeros.

     Let's find conditions when  decoupling  zeros  are simultaneously
invariant zeros. Rank equalities (5.21),(5.23) demonstrate that invariant zeros, defined via the matrix  $P(s)$,
contain transmission zeros, defined via the matrix $G(s)$, and decoupling zeros, defined  either via the matrix
$P_o(s)$ ($l>r$) or via the matrix $P_i(s)$ ($l<r$).  Let  for $l>r$  the matrix $G(s)$ has a full column rank.
Then if the rank of  $P_o(s)$ reduces then the rank of  $P(s)$ also reduces. Similarly, for $l<r$ if the matrix
$G(s)$ has a full row rank then reducing the rank of $P_i(s)$ involves reducing the rank  of $P(s)$. Therefore,
the following assertion is held.

\it{ASSERTION 5.3.} \rm  Let system (1.1),(1.2) has the matrix
$G(s)$ of the full rank. Then if this system has more outputs than
inputs ($l>r$) then every output decoupling zero is  an invariant
zero, i.e.
$$         \{ o.d. \}  \subset   \{ i \} \eqno(5.25)   $$
If the system has more inputs than outputs ($l<r$ ) then every
input decoupling zero is  an invariant zero, i.e.
$$                   \{ i.d. \} \subset   \{ i \}    \eqno(5.26) $$

     We may show also that the
controllability/observability properties are closely connected to
the structure of $\{ i \}$.
      Indeed, if system (1.1),(1.2) with $r>l$ is uncontrollable then
the rank of $P_i(s)$ is reduced. Let $rankP_i(s) = n-q$ then the
system has  $q$  uncontrollable poles and  the set $\{ i \}$
differs from the set $\{ p \}$ by  existence of $q$ input
decoupling zeros. From (5.24), (5.26) we obtain the following
inclusion
$$                  \{ i \} \supseteq \{ p \} + \{ i.d. \}
 \eqno(5.27) $$
Similarly, if system (1.1),(1.2) with $l>r$ is unobservable then
the rank of $P_o(s)$ is reduced. Let $rankP_o(s) = n-q$, then the
system has $q$  unobservable poles and  the following inclusion
takes place
$$       \{ i \} \supseteq  \{ p \} + \{ o.d. \}
 \eqno(5.28)  $$

     If system  (1.1),(1.2)  is  completely  controllable  and
observable then $rankP_o(s) = rankP_i(s) =n$ for any  $s$. Hence, the normal rank of $P(s)$ is reduced if and
only if the normal rank of $G(s)$ is reduced. In this case  sets of  invariant  and transmission zeros coincide
$$                  \{ i \} \equiv \{ p \}  \eqno   (5.29) $$

                \it { CONCLUSION}   \rm

     In the general case the set $\{ i \}$ differs from the set $\{ p \}$
by  existence of decoupling zeros. Inclusions  (5.27)  and (5.28)
represent the rough structure of the set $\{ i \}$. More exact
relations have been obtained in works [P6], [R3]. In [R3] it has
been shown that if $r=l$  then
$$            \{ i \} = \{ p \} + \{ o.d. \} + \{ i.d. \} - \{ i.o.d.\} $$
 if $l>r$  then
$$    \{ i \} = \{ p \} + \{ o.d.\} + \{ i.d.\} - \rm{some\; terms\; of\;}
        \it \{i.d.\}   $$
The calculation of those $\{ i.d. \}$ that is a part of $\{ i \}$ is represented in [P6].

     If a system is completely controllable  and  observable  then
$\{ i \}$ does not contain  decoupling zeros, therefore, the equality (5.29) takes place.

     The following example  illustrates the situation when
 the set $\{ i \}$ does not contain all decoupling zeros.

                  \it{EXAMPLE 5.4.} \rm

  We consider  the system  (5.6)  with  the matrix $P(s)$ (5.7).
This system has the only invariant  zero (-3) that coincides with
the output  decoupling zero because  the column rank
 of $P_o(s)$ is reduced at $s = -3$
$$ rankP_o(s)/_{s=3} \; = \; \left [ \begin{array}{rrr} -4 &  0 &  0\\
0 & -2 &  0 \\ 0 & 0 & 0 \\  1 & -1 &  0 \\ 0 & 2 & 0  \end{array} \right ]
 \; = \; 2 < 3  $$
Therefore, the zero $s = -3$ is equal to the output decoupling
zero and the invariant zero simultaneously.
     Moreover, this system is uncontrollable at  $s = 1$ (see Example 5.3).
 Hence,  $s = 1$ is the input decoupling zero. But   the
matrix $P(s)$  has  the complete rank at  $s = 1$  because  there
exists the following nonzero minor
$$ P(s)^{2,3,4,5}_{1,2,3,4}/_{s=1} \;=\; det \left [ \begin{array}{crcc}
   0 & 2 & 0 & 1\\ 0 & 0 & 4 & 1 \\  1 & -1 &  0 &  0 \\ 0 & 2 & 0 & 0
   \end{array} \right ] \;=\; -8 \neq 0 $$
     So, the input decoupling zero $s =  1$  is not  the  invariant
zero. That is why,  the set $\{ i \}$  does  not contain  all
decoupling zeros.

     In the next subsection  we show that only the set of system zeros
contains  all decoupling zeros.

\subsection[General structure of system zeros]
           {General structure of system zeros}
     To reveal the structure of the set of system zeros  $\{ n \}$
we will use  definitions of  system and transmission zeros from Sections 5.2 and  4.2  respectively. Let the
matrix $P(s)$ has the normal rank $n+ \delta (\delta \leq min(r,l)$. We consider all nonzero minors of the
$n+\delta$ order of the matrix $P(s)$ which are formed according to  the relation (5.9). By  the block structure
of $P(s)$ we can write [G1]
$$  P(s)^{1,2,\ldots,n,n+i_1,\ldots, n+ i_{\delta}}_{1,2,\ldots,n, n+j_1,\ldots,
n+j_{\delta}} \; = \; det(sI_n - A) det [ C^{i_1, i_2,\ldots,i_{\delta}}
 (sI_n - A)^{-1} B_{j_1, j_2,\ldots,j_{\delta}}  ]  $$
where $i_1,\ldots,i_{\delta}$ and $j_1,\ldots,j_{\delta}$ are rows and columns  of the matrices $C$  and $B$
respectively. Using the relation (4.5) and the last one we can express  minors of the transfer  function matrix
$G(s) =C(sI_n-a)^{-1}B$, which are formed  by deleting all rows expect $i_1,\ldots,i_{\delta}$  and all columns
expect $j_1,\ldots,j_{\delta}$, as follows
$$ G(s)^{i_1, i_2,\ldots,i_{\delta}}_{j_1, j_2,\ldots,j_{\delta}} \; = \;
 det [ C^{i_1, i_2,\ldots,i_{\delta}} (sI_n - A)^{-1}
 B_{j_1, j_2,\ldots,j_{\delta}}  ] \; = \;
\frac{P(s)^{1,2,\ldots,n,n+i_1,\ldots,
n+i_{\delta}}_{1,2,\ldots,n,n+j_1,\ldots,
n+j_{\delta}}}{det(sI_n-A)} \eqno (5.30) $$
On the other hand
minors of $G(s)$ may be represented as
$$ G(s)^{i_1, i_2,\ldots,i_{\delta}}_{j_1, j_2,\ldots,j_{\delta}} \; = \;
    \frac{Z(s)^{i_1, i_2,\ldots,i_{\delta}}_{j_1, j_2,\ldots,j_{\delta}}}
     {p(s)}  \eqno(5.31) $$
where $p(s)$ is the least common denominator of  minors $
G(s)^{i_1, i_2,\ldots,i_{\delta}}_{j_1, j_2,\ldots,j_{\delta}} $.
 In (5.31)  $ Z(s)^{i_1, i_2,\ldots,i_{\delta}}_{j_1,
j_2,\ldots,j_{\delta}} $ are polynomials, which are constructed from  numerators of $ G(s)^{i_1,
i_2,\ldots,i_{\delta}}_{j_1, j_2,\ldots,j_{\delta}} $ such a way that the new minors  of $G(s)$ have  the
polynomial $p(s)$ as  the common denominator.

     By Definitions  2.1,  2,2,   zeros  of  the
polynomial $p(s)$, which  are poles of  TFM $G(s)$, form a subset of  eigenvalues of the matrix $A$ because some
eigenvalues of $A$ may coincide with uncontrollable or/and unobservable poles ( decoupling zeros), which are
cancelled in the transfer function matrix $G(s)$. Hence, the following equality takes  place
$$             det(sI_n-A) \; = \; p(s)p_d(s)   \eqno(5.32) $$
where the polynomial $p_d(s)$ has  zeros that are unobservable or/and uncontrollable poles ( or decoupling
zeros). Substituting (5.32) into (5.30) and equating the right-hand sides of (5.30) and
 (5.31) we obtain
$$  P(s)^{1,2,\ldots,n,n+i_1,\ldots, n+ i_{\delta}}_{1,2,\ldots,n, n+j_1,\ldots,
n+j_{\delta}} \; = \; p_d(s)Z(s)^{i_1, i_2,\ldots,i_{\delta}}_{j_1,
 j_2,\ldots,j_{\delta}}
 \eqno (5.33)  $$
Thus, by Definitions 5.2, 4.4 and the relation (5.33) we have

\it{ASSERTION 5.4.} \rm   The set of  system  zeros is formed by sets of transmission and decoupling zeros.

     It  follows from Assertions 5.1, 5.4 and  results of Section
5.4.2  that   system  zeros  contains all decoupling zeros.
     This property of system zeros has been illustrated
in Example  5.4.  The set  of invariant  zeros  contains  the only zero  ($-3$) that is the output decoupling
zero. Two decoupling zeros, namely input decoupling  zero  ($1$)  and output decoupling zero ($-3$), are
contained in  the  set  of system zeros  $\{ n \} = \{ 1 , -3 \}$  calculated in Example 5.2.

     Since the complete set of decoupling zeros is formed by the
following sum
$$     \{ o.d. \} + \{ i.d. \} - \{ i.o.d. \}    $$
then using Assertion 5.4 we can write the structure of $\{ n \}$
as follows
$$  \{ n \} = \{ p \} + \{ o.d. \} + \{ i.d. \} - \{ i.o.d. \}
\eqno (5.34) $$
     If a system is controllable and observable then
$$  \{ o.d. \} + \{ i.d. \} - \{ i.o.d. \}= \emptyset $$
and
   $$ \{ n \} = \{ p \}   $$
     Using inclusion (5.24) and the relation  between  $\{ n \}$ and  $\{ i \}$:
 $\{ n \} \; \supseteq \; \{ i \}$, which has  been  obtained
in Section 5.2 ( see Assertion 5.1 ), we can write
$$  \{ n \} \supseteq \{ i \} \supseteq  \{ p \} \eqno
 (5.35) $$

     If a system is controllable and observable then the following
equalities take place
$$           \{ n \} \equiv  \{ i \} \equiv  \{ p \}   \eqno(5.36) $$

                       \it{EXAMPLE 5.5.}  \rm

We find zeros of different type  for the following system with $n
= 4$, $r = 1$,
 $l = 2$
$$ \dot{x} = \left[ \begin{array}{crrc} 1 & 0 & 0 & 0 \\ 0 & -1 & 0 & 0 \\
 0 & 0 & -5 & 0 \\ 0 & 0 & 0 & 7 \end{array} \right ]x +
\left[ \begin{array}{c}  0 \\ -1 \\ -1 \\ -1\end{array} \right ]u, \;
y =  \left[ \begin{array}{cccc} 1 & 0 & 2 & 1 \\ 0 & 0 & 2 & 1\end{array}
\right ]x  \eqno (5.37) $$
 Let's form the system matrix $P(s)$
$$   P(s) \; = \;\left[ \begin{array}{ccccc} s-1 &  0 & 0 & 0 & 0\\
 0 &  s+1 &  0  &  0  &  1 \\0 & 0 &  s+5 &  0 & 1 \\0 & 0 & 0 & s-7 &  1 \\
1 & 0 & 2 & 1 & 0 \\ 0 & 0 & 2 & 1 & 0 \end{array} \right ]
\eqno(5.38) $$
 and construct two minors of the form (5.9)
$$    P(s)^{1,2,3,4,5}_{1,2,3,4,5} = -s(s+1)(s-1)(s-3), \qquad
      P(s)^{1,2,3,4,6}_{1,2,3,4,5} = -3(s+1)(s-1)(s-3)
   \eqno(5.39)  $$
The  monic  greatest  common  divisor  of  these  minors  is   the
polynomial $\psi_c(s) = (s-1)(s+1)(s-3)$. Therefore, the system has three
system zeros: $1, -1, 3$.

     To find invariant  zeros  we  calculate   other
four minors of $P(s)$   of the order  5
$$    P(s)^{2,3,4,5,6}_{1,2,3,4,5} = -s(s+1)(s-3), \qquad
      P(s)^{1,3,4,5,6}_{1,2,3,4,5} \;= \; P(s)^{1,2,4,5,6}_{1,2,3,4,5} \; =
      P(s)^{1,2,3,5,6}_{1,2,3,4,5} = 0
  \eqno(5.40)  $$
and determine the monic greatest common divisor $\psi_I(s)$ of minors (5.39) and  nonzero  minors  (5.40).  We
obtain $\psi_I(s) = (s+1)(s-1)$.  Hence, the system has two invariant zeros : $1, -1$.

     To find transmission zeros we calculate  $G(s) =
C(sI_n-A)^{-1}B$
$$     G(s) \; = \; \left[ \begin{array}{cccc}  1 & 0 & 2 & 1\\
 0 & 0 & 2 & 1 \end{array} \right ] \left[ \begin{array}{cccc}
(s-1)^{-1} & 0 & 0 & 0 \\ 0 & (s+1)^{-1} & 0 & 0 \\ 0 & 0 & (s+5)^{-1} & 0 \\
0 & 0 & 0 & (s-7)^{-1}\end{array} \right ] \left[ \begin{array}{r} 0 \\ -1 \\
-1 \\ -1 \end{array} \right ] \; = $$
$$ = \; \frac{1}{(s+5)(s-7)} \left[ \begin{array}
{c} -3(s-3) \\ -3(s-3) \end{array} \right ]   $$ It is clear that
the system has the only transmission zero: $3$.

From  analysis of the matrices
$$   P_o(s) \; = \;\left[ \begin{array}{cccc} s-1 &  0 & 0 & 0 \\
 0 &  s+1 &  0  &  0 \\0 & 0 &  s+5 &  0 \\0 & 0 & 0 & s-7 \\ 1 & 0 & 2 & 1 \\
0 & 0 & 2 & 1 \end{array} \right ] , \qquad  P_i(s) \; = \;\left[ \begin{array}
{ccccc} s-1 & 0 & 0 & 0 & 0\\ 0 & s+1 &  0  &  0 & -1\\0 & 0 &  s+5 & 0 & -1\\
0 & 0 & 0 & s-7 & -1 \end{array} \right ]  $$ we find that input
and output decoupling zeros are $1$ and $-1$ respectively. So
$$   \begin{array}{ccl} \{ n \} &=& \{ 1, -1, 3 \} \\ \{ i \}& = &\{ 1, -1 \} \\
  \{ p \} &= &\{ 3 \} \\ \{ i.d. \} &= &\{ 1 \} \\ \{ o.d. \} &= &\{ -1 \}
\end{array} \eqno (5.41) $$
      These sets corroborate  the  equalities  and
inclusions, which have been obtained in the present chapter.

\section[Summary conclusions from chapters 3 - 5]
            {Summary conclusions from chapters 3 - 5}

     It has been studied four types of zeros. They are

     1. \it{TRANSMISSION ZEROS}: \rm They are defined via the transfer function
matrix  $G(s)$. They are physically  associated  with transmission-blocking  properties of a system, namely,
with the transmission (or blocking) of a steady signal through a system.

     2. \it{INVARIANT ZEROS}: \rm They are defined via  the system matrix $P(s)$.
They  are  physically associated with the zero-output behavior of
a system, namely, with the transmission (or  blocking) of all
parts of a signal ( free and forced ) through a system.

     3. \it{DECOUPLING ZEROS}: \rm  They are  defined by matrices
$P_i(s) = [sI_n-A,B], \; P_o(s)^T = [sI_n-A^T,C^T]$. They  are associated with existence of system modes that
are decoupled with an input or output of a system. These  modes are complete uncontrollable or unobservable
respectively.

     4. \it{ SYSTEM ZEROS} : \rm They form the set of zeros including all
transmission and decoupling zeros. System zeros  are  defined via
special formed minors of $P(s)$  (5.9).

\chapter[Property of zeros]
        {Property of zeros}
In this chapter we consider main properties, which are inherent
to all type of zeros. To study we will apply  elementary block row
and column operations on a polynomial matrix. There are

     1. interchange any two block rows (columns),

     2. premultiplication (postmultiplication) any block row
        (column) by a non singular  matrix,

     3. replacement of a block row (column) by itself plus any
        other row (column) premultiplicated (postmultiplicated)
        by any polynomial ( or constant) matrix.

     These elementary block operations correspond  to  usual
elementary operations  fulfilled on a group of rows
(columns) and do not change  a normal rank of a
polynomial matrix.

\section[Invariance of zeros]
        {Invariance of zeros}

     The important property of different type zeros is  invariance
 under nonsingular transformations of a
state and/or  inputs/outputs and  also under a state  and/or  an
output feedback control.
     We consider this property more detail.

     Denoting  a  set of any type zeros of system (1.1), (1.2)
 by $\Omega(A,B,C)$ we  study the following transformations.

 \it{1.  NONSINGULAR TRANSFORMATION OF THE STATE VECTOR }: $\hat{x} =Nx$  \rm
where  $\hat{x}$ is a new  state vector, $N$ is a nonsingular $
n\times n$  matrix. Matrices of  the transformed system
are defined as follows : $\hat{A} =
NAN^{-1}$, $\hat{B} = NB$, $\hat{C} = CN^{-1}$. Zeros of the
transformed system are calculated via the following system matrix
 $$  \hat{P}(s) \; = \; \left[ \begin{array}{cc} sI_n-\hat{A} & -\hat{B} \\
 \hat{C} & O \end{array} \right ] \; = \; \left[ \begin{array}{cc}
 sI_n-NAN^{-1} & -NB \\  CN^{-1} & O \end{array} \right ] $$
Applying elementary block  operations to $\hat{P}(s)$ we obtain
series of  rank equalities
 $$  \hat{P}(s) \; = \;rank \left[ \begin{array}{ccc} \left( \begin{array}{c}
 sI_n-NAN^{-1} \\ CN^{-1} \end{array} \right )N & \vdots & \begin{array}{c}
 -NB \\ O \end{array} \end{array} \right ] \;=\; rank \left[ \begin{array}{ccc}
\begin{array}{c} sN-NA \\ C \end{array} & \vdots & \begin{array}{c}
-NB \\ O \end{array} \end{array} \right ] \; = $$
$$ rank \left[ \begin{array}{cc} N^{-1}(sN-NA) & -NB \\ \dotfill & \dotfill\\
C & O \end{array} \right ] \; = \;rank \left[ \begin{array}{cc} sI_n-A & -B\\
C & O\end{array} \right ] \; = \; rank P(s) $$
 Hence, if   the
rank of $P(s)$  is  locally  reduced  below a normal one at $s=z$
then  this property  possesses the matrix $\hat{P}(s)$. We obtain
the following property.

\it{PROPERTY 6.1.} \rm   Zeros  are  invariant  under the
nonsingular transformation of state variables
$$  \Omega(A,B,C) \; = \; \Omega(NAN^{-1}, NB,CN^{-1})
 \eqno(6.1)    $$

\it{2. NONSINGULAR TRANSFORMATION OF THE INPUT VECTOR}: $\hat{u}
=Mu$  \rm where  $\hat{u}$ is a new input, $M$  is a nonsingular
$r\times r$  matrix. A new input matrix is defined as
$\hat{B} =BM^{-1}$. Calculating a rank of the transformed matrix
$\hat{P}(s)$ we obtain
$$    rank \hat{P}(s) = rank\left[ \begin{array}{cc}
 sI_n-A & -BM^{-1} \\ C & O \end{array} \right ] =
 rank \left[ \protect \begin{array} {ccc} \begin{array}{c}
sI_n-A\\C
\end{array}  & \vdots & \left ( \begin{array}{c} -BM^{-1}\\O \end{array}
 \right )M \end{array} \right ]  = rankP(s)  $$
Thus, the following property  takes place

\it{PROPERTY 6.2.} \rm Zeros are invariant  under  the nonsingular
transformation of input variables
$$  \Omega(A,B,C) \; = \; \Omega(A, BM^{-1},C)
 \eqno(6.2)    $$

  \it{3. NONSINGULAR TRANSFORMATION OF THE OUTPUT VECTOR}: \rm  $\hat{y} = Ty$
where $\hat{y}$ is a new output, $T$ is a nonsingular $l\times l$
matrix. A new output matrix is defined as  $\hat{C} =TC$.
Applying the following  elementary  block  operations we
transform the matrix $\hat{P}(s)$  as
$$  rank \hat{P}(s) \; = \; rank\left[ \begin{array}{cc}
sI_n-A & -B \\ TC & O \end{array} \right ] \; =\; rank \left [ \protect
\begin{array}{c} \begin{array}{cc} sI_n-A & -B \end{array} \\ \dotfill \\
 T^{-1}(TC,\;\; O) \end{array} \right ] \;=\;rankP(s) $$
and formulate the following property.

\it{PROPERTY 6.3.} \rm  Zeros are invariant under the nonsingular
transformation of output variables
$$  \Omega(A,B,C) \; = \; \Omega(A, B,TC)
 \eqno(6.3)    $$

We unite Properties 6.1-6.3 as follows
$$  \Omega(A,B,C) \; = \; \Omega(NAN^{-1}, NBM^{-1},TCN^{-1})
 \eqno(6.4)    $$

    \it{ 4.  STATE AND OUTPUT PROPORTIONAL FEEDBACK}.  \rm
 Let us inset  a linear proportional state feedback to system  (1.1),(1.2)
$$ u = Kx + v
\eqno(6.5) $$ where $v=v(t)$  is a new external reference input.
The closed-loop system is described by the equation
$$ \dot{x} = (A +BC)x +Bv
\eqno(6.6)   $$
 with the  output (1.2). To find a rank of the
system matrix $P_c(s)$ of  Eqns. (6.6), (1.2) we use the following
elementary block operations
$$    rank P_c(s) \; = \; rank \left[ \protect \begin{array} {ccc}
\begin{array}{c} sI_n-(A+BK)\\C \end{array}  & \vdots & \begin{array}{c}
 -B\\O \end{array} \end{array} \right ]  \; = \; $$
$$ = \; rank \left[ \protect \begin{array} {ccccc} \left (
 \begin{array}{c} sI_n-(A+BK)\\C \end{array}  \right )  & - & \left (
 \begin{array}{c} -B \\ 0 \end{array} \right )K & \vdots & \begin{array}{c}
 -B\\O \end{array} \end{array} \right ]  \; = \; rankP(s)
 \eqno (6.7)   $$
     We have the similar result if use a linear  proportional
output feedback $ u = \tilde{K}y + v = \tilde{K}Cx + v $ when
 the matrix $\tilde{K}$ is changed  by  $\tilde{K}C$ in (6.7).

    Therefore, we  deduce the following property.

\it{ PROPERTY 6.4.} \rm  Zeros are invariant under the
proportional state and output feedback
$$  \Omega(A,B,C) \; = \; \Omega(A+BK, B,C) \; = \; \Omega(A+BKC, B,C)
 \eqno(6.8)    $$
Uniting  (6.4) and (6.8) we obtain the  general formula of
zero invariance
$$  \Omega(A,B,C) \; = \; \Omega(N(A+BKC)N^{-1}, NBM^{-1},TCN^{-1})
 \eqno(6.9)    $$

\section[Squaring down operation]
        {Squaring down operation}

     Let  system (1.1), (1.2) has more outputs than inputs $(l>r)$.
To get a new  system  with  equal number of inputs and outputs
we combine output variables to replace the  $l$ vector $y=Cx$ by a
new output $r$ vector
$$      \tilde{y} \;=\; Ly \;=\; LCx    \eqno(6.10)     $$
where $L$ is an  $r\times l$ matrix of a full rank. The mentioned
operation is refereed as  'squaring down' [M1].

     If  we add extra input variables to form a new input $l$ vector
$\tilde{u}$ by the rule
$$ u \; = \; D\tilde{u}
 \eqno (6.11) $$
where  $D$ is an  $r\times l$ matrix then this operation is
refereed as 'squaring up' [M1]. Later we study in detail the
squaring down operation because its practical applicability.

    At first we consider the following important property of the squaring down
    operation.

 \it{ASSERTION 6.1.} \rm  A zero set of  system (1.1),(1.2) with $l>r$
is a subset of  zeros of the squared down system  (1.1),(6.10)
$$  \Omega(A,B,C) \; \subseteq \; \Omega(A,B,LC)
 \eqno(6.12)    $$
 but vice versa  of the relation (6.12) is not held.

     \it{PROOF}. \rm For definiteness we assume that the $r\times l$ matrix $L$ of
the rank $r$ has the form
$$                L \; = \; [L_1, L_2]
\eqno(6.13)  $$
 where $L_1$ is a nonsingular $r\times r$  matrix.
At first we consider a particular case when  $L = [I_r, O]$. The system matrix
$\tilde{P}(s)$ of the squared down system  has the following
structure
 $$  \tilde{P}(s) \; = \; \left[ \begin{array}{cc} sI_n-A & -B \\ C_1 & O
\end{array} \right ]
\eqno(6.14) $$ where $C_1$ is the $r\times n$ block row of the
$l\times n$ output matrix $C = \left[ \begin{array}{c} C_1 \\ C_2
\end{array} \right ]$. Zeros of the  squared down system
coincides  with  zeros of the greatest common divisor of minors of
$\tilde{P}(s)$ of a maximal order.

Then writing the system matrix for  Eqns. (1.1), (1.2) with $C^T
=[C^T_1,C^T_2]$
$$    P(s) \; = \; \left[ \begin{array}{cc} sI_n-A & -B \\ C_1 & O \\ C_2 & O
\end{array} \right ]
\eqno(6.15) $$
 we can see that zeros of system (1.1),(1.2) coincides
with  zeros of a  greatest common divisor of  a maximal order
minors of  the $(n+l)\times (n+r)$ matrix (6.15). It is evident that
 the set of minors of  $P(s)$ includes the set of minors of
$\tilde{P}$ for $l >r$. Hence, we have been proved the assertion
for $L = [I_r, O]$, i.e.
$$  \Omega(A,B,C) \; \subseteq \; \Omega(A,B,[I_r,O]C)
 \eqno(6.16)    $$

     Now we consider the general case of $L$ (6.13) and define a nonsingular
$l\times l$  matrix
$$    T \; = \; \left[ \begin{array}{cc} L^{-1}_1 & -L^{-1}_1L_2 \\ O & I_{l-r}
\end{array} \right ] $$
Such the  matrix exists because $rankL_1 = r$. It is clear that
$$  LT = [I_r,O]  \eqno(6.17) $$
From (6.16) and (6.17) we obtain the  following inclusions
for any $l\times n$ matrix $C^*$ of a full rank
$$  \Omega(A,B,C^*) \; \subseteq \; \Omega(A,B,[I_r,O]C^*) \subseteq \;
\Omega(A,B,LTC^*)  $$
or
$$  \Omega(A,B,C^*) \; \subseteq \; \Omega(A,B,LTC^*)
\eqno(6.18) $$
 On the other hand, since $T$ is the nonsingular
$l\times l$ matrix, then  the set $\Omega(A,B,C^*)$ becomes
by  Property  6.3
$$  \Omega(A,B,C^*) \; = \; \Omega(A,B,TC^*)
\eqno(6.19) $$
Substituting the right-hand side of (6.19) into the left-hand side
of (6.18) we get
$$  \Omega(A,B,TC^*) \; \subseteq \; \Omega(A,B,LTC^*)
\eqno(6.20) $$
 Then defining   $C^* =T^{-1}C$  we present (6.20) as

$$  \Omega(A,B,TT^{-1}C) \; \subseteq \; \Omega(A,B,LTT^{-1}C) $$
or
$$  \Omega(A,B,C) \; \subseteq \; \Omega(A,B,LC)    $$
The inclusion obtained completes  the proof.

     Let system (1.1),(1.2) has  $r>l$. We form a new input
$l$ vector $\tilde{u}$ by rule (6.11). This is the squared down operation for
inputs because the number of inputs are decreased from $r$ to $l$.
Similarly to  Assertion 6.1 we can prove

\it{ ASSERTION 6.2.} \rm  Any set of zeros of system (1.1), (1.2)
with $r>l$ is a subset of  zeros of the squared down system
(1.1),(1.2),(6.11)
$$  \Omega(A,B,C) \; \subseteq \; \Omega(A,BD,C)
\eqno(6.21)    $$ but vice versa   does not true.

                   \it{EXAMPLE 6.1.} \rm

 To illustrate  the  result we consider   the
 following system with  $n = 3, r = 1, l = 2 $
$$ \dot{x} = \left[ \begin{array}{rrr} 1 & 0 & 0 \\ 0 & -1 & -1 \\
1 & 0 & -1 \end{array} \right ]x + \left[ \begin{array}{c} -1 \\ 0 \\ 0
\end{array} \right ]u, \qquad y = \left[ \begin{array}{ccc} 1 & 0 & 0 \\
0 & 2 & 0\end{array} \right ]x
 \eqno(6.22) $$
   At first we construct the system matrix
$$   P(s) = \left[ \begin{array}{cccc} s-1 & 0 & 0 & 1 \\ 0  & s+1 & 1 & 0 \\
-1 & 0 & s+1 & 0 \\ 1  & 0 & 0 & 0 \\ 0 & 2 & 0 & 0 \end{array} \right ]
\eqno(6.23)   $$
and calculate two  minors of the order 4 :
$$    P(s)^{1,2,3,4}_{1,2,3,4} = (s+1)^2, \qquad
      P(s)^{1,2,3,5}_{1,2,3,4} = 2
   \eqno(6.24)  $$
The monic greatest common divisor of these minors is equal to
$\psi_c(s) =1$.  Hence, the  system has no zeros. Let's
combine  output variables to form the new scalar output $\tilde{y}$
$$  \tilde{y} = [\; 1\; 1\; ]y = [\; 1\; 1\; ] \left [ \begin{array}{c}
 y_1 \\y_2 \end{array} \right ] \; = \; y_1 + y_2
\eqno (6.25)  $$
In this case the squared down compensator  $L$ is
equal to $[\;1\; 1\;]$ and the new output matrix becomes
$$ \tilde{C} = LC = [\;1\; 1\;] \left [ \begin{array}{ccc} 1 & 0 & 0 \\ 0 & 2 & 0
\end{array} \right ] = [\; 1\; 2\; 0 \;]   $$
We build the system matrix of the new system
$$   \tilde{P}(s) = \left[ \begin{array}{cccc} s-1 & 0 & 0 & 1 \\
0 & s+1 & 1 & 0 \\ -1 & 0 & s+1 & 0 \\ 1  & 2 & 0 & 0 \end{array}
\right ] $$ and calculate the only minor of $\tilde{P}(s)$ that is
$\psi_c(s) = -(s^2+2s-1)$. Zeros of $\psi_c(s)$ are $ z_1 =
 -1+\sqrt{2}$,  $ z_2 = -1-\sqrt{2}$. Thus, we  see that the
squaring down operation introduces new zeros into the system. This
property must be taken into account when squared down
compensators $L$ and/or $D$ are used. For example, such a problem
inevitably appears in a cascade connection of  systems.

\section[Zeros of cascade system]
        {Zeros of cascade system}
     Let us consider two systems  $S_1$ and  $S_2$

$$ S_1 : \qquad \dot{x}_1 = A_1x_1 +  B_1u ,
\qquad y_1 = C_1x_1 $$
$$ S_2 : \qquad \dot{x}_2 = A_2x_2 +  B_2u ,
\qquad y_2 = C_2x_2 $$
where a number of outputs of the first system differs from a number
of inputs of the second one. In above equations
vectors $x_1$, $x_2$, $u_1$, $u_2$,
$y_1$,  $y_2$   have dimensions $n_1\times1$, $n_2\times 1$,
$r_1\times1$, $r_2\times 1$, $l_1\times1$, $l_2\times 1$
respectively.

The cascade connection of $S_1$ and  $S_2$  is as follows: we
insert a linear combination of  variables of the output $y_1$ to the input
$u_2$ , i.e. we use the connection
$$         u_2 \;=\; Gy_1   \eqno(6.26)    $$
with an $r_2\times l_1$  compensator $G$ .

Let every $S_i$, $i=1,2$ has a zero set $\Omega_i$, $i=1,2$. To find
a zero set of the augmented system we substitute the relation (6.26)
in $S_2$ and write the augmented system
$$  \left [ \begin{array}{c} \dot{x}_1 \\ \dot{x}_2 \end{array} \right ] \; =
\;  \left [ \begin{array}{cc} A_1 & O \\ B_2GC_1 & A_2 \end{array} \right ]
\left [ \begin{array}{c} x_1 \\ x_2 \end{array} \right ]  +  \left [
\begin{array}{c} B_1 \\ O \end{array} \right ] u_1, $$
$$ y = \left[ \begin{array}{cc}  O & C_2 \end{array} \right ] \left [
\begin{array}{c} x_1 \\ x_2 \end{array} \right ]
\eqno (6.27)  $$
 Let us find the rank of the system matrix $P_{s_1+s_2}$
of  system  (6.27) by using  the following rank
equalities
$$   rankP_{s_1+s_2} \; = \; rank \left [ \begin{array}{ccc}sI_{n_1}- A_1 &
O & -B_1 \\ B_2GC_1 & sI_{n_2}-A_2 & O \\ O & C_2 & O \end{array} \right ] \;=$$
$$ = \; rank \left ( \left [ \begin{array}{ccc}I_{n_1} & O & O \\
O & sI_{n_2}-A_2 & -B_2 \\ O & C_2 & O \end{array} \right ] \left [
\begin{array}{ccc} sI_{n_1}- A_1 & O & -B_1 \\ O & I_{n_2} & O \\ GC_1 & O & O
\end{array} \right ] \right )
 \eqno(6.28)   $$
On the other hand
$$   rankP_{s_1+s_2} \; = \; rank \left ( \left [ \begin{array}{ccc}
I_{n_1} & O & O \\ O & sI_{n_2}-A_2 & -B_2G \\ O & C_2 & O \end{array} \right ]
 \left [ \begin{array}{ccc} sI_{n_1}- A_1 & O & -B_1 \\ O & I_{n_2} & O \\
 C_1 & O & O \end{array} \right ] \right )
 \eqno(6.29)   $$
We limit our study by systems $S_1$ and $S_2$  having such a  number
of inputs and outputs that provides the product of square matrices in the right-hand  sides
 of (6.28), (6.29). Hence, we have  two cases:

     1. $l_1 \ge r_1 = r_2 = l_2$. Using Eqn.(6.28) we need to evaluate
$$  detP_{s_1+s_2} \; = \; det \left [ \begin{array}{ccc}I_{n_1} & O & O \\
O & sI_{n_2}-A_2 & -B_2 \\ O & C_2 & O \end{array} \right ] det\left [
\begin{array}{ccc} sI_{n_1}- A_1 & O & -B_1 \\ O & I_{n_2} & O \\ GC_1 & O & O
\end{array} \right ]   $$
At first we  decrease the  dimensions of block
matrices  by expanding unity blocks and represent $detP_{s_1+s_2} $ as
$$  detP_{s_1+s_2} \; = \; det \left [ \begin{array}{cc} sI_{n_2}-A_2 & -B_2 \\
C_2 & O \end{array} \right ] det\left [ \begin{array}{cc} sI_{n_1}- A_1 & -B_1 \\
 GC_1 & O \end{array} \right ]
\eqno (6.30)  $$
 It is evident that a rank of the matrix $P_{s_1+s_2}$
is reduced if and only if ranks of the system  matrix  of
 $S_2$ or the following squared down system
$$ S_1^* : \qquad \dot{x}_1 = A_1x_1 +  B_1u_1 , \qquad \tilde{y}_1 = GC_1x_1
\eqno (6.31)  $$
are reduced.

     Let's denote zero sets of
 systems (6.27) and (6.31) by $\Omega_{s_1+s_2}$  and $\Omega_{s_1}^*$ respectively.
 It follows from the equality (6.30)

$$ \Omega_{s_1+s_2} \; = \; \Omega_{s_2} \cup  \Omega_{s_1}^*
\eqno (6.32) $$
     2. $r_2 \ge l_2 =l_1 =r_1$.  Using similar way we obtain  from (6.29)
the following equalities
$$  detP_{s_1+s_2} \; = \; det \left [ \begin{array}{ccc}I_{n_1} & O & O \\
O & sI_{n_2}-A_2 & -B_2G \\ O & C_2 & O \end{array} \right ] det\left [
\begin{array}{ccc} sI_{n_1}- A_1 & O & -B_1 \\ O & I_{n_2} & O \\ C_1 & O & O
\end{array} \right ] \; =  $$
$$ = \; det \left [ \begin{array}{cc} sI_{n_2}-A_2 & -B_2G \\C_2 & O
\end{array} \right ] det\left [ \begin{array}{cc} sI_{n_1}- A_1 & -B_1 \\
C_1 & O \end{array} \right ]   $$
So, a rank of  the matrix $P_{s_1+s_2}$ is reduced if and  only  if
ranks of the system matrix of  $S_1$ or  the squared down system
$$ S_2^* : \qquad \dot{x}_2 = A_2x_2 +  B_2G\tilde{u}_2 , \qquad
\tilde{y}_2 = C_2x_2
\eqno (6.33)  $$
are reduced.
     Denoting a zero set of system (6.33)  by $\Omega_{s_2}^*$   we  can
write the following equality
$$ \Omega_{s_1+s_2} \; = \; \Omega_{s_1} \cup  \Omega_{s_2}^*
\eqno (6.34) $$

Now we analyze relations (6.32),(6.34). In the first case the system
$S_1^*$ is obtained from $S_1$  by squaring down its outputs.
In  the second case the system $S_2^*$ is obtained from $S_2$ by
squaring down its inputs. Above we have shown  that  the
squaring  down operation introduces new zeros into a system.
Denoting the set of  introducing zeros by $\Omega_{sq}$ we
represent sets $\Omega_{s_1}^*$  and $\Omega_{s_2}^* $ as
$$ \Omega_{s_i}^* \; = \; \Omega_{s_i}  \cup \Omega_{sq}, \qquad i=1,2
\eqno (6.35) $$
and   rewrite (6.32) or (6.34) using (6.35) as the only sum
$$ \Omega_{s_1+s_2} \; = \; \Omega_{s_1}  \cup \Omega_{s_2} \cup \Omega_{sq}
\eqno (6.36) $$
     So,  it  has  been  shown:  The  cascade
connection of systems with different  numbers  of  inputs  and
outputs may introduce additional zeros into an augmented
system. Hence, it is necessary to choose the  matrix $G$ to  shift these zeros
 to  the left-hand side of the complex
plan.

     Let's consider the important particular case when systems  $S_1$ and
$S_2$ have same numbers of inputs and outputs : $r_1 = l_1 = r_2 = l_2$.
In this case   $G$ is a square nonsingular matrix which transforms
outputs of $S_1$ or inputs  $S_2$. According Properties
6.2, 6.3 we have
$$ \Omega_{s_1}^* \; = \; \Omega_{s_1} ,\qquad \Omega_{s_2}^* \; =
\; \Omega_{s_2}, \qquad \Omega_{sq} = \emptyset  $$
and the  zero set of the augmented system  is
$$ \Omega_{s_1+s_2} \; = \; \Omega_{s_1} \cup  \Omega_{s_2}
\eqno (6.37) $$
The next assertion is the direct corollary of the equality (6.37).

\it{ASSERTION 6.3.} \rm  The set of zeros  of the  cascade connection of systems
$S_1$,\ldots,$S_k$ having  equal numbers of inputs and outputs is
defined as
$$ \Omega_{s_1+s_2+\cdots+s_k} \; = \; \Omega_{s_1} \cup \Omega_{s_2} \cup
   \cdots \cup \Omega_{s_k} $$

                \it{EXAMPLE} \rm 6.2.

Let us consider two systems
$$   \dot{x}_1 = \left[ \begin{array}{cc} 1 & 0 \\ 0 & 2 \end{array}
 \right ]x_1 + \left[ \begin{array}{c} 1 \\ 1 \end{array} \right ]u_1,
\qquad y_1 = x_1
\eqno(6.38)$$
$$ \dot{x}_2 = \left[ \begin{array}{cc} 2 & 0 \\ 1 & 1 \end{array} \right ]x_2
 + \left[ \begin{array}{c} 1 \\ 0 \end{array} \right ]u_2, \qquad
 y_2 = \left[ \begin{array}{cc} 1 & 2 \end{array} \right ]x_2
\eqno(6.39)$$
We insert  a linear combination of variables
of $y_1$ to the input $u_2$  of the  second  system, i.e. we  use
the connection :
$$    u_2 = Gy_1   \eqno(6.40)   $$
Substituting (6.40) into (6.39) and using (6.38) we can write  the
following augmented system
$$  \left [ \begin{array}{c} \dot{x}_1 \\ \dot{x}_2 \end{array} \right ] \; =
\;  \left [ \begin{array}{cc} A_1 & O \\ B_2GC_1 & A_2 \end{array} \right ]
\left [ \begin{array}{c} x_1 \\ x_2 \end{array} \right ]  +  \left [
\begin{array}{c} B_1 \\ O \end{array} \right ] u_1, $$
$$ y = \left[ \begin{array}{cc}  O & C_2 \end{array} \right ] \left [
\begin{array}{c} x_1 \\ x_2 \end{array} \right ]   $$
where
$$  A_1 =   \left[ \begin{array}{cc} 1 & 0 \\ 0 & 2 \end{array} \right ],\;
    A_2 =   \left[ \begin{array}{cc} 2 & 0 \\ 1 & 1 \end{array} \right ], \;
    B_1 =   \left[ \begin{array}{c} 1 \\ 1 \end{array}  \right ],\;
    C_1 =   \left[ \begin{array}{cc} 1 & 0 \\ 0 & 1 \end{array}  \right ], \;
    B_2 =   \left[ \begin{array}{c} 1 \\ 0  \end{array}  \right ] ,\;
    C_2 =  \left[ \begin{array}{cc} 1 & 2 \end{array}  \right ]  $$
Zeros of systems (6.38), (6.39)  are respectively
$\Omega_{s_1} = \{ \emptyset\}$,  $\Omega_{s_2} = \{ -1 \}$.

     Let's assign
$$                    G =[\;1 \; 1\;]  $$
and find zeros of the squared down system $ S_1^* : \qquad \dot{x}_1 =
 A_1x_1 +  B_1u_1 , \qquad \tilde{y}_1 = Gy_1 =GC_1x_1 $. We obtain
$\Omega_{s_1}^* = \{ 1.5 \}$. Therefore, the squared  down operation
introduces the only zero  ($z = 1.5$), i.e.  $\Omega_{sq} = \{ 1.5 \}$.
Using formula  (6.36)  we can determine the set of zeros of
the overall cascade connection: $ \Omega_{s_1+s_2}
= \{ \emptyset\} \cup  \{-1 \} \cup \{ 1.5 \} = \{ -1, 1.5 \}$.

    For  checking we calculate  zeros of the augmented system
with $G = [\;1\;1\;]$. The system matrix $P_{s_1+s_2}$ of the cascade connection  system is

$$ P_{s_1+s_2} \;=\; \left[ \begin{array}{ccccr} s-1 &  0 & 0 & 0 & -1\\
 0 &  s-2 &  0  &  0  &  -1 \\ 1 & 1 & s-2 &  0 & 0 \\0 & 0 & -1 & s-1 & 0 \\
0 & 0 & 1 & 2 & 0  \end{array} \right ]   $$
Determining $detP_{s_1+s_2} = -(2s-3)(s+1)$ we obtain  $ \Omega_{s_1+s_2} =
\{ 1.5, -1 \}$.

\section[Dynamic output feedback]
        {Dynamic output feedback}
     Now we study  the effect of a dynamic  regulator  (dynamic  output
feedback) on system zeros. To this point we insert the following linear  dynamic  output
feedback
$$ \dot{z} \; = \; Fz +Qy
  \eqno(6.41) $$
$$   u \;=\; v - K_1z
    \eqno (6.42) $$
into system (1.1), (1.2). Here $z$
 the $p\times 1$ state vector of the dynamic regulator, $v$ is the $r\times 1$
reference input vector and  constant matrices  $F$, $Q$, $K_1$
have the corresponding sizes.


Substituting (6.41), (6.42) into (1.1), (1.2) and denoting the
new state vector as  $[ x^T ,z^T ]$ gives the following closed-loop
augmented system
$$  \left [ \begin{array}{c} \dot{x} \\ \dot{z} \end{array} \right ] \; =
\;  \left [ \begin{array}{cc} A & -BK_1 \\ QC & F \end{array} \right ]
\left [ \begin{array}{c} x \\ z \end{array} \right ]  +  \left [
\begin{array}{c} B \\ O \end{array} \right ]v , $$
$$ \tilde{y} = \left[ \begin{array}{cc}  C & O \end{array} \right ] \left [
\begin{array}{c} x \\ z \end{array} \right ]
\eqno(6.43) $$
with the input  $r$ vector $v$ and  the output $l$ vector $\tilde{y}$. To find zeros
of this system we  need to  analyze its
the system matrix
$$  P(s) = \left [ \begin{array}{ccc} sI_{n+p} - \left [ \begin{array}{cc}
            A & -BK_1 \\ QC & F \end{array} \right ]  & \vdots & -\left [
           \begin{array}{c} B \\ O \end{array} \right ]  \\
           \dotfill & \dotfill & \dotfill \\  \left [ \begin{array}{cc}
           C & O \end{array} \right ] & \vdots & O \end{array} \right ] \; =
           \left [ \begin{array}{ccc} sI_n - A & -BK_1 & -B \\
           -QC & sI_p-F &  O \\ C & O & O \end{array} \right ] $$
Let's carry out several elementary block operations on the  matrix $P(s)$:
 We fulfill left and right multiplications of  $P(s)$   by   unimodular  matrices  and  then
interchange the second and the third block rows and the appropriate
columns. We result in  following  rank equalities
$$   rankP(s) = rank \left \{ \left [ \begin{array}{ccc} I_n & O & O\\
                O & I_p & O \\ O & O & I_l \end{array} \right ]
                \left [ \begin{array}{ccc} sI_n - A & -BK_1 & -B \\
               -QC & sI_p-F &  O \\ C & O & O \end{array} \right ]
                \left [ \begin{array}{ccc} I_n & O & O\\
           O & I_p & O \\ O & K_1 & I_l \end{array} \right ] \right \} \; = $$
$$           \; = rank \left [ \begin{array}{ccc} sI_n-A & O & -B\\
                O & sI_p-F & O \\ C & O & O \end{array} \right ] \; = \;
             \left [ \begin{array}{cccc} sI_n- A & -B &\vdots & O \\
              C & O & \vdots & O \\ \dotfill & \dotfill & \dotfill &\dotfill \\
              O & O & \vdots & sI_p-F \end{array} \right ] \eqno(6.44)$$
Hence, the rank of $P(s)$  is locally reduced at $s = s^*$ if  and
only if $s^*$ coincides with a zero of  system (1.1), (1.2) or
with an eigenvalue of  the matrix  $F$.

     The similar result may be obtain for the  feedback  regulator
of the  general structure
$$ \dot{z} \; = \; Fz +Qy ,  $$
$$   u \;=\; v - K_1z -K_2y = v - K_1z - K_2Cx
    \eqno (6.45) $$
where $K_2$ is a constant  $r\times l$ matrix.
The system  matrix  for
system (1.1),(1.2) with regulator (6.45) is
$$  P(s) \; = \; \left [ \begin{array}{ccc} sI_n - (A-BK_2C) & -BK_1 & -B \\
                 -QC & sI_p-F &  O \\ C & O & O \end{array} \right ] $$
Executing  elementary block operations on $P(s)$  we can show
that
$$ rankP(s) = rank \left [ \begin{array}{cccc} sI_n-(A-BK_2C) & -B &
           \vdots & O \\ C & O & \vdots & O \\ \dotfill & \dotfill & \dotfill &
           \dotfill \\ O & O & \vdots & sI_p-F \end{array} \right ] =
        rank\left [ \begin{array}{cccc} sI_n- A & -B &\vdots & O \\
            C & O & \vdots & O \\ \dotfill & \dotfill & \dotfill &\dotfill \\
            O & O & \vdots & sI_p-F \end{array} \right ]
            \eqno (6.46) $$
Thus, it follows from (6.44) and (6.46)

\it{ASSERTION 6.4.} \rm  The set of zeros of  the  augmented  system  with
the dynamic regulator (6.41), (6.42) or (6.45) consists of all  zeros
of system (1.1),(1.2) and all eigenvalues of
the matrix dynamics $F$ of the regulator.

     Therefore, we conclude that a dynamic  feedback  introduces
additional zeros in any system. This result generalizes the
similar property of the classic single-input/ single-output system,
namely, zeros of any closed-loop transfer function include
zeros of an open-loop transfer function and poles of a
compensator transfer function.

     Let  us consider the important case of a dynamic regulator, namely,
the proportional-integral (PI) regulator
$$ \dot{z} \; = \; y , \qquad  u \;=\; v - K_1z -\tilde{K}_2x
\eqno (6.47) $$
PI-regulator (6.47) is the particular  case  of the
dynamic regulator (6.45) with $p=l,\; F=O,\; Q=I_l,\;K_2C = \tilde{K}_2$.
We have from Assertion 6.4

\it{COROLLARY 6.1. } \rm  Any PI-regulator of the order $l$ introduces $l$ zeros in origin.

                    \it{EXAMPLE 6.3.}  \rm

To study the affect of the dynamic  feedback we consider  the following
simple system
$$ \dot{x} = \left[ \begin{array}{rc} -1 & 0 \\ 1 & 2 \end{array} \right ]x
 + \left[ \begin{array}{c} 1 \\ 0 \end{array} \right ]u, \qquad
 y = \left[ \begin{array}{cc} 1 & 1 \end{array} \right ]x
\eqno(6.48)$$
and the  dynamic  regulator of the structure (41),(42)
$$ \dot{z} \; = \;2z + y \;= \; 2z + [\; 1 \; 1 \;]x
\eqno(6.49) $$
$$   u \;=\; v - z - y \; = \; v - z - [\; 1 \; 1 \;]x
\eqno (6.50) $$
Here  $p=1,\; F=2,\; Q=1,\;K_1 = 1,\;K_2=1$.
Substituting  (6.50) in (6.48)
$$ \dot{x} = \left[ \begin{array}{rc} -1 & 0 \\ 1 & 2 \end{array} \right ]x
+ \left[ \begin{array}{c} 1 \\ 0 \end{array} \right ][ \; 1 \; 1 \;]x -
   \left[ \begin{array}{c} 1 \\ 0 \end{array} \right ]z  +
   \left[ \begin{array}{c} 1 \\ 0 \end{array} \right ]v  \; =
    \; \left[ \begin{array}{rc} -2 & -1 \\ 1 & 2 \end{array} \right ]x -
     \left[ \begin{array}{c} 1 \\ 0 \end{array} \right ]z  +
      \left[ \begin{array}{c} 1 \\ 0 \end{array} \right ]v  $$
and uniting this equation with (6.49) we  obtain  the  augmented
system
$$  \left [ \begin{array}{c} \dot{x} \\ \dot{z} \end{array} \right ] \; =
     \;  \left [ \begin{array}{rrr}-2 & -1 & -1\\ 1 & 2 & 0 \\ 1 & 1 & 2
      \end{array}\right ] \left [ \begin{array}{c} x \\ z \end{array} \right ]
    +  \left [ \begin{array}{c} 1\\0\\0 \end{array} \right ]v , $$
$$    y = \left[ \begin{array}{ccc}  1 & 1 & 0 \end{array} \right ]
          \left [ \begin{array}{c} x \\ z \end{array} \right ]
\eqno(6.51) $$
To find zeros of (6.51) we construct
$$   P(s) = \left[ \begin{array}{cccr} s+2 & 1 & 1 & -1 \\ -1  & s-2 & 0 & 0 \\
-1 & -1 & s-2 & 0 \\ 1  & 1 & 0 & 0  \end{array} \right ]  $$
and find the zero polynomial $\psi(s)=(s-2)(s-1)$. Therefore,
the closed-loop system (6.51) has two zeros : $s_1=2, s_2=1$.

     For testing we calculate  zeros of  system (6.48) and obtain  the only zero
($s_1=1$). Since one eigenvalue of the
matrix dynamics of (6.49) is  equal  to  2  then  we  have
obtained the corroboration of Assertion 6.4.

\section[Transmission zeros and high output feedback]
         {Transmission zeros and high output feedback}

Let a linear  negative  proportional  output feedback
$$     u \; = \; -Ky
\eqno(6.52) $$
is applied into  completely controllable and observable system
(1.1),(1.2) having equal numbers of inputs and outputs ($r=l$).

     We will investigate asymptotic behavior  of  eigenvalues  of
the dynamics matrix $A-BKC$  of the closed-loop system when elements of the gain matrix
$K$ unlimited increase. For this purpose we represent the matrix $K$
as
$$          K = k\tilde{K}     $$
where $\tilde{K}$ is a constant $r\times r$ matrix of a full  rank  with  bounded
elements and $k$ is a scalar value that increases to infinity.

     For  $\phi(s) = det(sI_n-A)$ and  $\phi_c(s) = det(sI_n-(A-BKC))$,
characteristic polynomials of the open-loop and closed-loop systems respectively,
we prove  the following assertion.

\it{ASSERTION 6.5. } \rm  [H1].
$$   \frac{\phi_c(s)}{\phi(s)}\; = \; det(I_r + KG(s))
\eqno(6.53) $$

\it{PROOF}.$\;$ \rm At first we express $\phi_c(s)$ via  $\phi(s)$
$$ \phi_c(s) = det(sI_n-(A-BKC)) = det \{ (sI_n-A)(I_n+ (sI_n-A)^{-1}BKC) \} =$$
$$ = det(sI_n-A)det(I_n+(sI_n-A)^{-1}BKC) = \phi(s)det(I_n+ (sI_n-A)^{-1}BKC)$$
Denoting $N = (sI_n-A)^{-1}B$, $\;M = KC$ and using the equality from
[K5, lemma 1.1]: $det(I_n +NM) = det(I_r + MN)$  where matrices  $M$, $N$
of dimensions $r\times n$  and  $n\times r$ respectively, we transform the
expression in the right-hand  side  of the last equality as
$$       \phi_c(s) = \phi(s)det(sI_n-A)^{-1}BKC) = \phi(s)det(I_r
         + KC(sI_n-A)^{-1}B) $$
This proves the  assertion.

    Rewriting (6.53) with $K=k\tilde{K}$ ($k \neq 0$)
$$       \phi_c(s) = \phi(s)det(I_r + k\tilde{K}G(s)) =  \phi(s)det \{
          k(\frac{1}{k}I_r + \tilde{K}G(s))\}   $$
and taking out $k$ from the determinant
$$       \phi_c(s) = \phi(s)k^r det(\frac{1}{k}I_r + \tilde{K}G(s))
\eqno(6.54) $$
we analyze  the relation (6.54) as $k \to \infty $
$$     \begin{array}{c} \\ lim \\ \scriptstyle k \rightarrow \infty \end{array}
        \displaystyle \phi_c(s) = \phi(s)k^r det(\tilde{K}G(s)) =
         k^r \phi(s)det\tilde{K}detG(s)
\eqno(6.55) $$
As it has been shown in Section 5.2 the following polynomial
$$ \phi(s)detG(s) \; = \; det(sI_n-A)det(C(sI_n-A)^{-1}B)   $$
is the zero polynomial of system (1.1),(1.2) with  $r=l$.
Denoting $\;\psi(s) = \phi(s)detG(s) $, $\;d=det\tilde{K}$  we can rewrite (6.55)
 as follows
$$     \begin{array}{c} \\ lim \\ \scriptstyle k \rightarrow \infty \end{array}
        \displaystyle \phi_c(s) = dk^r \psi(s)
\eqno(6.56)  $$
Hence, as $k \to \infty$  $\;n-r$ eigenvalues of the matrix $A-Bk\tilde{K}C$ will
asymptotically achieve  zero locations while the remainder $r$
eigenvalues will tend to infinity.

    \it{ REMARK 6.1.} \rm  The result obtained  extends the known
classic root-locus method to  multivariable systems.

\chapter[System zeros and matrix polynomial]
        {System zeros and matrix polynomial}

     In this chapter we will study a definition of system zeros
via an $l\times r$  matrix polynomial  of a degree $\nu-1$ where $\nu$  is the controllability index of the pair
$(A,B)$.  Above in   Section  4.3  we have already introduced the similar definition of transmission zeros for a
system with $n= r\nu$  where $n$, $r$ are  an order and number of inputs  (see Corollary 4.2). Now we consider
the general case $n \neq r\nu$. This definition was introduced by Smagina [S4] in 1981 and will use for study
important properties  of zeros such as a maximal number of zeros and its relations with the Markov parameter
matrices $CB, CAB,\ldots$.

\section[Zero definition via matrix polynomial]
        {Zero definition via matrix polynomial}

     Using the  nonsingular  transformation  of  state  and  input
variables
$$ z=Nx, \; \; v = M^{-1}u    $$
we  reduce  completely controllable system (1.1), (1.2) to Yokoyama's canonical form (1.61)
$$ \dot{z} = Fz +Gv , \qquad y = CN^{-1}z
\eqno(7.1)  $$
where
$$   F =NAN^{-1} = \left [ \begin{array}{c} F_1 \\ \dotfill\\ F_{\nu 1},
              F_{\nu 2},\ldots, F_{\nu \nu} \end{array} \right ], \qquad
     G =NBM = \left [ \begin{array}{c} O \\ G_{\nu} \end{array} \right ]
\eqno(7.2)  $$
 The structure of  nonzero blocks  $F_1$ (1.62), $G_{\nu}$ (1.65)  of dimensions  $(n-r)\times n$
and $r\times r$ respectively are depended on integers $\nu$ and $l_1,l_2,\ldots, l_{\nu}$ (1.59), (1.60).

     As it has been shown in  Section 6.1 system zeros  are
invariant  under state   and  input  nonsingular transformations. Therefore, system zeros of (1.1),(1.2)  are
equal to  system zeros of system (7.1) and defined via the following system matrix
$$    \bar{P}(s) \; = \; \left[ \begin{array}{cc} sI_n-F & -G \\
       CN^{-1} & O \end{array} \right ]
\eqno (7.3)   $$ Let's partition  the matrix $CN^{-1}$
$$ CN^{-1} = [\; C_1, C_2, \ldots, C_{\nu} \;] \eqno
(7.4)  $$
 where $C_i$ are  $l\times l_i$ blocks and construct  the following  $l\times r$ matrix polynomial
\footnote { The definition of a  matrix polynomial has been introduced in Sec.1.2.1}  of the degree $\nu -1$
$$
 \tilde{C}(s) = [O,C_1] + [O,C_2]s + \cdots +  [O,C_{\nu -1}]s^{\nu -2}
                +C_{\nu}s^{\nu -1}
\eqno(7.5)  $$
 Now we show that system zeros of controllable system (1.1), (1.2) are defined in terms of
the matrix polynomial (7.5).

At first we consider the particular case $r=l$.

\it{THEOREM 7.1.} \rm  System  zeros  of  controllable system (1.1),(1.2) with a similar number of inputs and
outputs are defined as zeros of the following polynomial
$$ \psi(s) = s^{n-r\nu}det\tilde{C}(s)
\eqno(7.6)   $$ where integer $\nu$ is the controllability index of the pair $(A,B)$ (see  (1.45)).

\it{PROOF}. $\;$\rm  In this case the system matrix (7.3) has the only minor of the maximal order $n+r$
$$    \bar{P}(s)^{1,2,\ldots,n,n+1,\ldots,n+r}_{1,2,\ldots,n, n+1,\ldots,n+r} =
     det\left[ \begin{array}{cc} sI_n-F & -G \\ CN^{-1} & O \end{array} \right ]
\eqno(7.7) $$
 Zeros of the minor (7.7) (that is a polynomial in $s$)
coincide with system zeros (see Definition 5.2). To find the determinant in (7.7) we partition the matrix
$sI_n-F$ with $F$ from (1.62)  into four blocks
$$ sI_n-F \;=\;\left [ \begin{array}{ccc} sI_{n-r}-F_{11} & \vdots & -F_{12}\\
                   \dotfill & \dotfill & \dotfill \\
                  -F_{\nu 1},-F_{\nu 2},\ldots, -F_{\nu, \nu -1} & \vdots &
                         sI_r-F_{\nu \nu} \end{array} \right ]
\eqno(7.8)$$
where $F_{\nu i}$ are $r\times l_i$  submatrices
($i=1,2,\ldots,\nu$), matrices $sI_{n-r}-F_{11}$  and $F_{12}$
have dimensions  $(n-r)\times
 (n-r)$ and $(n-r)\times r$ respectively. Substituting (7.8) and  (7.4)
into  the  right-hand  side of (7.7) and using the structure of
 $G$ (1.64) we can present $det\bar{P}(s)$ as
$$ det\bar{P}(s) \;=\;det\left [ \begin{array}{ccccc}
          sI_{n-r}-F_{11} & \vdots & -F_{12} & \vdots & O \\
            -F_{\nu 1},-F_{\nu 2},\ldots, -F_{\nu ,\nu -1} & \vdots &
           sI_r-F_{\nu \nu} & \vdots & -G_{\nu} \\
           C_1, C_2,\ldots, C_{\nu -1} & \vdots & C_{\nu} & \vdots & O
            \end{array} \right ] \; = $$
$$        det\left [ \begin{array}{ccccc}
           sI_{n-r}-F_{11} & \vdots & -F_{12} & \vdots & O \\
           C_1, C_2,\ldots, C_{\nu -1} & \vdots & C_{\nu} & \vdots & O \\
             F_{\nu 1},F_{\nu 2},\ldots, F_{\nu ,\nu -1} & \vdots &
            -sI_r+F_{\nu \nu} & \vdots & G_{\nu}
             \end{array} \right ] =
            det\left [ \begin{array}{ccc}
           sI_{n-r}-F_{11} & \vdots & -F_{12} \\
           \dotfill & \dotfill & \dotfill \\
           C_1, C_2,\ldots, C_{\nu -1} & \vdots & C_{\nu}
           \end{array} \right ]det(G_{\nu})
\eqno(7.9) $$
To calculate the determinant of the block matrix in
the right hand-side of (7.9) we  use the formula from [G1](
assuming $s \neq 0$)
$$     det\left [ \begin{array}{ccc} sI_{n-r}-F_{11} & \vdots & -F_{12} \\
     \dotfill & \dotfill & \dotfill \\
    C_1, C_2,\ldots, C_{\nu -1} & \vdots & C_{\nu} \end{array}
     \right ]det(G_{\nu})\; =$$
$$      \; det(sI_{n-r}-F_{11})(det(C_{\nu} +
     [C_1, C_2,\ldots, C_{\nu -1}](sI_{n-r}-F_{11})^{-1}F_{12})
\eqno  (7.10) $$
 Since  structures of  matrices $F_{11}$ and $F_{12}$ coincide with ones of  $P_{11}$ and
$P_{12}$ respectively (see (1.27) with $p=\nu$)  then according  results of  Section 1.2.1 we  find
$$
    det(sI_{n-r}-F_{11}) = s^{l_1+l_2+\cdots +l_{\nu -1}} = s^{n-r}
\eqno(7.11)$$
 $$          (sI_{n-r}-F_{11})^{-1}F_{12} \; = \;
           \left[ \begin{array}{c} s^{1-\nu }[O,I_{l_1}] \\
           s^{2- \nu }[O,I_{l_2}] \\ \vdots \\s^{-1}[O,I_{l_{\nu -1}}]
           \end{array} \right]
\eqno(7.12)$$
 where $[O,I_{l_i}]$ are $l_i\times r$  matrices. Substituting (7.11) and (7.12) into (7.10)
$$        det\left [ \begin{array}{ccc} sI_{n-r}-F_{11} & \vdots & -F_{12} \\
        \dotfill & \dotfill & \dotfill \\
       C_1, C_2,\ldots, C_{\nu -1} & \vdots & C_{\nu} \end{array}
        \right ] \; = \; s^{n-r}det(C_{\nu} + [C_1, C_2,\ldots, C_{\nu -1}]
        \left[ \begin{array}{c} [O,I_{l_1}] \\ s[O,I_{l_2}] \\ \vdots
        \\s^{\nu -2}[O,I_{l_{\nu -1}}] \end{array} \right]s^{1-\nu}) \; = \; $$
$$       = \; s^{n-r}s^{r-r\nu}det[C_1, C_2,\ldots, C_{\nu }]
        \left[ \begin{array}{c} [O,I_{l_1}] \\ s[O,I_{l_2}] \\ \vdots
      \\s^{\nu -2}[O,I_{l_{\nu -1}}] \\ s^{\nu -1}I_r \end{array} \right]) \;=
      s^{n-r\nu}det([O,C_1] + [O,C_2]s + \cdots + C_{\nu}s^{\nu -1}) \;=\; $$
$$       =\;   s^{n-r\nu}det\tilde{C}(s)
\eqno (7.13) $$ and  the right-hand side of (7.13) into (7.9) we get the final expression for $det\bar{P}(s)$
$$ det\bar{P}(s) \;= \;  s^{n-r\nu}det\tilde{C}(s)det(G_{\nu})
\eqno (7.14) $$
 Since  the $r\times r$ matrix $G_{\nu}$ is nonsingular one then zeros of $det\bar{P}(s)$ coincide with
zeros of the polynomial $\psi(s) = s^{n-r\nu}det\tilde{C}(s)$. This proves the theorem.

\it{REMARK 7.1.} \rm  The formula (7.6) is also true for $s= 0$. Indeed, calculating  the determinant of the
matrix $\bar{P}(s)$ (7.9) at $s= 0$ we obtain
$$ det\bar{P}(0) \; = \; det\left [ \begin{array}{ccc}
            -F_{11} & \vdots & -F_{12} \\ \dotfill & \dotfill & \dotfill \\
             C_1, C_2,\ldots, C_{\nu -1} & \vdots & C_{\nu}
             \end{array} \right ]det(G_{\nu})\; = \; (-1)^{n-r}
            det[C_1, C_{21},\ldots, C_{\nu 1}]det(G_{\nu})     $$
where $C_{i1}$ are $r\times (l_i-l_{i-1})$  submatrices of
matrices $C_i =[C_{i1},C_{i2}]$, $i=2,3,\ldots, \nu$. On the other
hand using the structure of (7.5) we can calculate
$$ s^{n-r\nu}det\tilde{C}(s)/_{s=0} \; =$$
$$ \; s^{n-r\nu }det\{ [O,C_1] +
       [O,C_2diag(I_{l_2}s)] +[O,C_3diag(I_{l_3}s^2)] + \cdots +
      [O,C_{\nu -1}diag(I_{l_\nu -1}s^{\nu -2})] + $$
$$    +  C_{\nu}diag(I_rs^{\nu -1}) \}/_{s=0} \; =
      det\{ [O,C_1] + [O,(C_{21},C_{22})diag(I_{l_2 -l_1}, I_{l_1}s] + $$
$$   +  [O,(C_{31},C_{32},C_{33})diag(I_{l_3 -l_2},I_{l_2 -l_1}s, I_{l_1}s^2]
     + \cdots +  $$
$$     + [C_{\nu 1},C_{\nu 2},\cdots, C_{\nu \nu}]diag(I_{r -l_{\nu -1}},
      I_{l_{\nu -1} -l_{\nu -2}}s,\ldots  I_{l_1}s^{\nu -1}) \}/_{s=0} \; = $$
$$    =\;  det(C_{\nu1},C_{\nu2}, \ldots,C_{21},C_{1})   $$
Thus, $det\bar{P}(s) /_{s=0} = 0$ if and only if $ s^{n-r\nu}det\tilde{C}(s)/_{s=0}= 0$.

\it{REMARK 7.2.} \rm  Consider  the particular  case  when system (1.1),(1.2) has  $n=r\nu, \;\; l_1 = l_2 =
\cdots = l_{\nu} = r$. Such system is reduced to Asseo's canonical form and the matrix $CN^{-1}$  is partitioned
into $(l\times r$)  blocks $C_i$. As a result the
 matrix polynomial (7.5) becomes  the simplest structure
$$ \tilde{C}(s) \;=\; C_1 + C_2s + \cdots + C_{\nu}s^{\nu -1}
\eqno  (7.15)  $$
 The zero polynomial is defined as
$$ \psi(s) \;=\; det\tilde{C}(s)
\eqno (7.16) $$

                 \it{EXAMPLE 7.1.} \rm

 To illustrate  the method we consider
system (1.1),(1.2) with  $n = 4$, $r = l = 2 $ and the following state-space model matrices
$$   A \; = \;  \left[ \begin{array}{cccc} 2 & 1 & 0 & 0 \\ 0 & 1 & 0 & 1 \\
     0 & 2 & 0 & 0 \\ 1 & 1 & 0 & 0 \end{array} \right ], \qquad
     B\; = \; \left[ \begin{array}{cc} 1 & 0 \\ 0 & 0 \\ 0 & 0 \\ 0 & 1
            \end{array} \right ], \qquad
     C \; = \;  \left[ \begin{array}{cccc} 1 & -1 & 1 & 0 \\ 1 & 1 & 0 & 1
             \end{array} \right ]
\eqno (7.17) $$

     As it has been shown in Sect.1.2.3.  (Example  1.3)  this system has
$\nu= 3$, $l_1=l_2=1$, $l_3 =2$ and  the following transformation matrix that reduces the system to Yokoyama's
form
$$ N = \left[ \begin{array}{cccc}   0 &  0 &  0.5 &  0 \\   0 &  1 & 0 &  0 \\
 1 &  0 &  0 &  0 \\ 0 & 1 &  0 &  1  \end{array} \right ] $$
Calculating
$$ N^{-1} = \left[ \begin{array}{crcc}   0 &  0 &  1 &  0 \\ 0 & 1 & 0 & 0 \\
 2 &  0 & 0 & 0 \\ 0 & -1 &  0 &  1  \end{array} \right ] $$
and
$$  CN^{-1} \;=\; \left [ \begin{array}{crcc} 1 &-1 & 1 & 0 \\2 &-1 & 1 & 0
     \end{array} \right ] N^{-1} \;=\; \left [ \begin{array}{crcc}
     2 & -1 & 1 & 0 \\ 0 & 0  & 1 & 1 \end{array} \right ]    $$
we can find the  matrix polynomial (7.5)
$$ \tilde{C}(s) \; =\; \left [ \begin{array}{cr} 0 & 2 \\ 0 & 0
       \end{array} \right ] + \left [ \begin{array}{cr} 0 & -1 \\ 0 & 0
       \end{array} \right ]s + \left [ \begin{array}{cc} 1 & 0 \\ 1 & 1
       \end{array} \right ]s^2 \; = \;\left [ \begin{array}{cc} s^2 & 2-s \\
       s^2 & s^2\end{array} \right ]   $$
and using   (7.6) determine the zero polynomial
$$   \psi(s) \;=\; s^{4-6}det\left [ \begin{array}{cc} s^2 & 2-s \\
                 s^2 & s^2\end{array} \right ] \;=\; s^2 + s -2  $$
For testing we calculate
$$ det\bar{P}(s) \; =\; det\left [ \begin{array}{cccrrr}
         s-2 & -1  & 0 & 0 & -1 & 0 \\ 0 & s-1 & 0 & -1 & 0 & 0 \\
          0 & -2 & s & 0 & 0 & 0 \\-1 & -1 & 0 & s & 0 &-1\\
           1 & -1 & 1 & 0 & 0 & 0 \\ 1 & 1 & 0 & 1 & 0 & 0
          \end{array} \right ] \;=\; s^2 +s - 2  $$

     Now we consider the general case  $l>r$.

\it{THEOREM 7.2.} \rm  System zeros  of  system  (1.1), (1.2) having more outputs than inputs ($l>r$) coincide
with zeros of the polynomial $\psi(s)$ that is the greatest common divisor of all non identically  zero  minors
of the $l\times r$ polynomial matrix $ s^{n-r\nu}det\tilde{C}(s)$ of the order $r$.

\it{PROOF} \rm[S4].$\;$ For  system (7.1) we construct
$(n+l)\times (n+r)$ system matrix $\bar{P}(s)$ (7.3) having the
normal rank  $\rho =n+min(r,l)=n+r$ and consider all its non
identically zero minors  of the form
$$  \bar{P}(s)^{1,2,\ldots,n,n+i_1,\ldots, n+ i_r}_{1,2,\ldots,n,n+1,
         \ldots,n+r}, \qquad i_k \in \{1,2,\ldots.l \}, \; k=1,2,\ldots,r $$
Let  $\psi(s)$ is the greatest common divisor of these minors. System zeros are zeros of $\psi(s)$ by Definition
5.2 and the invariance property of zeros. To calculate these minors we
     represent
$$  \bar{P}(s)^{1,2,\ldots,n,n+i_1,\ldots, n+ i_r}_{1,2,\ldots,n,n+1,
       \ldots,n+r} \; = \; det \left [ \begin{array}{cc} sI-F & -G \\
       \bar{C}(s)^{i_1,\ldots, i_r}  & 0 \end{array} \right ]
\eqno(7.18)  $$
 where  $ 1 \le i_1 \le i_2 \le \cdots \le i_r \le
l$ and the $r\times n$ matrix $\bar{C}(s)^{i_1,\ldots, i_r}$ is constructed from  the $l\times n$ matrix
$CN^{-1}$ by deleting all  rows except  $i_1, i_2, \ldots, i_r$. By using formulas (7.9)-(7.14) we  calculate
$$  \bar{P}(s)^{1,2,\ldots,n,n+i_1,\ldots, n+ i_r}_{1,2,\ldots,n,n+1,
       \ldots,n+r} \; = \; s^{n-r\nu}det \{[O,\bar{C}_1] + [O,\bar{C}_2]s +
      \cdots + \bar{C}_{\nu}s^{\nu -1} \} detG_{\nu}
\eqno (7.19) $$
 where $\bar{C_i}$, $i=1,2,\ldots,\nu$ are $r\times
l_i$ blocks of  the matrix $\bar{C}(s)^{i_1,i_2,\ldots, i_r} = [\bar{C}_1,\bar{C}_2, \ldots,\bar{C}_{\nu}]$.

     On the other hand  maximal order ($r$)  minors of the $l\times r$
polynomial matrix $ s^{n-r\nu}det\tilde{C}(s)$  constructed by deleting all rows except rows $i_1, i_2, \ldots,
i_r$ are
$$ s^{n-r\nu}\tilde{C}(s)^{i_1,i_2,\ldots, i_r} \; = \; s^{n-r\nu}
     det \{[O,\bar{C}_1] + [O,\bar{C}_2]s +\cdots + \bar{C}_{\nu}s^{\nu -1} \}
\eqno(7.20) $$ where $\bar{C_i}$, $i=1,2,\ldots,\nu$ are $r\times l_i$ blocks, which have been defined above.

     Substituting the left-hand side of (7.20) into the right-hand
side of (7.19) we obtain
$$  \bar{P}(s)^{1,2,\ldots,n,n+i_1,\ldots, n+ i_r}_{1,2,\ldots,n,n+1,
       \ldots,n+r} \; = \;s^{n-r\nu}\tilde{C}(s)^{i_1,i_2,\ldots, i_r}
        detG_{\nu}
\eqno (7.21)  $$
 Since $detG_{\nu} \neq 0$ then the greatest
common divisor of minors $ \bar{P}(s)^{1,2,\ldots,n,n+i_1,\ldots,
n+ i_r}_{1,2,\ldots,n,n+1, \ldots,n+r}$, which is equal to
$\psi(s)$, coincides with the greatest common divisor of minors
$s^{n-r\nu}\tilde{C}(s)^{i_1,i_2,\ldots, i_r}$. This proves the
theorem.

     Now we consider the case $l<r$. If  the pair of matrices $(A,C)$
is completely observable then the  pair of $(A^T,C^T)$ is  completely controllable.  Thus,  we  can find the
index observability $\alpha$, integers $\bar{l}_1 \le \bar{l}_2 \le \cdots \le \bar{l}_{\alpha} = l$  and the
nonsingular $n\times n$  matrix $N^*$ that reduces the pair $(A^T,C^T)$ to Yokoyama's canonical
 form. Calculating the $r\times n$ matrix $B^TN^{*-1}$,  partitioning its into
 $r\times l_i$ blocks $\bar{B}_i$ ($i=1,2,\ldots, \alpha$)
$$   B^TN^{*-1} \; =\; [\bar{B}_1,\bar{B}_2,\ldots,\bar{B}_{\alpha}]  $$
we can construct the following $r\times l$  matrix polynomial of the order $\alpha -1$
$$      \tilde{B}(s) = [O,\bar{B}_1] + [O,\bar{B}_2]s + \cdots +
    [O,\bar{B}_{\alpha -1}]s^{\alpha -2} + \bar{B}_{\alpha}s^{\alpha -1} $$
We obtain the dual theorem.

\it{THEOREM 7.3.} \rm System zeros  of system  (1.1), (1.2) having more inputs than outputs $(r>l)$ coincide
with zeros of the polynomial $\bar{\psi}(s)$ that is a greatest common divisor of all $l$ order non identically
zero  minors of the  $r\times l$ polynomial matrix $s^{n-l\alpha}\tilde{B}(s)$ of the order  $l$.

\it{COROLLARY 7.1.} \rm    Invariant  zeros   of    controllable
or observable system (1.1),(1.2) with $l>r$ or $r>l$  coincides
with zeros of the polynomial
$$\psi_I(s) = \epsilon_1(s) \epsilon_2(s) \cdots \epsilon_{\rho}(s)
\eqno(7.22)  $$
 where  $\rho = min(r,l), \;\epsilon_i(s)$ are
invariant  polynomials of   matrices $s^{n-r\nu}\tilde{C}(s)$ or
$s^{n-l\alpha}\tilde{B}(s)$ respectively.

\it{COROLLARY 7.2.} \rm  Invariant zeros of a controllable
(observable) system with $l \ge r$ $(l \le r)$ and $ n=r\nu \; (
n=l\alpha)$  coincide with  zeros of all invariant polynomials of
$\tilde{C}(s) (\tilde{B}(s))$, taken  all together.

\it{COROLLARY 7.3.} \rm  Transmission  zeros of a controllable and
observable system with $l \ge r$ $(l \le r)$ and $ n=r\nu \; (
n=l\alpha)$ coincide with zeros of all invariant polynomials of
$\tilde{C}(s)(\tilde{B}(s))$,
 taken  all together.

    The last result (Corollary 7.3) has been obtain  in Section 4.3 by the alternative way.

\section[ Markov's parameter matrices]
      {Markov's parameter matrices}

     In Section 7.1 we used the matrix polynomial
$\tilde{C}(s)$ (7.5) for zeros definition. Let's scrutinize block coefficients $C_1,\ldots, C_{\nu}$ of the
matrix polynomial $\tilde{C}(s)$. We may show [S10] that the mentioned coefficients are directly expressed via
matrices $A$, $B$, $C$ of system (1.1), (1.2). At first study  the case $l_1=l_2= \ldots=l_{\nu}=r, n=r\nu$.

\it{1. ASSEO'S FORM}.\rm $\;$  Let's partition  the matrix $N^{-1}$ on $n\times r$ blocks $R_i$ ($i=1,2,\ldots,
\nu$)
$$ N^{-1}=[R_1,R_2,\ldots,R_{\nu}]
\eqno(7.23) $$
 We will seek a structure of blocks $R_i$ based on
the relation  (7.2) where submatrices $F_{\nu 1},F_{\nu 2},\ldots,F_{\nu \nu}$ are known ones and $G_{\nu} =
I_r$.

     Since here  $M = I_r$ then the formula  $G =NB$ (7.2) may be used
to express the  matrix $B$  as
$$  B = N^{-1}G = [R_1, R_2, \ldots,R_{\nu} ] \left [ \begin{array}{c}
 O \\ I_r \end{array} \right ] = R_{\nu}
\eqno (7.24)$$
 Hence, the last block in (7.23) is
$$        R_{\nu} = B
\eqno (7.25) $$
 To find others blocks $R_1,R_2,\ldots,R_{\nu}$ we
 use the  relationship $AN^{-1} = N^{-1}F$ (see Eqn.(7.2)), which
is rewritten in the form
$$  A[R_1,R_2,\ldots,R_{\nu}] \; = \; [R_1,R_2,\ldots,R_{\nu}]F \eqno (7.26) $$
 At
first we find a structure of the product $[R_1,R_2,\ldots,R_{\nu}]F$. Since the matrix $F$ is in the form of the
block companion matrix  (1.20) with $p=\nu, \;-T_p =F_{\nu 1},\ldots, -T_1 =F_{\nu \nu}$ then
$$    [R_1,R_2,\ldots,R_{\nu}]F = [ R_{\nu}F_{\nu 1}, \; R_1 + R_{\nu}F_{\nu 2},
         \ldots, \; R_{\nu -1} + R_{\nu}F_{\nu \nu}]
\eqno (7.27)  $$ Using relations (7.26) and (7.27) we can express
blocks $AR_i$ via  $R_{i-1}$ and  $R_{\nu}$ as
$$  AR_i = R_{i-1} + R_{\nu}F_{\nu i}, \qquad i=\nu,\nu -1,\ldots, 2
\eqno  (7.28) $$
Thus, the following recurrent formula follows for
$R_{i-1}$
$$  R_{i-1} = AR_i - R_{\nu}F_{\nu i}, \qquad i=\nu,\nu -1,\ldots, 2
\eqno (7.29) $$
 Since $R_{\nu} = B $ (7.25) then using (7.29)  we
can successively  calculate
$$  \begin{array}{ccl}
    R_{\nu -1}  & = & AR_{\nu} - R_{\nu}F_{\nu \nu} = AB - BF_{\nu \nu},\\
 R_{\nu -2} & = & AR_{\nu -1} - R_{\nu}F_{\nu, \nu -1} = A(AB - BF_{\nu \nu})
               -BF_{\nu, \nu -1} =  A^2B - ABF_{\nu \nu} - BF_{\nu, \nu -1}  \\
              & \vdots & \\
    R_{\nu -i} & = & AR_{\nu -i+1} - R_{\nu}F_{\nu, \nu -i+1} =
                      A^iB - A^{i-1}BF_{\nu \nu} - \cdots -BF_{\nu, \nu -i+1}\\
              & \vdots & \\
     R_1 & = & A^{\nu -1}B - A^{\nu -2}BF_{\nu \nu} - \cdots - ABF_{\nu 3}
               -BF_{\nu 2}
    \end{array}
\eqno  (7.30)   $$
     Substituting the matrix  $N^{-1}$ (7.23) in (7.4) we  present
blocks $C_i$ $(i=1,2, \cdots, \nu )$ as
 $$ C_i =CR_i
 \eqno(7.31)$$
So, the matrix polynomial (7.5) becomes for case $l_1=l_2=\cdots=l_{\nu}=r$
$$
 \tilde{C}(s) = CR_1 + CR_2s + \cdots + CR_{\nu -1}s^{\nu -2} +
       CR_{\nu}s^{\nu -1}
\eqno(7.32)  $$
 Then substituting the right-hand side of (7.30) into
(7.32) we obtain  the matrix polynomial $\tilde{C}(s)$ in the final form
$$  \tilde{C}(s) = (CA^{\nu -1}B - CA^{\nu -2}BF_{\nu \nu} - \cdots
- CABF_{\nu 3} -CBF_{\nu 2} ) + (CA^{\nu -2}B -
$$
$$ -CA^{\nu -3}BF_{\nu \nu} - \cdots - CBF_{\nu 3})s + \cdots + (CA^2B - CABF_{\nu
\nu}-
$$
$$- CBF_{\nu, \nu-1})s^{\nu -3} +
          (CAB - CBF_{\nu \nu})s^{\nu -2} + CBs^{\nu -1}
 \eqno (7.33)    $$
     Hence, block coefficients of the matrix polynomial $\tilde{C}(s)$ of
 controllable system (1.1),(1.2)  with $n=r\nu$ are expressed  via
the matrices $CB, CAB,\ldots$  that are  blocks  of  the 'so-called'  output controllable matrix $[ CB, CAB,
\ldots , CA^{n-1}B ]$ [D1]. These matrices are known as Markov parameter matrices  (or Markov parameters in the
classic single-input/single output system).

If $r=1$ then $l_1=l_2=\cdots=l_{\nu}=r$, $\nu = n$ and the polynomial  $l$   vector $\tilde{C}(s)$ has the
following simple structure
$$  \tilde{C}(s) = (CA^{n -1}b - CA^{n -2}b\alpha_n - \cdots
      - CAb\alpha_3 -Cb\alpha_{2} ) + (CA^{n -2}b - CA^{n -3}b\alpha_n
       - \cdots - Cb\alpha_3)s + \cdots  $$
$$      + (CA^2b - CAb\alpha_n - Cb\alpha_{n-1})s^{n -3} +
          (CAb - Cb\alpha_n)s^{n -2} + Cbs^{n -1}
 \eqno (7.34)    $$
where $\alpha_2,\ldots, \alpha_n$ are  coefficients of the characteristic polynomial of  $A$: $ det (sI_n-A) =
s^n - \alpha_ns^{n-1} - \cdots - \alpha_2 s - \alpha_1$.

\it{2. YOKOYAMA'FORM}.\rm $\;$  Consider the general case: $ l_1 \le l_2 \le \cdots \le l_{\nu} = r$, $n <r\nu $
when the pair of matrices $(A,B)$  is reduced to Yokoyama's canonical form.
     We also apply the partition (7.23) with $\nu$  blocks $R_i$
of sizes $n\times l_i$, $i=1,2,\ldots,\nu$.

From  the formula  $G =NBM$ (7.2) we express the matrix $B$ as $B=N^{-1}GM^T$ and, using the special  structure
of the matrix $G =\left [ \begin{array}{c} O \\ G_{\nu} \end{array} \right ]$ and the
 partition (7.23), obtain
 $$   B = R_{\nu} G_{\nu}M^T
\eqno   (7.35)  $$
 From  (7.35) we get the last block $R_{\nu}$ of the  matrix  $N^{-1}$
$$  R_{\nu} = B\bar{G}_{\nu}^{-1}
\eqno (7.36) $$
 where $G_{\nu}M^T =\bar{G}_{\nu}$. For finding $n\times r$  blocks $[O,R_1],[O,R_2]
,\ldots,[O,R_{\nu -1}]$ in (7.31) we also  apply the equality  (7.26).  Using  the special structure of the
matrix $F$, which is in the form  of  the general block companion matrix (1.25) with $p=\nu,\; -\hat{T}_p
=F_{\nu 1}, \ldots, -\hat{T}_1 =F_{\nu \nu}$,  we can write
$$  AR_i = R_{i-1}[O,I_{l_{i-1}}] + R_{\nu}F_{\nu i},
\qquad i=\nu,\nu -1,\ldots, 2
\eqno  (7.37) $$
where $[O,I_{l_{i-1}}]$  are $l_{i-1}\times l_i$ matrices.

    From (7.37) and (7.36) we can express
blocks $R_{i-1}[O,I_{l_{i-1}}]$ in terms of  blocks $R_i$, $A$,
$R_{\nu}$ and $F_{\nu i}$ as follows
$$   R_{i-1}[O,I_{l_{i-1}}] =  AR_i - B\bar{G}_{\nu}^{-1}F_{\nu i},
      \qquad i=\nu,\nu -1,\ldots,2
\eqno  (7.38) $$
 The recurrent formula (7.38) is  used for finding
 $n\times r$ blocks  $[O,R_{i-1}] $ by varying $i$ from $\nu$ to
$2$.
     At first we determine  the $n\times r$ matrix $[O,R_{\nu-1}] =
 R_{\nu-1}[O,I_{l_{\nu-1}}]$. With $i=\nu$ the relation (7.38) becomes
$$   R_{\nu-1}[O,I_{l_{\nu-1}}] =  AR_{\nu} - B\bar{G}_{\nu}^{-1}F_{\nu \nu},
\eqno  (7.39) $$
 Taking into account the expression (7.36) we obtain
$$  [O,R_{\nu-1}] =  AB\bar{G}_{\nu}^{-1} - B\bar{G}_{\nu}^{-1}F_{\nu \nu},
  \eqno  (7.40) $$
Then with $i=\nu -1$ the relation (7.38) becomes
$$       R_{\nu-2}[O,I_{l_{\nu-2}}] =  AR_{\nu -1} -
         B\bar{G}_{\nu}^{-1}F_{\nu, \nu -1},
 \eqno  (7.41) $$
Let's postmultiply  both sides of (7.41) by the $l_{\nu -1}\times r$ matrix $[O,I_{l_{\nu -1}}]$. Since
$[O,I_{l_{\nu -2}}][O,I_{l_{\nu -1}}] = [O,I_{l_{\nu -2}}]^*$ where $[O,I_{l_{\nu -2}}]^*$ is the matrix of
sizes $l_{\nu -2}\times l_{\nu} = l_{\nu -2}\times r $ then  (7.41)  takes the form
$$       R_{\nu-2}[O,I_{l_{\nu-2}}]^* =  AR_{\nu -1}[O,I_{l_{\nu -1}}] -
         B\bar{G}_{\nu}^{-1}F_{\nu, \nu -1}[O,I_{l_{\nu -1}}],
 \eqno  (7.42) $$
Matrices $R_{\nu-2}[O,I_{l_{\nu-2}}]^* = [O,R_{\nu-2}]$ and $\;
 F_{\nu, \nu -1}[O,I_{l_{\nu -1}}]  = [O,F_{\nu, \nu -1}]$ are $n\times r$
 and $r\times r$ matrices. Using formulas (7.39), (7.40) we obtain
  the  $r\times r$ matrix  $[O,R_{\nu-2}]$
$$   [O,R_{\nu-2}] = A^2B\bar{G}_{\nu}^{-1} - AB\bar{G}_{\nu}^{-1}F_{\nu \nu}
     - B\bar{G}_{\nu}^{-1}[O,F_{\nu, \nu -1}]
\eqno  (7.43)  $$
  Continuing these reasonings we  determine $n\times r$ matrices $[O,R_{\nu-3}] ,\ldots,
[O,R_1]$
$$ \begin{array}{ccl}
 [O,R_{\nu-3}] & = & A^3B\bar{G}_{\nu}^{-1} - A^2B\bar{G}_{\nu}^{-1}F_{\nu \nu}-
 B\bar{G}_{\nu}^{-1}[O,F_{\nu,\nu-1}] - B\bar{G}_{\nu}^{-1}[O,F_{\nu, \nu -2}]\\
         & \vdots & \\   $$
$$  [O,R_1]  &  =  & A^{\nu -1}B\bar{G}_{\nu}^{-1} -
     A^{\nu -2}B\bar{G}_{\nu}^{-1}F_{\nu \nu} - \cdots -
     AB\bar{G}_{\nu}^{-1}[O,F_{\nu3}] - B\bar{G}_{\nu}^{-1}[O,F_{\nu2}]
\end{array}
\eqno(7.44) $$
 and substituting  $R_{\nu}$ (7.36) and $[O,R_{\nu-1}], [O,R_{\nu-2}], \ldots, [O,R_1]$   in
(7.31) find the matrix polynomial $\tilde{C}(s)$ of the general structure
$$  \tilde{C}(s) = (CA^{\nu -1}B\bar{G}_{\nu}^{-1} -
               CA^{\nu -2}B\bar{G}_{\nu}^{-1}F_{\nu \nu} - \cdots
     - CAB\bar{G}_{\nu}^{-1}[O,F_{\nu 3}] - CB\bar{G}_{\nu}^{-1}[O,F_{\nu 2}])
     + (CA^{\nu -2}B\bar{G}_{\nu}^{-1} - $$
$$   - CA^{\nu -3}BF_{\nu \nu} - \cdots - CB\bar{G}_{\nu}^{-1}[O,F_{\nu 3}])s +
       \cdots + (CA^2B\bar{G}_{\nu}^{-1} - CAB\bar{G}_{\nu}^{-1}F_{\nu \nu} -
        CB\bar{G}_{\nu}^{-1}[O,F_{\nu, \nu-1}])s^{\nu -3} +  $$
$$      (CAB\bar{G}_{\nu}^{-1} - CB\bar{G}_{\nu}^{-1}F_{\nu \nu})s^{\nu -2}
         + CB\bar{G}_{\nu}^{-1}s^{\nu -1}
 \eqno (7.33)    $$
where $[O,F_{\nu i}]$  are $r\times r$ matrices.

So,   here  we also  reveal the dependence of block coefficients of  $\tilde{C}(s)$ upon
 matrices $CB$, $CAB$,$\ldots$. In contrast to the first case,
this connection has the more complicated form.

    To improve the computation
accuracy it is desirable to use $\tilde{C}(s)$  in the form (7.5) [S10]
$$  \tilde{C}(s) = [O,(CA^{\nu -1}B\bar{G}_{\nu}^{-1} -
               CA^{\nu -2}B\bar{G}_{\nu}^{-1}F_{\nu \nu} - \cdots
     - CAB\bar{G}_{\nu}^{-1}[O,F_{\nu 3}] - CB\bar{G}_{\nu}^{-1}[O,F_{\nu 2}])
      \left [ \begin{array}{c} O\\ I_{l_1} \end{array} \right ] ] $$
$$     + [O,(CA^{\nu -2}B\bar{G}_{\nu}^{-1} - CA^{\nu -3}BF_{\nu \nu} -
       \cdots - CB\bar{G}_{\nu}^{-1}[O,F_{\nu 3}])
       \left [ \begin{array}{c} O\\ I_{l_2} \end{array} \right ]]s+ $$
$$        +\cdots + [O,(CAB\bar{G}_{\nu}^{-1} - CB\bar{G}_{\nu}^{-1}F_{\nu \nu})
      \left [ \begin{array}{c} O\\ I_{l_{\nu -1}}\end{array} \right ]]s^{\nu-2}
         + CB\bar{G}_{\nu}^{-1}s^{\nu -1}
 \eqno (7.46)  $$
 where  matrices $[O,I_{l_i}]^T$ have  sizes  $r\times l_i$.

     Let's consider several  examples.

                      \it{EXAMPLE 7.2.} \rm

 At first we  find  a  zero  polynomial of system (1.1),(1.2) with  $n = 4,
r = l = 2$ and the state space model matrices
$$   A \; = \;  \left[ \begin{array}{cccc} 2 & 1 & 0 & 1 \\ 1 & 0 & 1 & 1 \\
     1 & 1 & 0 & 0 \\ 0 & 0 & 1 & 0 \end{array} \right ], \qquad
     B\; = \; \left[ \begin{array}{cc} 0 & 0 \\ 1 & 0 \\ 0 & 1 \\ 0 & 0
            \end{array} \right ], \qquad
     C \; = \;  \left[ \begin{array}{cccc} 1 & 1 & 0 & 0 \\ 0 & 0 & 1 & 1
             \end{array} \right ]
\eqno (7.47) $$
Since
$$   rank [ B, AB, A^2B,A^3B ] \; = \; rank [ B, AB ] \; = \; 4     $$
then $\nu = 2$,  $n=r\nu $. and the system is reduced to Asseo's canonical form (7.1) with the matrices $F$ and
$G$ of the following structure
$$   F \; = \;  \left[ \begin{array}{cc} O & I_2\\ F_{21} & F_{22}
          \end{array} \right ], \qquad G \; = \;  \left[ \begin{array}{c}
          O \\ I_2 \end{array} \right ]
\eqno (7.48)  $$ where $I_2$ is  the $2\times 2$ unit matrix, $F_{21} , F_{22} $ are $2\times 2$
 submatrices. Using  (7.33) we  write
$$   \tilde{C}(s) = (CAB -CBF_{22}) +CBs
\eqno (7.49) $$
 To  determine $F_{22}$ we need at first to find
the $4\times 4$ transformation matrix $N$ that  reduces the system to Asseo's canonical form. By formulas of
Section 1.2.3 we find
$$    N \; = \;  \left[ \begin{array}{c}
          N_1 \\ N_2 \end{array} \right ]
\eqno (7.50) $$
 where an $2\times 4$ matrix $N_2$ is calculated from
 (1.54)
$$   N_2 \;=\;[ O, I_2 ][B,AB]^{-1}
\eqno  (7.51) $$ and an $2\times 4$ matrix $N_1$  is found from (1.49)
$$      N_1 \; = \; N_2A
\eqno(7.52)  $$
Determining
$$  [B,AB]^{-1} \; =\; \left[ \begin{array}{rrrr} 0 & 1 & 0 & -1 \\
                    -1 & 0 & 1 & 0 \\ 1 & 0 & 0 & 0 \\ 0 & 0 & 0 & 1
                    \end{array} \right ]    $$
$$  N_2\;=\; \left [ \begin{array}{crcc} 1 & 0 & 0 & 0 \\ 0 & 0 & 0 & 1
     \end{array} \right ], \qquad  N_1 =N_2A = \left [ \begin{array}{crcc}
     2 & 1 & 0 & 1 \\ 0 & 0  & 1 & 0 \end{array} \right ]    $$
we obtain from (7.50)
$$  N \; =\; \left[ \begin{array}{rrrr} 1 & 0 & 0 & 0 \\
                    0 & 0 & 0 & 1 \\ 2 & 1 & 0 & 1 \\ 0 & 0 & 1 & 0
                    \end{array} \right ]    $$
Then using the relation $F = NAN^{-1}$  we calculate
$$  F_{21}\;=\; \left [ \begin{array}{rr} 1 & 1\\ -1 & -1\end{array} \right ],
       \qquad  F_{22} \;=\;  \left [ \begin{array}{cc} 2 & 2\\ 1 & 0
         \end{array} \right ]    $$
Substituting $F_{22}$ into (7.49) and calculating $CB$, $CAB$  yields
the following matrix polynomial
$$ \tilde{C}(s) \; =\; \left(\left [ \begin{array}{cr} 1 & 1 \\ 1 & 1
       \end{array} \right ] - \left [ \begin{array}{cr} 2 & 2 \\ 0 & 1
    \end{array} \right ] \right ) + \left [ \begin{array}{cc} 1 & 0 \\ 0 & 1
       \end{array} \right ]s \; = \;\left [ \begin{array}{rr} -1 & -1 \\
       0 & 1\end{array} \right ] +  \left [ \begin{array}{cc} 1 & 0 \\ 0 & 1
          \end{array} \right ]s
\eqno(7.53) $$
       Thus  the zero polynomial is
$$ \psi(s) \;=\; det\tilde{C}(s) \; =\; s^2 -1
\eqno (7.54)  $$ To check the results obtained we calculate the
determinant of  the system matrix
$$ detP(s) \; =\; det\left [ \begin{array}{crrrrr}
         s-2 & -1  & 0 & -1 & 0 & 0 \\ -1 & s & -1 & -1 & -1 & 0 \\
          -1 & -1 & s & 0 & 0 & -1 \\0 & 0 & -1 & s & 0 & 0 \\
           1 & 1 & 0 & 0 & 0 & 0 \\ 0 & 0 & 1 & 1 & 0 & 0
          \end{array} \right ]  = s^2 -1  $$

                   \it{EXAMPLE 7.3.} \rm

  Let's determine a zero polynomial of  the system  from Example  7.1. Since
this system  has $n = 4, \; r = l= 2, \; \nu  = 3,\; l_1 = l_2 = 1, \; l_3 = 2 $ then using  (7.45) we obtain
the general structure of the matrix $\tilde{C}(s)$
$$  \tilde{C}(s) = (CA^2B\bar{G}_{\nu}^{-1} - CAB\bar{G}_{\nu}^{-1}[O,F_{33}] -
                    CB\bar{G}_{\nu}^{-1}[O,F_{32}]) +
                 (CAB\bar{G}_{\nu}^{-1} - CB\bar{G}_{\nu}^{-1}F_{33})s
                  + CB\bar{G}_{\nu}^{-1}s^2
 \eqno (7.55)    $$
where $F_{32}$, $F_{33}$ are $2\times 1$  and  $2\times 2$
matrices respectively.  For calculating $F_{32}$, $F_{33}$ and
$\bar{G}_{\nu}^{-1}$ we  use  the formulas from (7.2) ($ F =
NAN^{-1}$   and $G = NBM$)  with $N$ (1.93) and $M$ (1.90) (see
Example 1.3). We result in
$$  F_{32} \;=\; \left [ \begin{array}{c} 1 \\ 1 \end{array} \right ], \qquad
    F_{33}\;=\; \left [ \begin{array}{rr} 2 & 0\\ 1 & 1\end{array} \right ],
       \qquad  G_{\nu} \;=\;  \left [ \begin{array}{cc} 1 & 0\\ 0 & 1
         \end{array} \right ]
\eqno(7.56)  $$
 Since  $M =I_2$  then $\bar{G}_{\nu} = G_{\nu}M^T
= G_{\nu}$. Calculating
$$ CB \;=\; \left [ \begin{array}{cc} 1 & 0 \\ 1 & 1 \end{array} \right ],
     \qquad CAB\;=\;\left [ \begin{array}{rr} 2 & -1\\ 3 & 1\end{array} \right ],
       \qquad  CA^2B\;=\;  \left [ \begin{array}{cc} 3 & 2\\ 7 & 3
         \end{array} \right ]    $$
and substituting these matrices and matrices (7.56) in (7.55) we
obtain
$$ \tilde{C}(s) \; =\; \left [ \begin{array}{cc} 0 & 2 \\ 0 & 0
       \end{array} \right ] + \left [ \begin{array}{cr} 0 & -1 \\ 0 & 0
       \end{array} \right ]s  +  \left [ \begin{array}{cc} 1 & 0 \\ 0 & 1
          \end{array} \right ]s^2
\eqno(7.57) $$
 The matrix polynomial (7.57) coincides  with  one
 obtained in Example 7.1.

\section[A number of zeros]
        {A number of zeros}

    In this section  we obtain several upper
bounds of  a  system  zeros number   in terms of the  matrices $CB$, $CAB$,... and the controllability
characteristics of the pair $(A,B)$.

     At first we consider controllable  system  (1.1),(1.2)  with
equal number of independent inputs and outputs ( $r = l$, $rank B = rank C = r$). System zeros of this system
coincide with zeros of the polynomial $ \psi(s) = s^{n-r\nu}det(\tilde{C}(s))$ where the $r\times r$ matrix
$\tilde{C}(s)$
$$
 \tilde{C}(s) = [O,C_1] + [O,C_2]s + \cdots +  [O,C_{\nu -1}]s^{\nu -2}
                +C_{\nu}s^{\nu -1}
\eqno(7.58)  $$
 is the nonmonic matrix polynomial of the degree
$\nu-1$ with $r\times r$ matrices $[O,C_i]$, $i=1,2,\ldots,\nu $ of dimensions $r\times r$. Here  $r\times l_i$
submatrices $C_i$, $i=1,2,\ldots, \nu$   are defined from the relation (7.46) as
$$ \begin{array}{ccl} C_{\nu} & = & CB\bar{G}_{\nu}^{-1}\\
   C_{\nu-1} & = & (CAB\bar{G}_{\nu}^{-1} - CB\bar{G}_{\nu}^{-1}F_{\nu \nu})
 \left [ \begin{array}{c} O\\ I_{l_{\nu -1}}\end{array} \right ] \\
 & \vdots & \\
   C_1 & = & (CA^{\nu -1}B\bar{G}_{\nu}^{-1} -CA^{\nu -2}B\bar{G}_{\nu}^{-1}
           F_{\nu \nu} - \cdots - CB\bar{G}_{\nu}^{-1}[O,F_{\nu 2}])
         \left [ \begin{array}{c} O\\ I_{l_1} \end{array} \right ] \end{array}
\eqno (7.59) $$

We can see from (7.58) that  a number of zeros of the polynomial $ \psi(s) = s^{n-r\nu}det(\tilde{C}(s))$
depends on a rank of  the matrix $C_{\nu}  = CB\bar{G}_{\nu}^{-1}$.

\it{ASSERTION 7.1.} \rm  If  the $r\times r$ matrix $C_{\nu }$ has the full rank $(rank ( C_{\nu }) = r)$ then
system (1.1),(1.2) has exactly $n-r$  zeros.

\it{PROOF}.$\;$ \rm If  $rank (C_{\nu }) = r$ then there exists  the matrix $C_{\nu }^{-1}$ and the polynomial
$\psi(s)$  may be represent as
$$ \psi(s) \;=\; s^{n-r\nu}det \{C_{\nu}([O,Q_1] + [O,Q_2]s + \cdots +
             [O,Q_{\nu -1}]s^{\nu -2} + I_rs^{\nu -1} ) \} \; = \;
              det(C_{\nu})s^{n-r\nu}\bar{\psi}(s)
\eqno  (7.60)  $$
where
$$          Q_i = C_{\nu }^{-1}C_i, \;\; i=1,2,\ldots,\nu -1
\eqno (7.61)   $$
$$ \bar{\psi}(s) \;=\; det ([O,Q_1] + [O,Q_2]s + \cdots +
             [O,Q_{\nu -1}]s^{\nu -2} + I_rs^{\nu -1})
\eqno (7.62)  $$
 Since $C_{\nu}$ is the constant and nonsingular $r\times r$ matrix then  zeros of $\psi(s)$
coincide with zeros of the polynomial $s^{n-r\nu}\bar{\psi}(s)$ of the degree $r(\nu -1)$ because
  $[O,Q_1] + [O,Q_2]s + \cdots + [O,Q_{\nu
-1}]s^{\nu -2} + I_rs^{\nu -1}$ is the  monic $r\times r$ matrix polynomial of the degree  $\nu -1$. Calculating
the degree ($ \xi $) of the  polynomial $\psi(s) = s^{n-r\nu}\bar{\psi}(s)$ we obtain that $ \xi =
n-r\nu+r(\nu-1) = n-r$. Therefore,  the polynomial $\psi(s)$ has exactly $n-r$ zeros.

\it{ASSERTION 7.2.} \rm  If  matrix  $C_{\nu}$ has the rank
deficiency ($d$) then system (1.1),(1.2) has no more than $n-r-d$
system zeros.

\it{PROOF}.$\;$ \rm It is known that  any  constant $r\times r$ matrix $C_{\nu}$ of  rank deficiency $d < r$ may
be reduced by the series of elementary row and column operations  to the form
$$  U_LC_\nu U_R \; = \; \left [ \begin{array}{rr} I_{r-d} & O \\ O & O
  \end{array} \right ]  $$
where $U_L$, $U_R$   are unimodular  $r\times r$  matrices.

     Let's consider the polynomial $\tilde{\psi}(s)$ $=$
$s^{n-r\nu}det(U_L\tilde{C}(s)U_R)$ where $\tilde{C}(s)$ has form (7.58). It is evident that polynomials
$\psi(s)$ and $\tilde{\psi}(s)$ have  same zeros. Representing $\tilde{\psi}(s)$ as
$$   \tilde{\psi}(s) \; = \; s^{n-r\nu }det \{(U_L[O,C_1] + [O,C_2]s + \cdots +
      [O,C_{\nu -1}]s^{\nu -2} + C_{\nu}s^{\nu -1}) U_L\}
\eqno (7.63)
 $$ and performing the multiplications we get
$$   \tilde{\psi}(s) \; = \; s^{n-r\nu }det \{U_L[O,C_1]U_R + U_L[O,C]U_Rs + \cdots +
      U_L[O,C_{\nu -1}]U_Rs^{\nu -2} + U_LC_{\nu}U_Rs^{\nu -1} \}
\eqno (7.64)
$$ The analysis of the $r\times r$ matrix polynomial in the above expression  shows that its
diagonal elements  include $r-d$ polynomials of degrees $\nu -1 $ and $d$ polynomials of degrees that less the
value  $\nu -2 $. As a result the maximal degree of $\tilde{\psi}(s)$ is equal to $n-r\nu +(r-d)(\nu -1) + d(\nu
-2) = n-r-d \;$.

     From Assertions 7.1 and 7.2 we obtain

\it{ASSERTION 7.3.} \rm  If  the matrix $C_{\nu}$  is equal to the zero matrix i.e. $C_{\nu}=O$ (or $C_{\nu}$
has the rank deficiency $r$) then  we will have  following  variants:

 1. $\;l_{\nu -1} =r$, $rank(C_{\nu -1})=r$. The system has exactly $n-2r$ zeros.

 2. $\;l_{\nu -1} =r$, $rank(C_{\nu -1}) < r$. The system has no more than
$n-2r-\bar{d}$  zeros where $\bar{d} <r$ is a rank deficiency of $C_{\nu -1}$.

 3. $\;l_{\nu -1} < r$, $rank(C_{\nu -1}) = l_{\nu -1}$. The system has no more
than $n-2r-(r-l_{\nu -1}) = n-3r + l_{\nu -1}$ zeros.

 4. $\;l_{\nu -1} < r$, $rank(C_{\nu -1}) < l_{\nu -1}$. The system has no more
than $n-2r-(r-l_{\nu -1}) -\bar{d} = n-3r + l_{\nu -1} - \bar{d}$ zeros where $\bar{d} <r$ is a rank deficiency
of $C_{\nu -1}$.

     Let us study  the  relationship  between ranks of matrices $C_{\nu}$ and  $CB$.

\it{ASSERTION 7.4.} \rm  For  controllable system (1.1), (1.2)
with $l\ge r$ the following rank equality is true
$$           rankC_{\nu} \;=\; rank (CB)
\eqno (7.65) $$

     The proof follows from  the equality  $C_{\nu} \; = \; CB\bar{G}_{\nu}^{-1}$
 because  the $r\times r$ matrix $\bar{G}_{\nu} = G_{\nu}M^T$ is nonsingular one.

     Using  Assertion  7.4  and  formula   (7.59)   we   can
reformulate Assertions 7.1 - 7.3 in terms  of $CB,\; CAB,\;\ldots
$.

\it{THEOREM 7.5.} \rm  A number of system zeros of controllable system (1.1), (1.2) with $r=l$ is defined via
ranks of matrices $CB,\; CAB,\;\ldots $  as follows:

 1. if $rank(CB) = r$  then the system has exactly $n-r$ zeros.

 2. if  the matrix $CB$ has a rank deficiency $d$ then the
    system has no more than  $n-r-d$  zeros.

 3. if  the  matrix $CB$ is the zero matrix then:

   $\;$ 3.1. if $l_{\nu -1} =r$ and $rank(CAB) = r$ then the system has
         exactly $n-2r$   zeros,

   $\;$ 3.2. if $l_{\nu -1} =r$ and $rank(CAB) < r$ with a rank deficiency
         $\bar{d} < r$ then the system has no more than $n-2r-\bar{d}$ zeros,

 $\;$ 3.3. if $l_{\nu -1} < r$ and $rank(CAB) = l_{\nu -1}$ then the system has
         no more than $n-3r+l_{\nu-1}$ zeros,

   $\;$ 3.4. if if $l_{\nu -1} < r$ and $rank(CAB) < l_{\nu -1}$ with a rank
         deficiency $\bar{d} < r$ then the system has no more than
         $n-3r + l_{\nu -1} - \bar{d}$  zeros.
and etc.

     Let's consider the system with an unequal number of inputs and
outputs. For definiteness  we assume  $l>r$. By Theorem 7.2  system zeros coincide with zeros of a greatest
common divisor of all nonzero minors of the order $r$ of the  $l\times r$ polynomial matrix
$s^{n-r\nu}det\tilde{C}(s)$. As there are several such  minors then we may get only upper bounds on a number of
zeros that follow from Assertions 7.1 - 7.3 and Theorem 7.5. We select the most important cases.

\it{COROLLARY 7.4. } \rm There are following upper bounds on a  number of system zeros of controllable system
(1.1),(1.2) with $l > r$ :

 1. if $rank(CB) = r$  then the system has no more than $n-r$  zeros.

 2. if the matrix $CB$  has a rank deficiency  $d$ then the
    system has no more than  $n-r-d$    zeros.

 3. if $CB$ is equal to the zero matrix then:

  $\;$  3.1. if $l_{\nu -1} =r$ and $rank(CAB) = r$  then the system has no
          more than  $n-2r$    zeros,

    The next cases coincide with the corresponding points  of Theorem 7.5.

    In conclusion we consider  conditions when   system
(1.1),(1.2) with  $n=r\nu $  has no  zeros.  It follows from  $C_i$ (7.59) and Corollary 7.4 that zeros are
absent if the matrix $CA^{\nu -1}B$   has a full rank and matrices $CB, CAB, \ldots$
 are  zero matrices. The last condition  is true if
$$    C[B, AB,\ldots, A^{\nu -2}B] \; = \; O  $$
 Denoting  subspaces, which are formed from  linear independent rows of $C$ and  columns of the  $n\times
(n-r)$ matrix $[B, AB,\ldots, A^{\nu -2}B]$ by $\mbox{R}_c$  and $\mbox{R}_{\nu -2}$ respectively we obtain from
the last reasonings.

\it{ASSERTION 7.5.} \rm  If  the subspace $\mbox{R}_c$ is orthogonal to the subspace $\mbox{R}_{\nu -2}$   and
$rank(CA^{\nu -1}B) = r$ then  system  (1.1), (1.2) with $l \ge r$ and $ n=r\nu $ has no system zeros.

\it{REMARK 7.4.} \rm  If all matrices $CA^iB$, $i=0,1,\ldots, \nu -1$, are  zero  matrices in the mentioned
system  then  the system has zeros everywhere on the complex plan. Such a  system is known as the degenerate
system [D4] \footnote{ In this case a zero polynomial is identically equal to zero. That implies that system
zeros coincide with the whole complex plane.}.

     The completely controllable system with  $C \neq O$ and $r=1$, $\nu =n$ is
always the nondegenerate system. Indeed,  the subspace $\mbox{R}_{\nu -1}$, which  is formed from  columns of
matrix $[b,Ab,\ldots,A^{n-1}b]$, coincides with the complete state space. Therefore, an  intersection of
$\mbox{R}_c$ and $\mbox{R}_{\nu -1}$ is not the empty subspace.

\it{REMARK 7.5.} \rm  A similarly way may be used to derive  on a
number of zeros of  system (1.1), (1.2) with $l<r$.

\section[Zero determination via lower order
                           matrix pencil]
        {Zero determination via lower order\\
                           matrix pencil}

     In this section we reduce the problem of the zero calculation
to  eigenvalues  problem for a  matrix pencil of order  $n-r$ [S8]. We restrict our study by  system (1.1),
(1.2) with equal number of inputs and outputs. The method is based on constructing the generalized block
companion matrix of the structure (1.25) for the matrix polynomial $\tilde{C}(s)$ (7.5).

      We will consider two cases.

\it{CASE 1. $\;$ $CB$ is nonsingular matrix}. \rm Then the $r\times r$ matrix $C_{\nu}  = CB\bar{G}_{\nu}^{-1} $
also has a full rank and there exists the matrix $C_{\nu}^{-1}$. The matrix polynomial $C(s)$ can  be
represented as follows
$$ \tilde{C}(s) = C_{\nu}\{[O,\hat{T}_{\nu -1}] + [O,\hat{T}_{\nu -2}]s +
     \cdots +  [O,\hat{T}_1]s^{\nu -2} + I_rs^{\nu -1}\}
\eqno (7.66) $$
where $r\times l_i$ submatrices $\hat{T}_{\nu -i}$ are defined by the formula
$$   \hat{T}_{\nu -i}\; = \; C_{\nu}^{-1}C_i, \qquad i=1,2, \ldots, \nu -1
\eqno(7.67) $$
 The expression in brackets in (7.66) coincides with the matrix polynomial (1.19) with $p = \nu
-1$ and $T_i = [O,\hat{T}_i]$, $i=1,2,\ldots,\nu -1$.

Let's find a generalized block companion matrix for
$$     \bar{\Phi}(s) = [O,\hat{T}_{\nu -1}] + [O,\hat{T}_{\nu -2}]s +
        \cdots +  [O,\hat{T}_1]s^{\nu -2} + I_rs^{\nu -1}
\eqno(7.68)  $$
 We need to emphasis that here $l_{\nu -1} =r$,
hence the $r\times r$ block $[O,\hat{T}_1]$ coincides with the $r\times r$ matrix $\hat{T}_1$ and the structure
of (7.68) completely coincides with (1.19). Otherwise, ($l_{\nu -1} < r$) the structure of (7.68) is
distinguished from   (1.19)  by  the block $[O,\hat{T}_1]$. We consider these cases separately.

     If  $l_{\nu -1} =r$ and $[O,\hat{T}_1] = \hat{T}_1$ then   we  get
 the generalized block companion
$\bar{n}\times \bar{n}$  matrix of the form (1.25) with $p=\nu -1$ and $\bar{n} = l_1 + l_2 + \cdots +l_{\nu -1}
= n-r$ (see Sect 1.2.1). Let's denote this matrix as $\bar{P}$
$$   \bar{P} = \left [ \begin{array}{c} E \\ \dotfill\\ -\hat{T}_{\nu -1},
              -\hat{T}_{\nu -2},\ldots, -\hat{T}_1 \end{array} \right ]
\eqno (7.69)  $$ where  the $(n-r-l_{\nu -1})\times (n-r)$  block $E$ has the following form
$$
    E = \left[ \begin{array}{ccccc} O & E_{1,2} & O & \cdots & O \\
       O& O& E_{2,3}& \cdots & O\\\vdots & \vdots & \vdots & \ddots & \vdots \\
        O& O& O& \cdots & E_{\nu -2,\nu -1} \end{array} \right]
\eqno (7.70)   $$
with $l_i\times l_{i+1}$ submatrices  $E_{i,i+1} =[O, I_{l_i}]$.

     From Assertion 1.2 and replacing  $\bar{n}$ by $n-r$ in (1.26) and
 $p$ by $\nu -1$ we get the following equality
 $$     det(sI_{n-r} -\bar{P}) = s^{n-r-r(\nu -1)}det\bar{\Phi}(s) =
       s^{n-r\nu }det\bar{\Phi}(s)$$
Since $\bar{\Phi}(s) = C_{\nu}^{-1}\tilde{C}(s)$  then
 $$     det(sI_{n-r} -\bar{P}) = s^{n-r\nu }detC_{\nu}^{-1}det\tilde{C}(s)
\eqno (7.71)   $$
 This result may be formally expressed as

\it{ASSERTION 7.6.} \rm  System zeros of system  (1.1), (1.2) with an equal number of inputs and outputs having
$l_{\nu -1} =r$ and $det(CB)\neq 0 $ are defined as  eigenvalues of the $(n-r)\times (n-r)$
 matrix (7.69).

     Let $l_{\nu -1} < r$.  We introduce  the $l_{\nu -1}\times (n-r)$  matrix
$$  \hat{T} = [O,I_{l_{\nu -1}}][-\hat{T}_{\nu -1},-\hat{T}_{\nu -2},\ldots,
              -\hat{T}_1]
\eqno (7.72)  $$ and prove validity of the following lemma.

\it{LEMMA 7.1.} \rm
$$ det\tilde{C}(s) \; = \; detC_{\nu } s^{r\nu -n}det(sI_{n-r}-
             \left [ \begin{array}{c} E \\ \hat{T}\end{array} \right ] )
\eqno  (7.73) $$
 \it{PROOF}. \rm $\;$  At first we construct the $r\times r$ submatrix
$$ \hat{T}_1^* =   [O,\hat{T}_1]
\eqno (7.74)   $$
 and the $l_{\nu -2}\times r$   submatrix
$$       E_{\nu -2,\nu -1}^o \; = \; [O,  E_{\nu -2,\nu -1}]
\eqno  (7.75) $$
 from blocks  $\hat{T}_1$ and  $E_{\nu -2,\nu -1}$
of matrices (7.69) and (7.70). Then we form the  following $(n-l_{\nu -1})\times (n-l_{\nu -1})$ matrix $P^*$
from $\tilde{P}$  with blocks (7.74), (7.75) instead of $\hat{T}_1 $ and $E_{\nu -2,\nu -1} $ respectively
$$  P^* \; = \; \left[ \begin{array}{cccccc} O & E_{1,2} & O & \cdots & O & O \\
       O& O& E_{2,3}& \cdots & O & O \\
       \vdots & \vdots & \vdots & \ddots & \vdots &  \vdots\\
             O& O& O& \cdots & O & E_{\nu -2,\nu -1}^o \\
    \hat{T}_{\nu -1} & -\hat{T}_{\nu -2} & -\hat{T}_{\nu -3} & \ldots &
       -\hat{T}_2 &  -\hat{T}_1^*  \end{array} \right ]
\eqno(7.76)     $$
 To calculate the determinant of  the matrix $
(sI_{n-l_{\nu -1}} -P^*) $ we use Assertion 1.2 with changing $n$ by $n-l_{\nu -1}$, $p$ by $\nu -1$,
$\hat{T}_1$ by $\hat{T}_1^*$  and  $E_{p -1,p}$ by $E_{\nu -2,\nu -1}^o $. Thus
$$  det(sI_{n-l_{\nu -1}} -P^*)
\; = \;
      det\left[ \begin{array}{cccccc} sI_{l_1} & -E_{1,2} & O & \cdots & O & O \\
             O & sI_{l_2} & -E_{2,3}& \cdots & O & O \\
             \vdots & \vdots & \vdots & \ddots & \vdots & \vdots \\
           O& O& O& \cdots & sI_{l_{\nu -2}} & -E_{\nu -2,\nu -1}^o \\
           \hat{T}_{\nu -1} & \hat{T}_{\nu -2} & \hat{T}_{\nu -3} & \ldots &
           \hat{T}_2 & sI_r + \hat{T}_1^*       \end{array} \right ] \; =$$
$$   = \;s^{n-l_{\nu -1}-r(\nu -1)}det(I_rs^{\nu -1} + \hat{T}_1^*s^{\nu -2} +
        [O,\hat{T}_2]s^{\nu -3} + \cdots + [O,\hat{T}_{\nu -1}])\; = \;  $$
$$   = \; s^{n-l_{\nu -1}-r(\nu -1)}det(I_rs^{\nu -1} + [O,\hat{T}_1]s^{\nu -2} +
        [O,\hat{T}_2]s^{\nu -3} + \cdots + [O,\hat{T}_{\nu -1}]) $$
The analysis of the matrix polynomial in the last expression  and using (7.66) gives the following equality
$$  det(sI_{n-l_{\nu -1}} -P^*) \; = \; s^{n-l_{\nu -1}-r(\nu -1)}
     detC_{\nu}^{-1}det\tilde{C}(s)
\eqno (7.77)    $$
     On the other hand substituting (7.74)  and  (7.75)  into  the
right-hand side of (7.76) gives
$$  det(sI_{n-l_{\nu -1}} -P^*) \; = \;
      \left[ \begin{array}{cccccc} sI_{l_1} & -E_{1,2} & O & \cdots & O & O \\
             O & sI_{l_2} & -E_{2,3} & \cdots & O & O \\
             \vdots & \vdots & \vdots & \ddots & \vdots & \vdots \\
           O& O& O& \cdots & sI_{l_{\nu -2}} & -[O,E_{\nu -2,\nu -1}] \\
           \hat{T}_{\nu -1} & \hat{T}_{\nu -2} & \hat{T}_{\nu -3} & \ldots &
           \hat{T}_2 & sI_r + [O,\hat{T}_1]       \end{array} \right ] \
\eqno   (7.78)     $$
 Partitioning  submatrices $\hat{T}_i$, $i=1,2,\ldots,\nu -1$  as
$$   \hat{T}_i  \; = \; \left[ \begin{array}{c} \hat{T}_i^1 \\   \\ \hat{T}_i^2
\end{array}  \right] \begin{array}{cl} \} & r -l_{\nu -1} \\ \\ \} & l_{\nu -1}
\end{array}
\eqno  (7.79)  $$
 where  blocks  $\hat{T}_i^1$ , $\hat{T}_i^2$
have $r -l_{\nu -1}$ and $l_{\nu -1} $ rows respectively we find the structure of the matrix $sI_r +
[O,\hat{T}_1] $
$$ sI_r + [O,\hat{T}_1] \; = \; sI_r + \left[ \begin{array}{cc}
          O & \hat{T}_i^1 \\   \\ O & \hat{T}_i^2 \end{array}  \right] \; = \;
          \left[ \begin{array}{cc}  sI_q & \hat{T}_i^1 \\   \\
           O & sI_{l_{\nu -1}} + \hat{T}_i^2 \end{array}  \right]
\eqno  (7.80)   $$
where  $q=r-l_{\nu -1}$. Substituting (7.79) and (7.80) into (7.78)
$$  det(sI_{n-l_{\nu -1}} -P^*) \; = \;
  \left[ \begin{array}{ccccccc} sI_{l_1} & -E_{1,2} & O & \cdots & O & O & O \\
         O & sI_{l_2} & -E_{2,3} & \cdots & O & O & O\\
             \vdots & \vdots & \vdots & \ddots & \vdots & \vdots & \vdots \\
           O& O& O& \cdots & sI_{l_{\nu -2}} & O & -E_{\nu -2,\nu -1} \\   \\
        \hat{T}_{\nu -1}^1 & \hat{T}_{\nu -2}^1 & \hat{T}_{\nu -3}^1 & \ldots &
           \hat{T}_2^1 & sI_q &  \hat{T}_1^1 \\  \\
        \hat{T}_{\nu -1}^2 & \hat{T}_{\nu -2}^2 & \hat{T}_{\nu -3}^2 & \ldots &
           \hat{T}_2^2 & O &  sI_{l_{\nu -1}} + \hat{T}_1^2
       \end{array} \right ]
\eqno (7.81)           $$
   and calculating the determinant of the block matrix we  get
$$  det(sI_{n-l_{\nu -1}} -P^*) \; = \; s^q det \left[ \begin{array}{cccccc}
         sI_{l_1} & -E_{1,2} & O & \cdots & O & O \\
         O & sI_{l_2} & -E_{2,3} & \cdots & O & O \\
             \vdots & \vdots & \vdots & \ddots & \vdots & \vdots \\
           O& O& O& \cdots & sI_{l_{\nu -2}} & -E_{\nu -2,\nu -1} \\   \\
          \hat{T}_{\nu -1}^2 & \hat{T}_{\nu -2}^2 & \hat{T}_{\nu -3}^2 & \ldots &
           \hat{T}_2^2 & sI_{l_{\nu -1}} + \hat{T}_1^2
       \end{array} \right ]  $$
Using notations (7.70), (7.72) and  relations: $q=r-l_{\nu -1}$,
$l_1+l_2+\cdots+l_{\nu -2}+l_{\nu -1} = n-r$ we can rewrite the
last equality as follows
$$  det(sI_{n-l_{\nu -1}} -P^*) \; = \;s^{r-l_{\nu -1}}det(sI_{n-r}-
         \left [ \begin{array}{c} E \\ \hat{T}\end{array} \right ] )
\eqno(7.82)  $$
Equating right-hand sides of (7.82)  and   (7.77)
yields
$$  s^{n-l_{\nu -1}-r(\nu -1)}detC_{\nu}^{-1}det\tilde{C}(s) \; = \;
     s^{r-l_{\nu -1}}det(sI_{n-r}- \left [ \begin{array}{c}
       E \\ \hat{T}  \end{array} \right ] )    $$
The last equality proves the lemma.

    The following result follows immediately from  Assertion 7.6,
Lemma 7.1 and formulas (7.67) and (7.72).

\it{THEOREM 7.6.} \rm System zeros of system (1.1), (1.2) with an equal number of inputs  and  outputs  and
$det(CB) \neq 0$ are defined as eigenvalues of the  $(n-r)\times (n-r)$ matrix $\left [
\begin{array}{c} E \\ \hat{T}\end{array} \right ] $ where the $l_{\nu
-1}\times (n-r)$ submatrix $\hat{T}$  satisfies to the following formula
$$ \hat{T} \; = \;  \left \{ \begin{array}{rcc} - C_{\nu}^{-1}[C_1,C_2,\ldots
                   C_{\nu -1}] & , & l_{\nu -1} = r \\  \\
                  -[O,I_{l_{\nu -1}}]C_{\nu}^{-1}[C_1,C_2,\ldots
                   C_{\nu -1}] & , & l_{\nu -1} < r \end{array} \right.
\eqno (7.83)   $$

\it{CASE 2. $\;$ $CB$ is singular matrix}. \rm  Hence, the $r\times r$ matrix $C_{\nu}  =  CB\bar{G}_{\nu}^{-1}
$ is also singular one. Let's introduce the $(n-l_{\nu -1})\times (n-l_{\nu -1}) $ matrix
$$  \tilde{I}_{n-l_{\nu -1}} \;= \; \left [ \begin{array}{cc} I_{\beta} & O \\
      O & C_{\nu} \end{array} \right ], \; \beta = l_1+l_2+\cdots+l_{\nu -2}
\eqno(7.84)  $$
 and  the $r\times (n-r)$ matrix
$$       C^* = [C_1, C_2,\ldots, C_{\nu -1}]
\eqno (7.85) $$

\it{LEMMA 7.2.} \rm  If $l_{\nu -1} = r$ then the following
equality is true
$$ det(s\tilde{I}_{n-l_{\nu -1}} - \left [ \begin{array}{c} E \\ -C^*\end{array}
     \right ] )\; = \; s^{n-l_{\nu -1}-r(\nu -1)} det ( [O,C_1] +
      [O,C_2]s + \cdots +  [O,C_{\nu -1}]s^{\nu -2} + C_{\nu}s^{\nu -1})
\eqno (7.86) $$ where the $(n-r-l_{\nu -1})\times (n-r) $ submatrix $E$ is defined from (7.70).

\it{PROOF}.$\;$ \rm To calculate the determinant in  the left-hand side  of (7.86) we partition  the matrix
$s\tilde{I}_{n-l_{\nu -1}} - \left [ \begin{array}{c} E \\
-C^*\end{array} \right ] $ into four blocks similarly to ones in (1.27), (1.28) with $p =\nu -1$ and  find a
determinant of the resulting block matrix
$$  det(\tilde{I}_{n-l_{\nu -1}} \;- \; \left [ \begin{array}{c} E \\ -C^*
     \end{array} \right ] ) \; = \;
      det\left[ \begin{array}{cccccc} sI_{l_1} & -E_{1,2} & O & \cdots & O & O \\
             O & sI_{l_2} & -E_{2,3}& \cdots & O & O \\
             \vdots & \vdots & \vdots & \ddots & \vdots & \vdots \\
           O& O& O& \cdots & sI_{l_{\nu -2}} & -E_{\nu -2,\nu -1} \\
           C_1 & C_2 & C_3 & \ldots & C_{\nu -1} & sC_{\nu} + C_{\nu -1}
             \end{array} \right ] \; =$$
$$     = \; s^{(n-l_{\nu -1})-r}det(sC_{\nu} + C_{\nu -1} + [C_1, C_2,\ldots,
   C_{\nu -2}] \left[ \begin{array}{cccc} sI_{l_1} & -E_{1,2} & \cdots & O \\
                 O & sI_{l_2} & \cdots & O \\
                 \vdots & \vdots & \ddots &  \vdots \\
                 O & O & \cdots & sI_{l_{\nu -2}} \\
                 \end{array} \right ] ^{-1} \left[ \begin{array}{c} O \\
                 O \\ \vdots \\ E_{\nu -2,\nu -1}  \end{array} \right ]) =$$
$$ = \; s^{(n-l_{\nu -1})-r}det(sC_{\nu} + C_{\nu -1} + [C_1, C_2,\ldots,
     C_{\nu -2}] \left[ \begin{array}{c} s^{1-\nu }[O,I_{l_1}] \\
                 s^{2- \nu }[O,I_{l_2}] \\ \vdots \\s^{-1}[O,I_{l_{\nu -1}}]
                  \end{array} \right] ) \; = $$
$$         \; s^{(n-l_{\nu -1})-r}det(sC_{\nu} + C_{\nu -1} +
     [O,C_{\nu -2}]s^{-1} + \cdots + [O,C_2]s^{3-\nu } +[O,C_1]s^{2-\nu }) = $$
$$         \; = s^{n-l_{\nu -1}-r(\nu -1)} det ( [O,C_1] + [O,C_2]s +
           \cdots +  [O,C_{\nu -1}]s^{\nu -2} + C_{\nu}s^{\nu -1})
$$
The lemma has been proved.

     Let $l_{\nu -1} <r$. We introduce  the $r\times r$ matrix
$$    \bar{C}_{\nu -1} = [O,C_{\nu -1}]
\eqno  (7.87) $$
 and the $l_{\nu -2} \times r$ submatrix $E_{\nu -2,\nu -1}^o$ of the form (7.75). The following
equality is true by Lemma 7.2
$$ det(s\tilde{I}_{n-l_{\nu -1}} - \left [ \begin{array}{c} \bar{E} \\ \dotfill\\
  -C_1, -C_2,\ldots, -C_{\nu -2}, -\bar{C}_{\nu -1} \end{array}\right ] )\; = $$
$$= \; s^{n-l_{\nu -1}-r(\nu -1)} det ( [O,C_1] +
         [O,C_2]s + \cdots +  [O,C_{\nu -1}]s^{\nu -2} + C_{\nu}s^{\nu -1})
   \eqno (7.88) $$
where the $\beta \times (n-l_{\nu -1})$  matrix $\bar{E} \;(\beta = l_1+l_2+ \cdots+l_{\nu -2})$ has the form
(7.69) with $E_{\nu -2,\nu -1} = E_{\nu -2,\nu -1}^o$.

     Recalling that  system  zeros  of    controllable  system
(1.1), (1.2) with $l = r$  are equal to  zeros of the polynomial
$$ \psi(s) \; = \; s^{n-r\nu} det ( [O,C_1] + [O,C_2]s + \cdots +
        [O,C_{\nu -1}]s^{\nu -2} + C_{\nu}s^{\nu -1})
\eqno (7.89)   $$
 and substituting $det ( [O,C_1] + [O,C_2]s + \cdots
+[O,C_{\nu -1}]s^{\nu -2} + C_{\nu}s^{\nu -1})$ from (7.88) in
(7.89) we get
$$ \psi(s) \; = \; s^{n-r\nu} s^{r(\nu -1)-n+l_{\nu -1}}
 det(s\tilde{I}_{n-l_{\nu -1}} - \left [ \begin{array}{c} \bar{E} \\ -\bar{C}
     \end{array}\right ] )\; = \; s^{l_{\nu -1}-r}det(s\tilde{I}_{n-l_{\nu -1}}
     - \left [ \begin{array}{c} \bar{E} \\ -\bar{C} \end{array}\right ] )$$
where $\bar{C}=[C_1, C_2,\ldots, C_{\nu -1}, \bar{C}_{\nu -1}]$.
This result can be formally stated as

\it{THEOREM 7.7.} \rm  System zeros of system (1.1), (1.2) with an equal number of inputs and outputs coincide
with generalized eigenvalues of the  $(n-l_{\nu -1})\times (n-l_{\nu -1}) $ regular matrix pencil
$$ Z(s) \; = \; s \left [ \begin{array}{cc} I_{\beta} & O \\
                   O & C_{\nu} \end{array} \right ] \; + \;
         \left [ \begin{array}{c} -\bar{E} \\ \dotfill\\
    C_1, C_2,\ldots, C_{\nu -2}, [O,C_{\nu -1}] \end{array}\right ]
\eqno (7.90)   $$
 without $r-l_{\nu -1}$ generalized eigenvalues
in the origin. In (7.90) the $\beta\times (n-l_{\nu -1})$  matrix $\bar{E}$  has  the form (7.70) with the
$l_{\nu -2}\times r$ submatrix $E_{\nu -2,\nu -1} \; = \; [O, E_{\nu -2,\nu -1}]$ (see (7.75)).

\it{COROLLARY 7.5.} \rm  If $l_{\nu -1} = r$ then  system zeros of (1.1), (1.2) with an equal number of inputs
and outputs are defined as  generalized eigenvalues of the  $(n-r)\times (n-r)$ regular matrix pencil
$$ Z(s) \; = \; s \left [ \begin{array}{cc} I_{\beta} & O \\
                   O & C_{\nu} \end{array} \right ] \; + \;
         \left [ \begin{array}{c} -E \\ \dotfill\\
    C_1, C_2,\ldots, C_{\nu -1} \end{array}\right ]
\eqno (7.91) $$ where the $\beta \times (n-r)$  matrix $E$ has the form (7.70).

     We illustrate the method by the following examples.

                 \it{EXAMPLE  7.3.} \rm
Let us find zeros of  system (1.1), (1.2) with  $n = 4, r = l = 2$
and state-space model matrices

$$   A \; = \;  \left[ \begin{array}{cccc} 2 & 1 & 0 & 0 \\ 0 & 1 & 0 & 1 \\
     0 & 2 & 0 & 0 \\ 1 & 1 & 0 & 0 \end{array} \right ], \qquad
     B\; = \; \left[ \begin{array}{cc} 1 & 0 \\ 0 & 0 \\ 0 & 0 \\ 0 & 1
            \end{array} \right ], \qquad
C \; = \;  \left[ \begin{array}{crcc} 1 & 0 & 0 & 0 \\ 0 & 0 & 1 & 1
        \end{array} \right ]
\eqno (7.92) $$ In Example 1.3 (see  Sect  1.2.3.) it has been shown that this system has    $\nu =3,\; l_1 =
l_2 = 1, \; l_3 =2$ and following matrices $N$  and $N^{-1}$
$$ N = \left[ \begin{array}{cccc}   0 &  0 &  0.5 &  0 \\   0 &  1 & 0 &  0 \\
 1 &  0 &  0 &  0 \\ 0 & 1 &  0 &  1  \end{array} \right ], \;
   N^{-1} = \left[ \begin{array}{crcc}   0 &  0 &  1 &  0 \\ 0 & 1 & 0 & 0 \\
 2 &  0 & 0 & 0 \\ 0 & -1 &  0 &  1  \end{array} \right ] $$
We calculate
$$  CN^{-1} \;=\; \left [ \begin{array}{crcc} 1 & 0 & 0 & 0 \\ 0 & 0 & 1 & 1
     \end{array} \right ] N^{-1} \;=\; \left [ \begin{array}{crcc}
    0 & 0 & 1 & 0 \\ 2 & -1 & 0 & 1 \end{array} \right ]
\eqno (7.93) $$
and partition this matrix on three blocks $C_1,
C_2, C_3$ of dimensions $r\times l_1 = 2\times 1$, $r\times l_2=
2\times 1$, $r\times l_3 = 2\times 2$ respectively
$$    CN^{-1}\; = \; [C_1, C_2, C_3] \; = \; \left [ \begin{array}{ccccc}
      \begin{array}{c} 0 \\ 2 \end{array} & \begin{array}{c} \vdots \\ \vdots
      \end{array} & \begin{array}{r} 0 \\ -1 \end{array} & \begin{array}{c}
      \vdots \\ \vdots \end{array} & \begin{array}{cc} 1 & 0 \\ 0 & 1
         \end{array} \end{array} \right ] $$
 Since  $ CB =\left [ \begin{array}{cc} 1 & 0 \\ 0 & 1
 \end{array} \right ]$   is the nonsingular matrix then
we need to use Theorem 7.6 to calculate system zeros. Using formula (7.70) we form  the $(n-r-l_{\nu -1})\times
(n-r) $ submatrix $E$. Since $\nu =3,\; l_1 = l_2 = 1, \; l_3 =2$  then $n-r-l_{\nu -1} =1$,  $n-r =2$, and
$E_{\nu -2,\nu -1} = E_{12}$ is the $1\times 1$ submatrix that equals to 1. Thus, we get
$$           E \; = \; [ 0 \; 1 ]      $$
and by formula (7.83) calculate
$$ \hat{T} \;=\;-[O,I_{l_2}]C_3^{-1}[C_1,C_2] = -[ 0 \;\;1 ] \left [
    \begin{array}{cr} 0 & 0 \\ 2 & -1 \end{array} \right ] =
     [ -2 \;\; 1 ]    $$
The matrix $\left [ \begin{array}{c} E \\ \hat{T} \end{array} \right ] $ of the order $n-r = 2$ becomes
$$ \left [ \begin{array}{c} E \\ \hat{T}\end{array} \right ]  \; = \;
   \left [ \begin{array}{rc} 0 & 1 \\ -2 & 1\end{array} \right ] $$
It has the characteristic polynomial
$$ \psi(s)\; =\;det\left [ \begin{array}{cc} s & -1 \\ 2 & s-1\end{array}
      \right ]\;=\;s^2 -s +2        $$
 that is equal  to  the zero polynomial of the present system.

     For testing we compute the only minor of the system matrix
 $P(s)$
$$ detP(s) \; =\; det\left [ \begin{array}{cccrrr}
         s-2 & -1  & 0 & 0 & -1 & 0 \\ 0 & s-1 & 0 & -1 & 0 & 0 \\
          0 & -2 & s & 0 & 0 & 0 \\-1 & -1 & 0 & s & 0 &-1\\
           1 & 0 & 0 & 0 & 0 & 0 \\ 0 & 0 & 1 & 1 & 0 & 0
          \end{array} \right ]   $$
We have been obtained the same result.

                  \it{EXAMPLE 7.4.} \rm

 Let the model (7.92) has the following output matrix
$$C \; = \;  \left[ \begin{array}{crcc} 1 & -1 & 1 & 0 \\ 1 & 0 & 0 & 0
        \end{array} \right ]     $$
        Since in this case the  matrix $N$ is same as above
then
$$  CN^{-1} \;=\; \left [ \begin{array}{crcc} 1 & -1 & 1 & 0 \\ 1 & 0 & 0 & 0
     \end{array} \right ] N^{-1} \;=\; \left [ \begin{array}{crcc}
    2 & -1 & 1 & 0 \\ 0 & 0 & 1 & 0 \end{array} \right ]
\eqno (7.94) $$
 Using above calculated  controllability characteristics: $\nu =3,\; l_1 = l_2 = 1,\;\; l_3 =2$
we  partition the matrix (7.94) into three blocks $C_1, C_2, C_3$ of the sizes $2\times 1$, $2\times 1$,
$2\times 2$ respectively
$$    CN^{-1}\; = \; [C_1, C_2, C_3] \; = \; \left [ \begin{array}{ccccc}
      \begin{array}{c} 2 \\ 0 \end{array} & \begin{array}{c} \vdots \\ \vdots
      \end{array} & \begin{array}{r} -1 \\ 0 \end{array} & \begin{array}{c}
      \vdots \\ \vdots \end{array} & \begin{array}{cc} 1 & 0 \\ 1 & 0
         \end{array} \end{array} \right ]
\eqno   (7.95) $$
  In this case  the matrix  $ CB =\left [
\begin{array}{cc} 1 & 0 \\ 1 & 0
\end{array} \right ]$ is  singular  one. Thus, to calculate  system zeros we
need to use  Theorem 7.7. We form the $l_{\nu -2}\times r$ submatrix $E_{\nu -2,\nu -1}^o \; = \; [O,  E_{\nu
-2,\nu -1}]$ and the $r\times r$ submatrix  $[O,C_{\nu -1}]$. Since $r=3,\; l_{\nu -1}=l_1 = 1,\; \nu -2
=1,\;\nu -1 =2$ then by the formula (7.75)  we find
$$     E_{1,2}^o \; = \; \left [ \begin{array}{cc} O & E_{1,2}\end{array}
       \right ] \;=\; \left [ \begin{array}{cc} 0 & 1 \end{array}\right ]
\eqno (7.96)    $$
 The $2\times 2$ matrix $[O,C_{\nu -1}] \;=\; [O,C_2]$  with the
 $2\times 1$ block  $C_2$  from (7.95) is
$$    [O,C_{\nu -1}] \;=\; [O,C_2]\;=\; \left [ \begin{array}{ccc}
          \begin{array}{c} 0 \\ 0 \end{array} & \begin{array}{c} \vdots \\
      \vdots \end{array} & \begin{array}{r} -1 \\ 0 \end{array} \end{array}
         \right ] $$
To form the $\beta\times (n-l_{\nu -1})$ matrix $\bar{E}$ and the $(n-l_{\nu -1})\times (n-l_{\nu -1})$ matrix
$\tilde{I}_{n-l_{\nu -1}}$ we use (7.96) and (7.84). Since $\nu =3, \; \beta = l_1+\cdots +l_{\nu -2} = l_1
=1,\; n-l_{\nu -1} = 4-1=3 $ then
$$         \bar{E}\;=\; [O, E_{1,2}^o] \;=\; [ 0 \; 0\;  1 ]   $$
$$  \tilde{I}_{n-l_{\nu -1}} \;=\; \tilde{I}_{n-l_2}\;=\;\tilde{I}_3 \;=\;
    \left [ \begin{array}{cc} I_1 & O \\ O & C_3 \end{array} \right ]  \;=\;
    \left [ \begin{array}{ccc} 1 & 0 & 0 \\ 0 & 1 & 0 \\ 0 & 1 & 0
     \end{array} \right ] $$
We result in the following matrix pencil
$$ Z(s) \; = \; s \left [ \begin{array}{cc} I_1 & O \\
                   O & C_3 \end{array} \right ] \; + \;
         \left [ \begin{array}{c} -\bar{E} \\ \dotfill\\
          C_1, [O,C_2] \end{array}\right ]\;=\; s \left [ \begin{array}{ccc}
       1 & 0 & 0\\ 0 & 1 & 0\\0 & 1 & 0 \end{array} \right ] \; + \;
        \left [ \begin{array}{ccr} 0 & 0 & -1\\
         2 & 0 & -1\\0 & 0 & 0 \end{array} \right ]      $$
having the following characteristic polynomial
$$  detZ(s) \;=\; det(s \left [ \begin{array}{ccc} 1 & 0 & 0\\ 0 & 1 & 0\\
             0 & 1 & 0 \end{array} \right ] \; + \;
          \left [ \begin{array}{ccr} 0 & 0 & -1\\ 2 & 0 & -1\\0 & 0 & 0
            \end{array} \right ] ) \;=\; s^2 - 2s      $$
 According to Theorem 7.7  the zero
polynomial is
$$  \psi(s) \;=\; s^{l_{\nu -1}-r}detZ(s) \;=\;s^{l_2-r}detZ(s) \;=\;
                  s^{-1}detZ(s)\;=\; s-2 $$
Hence, the system has the only system zero  $s=2$.

    To test we find the determinant of  the system matrix
$P(s)$ that equal to  $s-2$. This result is in agreement with the
obtained one.

\chapter[Zero computation]
        {Zero computation}

In the present chapter we develop several computational techniques. The natural way to compute invariant and
transmission zeros is based on their definitions via the Smith and Smith-McMillan canonical forms for matrices
$P(s)$ and $G(s)$ respectively. Decoupling zeros may be calculated by using the Smith form for matrices $P_i(s)
= [sI_n-A,B]$ and $P_o(s)^T = [sI_n-A^T,C^T]$. This  approach is laborious because of  operations with
polynomial and/or rational matrices. Here we will study alternative  methods that are based on efficient
numerical procedures and have the simple computer-aided realization.

     For the most part we will consider a system with an equal  number
of inputs and outputs $(r = l)$. This restriction does not essential because otherwise we can recommend to
perform twice the squaring down  operation of outputs (if $l>r$) or inputs (if $l<r$) and to find  system zeros
as an intersection of sets of zeros of  squared down systems. Indeed, let for definiteness $l>r$. Using
different $r\times l$  matrices  $E_1$ and $E_2$ of full ranks we construct two squared down systems from (1.1),
(1.2):

$$     \begin{array}{ccc} S_1 & : & \dot{x} = Ax + Bu, \;\;y_1 =E_1Cx \\
                          S_2 & : & \dot{x} = Ax + Bu, \;\;y_2 =E_2Cx
         \end{array} $$
Zero sets $\Omega(E_1)$ and $\Omega(E_2)$ of systems $S_1$ and $S_2$  respectively are calculated via system
matrices $P_1(s)$ and $P_2(s)$   of the following form
 $$  P_i(s) \; = \; \left[ \begin{array}{cc} sI_n-A & -B \\ E_iC & O
\end{array} \right ],\; \; i=1,2
\eqno(8.1) $$
It is evident that the intersection
$$ \Omega\;=\;  \Omega(E_1) \cap \Omega(E_2)
\eqno (8.2) $$ 'almost always'\footnote{I.e. the class of systems  that don't possess such property is either
empty or lies on a hypersurface of the parameter space of $(A,B,C)$ [D4].} is equal to the zero  set  of
 system  (1.1), (1.2).

\section[Zero computation via matrix P(s)]
        {Zero computation via matrix P(s)}

Analysis of the  system matrix
 $$  P(s) \; = \; \left[ \begin{array}{cc} sI_n-A & -B \\ C & O
\end{array} \right ]
\eqno(8.3) $$
 shows that the complete set of  system zeros is
formed as the set of complex $s$  for which the normal rank  of (8.3) is locally reduced. For the system with $r
= l$   and non identically zero $detP(s)$ we can consider  the matrix (8.3) as the
 regular [G1] pencil of matrices
$$        P(s) \;=\; sD + L
\eqno (8.4) $$
where
$$  L \;=\;\left[ \begin{array}{cc} -A & -B \\ C & O \end{array} \right ] ,\;
    D \;=\;\left[ \begin{array}{cc} I_n & O \\ O& O \end{array} \right ]
\eqno  (8.5)    $$
 Thus, the problem is reduced to  calculating
 generalized eigenvalues of the regular matrix pencil (8.4).
These generalized eigenvalues coincide with  finite $s$    for which there exist a nontrivial solution  of the
equation $(sD + L)q \;=\; 0$ where $q$ is an  $n + r$  vector [P1]. To calculate generalized eigenvalues we can
apply  QZ algorithm [P1]. Moreover, the special modification [M6] of this algorithm  may be used for  a singular
matrix $D$. In [L3] this  modification was successfully applied to  calculating zeros via the matrix pencil
(8.4).

The advantage of the mentioned approach is its numerical stability because of  effectiveness of  QZ algorithm.
The main shortcoming is associated with separating decoupling zeros from transmission zeros.

To overcome  this problem we may use the method proposed by Porter [P7]. This approach proposes at first to
solve twice the QZ problem for
 the regular $(n+r+l)\times (n+r+l)$ matrix pencil
$$ s \left [ \begin{array}{ccc} I_n & O & O\\ O& O & O\\
   O & O & O \end{array} \right ] \; + \;
   \left [ \begin{array}{ccr} -A & O & B\\ -C & I_l & O \\O & K_i & I_r
      \end{array} \right ]
\eqno(8.6) $$
 where $K_i$, $i=1,2$ are some full rank matrices
with bounded elements.
     Indeed, using the formula Shura [G1] we can calculate
$$ det\left [ \begin{array}{ccr} sI_n-A & O & B\\ -C & I_l & O \\O & K_i & I_r
      \end{array} \right ]  \;=\; det(sI_n -A-BK_iC),\; i=1,2
\eqno(8.7)   $$
 Thus, applying  matrices $K_i$ $(i=1,2)$ in (8.6) is equivalent to introducing a proportional
output feedback $u = K_iy$ into the system: $\dot{x} = Ax + Bu, \;\;y =Cx $. Such feedback shifts only
controllable and observable eigenvalues of the matrix $A$. Uncontrollable and unobservable eigenvalues are the
invariants under the proportional output feedback and they are decoupling zeros of the system (see Sections 2.4,
5.3). Hence, those generalized eigenvalues of the problem (8.6), which are not changed for different matrices
$K_i$ $(i=1,2)$ coincide with decoupling zeros  of system (1.1), (1.2) and a set of decoupling zeros $\Omega_d$
is almost always computed from the intersection
$$ \Omega_d\;=\;  \Omega(K_1) \cap \Omega(K_2)
\eqno (8.8) $$
 where $\Omega(K_i), i=1,2 $ are sets of
generalized eigenvalues of matrix pencils (8.6).

If the whole set of system zeros $ \Omega $  has been computed by somehow method then from the union
$$ \Omega \;=\;  \Omega_d \cup \Omega_t
\eqno (8.9) $$ we can separate the set of transmission zeros $ \Omega_d $ from the set of  decoupling zeros.

For calculating the  set  $\Omega$ for a system with an unequal number of inputs and outputs $(r \neq l)$ we can
use the following strategy. At first the QZ problem is twice solved for  following regular matrix pencils of the
order $n + min (r,l)$
$$ s\left[ \begin{array}{cc} I_n & O \\ O& O \end{array} \right ] \; +
    \;  \left[ \begin{array}{cc} -A & -B \\ E_iC & O \end{array} \right ],
     \qquad   \rm{ if} \qquad l>r
\eqno (8.10)   $$
$$ s\left[ \begin{array}{cc} I_n & O \\ O& O \end{array} \right ] \; +
    \;  \left[ \begin{array}{cc} -A & -BE_i \\C & O \end{array} \right ],
     \qquad   \rm{ if} \qquad l<r
\eqno (8.11)   $$ where $E_i$, $i=1,2$ are  some  $r\times l$ matrices  of a full rank with bounded elements.
     In view of preceding reasoning these  pencils  correspond  to
 system matrices $P_1(s)$ and $P_2(s)$  of squared down
systems: $P_1(s)$ defines the system with squared down outputs and
$P_2(s)$ defines the system with squared down inputs. Let's denote
zero sets of  the first and  second  squared  down systems by
$\Omega_1(E_i)$ and $\Omega_1(E_i)$ respectively. It is evident
that  the  zero  set  of   original system (1.1), (1.2) can be
calculated from the  following intersections
$$ \Omega\;=\;  \Omega_1(E_1) \cap \Omega_1(E_2)  \qquad \rm{for} \qquad \it l>r
\eqno (8.12) $$
$$ \Omega\;=\;  \Omega_2(E_1) \cap \Omega_2(E_2)  \qquad \rm{for} \qquad \it l<r
\eqno (8.13) $$
 In  the work [7] it is presented  the following computational method for a system with
 $l \geq r $: to solve twice the QZ problem for the following regular $(n+r+l)\times (n+r+l)$ matrix pencil
$$ s \left [ \begin{array}{ccc} I_n & O & O\\ O& O & O\\
   O & O & O \end{array} \right ] \; + \;
   \left [ \begin{array}{ccr} -A & O & B\\ -C & I_l & O \\O & K_i & O
      \end{array} \right ]
\eqno(8.14) $$
 which is distinguished from the matrix pencil
(8.6)  by the only block. The set $\Omega$ of system zeros
 is calculated from the following intersection
$$ \Omega\;=\;  \Omega^*(K_1) \cap \Omega^*(K_2)
\eqno (8.15) $$
 where $\Omega^*(K_i)$ $(i=1,2)$ are  sets of generalized eigenvalues of  the matrix pencil
(8.14).

     Let us  show that the present method actually
computes  the full set of system zeros. Indeed, interchanging the second and third block rows  and  columns  of
the matrix pencil (8.14) and using the formula Shura [G1] we calculate the determinant of (8.14)
$$ det \left [ \begin{array}{ccr} sI_r-A & O & B\\ -C & I_l & O \\O & K_i & O
          \end{array} \right ]\;=\;det \left [ \begin{array}{ccr}
         sI_r-A & B & O\\ O & O & K_i \\-C & O & I_l\end{array} \right ]\;=\; $$
 $$ =\; det \left [ \begin{array}{cc} sI_r-A & B \\ K_iC & O
      \end{array} \right ] \;=\;
    (-1)^r det \left [ \begin{array}{cc} sI_r-A & -B \\ K_iC & O
    \end{array} \right ]
\eqno  (8.16)  $$
 It follows from (8.16)  that  generalized
eigenvalues of the regular pencil (8.14) coincide with ones of the regular pencil (8.10) with $K_i =E_i$
$(i=1,2)$. Therefore, intersection (8.15) gives the set of system zeros.

The similar procedure may be used when $r \ge l$. Here it is applied the following $(n+r+l)\times (n+r+l)$
matrix pencil
$$ s \left [ \begin{array}{ccc} I_n & O & O\\ O& O & O\\
   O & O & O \end{array} \right ] \; + \;
   \left [ \begin{array}{ccr} -A & O & B\\ -C & O & O \\O & K_i & I_r
      \end{array} \right ]
\eqno(8.17) $$ having the following determinant
$$ det \left [ \begin{array}{ccr} sI_r-A & O & B\\ -C & O & O \\O & K_i & I_r
         \end{array} \right ]\;=\;(-1)^l det \left [ \begin{array}{cc}
         sI_r-A & -BK_i \\ C & O \end{array} \right ]
\eqno  (8.18)  $$
 It follows  from (8.18) that  generalized
eigenvalues  of   the pencil (8.17) coincide with ones of the regular pencil (8.11) with $K_i =E_i$ $(i=1,2)$.
Therefore, all previous reasoning are held.

The advantage of this approach is  the applicability of the universal QZ algorithm for matrix pencils (8.6),
(8.14) or (8.17). Hence, similar computational algorithms and computer software can be used as for computing
zeros  as for separating different type zeros.  But  such approach increases considerably a dimension of the
problem. Therefore, when  $r \neq l$ it is  more preferable to use the QZ procedure for pencils (8.10) or (8.11)
of  the order $n + min (r,l)$.

\section[Zero computation based on matrix  A+BKC]
        {Zero computation based on matrix  A+BKC}

This approach uses  invariance of zeros under a high gain output feedback (see Sect.6.5). Let us remind this
property.  We consider a close-loop controllable and observable  system having the following dynamics matrix:
$A(K) = A + kBKC$ where $K$ is some arbitrary matrix with limited elements, $k$  is a real scalar. If $k$ goes
to infinity ($k \to \infty $) then $n-r$ eigenvalues of the matrix $A(k)$ approach positions of transmission
zeros of system (1.1), (1.2) and $r$ reminder eigenvalues tend to infinity. Therefore, to calculate zeros we can
use the following approach [D4]:

     1. Compute $n$ eigenvalues of  the matrix $A(K) = A + kBKC$ for any
arbitrary  matrix $K$ of a full rank and a large value of $k\approx 10^{15}$.

     2. Separate $n-r$  finite eigenvalues that are equal  to
        transmission zeros.

     This approach  may  be successfully  used  for  uncontrollable
and/or unobservable system. Indeed,  decoupling zeros coincide with limited eigenvalues of the dynamics matrix
of a closed-loop system that are invariant under a proportional output feedback and the  procedure above
calculates system zeros of an uncontrollable and/or unobservable system. But here it is necessary to separate
sets of transmission and decoupling zeros. For  this purpose we  can  use, for example, the approach of Section
8.1.

     Consequently, for  uncontrollable and/or unobservable
system (1.1), (1.2) with $r=l$  we propose the following general procedure:

     1. Find  finite eigenvalues of the matrix  $A+kBK_1C$
        with $k\approx 10^{15}$ and  $K_1$  being an arbitrary
        matrix of a full rank. Denote these eigenvalues
        by $\Omega(K_1)$.

     2. Repeat step 1 with another matrix $K_2$. Denote resulted
        eigenvalues as $\Omega(K_2)$ .

     3. From  intersection (8.8) calculate decoupling zeros $\Omega_d$.

     4. From union (8.9) calculate transmission zeros $\Omega_t$ .

The advantage of  this  approach  is its   simplicity . Moreover, in contrast to the approach of Sect.8.1, the
processed  matrices are of small sizes. But low computational accuracy (using a large number $k\approx 10^{15}$)
makes difficulties for applicability of this method.

\section[Zero computation via transfer
                             function matrix]
        {Zero computation via  transfer\\
                             function matrix}

Now  we consider the numerical method proposed by Samash in [S1].
This method  is  based on the  definition of system zeros of  a
square system  as zeros  of  the following polynomial $\psi(s)$
(see Sect 5.2, formula (5.11))
$$ \psi(s)\;=\; det(sI_n-A)det(C(sI_n-A)^{-1}B)
\eqno    (8.19)  $$
 Let the system has $\mu$ zeros, $\mu \le n-r$.
Then the zero polynomial $\psi(s)$
$$   \psi(s) \;=\; a_o
+ a_1s +\cdots+ a_{\mu}s^{\mu} \eqno  (8.20) $$
 with  unknown real
coefficients $a_i$ , $i=0,1,\ldots,\mu $  to be found. Substituting in the right-hand side of (8.19) $\mu +1$
different real  numbers $s_i$, $i=1,2,..., \mu+1$  that  differ from eigenvalues of  the matrix $A$  we result
in the following $\mu+1$ real numbers $b_i$, $i=1,2,..., \mu+1$
$$   b_i\; =\; \psi(s_i) \;=\; det(s_iI_n-A)det(C(s_iI_n-A)^{-1}B)
\eqno (8.21) $$
 Substituting same $s_i$ into (8.20) we  write
$\mu +1$  equations in the coefficients $a_i$ , $i=0,1,\ldots,\mu $
$$   \psi(s_i) \;=\; a_o + a_1s_i +\cdots+ a_{\mu}s^{\mu}_i \;=\;
   [1,s_i,\ldots, s_i^{\mu}] \left [ \begin{array}{c} a_o \\ a_1 \\
    \vdots \\ a_{\mu}     \end{array} \right ]
\eqno  (8.22)  $$
 Equating the left-hand side of (8.21) to the right-hand side of (8.22) for $i=1,2,\ldots, \mu
+1$ we get the following system of linear algebraic equations in unknown  $a_i$, $i=0,1,\ldots,\mu$
$$ \left [ \begin{array}{ccccc} 1 & s_1 & s_1^2 & \cdots &  s_1^{\mu} \\
     1 & s_2 & s_2^2 & \cdots &  s_2^{\mu} \\  \vdots & \vdots & \vdots &
 \cdots & \vdots \\ 1 & s_{\mu +1} & s_{\mu +1}^2 & \cdots & s_{\mu +1}^{\mu}
 \end{array} \right ]  \left [ \begin{array}{c} a_o \\ a_1 \\ \vdots \\
   a_{\mu} \end{array} \right ] \;=\; \left [ \begin{array}{c} b_1 \\ b_2 \\
   \vdots \\ b_{\mu +1} \end{array} \right ]
\eqno(8.23) $$
 The square matrix in (8.23)    is
nonsingular  one if $s_i \neq s_j $ $(i=1,2,\ldots, \mu +1)$ because it is the Vandermonde matrix. Thus, system
(8.23) has the only solution in $a_o, a_1, \ldots, a_{\mu}$.

     For  realization  of  this  method  it  is
necessary to know a number of zeros ($\mu$). For this purpose we can use  results of Section 7.3. If it is
difficult to find $\mu$  then we should change this value by an upper bound $\bar{\mu} \ge \mu $. Then we need
to separate actual system zeros from zeros of the polynomial $a_o + a_1s + \cdots+ a_{\bar{\mu}}s^{\bar{\mu}}$
by finding such $s_i$  that reduce the rank of the system matrix $P(s)$.

     Thus we can write the following algorithm:

     1. Evaluate a number $\mu $ or its an upper estimate $\bar{\mu}$ .

     2. Assign different $s_i \neq \lambda_j$ $(i=1,2,\ldots, \eta +1$,
         $j=1,2,\ldots,n )$  where $\lambda_j$ are eigenvalues of $A$, $\eta =
         \mu $ or $\eta = \bar{\mu}$.

3. Calculate  $b_i$ $(i=1,2,\ldots,\eta +1)$, $\eta = \mu $ or $\eta = \bar{\mu}$.

     4. Build the Vandermonde matrix  (8.23).

     5. Calculate  $a_o, a_1, \ldots, a_{\mu}$, $\eta = \mu $ or $\eta =
         \bar{\mu}$.

     6. If $ \eta = \mu$ then $a_o, a_1, \ldots, a_{\mu}$  are
        coefficients of the zero polynomial; otherwise  ($\eta =
        \bar{\mu}$), select among zeros of the polynomial $a_o + a_1s +
        \cdots+ a_{\bar{\mu}}s^{\bar{\mu}}$   such  $z_i^*$ that reduce
       the rank  of  the matrix $P(s)$ at $s =z_i^*$.

     This method is less laborious  than
the method of Sect.8.1 because it processes  $n\times n$ matrices. But its numerical accuracy depends on
accuracy of inverting
 $n\times n$  matrices, which  may be ill conditioned matrices.
Moreover, the separation of system zeros  from zeros of the polynomial $\psi(s)$ also influences on numerical
accuracy.

\section[Zero computation via matrix polynomial
                    and matrix pencil]
        {Zero computation via matrix \\polynomial
                   and  matrix pencil}

These approaches are based on results of Chapter 7.

\it{METHOD 1.} \rm  To calculate zeros of   controllable square system (1.1), (1.2) it is used the definition of
zeros via the following polynomial
$$   \psi(s) \; = \; s^{n-r\nu }det([O,C_1] + [O,C]s + \cdots +
      [O,C_{\nu -1}]s^{\nu -2} + C_{\nu}s^{\nu -1})
\eqno (8.24) $$
 where $C_i$ are $r\times l_i$ submatrices, $\nu,\; l_1,\; l_2,\ldots,
l_{\nu}$ are the controllability characteristics of the pair $(A,B)$ (see Sect. 7.1 for details).
     The  method is a modification  of  Samash's  method.  We
find a zero polynomial $\psi(s)$ in the form (8.24). Then assigning different real numbers $s_1,
s_2,\ldots,s_{\mu +1}$ that do not coincide with  eigenvalues of the matrix $A$ we  calculate real numbers $b_i$
$(i=1,2,\ldots, \mu +1)$
$$   b_i\; =\; \psi(s_i) \;=\;s_i^{n-r\nu }det([O,C_1] + [O,C]s_i +
            \cdots +[O,C_{\nu -1}]s_i^{\nu -2} + C_{\nu}s_i^{\nu -1})
\eqno (8.25)  $$
 Equating the right-hand side of (8.20) to  $b_i$, $i=1,2,\ldots, \mu +1$ gives the system of
linear equations in $a_i$ ($i=0,1,\ldots, \mu$) used for system zeros  calculation.

     Above  reasonings can be expressed in the following algorithm:

     1.  Determine   controllability  characteristics  of
       system (1.1), (1.2): $\nu$, $l_1, l_2,\ldots,l_\nu $ (see
       formulas (1.45), (1.59)) and the $n\times n$ transformation matrix
        $N$ that reduces the pair $(A,B)$ to Yokoyama's $(n \le r\nu)$
         or Asseo's $(n = r\nu)$ canonical form.

     2. Calculate $r\times l_i$ submatrices $C_i$ , $i=1,2,\ldots, \nu$
       from  the relation $[C_1, C_2, \ldots, C_\nu]\;=\; CN^{-1}$.

     3. Evaluate a zeros number $\mu $ or an upper bound  $\bar{\mu}$
        from  analysis of  blocks $C_\nu$,  $C_{\nu -1},\ldots$ (see
        Sect.7.3).

     4. Assign  different real numbers $s_i \neq \lambda_j$ $(i=1,2,\ldots,
       \eta +1$, $j=1,2,\ldots,n)$ where $\lambda_j$ are eigenvalues of $A$,
        $\eta = \mu $ or $\eta = \bar{\mu}$.

     5. Calculate  $b_i$ $(i=1,2,\ldots,\eta +1)$ from  formula (8.24),
        $\eta = \mu $ or $\eta = \bar{\mu}$.

Further see steps 4-6 of the algorithm of Sect 8.3.

     This algorithm has more steps than   the Samash's  method  but
computational difficulties are decreased because we operate with  matrices of the order $n-r+1$ and $\mu +1$.
Moreover, to calculate submatrices $C_i$ we may  use  formulas of Section 7.2: $C_{\nu} \; = \;
CB\bar{G}_{\nu}^{-1} \; \;$, $C_{\nu-1} \; = \; CAB\bar{G}_{\nu}^{-1} - CB\bar{G}_{\nu}^{-1}F_{\nu \nu}$ and so
on. Also the iterative procedure of the work [S10] may be used to  compute  $F_{\nu 1},\ldots,F_{\nu \nu}$
without any inverting an $n\times n$ matrix $N$.

 \it{METHOD 2.} \rm  It is based on Theorems  7.6  and  7.7  of
Sect.7.4.  Zeros are computed via square matrices of the order $n-r$ or $n - l_{\nu -1}$  as follows:

     1.  If  the  condition  $det(CB) \neq 0$ is satisfied then  system  zeros are
calculated as eigenvalues of the following  $(n-r)\times (n-r)$ matrix
$$ Z(s) \; = \; \left [ \begin{array}{c} \bar{E} \\ \dotfill\\ -\theta
   C_{\nu}^{-1}[C_1, C_2,\ldots, C_{\nu -2}, C_{\nu -1}] \end{array}\right ]
\eqno (8.26)   $$
 where $\theta = I_r$ for $l_{\nu -1}  = r$ or
$\theta = [O,I_{l_{\nu -1}}]$ for  $l_{\nu -1}  < r$, the $(n-2r)\times (n-r)$ matrix  $E$ has the form (7.70),
$[C_1, C_2, \ldots, C_\nu]\;=\; CN^{-1}$, $N$ is a $n\times n$ transformation matrix reducing pair $(A,B)$ to
Yokoyama's canonical form.

     2. If  $det(CB) = 0$  then  system
zeros  are      generalized  eigenvalues  of  the following matrix pencil of order $n-l_{\nu -1}$
$$ Z(s) \; = \; s \left [ \begin{array}{cc} I_{\beta} & O \\
                   O & C_{\nu} \end{array} \right ] \; + \;
         \left [ \begin{array}{c} -\bar{E} \\ \dotfill\\
    C_1, C_2,\ldots, C_{\nu -2}, [O,C_{\nu -1}] \end{array}\right ]
\eqno (8.27)   $$
  where $\beta = l_1+l_2+\cdots+l_{\nu -2} = n-r-
l_{\nu -1}$, the $\beta\times (n-l_{\nu -1})$  matrix $\bar{E}$ has the form (7.70) with $E_{\nu -2,\nu -1} \; =
\; [O, E_{\nu -2,\nu -1}]$ being the $l_{\nu -2}\times r$ submatrix, $[C_1, C_2, \ldots,C_\nu]\;=\; CN^{-1}$.

     The present approach  decreases computational
difficulties because it processes square matrices of the order  $n-r$ or $n-l_{\nu -1}$.
     To calculate submatrices $C_i$ $(i=1,2,\ldots,\nu )$ we
also can use formulas of Sect.7.2  that allow to avoid inverting  an $n\times n$  matrix $N$. For computing
eigenvalues or generalized eigenvalues of (8.26) or (8.27)  we can use  QZ algorithm.

\it{REMARK 8.1.} \rm These computational methods are applied to  the well conditional controllability matrix
$Y_{AB} = [B,AB, \ldots,A^{n-1}B]$.  Otherwise, we can recommend to make the following operations:

     1. Separate a well conditioned part of $Y_{AB}$,

     2. Decrease a dimension of  the controllability subspace,

     3. Calculate  system zeros of the system obtained. They
        form the set of invariant zeros.

     The rest of zeros (input decoupling zeros ) may be found by  another approach.

\chapter[Zero assignment]
        {Zero assignment}

     In this chapter we  consider the assignment of
system zeros  by choosing the output (input) matrix or by using the squaring down operation. This problem is
caused by large influence of  system zeros on dynamic behavior of any control system. It is known from the
classic control theory that right-half plan zeros create severe difficulties for the control design.
Multivariable systems have  similar properties. For example,   maximally achievable accuracy of optimal
regulators and/or filters is attainable if an open-loop system  has not  zeros in the right-half part of the
complex plane [K4], [K5]. In Section 10 we demonstrate that solvability conditions  of different type tracking
problems also contain limitations on system zero locations.

System zeros are invariant under state and/or output proportional feedbacks. They may be  shift  only  by  an
appropriate correction of output and/or input matrices.\footnote{Later  we will use only an output matrix for
the zero assignment.} We consider two different techniques:  the first one shifts an output matrix and second
one calculates a squaring down compensator.

\section[Zero assignment by selection of output  matrix]
        {Zero assignment by selection of output  matrix}

This approach is based  on    appropriate  choice  of  an output matrix $C$ of a system  and  may be recommended
when a freedom 'may still exist to choose that sets of variables are to be manipulated and what sets to be
measured for control purpose'[M1].  For example, an output matrix $C$ may be selected in an estimation system
where the whole state vector is available for measurement.

The   zero assignment  problem  by choosing an output matrix has  been  set up  by Rosenbrock in [R1,Theorem
4.1] and briefly is  formulated  as follows.

     Let  in    given  square  system  (1.1),  (1.2)  with  the
completely controllable pair $(A,B)$ we may choose  the $r\times n$ matrix $C$ in (1.2). Form the problem : it
is necessary to choose a matrix $C$ so that

     1. the pair $(A,B)$ is observable,

     2. The Smith-McMillan form  $M(s)$ of the $r\times r$ transfer  function  matrix $G(s) =C(sI_n-A)^{-1}B$
        has  assigned numerator polynomials  $\epsilon_i(s)$.

Rosenbrock has shown that such $C$ always can be chosen if degrees of assigned $\epsilon_i(s)$ satisfy certain
conditions.  But there are some restrictions on the choice of the output (measurable) matrix, which also should
take into  account, for example,  a fullness  rank  of the matrix $C$,  well-posed of the observability matrix
$Z_{AC}$  and others. Moreover, Rosenbrock's necessary conditions on degrees of polynomials $\epsilon_i(s)$ are
rather complicated .

Further we will  ensure \it {distinct} \rm  assigned zeros. This restriction on zero locations considerably
simplifies Rosenbrock's conditions. Also we will ensure fullness of a rank of  the matrix $C$ and  several other
requirements.  The  methods  are based on works [S6], [S11], which use zero definitions via a matrix polynomial
and the reduced (lower order) matrix (see Chapter 7).

\subsection[Iterative method of zero assignment]
           {Iterative method of zero assignment}

     We  consider completely controllable system (1.1)
with $r$ independent inputs  $(rankB = r)$ and  the completely measurable state vector $x$.

    Let us assign $n-r$ \footnote{The maximal number of zeros of a proper
square system (see sect.7.3).} distinct real numbers $\bar{s}_i$ $(i=1,2, \ldots,n-r)$ and denote by
$\bar{\psi}(s)$  the following zero polynomial
$$    \bar{\psi}(s)\;=\; \prod_{i=1}^{n-r}(s-\bar{s}_i)
\eqno (9.1)  $$

     We consider the following problem.

\it{PROBLEM  1.} \rm To choose  the $r\times n$ output matrix $C$ that assigns  zeros  of system (1.1) with the
output
$$       y \;=\; Cx
\eqno (9.2)  $$
 at desired positions  $\bar{s}_i$.
 Matrix $C$ must also satisfy the following requirements:

$$ (a).\;\; \rm{} \; \rm{the \;pair} \; \it (A,B)\; \rm{is}\;\rm{observable}, $$
$$ (b).\; \;\;rank C \;=\; r \qquad \qquad\eqno (9.3)$$

Now  we find conditions, which assure a  solution of the problem. These conditions may be considered as a
particular case of Rosenbrock's ones.

\it{THEOREM 9.1.} \rm  If the pair of matrices $(A,B)$ is controllable and  distinct zeros $\bar{s}_i$
($i=1,2,\ldots, n-r$)  do not coincide  with eigenvalues of the matrix $A$ then there exists a matrix  $C$ that
ensures both the assigned zero polynomial to  system  (1.1), (9.2) and the
 condition  (9.3a).

\it{PROOF} \rm [S6]. We  base on Theorem 4.1 from  [R1, p.186] where general zero assignment conditions are
defined.

At first we consider  values  $\mu_1$,$\mu_2$,\ldots,$\mu_q$ that are nonzero minimal indices of the singular
matrix pencil $(sI_n-A,B)$ with  the $n\times p$  matrix $B$ ($rank B=q \le p$), $1\le \mu_1 \le \mu_2  \le
\cdots \le \mu_q$, $\;\;\mu_1 + \mu_2 + \cdots + \mu_q = n$. As  it has been studied in [R1], the minimal
indices coincide with ordered numbers $\beta_i$ obtained from
 the sequence of linearly independent vectors
$$ b_1,Ab_1,\ldots A^{\beta_1-1}b_1, b_2,Ab_2,\ldots A^{\beta_2-1}b_2,
\ldots ,b_q,Ab_q,\ldots A^{\beta_q-1}b_q \eqno  (9.4) $$ which are
selected from  vectors of the controllability matrix
$$ Y_{AB} \;=\; [B,AB,\ldots,A^{n-1}B]\;=\; [b_1,b_2,\ldots,b_q,
Ab_1,Ab_2,\ldots,Ab_q, \ldots,A^{n-1}b_q]
 \eqno(9.5) $$
To form (9.4) we accept a such new vector from (9.5), which is not linearly depended from all previously
accepted vectors; otherwise we reject it. When $n$ vectors have been accepted then we arrange them in the order
(9.4). It is evident that a number of  nonzero minimal indices is equal to the rank of $B$.

    Then we consider invariant polynomials  $\epsilon_i(s)$ of  the  transfer
function matrix of the constructed system. Zeros of $\epsilon_i(s)$, taken all together, form the set of system
zeros of the completely controllable and observable system. Now, by Theorem 4.1 [R1], the degrees of   desired
invariant polynomials $\epsilon_i(s)$ for the controllable  pair $(A,B)$ with an $n\times p$ matrix $B$ of
$rank\;q \le p $ must satisfy conditions:

     (a). the number of non identically zero invariant polynomials
       is equal to  $r \le q$  where  $q$  is the number of nonzero
       minimal indices of the pencil $(sI_n-A,B)$,

     (b). $\epsilon_i(s)$ divides $\epsilon_{i+1}(s)$, $\;\;i=1,2,\ldots,r-1$,

     (c).  degrees $\delta(\epsilon_i)$ of nonzero $\epsilon_i(s)$
     satisfy the following inequalities
$$  \sum _{i=1}^{k}\delta(\epsilon_i) \le \sum _{i=1}^{k}(\mu_{p-r+i}-1), \;\;
      k=1,2,\ldots, r
\eqno (9.6)   $$

       (d). $r = q$; polynomials $\epsilon_r(s)$ and $\phi_1(s)$ are
        relatively prime where $\phi_1(s)$ is  a  minimal  polynomial
       [G1] of the matrix $A$.

     Now we need to show that the conditions of above Theorem 9.1 are a particular
case of  the conditions  (a)- (d)  for distinct assigned zeros and an  $n\times r$  matrix $B$ of rank $r$.

For distinct zeros polynomials $\epsilon_i(s)$ $(i=1,2,\ldots, r)$ become
$$  \epsilon_1(s) = \epsilon_2(s) =\cdots= \epsilon_{r-1}(s) =1,\;\;
      \epsilon_r(s) = \psi(s)
\eqno   (9.7)  $$
 where $\psi(s)$ is the zero  polynomial. Equalities (9.7) ensure the condition (b).
  The condition (a) and the first part of (d)  are assured by the assumption $rankB = r$. The condition (c) for
$p = r$ and $\epsilon_i(s)$  satisfied (9.7) may be rewritten as follows
$$  0 \le \sum _{i=1}^{t}(\mu_i-1), \;\; t=1,2,\ldots, r-1,\;\;\;
     n-r \le \sum _{i=1}^{r}(\mu_i-1)
\eqno (9.8)   $$
 Since $\mu_i \ge 1$ $(i=1,2,\ldots, r)$ and $
\mu_1 + \mu_2 + \cdots + \mu_r = n$   then  inequalities (9.8) always are true. The second part of the condition
(d) is also satisfied in accordance with the hypothesis  of Theorem 9.1 about controllability of the  pair
$(A,B)$.

Therefore, all conditions of Theorem 4.1 from [R1]  are carried
out. This completes the proof.

Now we consider a method for calculating the matrix $C$. In accordance with the structural restrictions on $C$
we study two cases.

\it{CASE 1. OUTPUT MATRIX WITHOUT STRUCTURAL RESTRICTIONS}. \rm We assume that all elements of the matrix $C$
may be  any bounded numbers. For square system (1.1), (9.2) we consider the following definition of the zero
polynomial from Sect. 7.1
$$ \psi(s)\; = \;s^{n-r\nu}det\tilde{C}(s)
\eqno (9.9) $$ where the $r\times r$ polynomial matrix $\tilde{C}(s)$ has the following structure
$$  \tilde{C}(s) = [O,C_1] + [O,C_2]s + \cdots +  [O,C_{\nu -1}]s^{\nu -2}
                +C_{\nu}s^{\nu -1}
\eqno(9.10)  $$
 Here $\nu $ is an index of controllability of
$(A,B)$ (see formula (1.45)) and $r\times l_i$ submatrices $C_i$ calculated from the expression
$$   CN^{-1} \;=\; [C_1, C_2,\ldots, C_\nu ]
\eqno (9.11)    $$
 where an  $n\times n$  matrix $N$ reduces the pair $(A,B)$
to Yokoyama's canonical form and integers $l_i$ $(i=1,2,\ldots,\nu )$ are defined from  formulas (1.59).

Since  zeros of the polynomial (9.9) differ from assigned numbers $\bar{s}_i$ $(i=1,2,\ldots, n-r)$ then we
obtain the following equalities at  points $s = \bar{s}_i$
$$ \psi(\bar{s}_i)\; = \;\bar{s}_i^{n-r\nu}det\tilde{C}(\bar{s}_i) \; \neq 0,
     \;\; i=1,2,\ldots n-r  $$
Thus, values $\psi(\bar{s}_i)$ are functions of elements $c_{kl}, k=1,2,\ldots,r, l=1,2,\ldots,n$ of the matrix
$C$. To find these elements we consider the minimization  of  the following performance criterion with respect
to elements of $C$
$$      J \;=\; = J_1 +qJ_2, \;\; q \ge 0
\eqno   (9.12) $$  where $q>0$  is a weight coefficient and
$$    J_1 \;=\; 0.5 \sum _{i=1}^{n-r} \psi(\bar{s}_i)^2,
\eqno   (9.13) $$
$$     J_2 \;=\; (det(CC^T))^{-1}
\eqno (9.14)    $$
 In (9.12) the first term ($J_1$) depends on
zero locations of system (1.1), (9.2) and the second term ($J_2$) depends on the  rank  of  the matrix $C$
because it is the inversion of Gram's determinant [G1] builded from rows of  $C$. Since Gram's determinant is a
nonnegative value then  (9.14) ensures the rank fullness of $C$. So, the minimization of (9.12)  with respect to
elements of matrix $C$  guarantees assigned zeros to system (1.1), (9.2) and conditions (9.3a,b).

\it{REMARK 9.1. } \rm  To improve the conditionality of the observability matrix $Z_{CA}$ we replace (9.14) by
the following modify criterion
$$  \bar{J}_2\;=\;det(Z_{\gamma}(C,A)^T Z_{\gamma}(C,A))^{-1},\;\;
     Z_{\gamma}(C,A)^T = [C^T,A^TC^T,\ldots,(A^T)^{\gamma -1}C^T]
\eqno (9.15)    $$
 In (9.15) $\gamma$  is  equal  to
the smallest integer from  $n/r$.

     The numerical minimization of  criterion
(9.12) is realized by the following simple iterative scheme
$$ c_{kl}^{(i+1)}\; =\; c_{kl}^{(i)} \;-\; \alpha \frac{\partial J}
{\partial c_{kl}}\mid^{(i)}, \qquad k=1,2,\ldots,r,\; l=1,2,\ldots,n \eqno  (9.16) $$ where $\alpha > 0$  is a
some constant, $dJ/dc_{kl}$ is a gradient of $J$ with the respect to elements $c_{kl}$.

     Let us find an analytic formula for  $dJ/dc_{kl}$. For this  we
 apply repeatedly  the following equality  from [A3]
$$\frac{\partial f(Z(X))}{\partial x_{ij}}\;=\;tr \{ \frac{\partial f(Z(X))}
{\partial Z}\;\frac{\partial Z^T}{\partial x_{ij}} \} \eqno (9.17)
$$ where $Z$ and $X$ are some rectangular matrices. At first we express $\partial J_1/
\partial c_{kl}$ as
$$\frac{\partial J_1}{\partial c_{kl}} \;=\; \sum _{i=1}^{n-r} \psi(\bar{s}_i)
   \frac{\partial \psi(\bar{s}_i)}{\partial c_{kl}}
\eqno  (9.18)$$
 To calculate $\partial \psi(\bar{s}_i)/\partial
c_{kl}$ we employ (9.17) to the expression (9.9)
$$  \frac{\partial \psi(\bar{s}_i)}{\partial c_{kl}} \;=\;\bar{s}_i^{n-r\nu} tr \{
     \frac{\partial (det\tilde{C}(\bar{s}_i)}{\partial \tilde{C}(\bar{s}_i)}\;
      \frac{\partial \tilde{C}(\bar{s}_i)^T}{\partial c_{kl}} \}   $$
Applying the following formula from [A3]
$$  \frac{\partial (detX)}{\partial X} \;=\; detX(X^{-1})^T  $$
and using the property of  operation $tr<.>$ : $tr(XY)
=tr(Y^TX^T)$ we can rewrite the last relation as follows
$$\frac{\partial \psi(\bar{s}_i)}{\partial c_{kl}} \;=\;\bar{s}_i^{n-r\nu}
tr \{\frac{\partial \tilde{C}(\bar{s}_i)^T}{\partial c_{kl}}\; adj(\tilde{C}
(\bar{s}_i)) \}
\eqno  (9.19)  $$

To find $\partial \tilde{C}(\bar{s}_i)/\partial c_{kl}$ we  represent $\tilde{C}(\bar{s}_i)$ in the more
convenient form. Partitioning the $n\times n$ matrix $N^{-1}$ in (9.11) into $\nu $   blocks of sizes $n\times
l_i$
$$      N^{-1} \;=\; [P_1, P_2, \ldots, P_{\nu}]
\eqno  (9.20)  $$ and representing  $CN^{-1} \;=\; C[P_1, P_2, \ldots, P_{\nu}] $ we express blocks $C_i$ in
(9.10) as $C_i = CP_i \;\;(i=1,2,\ldots, \nu)$. Thus
$$  \tilde{C}(s) = C \{ [O,P_1] + [O,P_2]s + \cdots +  [O,P_{\nu -1}]s^{\nu -2}
                +P_{\nu}s^{\nu -1} \}
\eqno(9.21)  $$
  where $[O,P_i]$   are $n\times r$ matrices. Differentiating  (9.21)
$$ \partial \tilde{C}(\bar{s}_i)/\partial c_{kl} \;=\; E^{kl}_{r\times n}
(\sum_{t=1}^{\nu }[O,P_i]\bar{s}_i^{t-1}) \eqno  (9.22)  $$
 where
$E^{kl}_{r\times n}$ is the $r\times n$ matrix having the unit $kl$-th element and  zeros otherwise and
substituting  (9.22)  into (9.19) and   the result into (9.18) yields
$$ \frac{\partial J_1}{\partial c_{kl}} \;=\; \sum _{i=1}^{n-r} \psi(\bar{s}_i)
     \bar{s}_i^{n-r\nu} tr \{ E^{kl}_{r\times n} ( \sum_{t=1}^{\nu }
      [O,P_i]\bar{s}_i^{t-1})adj(\tilde{C}(\bar{s}_i)) \}
\eqno (9.23)  $$
     Then using  equality (9.17) and property of the operation $tr<.>$
we calculate $\partial J_2/\partial c_{kl}$
$$ \frac{\partial J_2}{\partial c_{kl}} \;=\; \frac{\partial (det(CC^T)^{-1})}
{\partial c_{kl}}\;=\; -det(CC^T)^{-2}tr \{ \frac{\partial (det(CC^T))}
{\partial (CC^T)}\; \frac{\partial (CC^T)^T)}{\partial c_{kl}}\}  \;=$$
$$\;=\; -det(CC^T)^{-2}tr \{ det(CC^T)(CC^T)^{-1T}\frac{\partial (CC^T)^T)}
{\partial c_{kl}} \} \;=\;
 -det(CC^T)^{-1}tr \{ \frac{\partial (CC^T)}{\partial c_{kl}} (CC^T)^{-1} \}
\eqno  (9.24)  $$
Since
$$  \frac{\partial (CC^T)}{\partial c_{kl}} \;=\; E^{kl}_{r\times n}C^T +
     C(E^{kl}_{r\times n})^T
\eqno (9.25) $$
then substituting (9.25) into the right-hand side
of (9.24) and denoting  $E^{lk}_{n\times r}\;=\; (E^{kl}_{r\times
n})^T$ yields the final expression for $\partial J_2/\partial
c_{kl} $
$$ \frac{\partial J_2}{\partial c_{kl}} \;=\; -det(CC^T)^{-1}tr
\{ (E^{kl}_{r\times n}C^T + CE^{lk}_{n\times r})(CC^T)^{-1} \}
\eqno (9.26) $$
Uniting (9.23) and (9.26) we result in the general
formula for $\partial J/
\partial c_{kl}$
$$   \frac{\partial J}{\partial c_{kl}} \;=\; \sum _{i=1}^{n-r} \psi(\bar{s}_i)
 \bar{s}_i^{n-r\nu} tr \{ E^{kl}_{r\times n} ( \sum_{t=1}^{\nu }
  [O,P_i]\bar{s}_i^{t-1})adj(\tilde{C}(\bar{s}_i)) \} \;-\;$$
$$   - q det(CC^T)^{-1}tr \{ (E^{kl}_{r\times n}C^T +
     CE^{lk}_{n\times r})(CC^T)^{-1} \}
\eqno (9.27) $$ To calculate $\partial \bar{J}_2/\partial c_{kl}$ we find the $n\times n$  matrix $Z_{\gamma}
\;=\;Z_{\gamma}(C,A)^T Z_{\gamma}(C,A)$
$$ Z_{\gamma} \;=\; C^TC + A^TC^TCA + \cdots +
(A^T)^{\gamma -1}C^TCA^{\gamma -1} \eqno  (9.28) $$ and carrying out the similar operations (see (9.24)) obtain
$$ \frac{\partial \bar{J}_2}{\partial c_{kl}} \;=\;
\frac{\partial (det(Z_{\gamma})^{-1})}{\partial c_{kl}} \;=\;
 -det(Z_{\gamma})^{-1} tr \{ \frac{\partial Z_{\gamma}}
 {\partial c_{kl}} Z_{\gamma}^{-1} \}
\eqno  (9.29) $$
 where $\partial Z_{\gamma}/\partial c_{kl}$ is calculated  as
$$  \frac{\partial Z_{\gamma}}{\partial c_{kl}}\;=\;
\sum_{t=0}^{\gamma -1}(A^T)^t (E^{kl}_{r\times n}C^T + CE^{lk}_{n\times r}) A^t
 $$
Substituting the last expression into (9.29) we get  the  final
formula for $\partial \bar{J}_2/\partial c_{kl}$
$$ \frac{\partial \bar{J}_2}{\partial c_{kl}} \;=\; -det(Z_{\gamma})^{-1} tr
    \{ (\sum_{t=0}^{\gamma -1}(A^T)^t(E^{lk}_{n\times r}C +
    C^TE^{kl}_{r\times n}) A^t) Z_{\gamma}^{-1} \}
\eqno (9.30)  $$
    In the final we summarize the results as the following
algorithm  for   zero assignment:

     1. Check  controllability of the pair  $(A,B)$.  If it
         is completely controllable then go to the next
        step; otherwise the problem has no solution.

     2. Assign  $n-r$ desirable  distinct real zeros $\bar{s}_i$ $(i=1,2,\ldots,
        n-r) $, which don't coincide with eigenvalues of $A$.

     3. Calculate controllability characteristics of the  pair
        $(A,B)$ : $ \nu, \;l_1,\; l_1,\;\ldots,\; l_{\nu }$ (see formulas (1.45),
        (1.59)) and the transformation  $n\times n$  matrix $N$. Determine
         $N^{-1}$  and partition it into $\nu $   blocks in according with
        (9.20).

     4.  Calculate $\partial J_1/\partial c_{kl}$ (9.23) by formulas (9.9), (9.10)

     5. In according with a chosen criterion $J_2$ or $\bar{J}_2$  calculate
$\partial J_1/\partial c_{kl} $ or  $\partial \bar{J}_2/\partial c_{kl} $
 by formulas (9.26) and (9.30) respectively.

     6. Calculate $C^{(i+1)}$ by the  recurrent scheme (9.16).

 7. If $ \parallel \partial J/\partial C \parallel = \sum_{k=1}^{r}
     \sum_{l=1}^{n} \mid \partial J/\partial c_{kl} \mid >
     \epsilon$  where $\epsilon >0 $ is a given real number, then go to step 4;
     otherwise the end of calculations.

     The  method may  be illustrate by the
following example [S6].

                \it{EXAMPLE 9.1.} \rm

Consider a completely controllable system  with two inputs and outputs
$$ \dot{x} \;=\; \left [ \begin{array}{rrr}  -2  & 0.9374 & -2.062 \\
           2 &  -0.4375 &  2.562 \\ -1 &  -1.563 &  -1.562  \end{array}
         \right ] x \;+\; \left [ \begin{array}{rr} 1 & 0 \\ -1 & 1 \\ 1 & 1
         \end{array} \right ] u
\eqno(9.31) $$
$$  y \;=\; \left [ \begin{array}{rrr} 1 &  0 & 0 \\ 1  &  1  &  0
     \end{array} \right ] x
\eqno(9.32) $$
 Since here $n-r\; =\; 3-2 \;= \;1$ then the  system has no more than one zero. To calculate a
zero polynomial we find  $\nu =2 ,\; l_1=1,\; l_2 =2 \;$ and the transformation matrix $N$ of order 3 which
reduces the pair $(A,B)$ to Yokoyama's canonical form
$$  N \;=\; \left [ \begin{array}{rrr} -0.5 & -0.25 & 0.25 \\ 1&1& 0 \\
        0.25 & -0.75 & 0 \end{array} \right ]  $$
Calculating
$$  CN^{-1} \;=\; \left [ \begin{array}{rrr} 0 & 0.75 & 1\\ 0 & 1 & 0
    \end{array} \right ] $$
with
$$    C_1 \;=\; \left [ \begin{array}{r} 0 \\ 0  \end{array} \right ], \;\;
       C_2 \;=\; \left [ \begin{array}{rr} 0.75 & 1 \\ 1 & 0
    \end{array} \right ]  $$
and using formulas (9.9),(9.10) we find  the zero polynomial $\psi(s)$
$$  \psi(s) \;=\; s^{3-4}det \{ \left [ \begin{array}{rr} 0 & 0\\ 0 & 0
      \end{array} \right ] \; + \left [ \begin{array}{rr} 0.75 & 1 \\ 1 & 0
      \end{array} \right ]s \} \;=\; -s $$
Therefore,  system (9.31),(9.32) has a zero in the origin.

     Assuming that all state variables are accessible  we try to find a
new output matrix $\bar{C}$  that shifts the zero to the value $-1$. Since  eigenvalues of $A$ are equal to
$-0.5, -1.5, -2$ then Theorem  9.1  is  satisfied. Calculating  $\gamma =1$ and forming the performance
criterion  $ J = J_1 +q\bar{J}_2 $ with $J_1$ from (9.13) and $\bar{J}_2 $ from (9.15) with for  $n-r = 1$,
$\bar{s}_i = \bar{s}_1 = -1$, $ l_1 = 1, l_2 = \nu = 2$ and $ q = 0.25$ we get
$$ J \;=\; 0.5(\psi(\bar{s}_1))^2 + 0.25(det(\bar{C}^T\bar{C}))^{-1} \;=\;
     0.5\bar{s}_1^{-2}\{det(\bar{C}([O,P_1] + P_2]\bar{s}_1)) \}^2 +
     0.25(det(\bar{C}^T\bar{C}))^{-1}  $$
where $P_1$, $P_2$  are respectively $3\times 1$, $3\times 2$ blocks of the matrix  $N^{-1} =[P_1, P_2]$. The
numerical minimization of  this  criterion by the recurrent scheme (9.16) is finished as
 $ \parallel \partial J/\partial \bar{C} \parallel \le 0.02$. We result in the
following matrix
$$ \bar{C} \;=\; \left [ \begin{array}{rrr} 1.242 & 1.098 & -0.2423 \\
    -0.027 & 1.343 & -0.3478 \end{array} \right ]
\eqno (9.33) $$ which creates  the system  having the  system zero $-1.000$. Rounding off elements of (9.33)
yields the matrix
$$ \bar{C}_r \;=\; \left [ \begin{array}{rrr} 1.2 & 1.1 & -0.2\\
         0 & 1.3 & -0.35 \end{array} \right ] $$
which ensures the zero $-0.9$.

  \it{CASE 2. OUTPUT MATRIX WITH STRUCTURAL RESTRICTIONS}. \rm  In most
practical situations only $m$ ($m<n$) components  of  a  state  vector are accessible. Without loss of
generality we may  assume that these components are the first $m$ elements $x_1, x_2 ,\ldots, x_m$ of $x$. Any
linearly independent  combinations of these components, which form an $r$-vector $y$, are realized by the
$r\times n$ output matrix $C$ of the structure
$$    C \;=\;  [\; C_m,\; O \;]
\eqno(9.34)   $$
 where $C_m$ is an $r\times m$ submatrix of full rank $r$. Let's consider  the zero assignment
problem by choosing  the structural restricted matrix $C$  (9.34).

   At first we show that this problem has no solution for
$m=r$. Indeed, let system (1.1), (1.2) with  $C\;=\;[\; C_m,\; O \;] $ has a zero polynomial  $\psi(s)$.
Changing the output vector (1.2)  by  $ \bar{y} \;=\; [\;\bar{C}_m,\;O\;]$  with $r$ independent components we
build a new system with a zero polynomial $\bar{\psi}(s) \neq \psi(s)$. But, since $r\times r$ matrices $C_m$
and $\bar{C}_m $ are nonsingular ones then the vector $\bar{y}$ is expressed via $y$ as follows
$$ \bar{y} \;=\;[\;\bar{C}_m,\;O\;]x\;=\; \bar{C}_m C_m^{-1}[\;C_m,\;O\;]x\;=\;
    \bar{C}_m C_m^{-1}y  $$
We can see that  the new output $\bar{y}$  is obtained from the old one by a nonsingular transformation of its
components, hence zero polynomials of these systems must be similar: $\bar{\psi}(s)\;= \; \psi(s)$ and system
zeros are not shifted.

     If $m>r$ then the  zero  assignment problem can  be solvable
because the output matrix has enough number of  free  variables to minimize the  criterion (9.12).
%
%
To solve the zero assignment problem we may use the  above recurrent scheme (9.16).

\it{REMARK 9.2.} \rm  If $q$ is a rank deficiency of the  matrix  $CB =[C_m,O]B$ then a new system will have
less than $n-r-q$ zeros (see Sect. 7.3). Thus, we need to   change  the value $n-r$ by $n-r-q$ in (9.13).

                  \it{EXAMPLE  9.2.} \rm

 For illustration we  consider the zero assignment problem
for the following controllable system [S6]
$$ \dot{x} \;=\; \left [ \begin{array}{rrrr} 14.39 & -62.43 & -30.81 & 10.33 \\
    3.752 & -19.39 & -10.0 & 2.997 \\ 1 & 1 & 1 & 0  \\
     -0.867 & -1.267 &-1.8 & -0.6  \end{array} \right ] x \;+\;
     \left [ \begin{array}{rr} 3 & -1 \\ 1 & 0 \\ 0 & 0 \\ 0 & -1
     \end{array} \right ] u
\eqno(9.35) $$
$$  y \;=\; \left [ \begin{array}{rrrr} 1 & 1 &  0 & 0 \\ 1  &  0  &  0 & 0
  \end{array} \right ] x
\eqno(9.36) $$
 At first using results of Sect.1.2.3 we calculate $\nu=2, l_1=l_2= 2$  and the
transformation $4\times 4$ matrix $N$ reducing the pair $(A,B)$ to Asseo's canonical form\footnote{This system
has Asseo's canonical form because of $\nu = n/r = 2$.}
$$ N \;=\; \left [ \begin{array}{rrrr} -0.2642 &  0.7926 &  0.6348 &  0.2642 \\
    0.660 & 0.1982 &  0.4062 &  0.066 \\ -0.4719 & 1.416 & 0.3622 & -0.521  \\
    0.132 & 0.6039 & 0.3408 & -0.132 \end{array} \right ]  $$
Using formulas (9.9)-(.11) we  find the zero polynomial $\psi(s) = s^2-s-2$ having zeros $s_1=1, s_2 = -2$.
Thus, system (9.35), (9.36) has the  right-half zero and we can consider the problem of zero shifting to
locations: $-1, -2$. Let us assume that only three first components $x_1, x_2, x_3$ of $x$ are accessible ($m =
3$). Thus, the problem with structural restricted matrix $C=[C_m,O]$  may have a solution because $ m=3 > r=2$.
Since conditions of Theorem 9.1 are satisfied then calculating $\gamma = n/r = 2$, forming  the performance
criterion $ J = J_1 +q\bar{J}_2 $ with $n-r = 2$, $\bar{s}_1 = -1, \bar{s}_2 = -2 $, $ l_1 =l_2 = 2, \;\;\nu =
2,\;\; q = 0.5$ and minimizing $J$ by the recurrent  scheme (4.16) for $k=1,2$; $l=1,2,3$ we get the following
structural restricted matrix
$$ C \;=\; \left [ \begin{array}{rrrr}-0.1252 & 0.3741 & 0.8442 &  0 \\
   0.1118 & 0.7959 & 0.723 & 0  \end{array} \right ]
\eqno (9.37) $$
  that assures system zeros at locations: $ -0.9985,
-2.000$.  The minimization process has been finished as $
\parallel \partial J/\partial \bar{C} \parallel \le 0.01$.

\subsection[Analytical zero assignment]
           {Analytical zero assignment}

    In this section we try to find an analytical solution of Problem 1 with one an
additional requirement
$$    rank (CB) \;=\; r
\eqno        (9.38)  $$
 This restriction ensures that a new system
has exact $n-r$ zeros. The method have been suggested in [S11]. Let's note that  (9.38) is contained in the
above requirement (9.1) on a number of assigned zeros.

     At first we note that the output matrix $C$ that
satisfies the condition (9.3a) always exists if assigned zeros are differ from  eigenvalues of the matrix $A$.
Then we ought to find the  matrix $C$  that ensures  conditions (9.3b) and (9.38).

We recall that a controllable system with  $r$ inputs and outputs that  satisfies the condition  $det(CB) \neq 0
$ has system zeros coinciding with eigenvalues of the following $(n-r)\times (n-r)$ matrix (see Theorem 7.6 from
Sect. 7.4)
$$  Z \; = \; \left[ \begin{array}{ccccc} O & [O,I_{l_1}] & O & \cdots & O \\
       O& O& [O,I_{l_2}]& \cdots & O \\
       \vdots & \vdots & \vdots & \ddots & \vdots\\
             O & O & O & \cdots & [O,I_{l_{\nu -2}}] \\
       \dotfill & \dotfill & \dotfill & \cdots & \dotfill \\
    -T_{\nu 1} & -T_{\nu 2} & -T_{\nu 3} & \ldots & -T_{\nu ,\nu -1}
  \end{array} \right ]
\eqno(9.39)     $$
 where $[O,I_{l_i}]$  are $l_i\times l_{i+1} $
submatrices, $I_{l_i}$ are $l_i\times l_{i} $ unity blocks,
$l_{\nu -1}\times l_i $ submatrices $T_{\nu i}$, $i=1,2,\ldots,\nu
-1$ have the form
$$ T_{\nu i} \; = \;  C_{\nu}^{-1}C_i ,\qquad  l_{\nu -1} = r
\eqno(9.40a) $$
$$ T_{\nu i} \; = \; [O,I_{l_{\nu -1}}]C_{\nu}^{-1}C_i, \qquad l_{\nu -1} < r
\eqno (9.40b)   $$
 In (9.40) $C_i$  are  $r\times l_i$  blocks of
the matrix (9.11). Since $rank C_{\nu }= rank (CB)$ (see Assertion 7.4 from Sect. 7.3), hence, if the $r\times
r$ submatrix $C_{\nu }$ is constructed as a full rank matrix then the problem of zero assignment can be
reformulated as follows: Find an $l_{\nu -1}\times (n-r)$  submatrix $T\;=\;[-T_{\nu 1}, -T_{\nu 2},
 -T_{\nu 3}, \ldots,-T_{\nu ,\nu -1}]$  which places eigenvalues of the
matrix (9.39) at desirable locations.
     Consequently, the zero assignment problem is reduced to  the
eigenvalue assignment problem. Now we show that this problem
always has a solution.

\it{ASSERTION 9.1.} \rm  For any given polynomial $\psi^*(s)$ of order $n-r$ there is  an  $l_{\nu -1}\times
(n-r)$ submatrix $T\;=\; [-T_{\nu 1}, -T_{\nu 2}, -T_{\nu 3}, \ldots,-T_{\nu ,\nu -1}]$ such that  zeros of
polynomials $det(sI_{n-r} -Z)$ and $\psi^*(s)$ are similar.

\it{PROOF}. $\;$ \rm At first we let $l_{\nu -1} =1$. Then, $l_1=l_2= \cdots= l_{\nu}= 1$ and  $T\; =\; q \;=\;
[-q_1, -q_2,\ldots,-q_{n-r}]$ is a vector-row. The matrix $Z$ becomes the following companion form
$$ Z \;=\; \left[ \begin{array}{ccccc} 0 & 1 & 0 & \cdots & 0 \\
      0 & 0 & 1 & \cdots & 0 \\ \vdots & \vdots & \vdots & \ddots & \vdots \\
     0 & 0 & 0 & \cdots & 1 \\ -q_1 & -q_2 & -q_3 & \cdots & -q_{n-r}
\end{array} \right]
\eqno(9.41) $$
It follows from the obviously equality
$$    det(sI_{n-r} -Z)\;=\; s^{n-r} + q_{n-r}s^{n-r-1} +\cdots +q_1 $$
that the vector-row $[-q_1 , -q_2 , -q_3 , \ldots , -q_{n-r}]$ always exists such that the polynomial in the
right-hand side of the last expression is the assigned polynomial $\psi^*(s)$.

 Now we consider the case $l_{\nu -1} >1$. We  construct first $l_{\nu -1} -1$
rows of $T$  in such a way that every row has the only unit
element and the rest elements are zeros. Moreover,  unit elements
are situated in such columns of $Z$  that its first $n-r-1$ rows
form a submatrix with the only unit element in the every column
except the first one. The rest elements of this submatrix are
zeros.
     Elements of the last row of $T$ are uncertainty ones.
We denote  them  by $ -q_1 , -q_2 , -q_3 , \ldots , -q_{n-r}$  and conclude that the matrix $Z$ is obtained from
the matrix (9.41) by appropriate permutations of all rows excluding  the last one. Thus, we can write
$$   det(sI_{n-r} -Z)\;=\; (-1)^{\alpha} (s^{n-r} + q_{n-r}s^{n-r-1}
 +\cdots + q_1 ) $$
where $\alpha$ is the number of row permutations. It  is  evident  that
 elements $q_i$ $(i=1,2,\ldots, n-r)$ can be assigned so that
$$  det(sI_{n-r} - Z)\;=\; \psi^*(s)
\eqno (9.42) $$

The proof is completed.

     Applying Assertion 9.1 we can always find a submatrix
$T$  that guarantees (9.42). Let us consider two cases.

\it{CASE 1.} \rm  $\;l_{\nu -1} =r$. From (9.40a) we get $C_i = C_{\nu}T_{\nu i},\;\; i=1,2,\ldots,\nu -1$. It
implies the following structure of the  matrix $CN^{-1} \;=\;[C_1, C_2,\ldots, C_{\nu }]$
$$ CN^{-1} \;=\;C_{\nu}[T_{\nu 1}, T_{\nu 2}, \ldots,-T_{\nu ,\nu -1}, I_r]
         \;=\;C_{\nu}[-T,\;I_r]
\eqno   (9.43)$$
 In (9.43)  the $r\times r$ submatrix $C_{\nu}$ is chosen in according the condition  $rank
C_{\nu}= r$. Thus, the output matrix $C$  of system (1.1), (9.2) becomes
$$         C \;=\; C_{\nu}[-T,\;I_r]N
\eqno   (9.44)  $$
This  $C$   satisfies also conditions (9.3a,b)
and (9.38).

\it{CASE 2.} \rm $\;l_{\nu -1} <r $. It  follows from the expression (9.40b) that  upper blocks $\bar{T}_{\nu i}
\;=\;$ $-[I_{r-l_{\nu -1}}, O]C_{\nu}^{-1}C_i$ of  submatrices $T_{\nu i}$ can  be arbitrary ones. Uniting
(9.40b) with the last expression we  represent  blocks $C_{\nu}^{-1}C_i$ as follows
$$  C_{\nu}^{-1}C_i \;=\; \left [ \begin{array}{c} \bar{T}_{\nu i} \\
                      T_{\nu i}\end{array} \right ] \;=\; T^*_{\nu i}  $$
Since $C_i \;=\; C_{\nu}T^*_{\nu i}$  then varying  $i$  from  1 to $\nu -1$
$$ [C_1, C_2,\ldots, C_{\nu -1}]\;=\; C_{\nu}[T_{\nu 1}^*, T_{\nu 2}^*,
 \ldots,-T_{\nu ,\nu -1}^*]   $$
and using the expression $ [C_1, C_2,\ldots, C_{\nu}] \;=\; CN^{-1}$ we obtain the matrix $CN^{-1}$ in the form
$$   CN^{-1}\;=\;  C_{\nu}[T_{\nu 1}^*, T_{\nu 2}^*,\ldots,
                -T_{\nu ,\nu -1}^*, I_r] \;=\;  C_{\nu}[-T^*,I_r]
\eqno (9.45) $$
 To calculate  $C$ we use the formula (9.44) with $T\;=\;T^*$.

At the final we summarize the algorithm for zero placement:

     1. Check  controllability of the  pair  $(A,B)$.  If it is completely controllable then go to step 2;
       otherwise the problem has no solution.

     2. Set  $n-r$ desirable distinct real zeros $\bar{s}_i,\;i=1,2,\ldots,
        n-r$, which don't coincide with eigenvalues of $A$.

     3. Define integers $\nu,\;l_1,\; l_2 ,\;\ldots,\; l_{\nu}$ (see
        formulas (1.45),(1.59)) and  the $n\times n$ matrix $N$.

     4.  Find the $l_{\nu -1}\times (n-r)$ submatrix
        $T$ from the condition (9.42); if $l_{\nu -1}< r$ then form the matrix $T^*$.

     5. Construct the submatrix $C_{\nu}$  from  condition  $rankC_{\nu} = r$.

     6. Calculate the matrix $C$ by formulas (9.43), (9.44) (if $l_{\nu -1}
        = r$) or by formulas (9.45), (9.44) (if  $l_{\nu -1}< r$).

\it{REMARK 9.3.} \rm  To satisfy condition (9.42) it is sufficient the only row  of  the submatrix $T$.
Therefore, if $r>1$  then $T$ has $(r-1)(n-r)$ free elements. These elements may be used to fulfil supplementary
requirements, for example, to minimize a performance criterion  $ J = trCC^T$ or to ensure structural
restrictions  on  the  matrix $C$.

     Consider some numerical examples.

                    \it{EXAMPLE 9.3.} \rm

Let's consider a system with completely  accessible state variables and the following matrices $A$ and  $B$
$$   A \; = \;  \left[ \begin{array}{cccc} 2 & 1 & 0 & 0 \\ 0 & 1 & 0 & 1 \\
     0 & 2 & 0 & 0 \\ 1 & 1 & 0 & 0 \end{array} \right ], \qquad
     B\; = \; \left[ \begin{array}{cc} 1 & 0 \\ 0 & 0 \\ 0 & 1 \\ 0 & 1
            \end{array} \right ]
\eqno (9.46)$$
 We assign two desired zeros  $\bar{s}_1 = -1, \; \bar{s}_2 = -2$ ($\psi^*(s) =s^2 + 3s +2$ ) and
 will find an appropriate matrix $C$.

     At first we check  conditions of  Theorem 9.1. We can see that
the  pair (9.46) is completely controllable and eigenvalues of matrix $A$ don't coincide with $\bar{s}_1,
\bar{s}_2 $. Since $rank [B,AB] = 4 $ then we get $\nu = 2, l_1= l_2 = 2$.

    Using   results of Section 1.2.3 we calculate  the
transformation matrix $N$ that reduces   pair  matrices (9.46) to   Asseo's canonical form
$$ N = \left[ \begin{array}{crrc} 0 & 0 & -1 &  0 \\   0 &  1 & 0 &  0 \\
 1 &  -1 &  0 &  0 \\ 0 & 1 &  0 &  1  \end{array} \right ]
\eqno(9.47)$$
 Then, since here  $n-r = 2, \;\;l_{\nu -1} = l_1 = r = 2$ then the upper  block in the matrix $Z$
is absent. Constructing $Z$ as   $Z \;=\; T\;=\; -T_1 \;=\; \left [
\begin{array}{rr} -t_{11} & -t_{12} \\ -t_{21} & -t_{22}
\end{array} \right ]$ we find elements $t_{ij}\;
(i=1,2; \; j=1,2)$ so that  equality (9.42)  be true. For $t_{11} = 1$, $t_{12} = 1$ a polynomial $det(sI_2-Z)$
becomes
$$  det(sI_2-Z)\;=\;  det(sI_2-T)\;=\; det \left [ \begin{array}{cc} s+1 & 1 \\
      t_{21} & s+t_{22} \end{array}\right ] \;=\; s^2 + s(1 + t_{22}) - t_{21}
+ t_{22}  $$
 Comparing the right-hand side of the last
expression with the assigned  zero polynomial $\psi^*(s) =s^2 + 3s +2$  we obtain equations $ 1 + t_{22} = 3, \;
- t_{21} + t_{22} = 2 $ having the solution $t_{21}= 0, \; t_{22} = 2$. Thus, we get
$$  T \;=\; \left [ \begin{array}{rr} -1 & -1 \\ 0 & -2\end{array}\right ]$$
Putting  $C_{\nu} =I_2$  and substituting these  $C_{\nu}$, $T$
and $N$ from (9.47) into (9.44) we obtain the final matrix
$$  C \;=\; \left [ \begin{array}{rrrr} 1 & 0 & -1 & 1\\ 0 & 3 & 0 & 1
\end{array}\right ]
\eqno(9.48)$$
 To check  we form the system matrix $P(s)$  with $A,B$
 (9.46) and $C$ (9.48)  and calculate $detP(s)$
$$ detP(s) \; =\; det\left [ \begin{array}{cccrrr}
         s-2 & -1  & 0 & 0 & -1 & 0 \\ 0 & s-1 & 0 & -1 & 0 & 0 \\
          0 & -2 & s & 0 & 0 & -1 \\-1 & -1 & 0 & s & 0 &-1\\
           1 & 0 & -1 & 1 & 0 & 0 \\ 0 & 3 & 0 & 1 & 0 & 0
          \end{array} \right ] \;=\; -(s^2 +3s +2)  $$
Hence, the zero polynomial coincides with the desirable one.

         \it{EXAMPLE 9.4.} \rm

 We consider the model from Example 9.3 but with the another input matrix
$$    B\; = \; \left[ \begin{array}{cc} 1 & 0 \\ 0 & 0 \\ 0 & 0 \\ 0 & 1
            \end{array} \right ] $$
Let the desired zero polynomial be the same as above one $(\psi^*(s) =s^2 + 3s +2)$. It is  easily verify that
conditions of Theorem 9.1 are also held for given system. But here  $rank [B,AB,A^2 B] = 4$, $\nu = 3$ and $l_1=
l_2 = 1, l_3 =2 $. Thus, $r\nu =6$ and the pair $(A,B)$ is reduced to  Yokoyama's canonical form. Here the
matrix $Z$ of the order  $n-r = 2$  with $l_1 =1,\; l_2 =l_{\nu -1} =1 (< r=2)$ should have  the following
structure
$$  Z \;=\; \left [ \begin{array}{rr} 0 & 1 \\ -t_{11} & -t_{12} \end{array}
\right ]  \qquad \rm{where} \qquad \it [-t_{11} ,\; -t_{12}] \;=\;T  $$
 Substituting $t_{11}$ and $t_{12}$ into
(9.42) yields the equation
$$  det(sI_2-Z)\;=\; det \left [ \begin{array}{cc} s & -1 \\
      t_{11} & s+t_{12} \end{array}\right ] \;=\; s^2 + s t_{12} + t_{11}
\;=\; \psi^*(s)  $$
 from which we  find: $t_{11} = 2$, $t_{12}= 3$. Thus,  $T= [-2 -3]$. Setting $\bar{T} = [-1
-1]$ we form the matrix $T^*$ as follows
$$  T^* \;=\; \left [ \begin{array}{rr} -1 & -1 \\ -2 & -3\end{array}\right ]$$
Putting $C_{\nu} =I_2$ and substituting these $C_{\nu}$ and $T^*$ into (9.44) with $T=T^*$  and  the  matrix
$N$ calculated earlier in Example 1.3 (see Sect.1.2.3) we get
$$ C\;=\; \left [ \begin{array}{rr} 1 & 0 \\ 0 & 1\end{array}\right ]
    \left [ \begin{array}{rrrr} 1 & 1 & 1 & 0\\ 2 & 3 & 0 & 1
      \end{array}\right ] \left[ \begin{array}{cccc}   0 &  0 &  0.5 &  0 \\
     0 &  1 & 0 &  0 \\ 1 &  0 &  0 &  0 \\ 0 & 1 & 0 & 1 \end{array} \right ]
    \;=\; \left [ \begin{array}{cccc} 1 & 1 & 0.5 & 0\\ 0 & 4 & 1 & 1
    \end{array}\right ]
\eqno (9.49)  $$

    For testing we calculate the determinant of the system matrix $P(s)$
with $C$  from (9.49)
$$ detP(s) \; =\; det\left [ \begin{array}{cccrrr}
         s-2 & -1  & 0 & 0 & -1 & 0 \\ 0 & s-1 & 0 & -1 & 0 & 0 \\
          0 & -2 & s & 0 & 0 & 0 \\-1 & -1 & 0 & s & 0 &-1\\
           1 & 1 & 0.5 & 0 & 0 & 0 \\ 0 & 4 & 1 & 1 & 0 & 0
          \end{array} \right ] \;=\; s^2 +3s +2  $$
This completes the verification.

\section[Zero assignment by squaring down operation]
         {Zero assignment by squaring down\\
                            operation}

     Loops in multivariable  feedback  systems
are often introduced between a selected set of accessible variables and an equal number  of independent control
inputs. Thus, the first stage of control design  in a system with $l>r$  contains the stage of combining all
output variables into a new output such that the resulting system has equal number of inputs  and outputs. As it
has been shown in Sect.6.2 this operation ('squaring down') introduces new zeros into the system. The similar
operation may be carried in a system having more inputs then outputs.  We will consider squaring down procedure
only for outputs. But all results may be easily extended for inputs.

     Let controllable and observable system (1.1),  (1.2)  has
more outputs than inputs $(l>r)$ and input and output matrices are of full rank, i.e. $rankB = r$, $rankC = l$.
Combining  output variables by means of a feedforward proportional post-compensator we get a new output $r$
vector $\tilde{y}$
$$    \tilde{y} \;= \;  Dy \;=\; DCx
\eqno   (9.50)   $$
 Let's suppose that system (1.1), (1.2)
possesses  $\mu < n-r$ system zeros: $s_1, s_2,\ldots,s_{\mu}$. As it has been shown above these zeros are not
affected  by any squaring down  operation.  But  this operation introduces new $\eta $ zeros ($ \mu + \eta \le
n-r$) into  squared down system (1.1), (9.50). Therefore, we can consider  the following problem:

\it{PROBLEM 2.} \rm Choose an $r\times l$ constant squaring down matrix $D$ in (9.50) to assign introducing
zeros. The  matrix $D$ must  satisfy in additional the following requirements
$$ (a)\;\; \rm{the} \; \rm{pair} \;\it (A,DC)\; \rm{is}\;\rm{observable} $$
$$ (b)\; \;\;rank DC \;=\; r \qquad \qquad\eqno (9.51)$$
      The  above zero placement problem was first  formulated  and  studied  in  [K3].  We
consider the approach [S5] that is the natural extension of results of Sect. 9.1.1.

     Let us assign distinct real numbers $\bar{s}_1, \bar{s}_2, \ldots,
\bar{s}_{\eta}$ ($\eta = n-r-\mu $)
$$  \bar{s}_i  \neq  \bar{s}_j, \;\;\; \bar{s}_i \neq  \lambda_k, \;\;\;
     i=1,2,\ldots, \eta, \;j=1,2,\ldots, \mu, \; k=1,2,\ldots,n
\eqno (9.52)   $$
where $\lambda_k $ is eigenvalues of the matrix
$A$ and denote by $\bar{\psi}(s)$ the polynomial having  numbers
$\bar{s}_i$ as its zeros
$$    \bar{\psi}(s)\;=\; \prod_{i=1}^{\eta }(s-\bar{s}_i)
\eqno (9.53)  $$

\it{ASSERTION 9.2.}  \rm If the pair $(A,B)$ is controllable, system (1.1), (1.2) has distinct  zeros $s_1,
s_2,\ldots, s_\mu $\footnote{This restriction does not severe  because a system with an unequal number of inputs
and outputs almost always has no zeros.} and  assigned zeros
 $\bar{s}_1, \bar{s}_2, \ldots, \bar{s}_{\eta}$ satisfy  the requirement
(9.52) then  there  is  a  matrix $DC$ that  ensures  both conditions: setting introducing  zeros $\bar{s}_i$
($i=1,2,\ldots,\eta$) and (9.51a).

\it{PROOF}.$\;$ \rm This assertion is the  direct corollary of
Theorem 9.1.

     Now we want to find the matrix $D$ that
satisfies the requirement (9.51b).
     At first we show that (9.51b) is equivalent to
the following one
$$                 rank D = r
\eqno  (9.54)  $$

 \it{ASSERTION 9.3.} \rm  If $r < l \le n$,  $rank
C = l$ and $rank D = r$ then $rank (DC) = r$.

\it{PROOF}. \rm Let's suppose the contrary property : $rank CD < r$. Thus, there is an $r$ vector-row $q^T$
providing the equality $q^TDC = 0$. Denoting  the $l$ vector-row $q^TD$ by $\tilde{q}^T$ we obtain
$$  \tilde{q}^TC \;=\;0
\eqno  (9.55) $$
 It is follows from (9.55)  that $rankC < l$. This
contradicts with the above proposition about fullness of the rank  $C$. This  proves the assertion.

     Now we consider a method for the numerically calculation of the matrix $D$.
The zero polynomial $\psi(s)$  of the  squared down system (1.1), (9.50) are defined from relations (9.9),
(9.10) with  matrices $C_i$ being blocks of the matrix
$$   \bar{C} \;=\; DCN^{-1} \;=\; [C_1, C_2,\ldots, C_\nu ]
\eqno (9.56)    $$
  It follows from inequality (9.52) that $ \psi(\bar{s}_i)\; =
\;\bar{s}_i^{n-r\nu}det\tilde{C}(\bar{s}_i) \; \neq 0$ for all $\bar{s}_i \;\; (i=1,2,\ldots, \eta) $. Let us
form the following criterion
$$    J_1 \;=\; 0.5 \sum _{i=1}^{\eta } \psi(\bar{s}_i)^2
\eqno   (9.57) $$
 that has a minimal value for $\bar{s}_i$ coinciding with zeros of the polynomial $\psi(s)$. The minimization of $J_1$ with
respect to elements $d_{ij}\;$ ($\;i=1,2,\ldots,r,\; j=1,2,\ldots, l$) enables to shift    zeros to
 desirable locations\footnote{It is evident that changing  $D$
 does't affect on zeros $s_1, s_2, \ldots, s_{\mu}$ of an original system  (before the squaring down operation).}.

     To ensure  the condition  (9.54)  we introduce  the following term
$$     J_2 \;=\; (det(DD^T))^{-1}
\eqno (9.58)    $$
 which is the inversion of Gram's determinant
for rows of  the $r\times l$ matrix $D$ ($ r <l $). Since $det(DD^T)
>0$ then the minimization of (9.58) ensures the rank fullness of
$D$ and, by virtue of Assertion 9.3, the rank fullness of $DC$.

     Thus, the minimization of the  sum
$$      J \;=\; J_1 +qJ_2, \;\; q \ge 0
\eqno   (9.59) $$
 where $q$ is a weight coefficient,  enables to find  the  matrix $D$  ensuring desirable
locations to zeros and the condition (9.51).

     We  execute the minimization  in    according   to the
iterative scheme
$$ d_{ij}^{(k+1)}\; =\; d_{ij}^{(k)} \;-\; \alpha \frac{\partial J}
{\partial (d_{ij})}\mid^{(k)}, \qquad i=1,2,\ldots,r,\;
j=1,2,\ldots,l \eqno  (9.60) $$
 where $\alpha >0 $ is a some
constant, $\partial J/\partial (d_{ij})$ is the gradient of $J$ with the respect to elements $d_{ij}$ of the
matrix $D$. For finding an analytic expression for $\partial J/\partial (d_{ij})$ we differentiate the
right-hand sides of (9.57) and (9.58) in a similar way as above in Sect.9.1.1.  We result in
$$ \frac{\partial J_1}{\partial (d_{ij})} \;=\; \sum _{i=1}^{\eta }
\psi(\bar{s}_i)\bar{s}_i^{n-r\nu} tr \{ E^{ij}_{r\times l}C( \sum_{t=1}^{\nu }
      [O,P_i]\bar{s}_i^{t-1})adj(\tilde{C}(\bar{s}_i)) \}
\eqno (9.61)  $$
$$ \frac{\partial J_2}{\partial (d_{ij})} \;=\; -det(DD^T)^{-1}tr
\{ (E^{ij}_{r\times l}D^T + DE^{ji}_{l\times r})(DD^T)^{-1} \} \eqno (9.62) $$ where $E^{ij}_{r\times l}$ is the
$r\times l$ matrix with unit  $ij$-th element and  zeros otherwise, $E^{ji}_{l\times r} = (E^{ij}_{r\times
l})^T$, $n\times l_i$ submatrices $P_i$  $\;(i=1,2,\ldots,\nu )$  defined from the formula (9.20).


           \it{ EXAMPLE 9.5.} \rm

 Let's consider   completely   controllable   and observable  system  (1.1),
 (1.2)  with  $n=3, \; r=2,\; l=3$, matrices  $A$  and
$B$  from  Example 9.1 and the following output
$$  y \;=\; \left [ \begin{array}{rrrr} 1 & 0 & 0 \\ 1  &  1 & 0 \\ 0 & 0 & 1
  \end{array} \right ] x
\eqno(9.63) $$ This  system  has  no  zeros.  The  squaring  down operation  may introduce a undesirable zero.
For example, if
$$  D \;=\; \left [ \begin{array}{rrr} 1 & 0 & 0 \\ 0 & 1 & 0
  \end{array} \right ]      $$
then
$$  DC \;=\; \left [ \begin{array}{rrr} 1 & 0 & 0 \\ 1 & 1 & 0
  \end{array} \right ]      $$
and the  squared down system obtained has the only zero ($0$). This zero may create undesirable difficulties for
control design. We will choose the matrix $D$ to assure the negative zero ($-1$).

     Since   eigenvalues $-0.5, -1.5, -2 $
of $A$ don't coincide with the assigned zero and the original
system has  no   zeros  then Assertion  9.2  is fulfilled and the
zero assignment problem  be to  have  a  solution. Taking into
account that $l_1  = 1, l_2 = 2, \nu = 2, \eta = 1, \bar{s}_1 =
-1$ and  using  formulas (9.59), (9.57), (9.58) we form the
criterion
$$ J \;=\; 0.5(\psi(\bar{s}_1))^2 + q(det(DD^T)^{-1} \;=\;
     0.5\bar{s}_1^{-2}\{det(DC([O,P_1] + P_2]\bar{s}_1)) \}^2 +
     q(det(DD^T)^{-1}
\eqno(9.64) $$
 where $P_1$, $P_2$  are $3\times 1$ and $3\times 2$ blocks of the matrix $N^{-1} =[P_1, P_2]\;$
and the $3\times 3$ matrix $N$ was calculated in  Example 9.1.

The numerical minimization of (9.64) by the recurrent scheme  (9.60) is  finished  as $ \parallel \partial J/
\partial\bar{C} \parallel \le 0.035$. We result in the following matrix
$$ D \;=\; \left [ \begin{array}{rrr}  2.334 & 0.0224 & -0.2653 \\
       0.0224 & 0.5319 & -0.3768 \end{array} \right ]
\eqno(9.65) $$
 which introduces the zero $-1.000$ into the squared down system.

\chapter[Using zeros in analysis and control design  ]
         {Using zeros in analysis and control design }

In  this chapter we consider several control  problems for multivariable
 systems where the notion of
system zeros is useful.

\section[Tracking for constant reference signal.
                           PI-regulator]
          {Tracking for constant reference signal.\\
                           PI-regulator}

From the classic control theory it is known that   steady output tracking for a reference step signal may be
occurred by using an integrator having an error between regulated and reference variables as  an input. We
 study this problem for multivariable systems  with  several  inputs  and outputs. Let's consider a linear
time-invariant multivariable model  of  a dynamical system in the state space
$$ \dot{x} \;=\; Ax + Bu + Ew
\eqno (10.1) $$
$$   z \; = \; Dx
\eqno (10.2)  $$
 where $x \in \mbox{R}^n$ is an state vector, $u \in \mbox{R}^r$ is an input vector, $z \in
\mbox{R}^d$ is an output vector (to be regulated), $w \in \mbox{R}^r$ is a vector of unmeasurable constant
disturbances satisfied to a linear dynamical system
$$ \dot{w}(t) \;=\; 0, \qquad w(t_o) = w_o
\eqno  (10.3)  $$
 with  an unknown  $r$-vector $w_o$. Thus,  vector $w$ is an unmeasurable step function. Matrices
$A$, $B$, $E$, $D$ are the constant ones of appropriate dimensions, $rank B = r$, $rank D =d$. It is assumed
that  the state vector $x$ is completely accessible one.

Let  an  $d$ - vector $z_r$  is a  desirable accessible reference signal described by the  dynamical system
$$ \dot{z}_r(t) \;=\; 0, \qquad z_r(t_o) = z_{ro}
\eqno  (10.4)  $$ with the known  $z_{ro}$.

\it{PROBLEM 1.} \rm  For  plant (10.1), (10.2) it is required to find a feedback regulator $u$  as a function of
$x$, $z$ and  $z_r$ : $u = u(x,z,z_r)$  such that   following output tracking
$$   z(t) \to z_r, \qquad t \to \infty
\eqno    (10.5)  $$
 is executed for all  disturbances  $w$ and  for arbitrary initial conditions.

We  will find the  feedback controller as a proportional-integral (PI) regulator having the error $z-z_r$ as an
input. This  regulator is described in the state-space by equations
$$ \dot{q} \;=\; z-z_r \eqno   (10.6) $$
$$ u \;=\;  K_1x + K_2q
\eqno (10.7)   $$ where $q \in \mbox{R}^d$ is a state vector of the regulator, $K_1$ and $K_2$ are constant
 gain matrices of dimensions $r\times n$  and $r\times d$  respectively.

     Let's note that  PI-regulator (10.6), (10.7) may be represented
in the alternative (classic)  form
$$ u \;=\; K_1x + K_2 \int_{t_o}^{t}(z-z_r)dt + K_2z_o \eqno   (10.8) $$
where  an $d$-vector $z_o$ usually is equal to zero.


     To unite (10.1), (10.2) with the dynamical equation of the regulator
(10.6) we introduce a new state vector $\bar{x}^T =[x^T,\;q^T]$
 and write the augmented differential equation with a
new $(n+d)$ state $\bar{x}$
$$  \dot{\bar{x}} \;=\; \left [ \begin{array}{cc} A & O \\ D & O \end{array}
     \right ]\bar{x} \;+\; \left [ \begin{array}{c} B \\ O
    \end{array} \right ] u \;+\; \left [ \begin{array}{c} E \\ O \end{array}
     \right ] w \;- \;\left [\begin{array}{c}O \\ I_d \end{array} \right ]z_r
\eqno(10.8)  $$
The feedback control (10.7) is rewritten  as follows
$$          u \;=\; [ K_1,\;K_2]\bar{x}
\eqno(10.9)     $$ and present  the  linear proportional state feedback introduced into the  open-loop system
(10.8).
     Problem 1 is reformulated in terms of the augmented system  as:
it is necessary to find matrices $K_1$ and $K_2$ of the proportional feedback regulator (10.9) that ensure
asymptotic stability to the following closed-loop system
$$  \dot{\bar{x}} = \left [ \begin{array}{cc} A +BK_1& BK_2 \\ D & O \end{array}
 \right ]\bar{x} \;+\;  \left [\begin{array}{c} E \\ O \end{array} \right ] w -
\left [\begin{array}{c}  O \\ I_d \end{array} \right ]z_r \eqno (10.10)  $$
 i.e. the dynamics matrix of (10.10)
must satisfy the following condition
$$  Re\lambda_i \left ( \begin{array}{cc} A +BK_1& BK_2 \\ D & O\end{array}
  \right ) < 0, \;\, i=1,2,\ldots,n+d
 \eqno (10.11) $$
where  $\lambda_i<.>$  is an eigenvalue of an matrix.

     Let us show that solvability of this problem guarantees simultaneously solvability  of
Problem 1. To this purpose we  study  asymptotic behavior of the vector $\dot{q} =z-z_r$. Differentiating the
both sides of (10.10) with respect to $t$ and denoting $ \tilde{x}^T = \dot{\bar{x}}^T =[\dot{x}^T,\;\dot{q}^T]$
we get the following linear homogeneous differential equation in $\tilde{x}$
$$  \dot{\tilde{x}} = \left [ \begin{array}{cc} A +BK_1& BK_2 \\ D & O
     \end{array} \right ]\tilde{x}
\eqno (10.12) $$
 If  the condition (10.11) is satisfied then (it follows from (10.12))   $\;\;\tilde{x} \to 0$  as $t
\to \infty$  or $\dot{x} \to 0$, $\dot{q} \to 0$ as $t \to \infty$. Consequently, if we have been found matrices
$K_1$ and $K_2$ that  ensure the condition (10.11)  then  PI-regulator (10.6), (10.7) with same gain matrices
$K_1$ and $K_2$ ensures  asymptotic steady output tracking in system (10.1), (10.2).
     Solvability  conditions of  Problem 1 are same as  the state feedback  problem, i.e.  a solution
exists  if and only if the pair of  matrices
$$   \tilde{A} \;=\;\left [ \begin{array}{cc} A & O \\ D & O \end{array}
        \right ], \qquad  \tilde{B} \;=\; \left [ \begin{array}{c} B \\ O
       \end{array} \right ]
\eqno (10.13)  $$ is stabilizable [W3]. Let us express stabilizability of matrices $\tilde{A}$, $\tilde{B}$ via
matrices $A$, $B$ and $D$.

\it{ASSERTION 10.1.} \rm  The pair ($\tilde{A}, \tilde{B}$) is stabilizable if and only if the following
conditions take place

     a. the pair $(A,B)$ is stabilizable,

     b.  $d \le r $ ,

     c. the system
$$       \dot{x} = Ax + Bu \;\; y = Dx
\eqno (10.14) $$ has no system zeros in the origin.

\it{PROOF}.$\;$ \rm It follows from the stabilizability criterion  [W3, Theorem 2.3]  that the pair of matrices
$(\bar{A},\bar{B})$  is stabilizable if and  only if $rank[ \lambda I_n - \bar{A}, \bar{B} ] = n$ where
$\lambda$ is an unstable eigenvalue of  the $n\times n$  matrix $\bar{A}$.

     Let's denote unstable eigenvalues of the matrix $\tilde{A}$ (10.13) by
     $\tilde{\lambda}_i^*$ $(i=1,2,\ldots, \mu; \;\; \mu \le n+d)$. One can  see  that  the set of
$\tilde{\lambda}_i^*$, $i=1,2,\ldots, \mu $ consists of unstable eigenvalues $\lambda_i^*$ $(i=1,2,\ldots, \eta;
\;\;\eta \le n)$  of  the  matrix $A$ and  $d$ eigenvalues that are equal to zero. Therefore, the
stabilizability criterion for  the pair (10.13) may be formulated as follows: the pair of matrices ($\tilde{A},
\tilde{B}$) is stabilizable if and only if the  $(n+d)\times (n+d+r)$ matrix
$$             Q(\lambda )\;=\; (\lambda I_{n+d} - \tilde{A}, \tilde{B}) \;
      =\; \left [ \begin{array}{ccc} \lambda I_n-A & O & B\\
           -D & \lambda I_d & O \end{array} \right ]
\eqno   (10.15)  $$ has the full rank  $n+d$  for $\lambda = \lambda_i^*$ $(i=1,2,\ldots, \eta ) $ and  $\lambda
=0$.

     Further we  separate two cases.

     1. $\lambda = \lambda_i^* \neq 0$.$\;$ Using  equivalent block
operations we write series of  rank equalities
$$ rankQ(\lambda_i^*)\;=\; rank\left [ \begin{array}{ccc}
        \lambda_I^* I_n-A & O & B\\-D & \lambda_i^* I_d & O \end{array}
      \right ] \;=\; rank\left [ \begin{array}{ccc} \lambda_I^* I_n-A & B & O\\
            -D & O & \lambda_i^* I_d \end{array} \right ] \;=\; $$
$$      rank\left [ \begin{array}{ccc} \lambda_I^* I_n-A & B & O\\
           O & O & \lambda_i^* I_d \end{array} \right ] \;=\;
           d + rank[\lambda_i^* I_n -A, B]
\eqno   (10.16)  $$
 Analysis of the right-hand side of (10.16)  reveals  that the  matrix $Q(\lambda_i^*)$ has
 the full rank $n+d$ if and only if
$$       rank[\lambda_i^* I_n -A, B] = n $$
The above rank condition is fulfilled if and only if the pair of  matrices $(A,B)$ is stabilizable.  Hence, we
prove the  condition (10.14a) of the assertion.

     2. $\lambda =0$. $\;$ In this case
$$ rankQ(0)\;=\; rank\left [ \begin{array}{ccc} -A & O & B\\
           -D & O & O \end{array}\right ] \;=\; rank\left [ \begin{array}{rr}
            -A & -B \\ D & O \end{array} \right ]
\eqno   (10.17)   $$
 Consequently, the matrix $Q(0)$ coincides with  the system matrix  $P(s)$ at  $s=0$ for the  system
$\dot{x} =Ax + Bu, \; y=Dx $ and the requirement of the rank fullness of $Q(0)$ is  equivalent to absence of
 system zeros in origin. Moreover,  $rank Q(0)= n+d$  if  and  only  if $d \le r$. Therefore,
we  validate   conditions (10.14 b,c) of the assertion. The proof is completed.

     Thus, we have shown that the  problem  of  asymptotic
steady-output tracking with simultaneously rejecting constant disturbances  has a solution if, apart from other
conditions, open-loop system zeros satisfy the requirement on zero locations.

\section[Using state estimator in PI-regulator]
         {Using state estimator in PI-regulator}

     Let  the  state vector $x$ in  system  (10.1)  does  not
completely accessible: only $l$  state variables form the $l$ output vector
$$         y \;=\; Hx +Fw
\eqno(10.18) $$
 which is accessible for a measurement. In the general case the regulated output $z$ (10.2)
distinguishes  from the measurable vector $y$, for example, the vector $z$ may be a part of $y$. In (10.18) $w
\in \mbox{R}^p$ is  the vector of unmeasurable constant disturbances describing by  equation (10.3) with an
unknown initial state $w_{o1} \neq w_o$; the matrix $F$ is a constant $l\times p$ matrix.

 \it{PROBLEM} \rm 2. For plant (1.1), (10.18), (10.2) it is
required to find a feedback regulator $u$  as a  function of  $\hat{x}$,  $z$  and  $z_r$ : $u
\;=\;u(\hat{x},z,z_r)$ such that steady output tracking (10.5) takes place for any constant disturbances $w$ and
for arbitrary initial conditions.

According to the approach [S12] we apply a feedback PI-regulator of the  following structure\footnote{ This
structure is similar to  one considered in [K5, p.477].}
$$        \dot{q} \;=\; z - z_r
\eqno(10.19)$$
$$            u\;=\; K_1 \hat{x} + K_2q
\eqno(10.20)  $$
where $q \in \mbox{R}^d$ is a state of the dynamic regulator, $K_1$, $K_2$ are constant
$r\times n$ and  $r\times d$   matrices  respectively, $\hat{x} \in \mbox{R}^n$ is an estimate  of vector $x$.

     To find $\hat{x} $  we use a full-order  state
observer [O1]. Apart from the vector $x$ we need also to estimate simultaneously the  disturbance vector $w$.

     Since the disturbance model coincides with the dynamical system (10.3) then
introducing  $(n+p)$  an state vector $[x^T, w^T]$ gives the following augmented differential equation
$$  \left [ \begin{array}{c} \dot{x} \\ \dot{w} \end{array} \right ] \; =
\;  \left [ \begin{array}{cc} A & E \\ O & O \end{array} \right ] \left [ \begin{array}{c} x \\ w \end{array}
\right ]  +  \left [
\begin{array}{c} B \\ O \end{array} \right ] u,
\eqno(10.21) $$
$$ y = \left[ \begin{array}{cc}  H & F \end{array} \right ] \left [
\begin{array}{c} x \\ w \end{array} \right ]
\eqno (10.22)  $$
The full order state observer that estimates  $[x^T, w^T]$
 has the following structure [O1]
$$  \left [ \begin{array}{c} \dot{\hat{x}} \\ \dot{\hat{w}}
     \end{array} \right ] \; = \; \left [ \begin{array}{cc} A & E \\ O & O
      \end{array} \right ] \left [ \begin{array}{c} \hat{x} \\ \hat{w}
      \end{array} \right ]  +  \left [ \begin{array}{c} B \\ O \end{array}
     \right ] u \;-\; L[ H, F ]\epsilon
\eqno(10.23) $$ where $L$ is an $(n+p)\times l$ constant matrix, $\epsilon ^T = (x^T, w^T) - (\hat{x}^T,
\hat{w}^T)$ is the  $(n+p)$ vector of the error. It easy to show based  on differential equations (10.23),
(10.21) that  the vector $\epsilon = \epsilon (t)$ satisfies to the following linear homogeneous differential
equation
$$      \dot{\epsilon} \;=\; \left ( \left [ \begin{array}{cc} A & E \\ O & O
        \end{array} \right ]  \;-\; L[ H, F ] \right )\epsilon
\eqno   (10.24)  $$
 and $\epsilon \to 0$ at  $t \to \infty$ if  the dynamic matrix of
(10.24) satisfies to the following condition
$$  Re\lambda_i \left ( \left [ \begin{array}{cc} A & E \\ O & O\end{array}
  \right ] - L[ H, F ] \right ) < 0, \;\, i=1,2,\ldots,n+p
 \eqno (10.25) $$
Thus, we need to find a constant matrix $L$ that guarantees the condition (10.25). This problem has a solution
if and only if the pair of  matrices
$$   \left ( \left [ \begin{array}{cc} A & E \\ O & O \end{array} \right ]^T,\; \left [ \begin{array}{c} H^T \\ F^T
       \end{array} \right ] \right )
\eqno (10.26)  $$ is stabilizable. The stabilizability of this pair is guaranteed by the following assertion.

\it{ASSERTION 10.2.} \rm  The  pair  of   matrices   (10.26)   is stabilizable if and only if

     a. the pair $(A^T,H^T)$ is stabilizable,

     b.  $l \ge p $ ,

     c. the system
$$       \dot{x} = Ax + Ew, \;\; y = Hx+Fw  $$
has no system zeros in origin.

     To prove we can use a similar way as in Assertion 10.1.

Thus, if  conditions of Assertion  10.2  are fulfilled then there always exists such a matrix $L$ that the
asymptotically  exact reconstruction of vectors $x$ and $w$ takes place
$$    \hat{x} \to x, \;  \hat{w} \to w \;\; \rm{as} \;\; t \to \infty
\eqno(10.27) $$

Then we find conditions when PI-regulator (10.20) exists. We unite equations (10.1), (10.19),(10.20) and (10.24)
by introducing a state vector $[x^T, q^T, \epsilon ^T]$ and express the vector $\hat{x}$ via $x$  and $\epsilon$
$$        \hat{x} \;=\; x - [I_n, O]\epsilon
\eqno  (10.28)  $$ The differential equation of the closed-loop system: object + PI-regulator + observer becomes
$$  \left [ \begin{array}{c} \dot{x} \\ \dot{w} \\ \dot{\epsilon} \end{array}
 \right ] \; = \; \left [ \begin{array}{cccc} A +BK_1 & BK_2 & \vdots &
             -B[K_1,O] \\ D & O & \vdots & O \\ \dotfill & \dotfill &
              \vdots & \dotfill \\ O & O & \vdots &
             \left [ \begin{array}{cc}A & E \\ O & O\end{array} \right ]
             - L\left [ \begin{array}{cc} H & F \end{array} \right ]
            \end{array} \right ] \left [ \begin{array}{c} x \\ w \\ \epsilon
        \end{array} \right ] \; +\;  \left [ \begin{array}{c} E \\ O \\ O \\ O
        \end{array} \right ] w \;-\; \left [ \begin{array}{c} O \\ I_d \\ O \\ O
        \end{array} \right ] z_r
\eqno    (10.29) $$
If we find  matrices $K_1$ and $K_2$  such as the following condition takes place
$$  Re\lambda_i \left ( \left [ \begin{array}{cc} A & O \\ D & O\end{array}
  \right ] \;+\; \left [ \begin{array}{cc} B \\ O \end{array} \right ]
       [ K_1, K_2 ] \right ) < 0, \;\, i=1,2,\ldots,n+d
 \eqno (10.30) $$
 then the problem has a solution.

 Now we  show that
 simultaneously fulfilment of conditions (10.30) and (10.25) guarantees the solution of  Problem 2. For
this purpose we investigate asymptotic behavior of the vector $\dot{q} = z - z_r$. Since  $\dot{w} =0,\;
\dot{z}_r =0$ then differentiating  both sides of (10.29) with respect to $t$ yields  the following linear
homogeneous differential equation in $\tilde{x}^T = [\dot{x}^T, \dot{q}^T, \dot{\epsilon}^T]$
$$  \dot{\tilde{x}} \;=\;\left [ \begin{array}{cccc} A +BK_1 & BK_2 & \vdots &
                    -B[K_1,O] \\ D & O & \vdots & O \\ \dotfill & \dotfill &
                     \vdots & \dotfill \\ O & O & \vdots &
                  \left [ \begin{array}{cc}A & E \\ O & O\end{array} \right ]
                    - L\left [ \begin{array}{cc} H & F \end{array} \right ]
                   \end{array} \right ] \tilde{x} $$
which is asymptotic stable ( $\tilde{x} \to 0$ as $t \to \infty$  or $\dot{q} \to 0$ as $t \to \infty$) if
conditions (10.30), (10.25) are held. It is evident that  matrices $K_1$ and $K_2$ exist to assure the condition
(10.30) if and only if  the pair of  matrices $ \left ( \left [
\begin{array}{cc}
 A & O \\ D & O\end{array} \right ],\; \left [ \begin{array}{cc} B \\ O
\end{array} \right ] \right )$ is stabilizable. The stabilizability
of this  pair  guarantees by Assertion 10.1 (see Sect.10.1).

     Uniting Assertions 10.1 and  10.2  we  obtain
total solvability  conditions of  Problem 2.

\it{THEOREM 10.1.} \rm Necessary and sufficient conditions for existence  of PI-regulator (10.19), (10.20) for
system (10.1), (10.2), (10.18), which ensures that $z \to z_r$  when $t \to \infty$ for all constant
unmeasurable disturbances $w$  (10.3) and for all constant reference signals $z_r$  (10.4) are follows

     a. the pair $(A,B)$ is to be stabilizable,

     b. the pair $(A,H)$ is to be detectable,

     c.  $r  \ge d, \;\;   l \ge p$,

     d.  system zeros of the following systems:
$$       \dot{x} = Ax + Ew, \;\; y = Hx + Fw
\eqno (10.31) $$ and
$$    \dot{x} = Ax + Bu, \;\; z = Dx
\eqno    (10.32) $$ are not equal to zero.


     Let's consider the general case of  the regulated output
$$                  z \;=\; Dx + Qw
\eqno(10.33)$$
 where $w \in \mbox{R}^p$ is a vector of  unmeasurable  constant disturbances satisfying the  differential
equation  (10.3), $Q$  is an  $d\times p$ matrix. We will study asymptotic steady tracking of $\hat{z}
=D\hat{x}$ for the reference signal $z_r(t)$
$$      \hat{z}(t) \;\to \; z_r(t), \qquad t \to \infty
\eqno   (10.34)$$
     In this case  PI-regulator becomes
$$            \dot{q} \;=\; D\hat{x} - z_r
\eqno   (10.35)  $$
 To deduce a differential equation for the closed-loop system: object + PI-regulator
+ observer we  express the  estimate $\hat{x}$ via $x$ and $\epsilon$ (see (10.28)) and substitute the result
into (10.35). We get
$$            \dot{q} \;=\; Dx - [D,\;O]\epsilon -z_r
\eqno   (10.36)  $$
 Uniting  equations  (10.1),  (10.36), (10.20)  and  (10.24)  by introducing a new state
vector $[x^T, q^T, \epsilon ^T]$  we get the  differential  equation of the closed-loop system
$$  \left [ \begin{array}{c} \dot{x} \\ \dot{w} \\ \dot{\epsilon} \end{array}
 \right ] \; = \; \left [ \begin{array}{cccc} A +BK_1 & BK_2 & \vdots &
             -B[K_1,O] \\ D & O & \vdots & -[D,\;O] \\ \dotfill & \dotfill &
              \vdots & \dotfill \\ O & O & \vdots &
             \left [ \begin{array}{cc}A & E \\ O & O\end{array} \right ]
             - L\left [ \begin{array}{cc} H & F \end{array} \right ]
            \end{array} \right ] \left [ \begin{array}{c} x \\ w \\ \epsilon
        \end{array} \right ] \; +\;  \left [ \begin{array}{c} E \\ O \\ O \\ O
        \end{array} \right ] w \;-\; \left [ \begin{array}{c} O \\ I_d \\ O \\ O
        \end{array} \right ] z_r
\eqno    (10.37) $$
 Dynamic  behavior  of the state vector is defined now by the dynamics matrix
of Eqn.(10.37) that has diagonal blocks coinciding with diagonal blocks of the dynamics matrix of Eqn.(10.29).
Thus,  solvability conditions are formulated here by Theorem 10.1.

                     \it{EXAMPLE 10.1.} \rm

 To illustrate the main results of this section
we consider control of an aerial antenna position. The model is [K5, Example 2.4]
$$   \dot{x} \;=\; \left [ \begin{array}{cr} 0 & 1 \\ 0 & -4.6 \end{array}
    \right ] x \;+\;\left [ \begin{array}{c} 0 \\ 0.787 \end{array} \right ] \mu
     \;+\; \left [ \begin{array}{c} 0 \\ 0.1 \end{array} \right ]\tau _o
\eqno  (10.38)  $$
 where  state variables of $x^T =[x_1,x_2]$  have the following sense: $x_1$ is a aerial antenna
angle, $x_2 = \dot{x}_1$ is an aerial antenna angle  speed. Here $\mu$ is  a  control  variable and $\tau _o$ is
a constant disturbance. In this system $z=x_1+x_2$ is the regulated output. It is proposed that the measurable
output is   $x_1$, i.e.  $y= x_1$.

We will design PI-regulator to  assure asymptotic tracking the  regulated output $z=x_1+x_2$ for the preassigned
value $z_r$. The last is a constant value during a long time interval but it may change unevenly in some
moments.

     Here we have  the system (10.1),(10.2), (10.18) with $n = 2, \;\;r = l  = d = p = 1$ and
$$   A \;=\; \left [ \begin{array}{cr} 0 & 1 \\ 0 & -4.6 \end{array}
\right ] , \;\; B\;=\;\left [ \begin{array}{c} 0 \\ 0.787 \end{array} \right ],
  \;\; E\;=\;\left [ \begin{array}{c} 0 \\ 0.1 \end{array} \right ], \;\;
     D\;=\; [1\;1], \;\; H\;=\;[1\;0],\;
     \; F=0 \eqno  (10.39)  $$
It is evident that the disturbance variable $\tau _o$ and the reference signal $z_r$   may be described by
differential equations
$$ \begin{array}{cc} \dot{\tau}_o \;=\; 0, &  \tau _o(0) = \bar{\tau}_o \\
  \dot{z}_r(t) \;=\; 0, &  z_r(t_o) = z_{ro} \end{array}
\eqno   (10.40) $$
with the known $z_{ro}$ and some an unknown $ \bar{\tau}_o $.

     Let us check  conditions  (a)- (d) of  Theorem 10.1. Conditions  (a) and (b) are fulfilled
$$ rank[B,AB] \;=\; rank \left [ \begin{array}{cc} 0 & 0.787\\ 0.787 & 3.62
    \end{array} \right ] \;=\;2 , \qquad rank[H^T,A^TH^T] \;=\;
    rank \left [ \begin{array}{cc} 1 & 0\\ 0 & 1  \end{array} \right ] \;=\;2
$$
The performance of  the condition  (c) is evident.  Testing  the   condition  (d) gives
$$ rank\left [ \begin{array}{rr} -A & -E\\ H & F  \end{array} \right ] \;=\;
       rank \left [ \begin{array}{crr} 0 & -1 & 0\\ 0 & 4.6 & -0.1 \\
        1 & 0 & 0  \end{array} \right ] \;=\; 3 $$
$$ rank\left [ \begin{array}{rr} -A & -B\\ D & O  \end{array} \right ] \;=\;
       rank \left [ \begin{array}{crr} 0 & -1 & 0\\ 0 & 4.6 & -0.787 \\
        1 & 1 & 0  \end{array} \right ] \;=\; 3 $$
Therefore, PI-regulator  (10.35),(10.20) exists\footnote{ Note, if   $y=x_2$, i.e. $H=[0\;\;1]$ then the
condition (10.31) is not fulfilled.}    and it has the following form
$$ \begin{array}{ccl} \dot{q}&  =&  \hat{x}_1 +x_2 - z_r,\\
    \mu & = & k_1\hat{x}_1 + k_2\hat{x}_2  + k_3q \end{array}
\eqno     (10.41) $$
where $k_1$, $k_2$, $k_3 $ are constant feedback gains, $\hat{x}_1$, $\hat{x}_2$   are
estimates  of   variables $x_1$ and $x_2$, which  are defined from the formulas
$$   \hat{x}_1  = x_1 -\epsilon _1, \qquad  \hat{x}_2 = x_2 -\epsilon _2 $$
Here $\epsilon _1$, $\epsilon _2$ are first variables of the error vector $\epsilon ^T =[\epsilon _1, \epsilon
_2, \epsilon _3]$   where
$$  \epsilon _1 = x_1- \hat{x}_1, \;\;  \epsilon _2 = x_2 -\hat{x}_2, \;\;
     \epsilon _3 = \tau_o - \hat{\tau}_o
\eqno  (10.42)     $$
     Substituting  concrete matrices $A$, $E$, $H$, $F$ in
(10.24) we obtain the differential equation for $\epsilon$
$$      \dot{\epsilon} \;=\; \left ( \left [ \begin{array}{crc} 0 & 1 & 0 \\
          0 & -4.6 & 0.1 \\ 0 & 0 & 0 \end{array} \right ]  \;-\; L[ 1\; 0\; 0 ]
         \right )\epsilon
\eqno   (10.43)  $$ where a row  vector  $L^T = [ l_1,\; l_2,\; l_3 ]$ should be chosen to assure asymptotic
stability of the system (10.43). Assigning observer poles as $( -1+j,\; -1-j,\; -2)$ and using the method of
modal control [P4] yields
$$ L^T = [-0.6,\;  8.76,\;  40 ]$$
Substituting this $L$ into (10.43) gives the following differential equation in $\epsilon _1, \epsilon _2,
\epsilon _3$
$$ \begin{array}{ccl} \dot{\epsilon} _1 & = & 0.6\epsilon _1 + \epsilon _2 \\
  \dot{\epsilon} _2 & = & -8.76\epsilon _1 - 4.6\epsilon _2 - 0.1\epsilon _3\\
      \dot{\epsilon} _3 & = & -40\epsilon _1 \end{array}
\eqno(10.44)     $$
 To calculate feedback gains $k_1$, $k_2$, $k_3$ we also  use the modal control method.
Assigning poles of the closed-loop system ((10.38),(10.41))  as $(-0.5, -1, -1.5)$ we find
$$         (k_1, k_2, k_3) \;=\; ( -1.906,  2.033,  -1.080 )$$
Thus,  PI-regulator (10.41) becomes
$$ \begin{array}{ccl} \dot{q}&  =&  \hat{x}_1+\hat{x}_2 -z_r,\\
    \mu & = & -1.906\hat{x}_1 + 2.033\hat{x}_2  -1.08q \end{array}
\eqno     (10.45) $$
 where $\hat{x}_1$ and $\hat{x}_2$  are defined from  (10.23) with
concrete matrices (10.39), $\hat{w} = \hat{\tau} _o$ and  $ u = \mu $ from (10.45).
In work [S21] is demonstrated the response of output $z=x_1+x_2$ of the closed-loop system for  $\tau _o =10,\;
\hat{x}_1(0) = \hat{x}_2(0) =0,\; \hat{\tau} _o(0) = q(0) =0,\; x_1(0) = x_2(0) =1$ and the reference signal
$$  z_r\;=\; \left \{ \begin{array} {ccc} 1 & , & 0 \le t \le 40\\
                3 & , & 40 < t \le 80 \end{array} \right.$$

\section[Tracking for polynomial reference signal]
        {Tracking for polynomial reference signal}

     Now we consider the general tracking  problem that  has similar solvability
conditions as in Problem 1. We will study tracking for a polynomial reference signal of the form
$$ z_{ref}(t) = \alpha _o + \alpha _1t + \alpha _2t^2 +\cdots+
    \alpha _{\eta -1}t^{\eta -1}
\eqno(10.47)    $$ where $ \alpha _o , \; \alpha _1. \; \ldots, \alpha _{\eta -1}$ are known $d$-vectors.

     Let us consider the following linear time-invariant
multivariable state-space model \footnote{For simplicity  we consider the model without disturbances.}
$$ \dot{x} \;=\; Ax + Bu
\eqno (10.48) $$
$$   z \; = \; Dx
\eqno (10.49)  $$ where  vectors $x \in \mbox{R}^n$,  $u \in \mbox{R}^r$ and $z \in \mbox{R}^d$ have the  same
sense as in Sect.10.1, $A$, $B$, $D$ are  $n\times n$, $n\times r$  and  $d\times n$ constant matrices
respectively. It is assumed that  the state vector $x$  is completely accessible.

\it{PROBLEM 3.} \rm  For plant (10.48), (10.49) it is required to find a feedback dynamic regulator $u = u(x, z,
z_{ref})$ such as following output tracking
$$   z(t) \to z_{ref}, \qquad t \to \infty
\eqno   (10.50)   $$ takes place in the closed-loop system.

For seeking a solution we  use the approach proposed by Porter, Bradshow [P5] where the following feedback
dynamic  regulator is used
$$ \begin{array}{ccl} \dot{q}_1 & = & z - z_{ref} \\
                      \dot{q}_2 & = &q_1\\
                        & \vdots & \\
                       \dot{q}_{\eta} & = &q_{\eta -1}\end{array}
\eqno(10.51)$$
$$             u \; = \; Kx + [K_1, \; K_2,\; \ldots,K_{\eta}] \bar{q}
\eqno(10.52)   $$
 In (10.51), (10.52)  $q_i, i=1,2,\ldots,\eta $ are $d$-vectors, which form  the state vector
$\bar{q}^T = [ \bar{q}_1^T, \bar{q}_2^T, \ldots, \bar{q}_{\eta}^T]$ of the dynamic regulator, $\bar{q} \in
\mbox{R}^{\eta d}$, $K$ is an $r\times n$ constant matrix and $K_i, i=1,2,\ldots, \eta $ are  $r\times d$
 constant matrices.

    Uniting  equations (10.48), (10.49), (10.51) we write the
augmented differential equation with respect to a new $n + d\eta $ vector $\bar{q}^T = [ x^T, \bar{q}_1^T,
\bar{q}_2^T, \ldots, \bar{q}_{\eta}^T]$
$$  \left [ \begin{array}{c} \dot{x} \\ \dot{q}_1 \\ \dot{q}_2\\ \vdots \\
  \dot{q}_{\eta} \end{array} \right ] \; = \; \left [ \begin{array}{cccccc}
        A  & O & O & \cdots & O & O \\ D & O & O & \cdots & O & O \\
        O  & I_d & O & \cdots & O & O \\ \vdots& \vdots & \vdots & \cdots &
              \vdots & \vdots   \\  O  & O & O & \cdots & I_d & O
              \end{array} \right ] \left [ \begin{array}{c} x \\ q_1 \\ q_2\\
         \vdots \\ q_{\eta} \end{array} \right ] \; +\;
 \left [ \begin{array}{c} B \\ O \\ O \\ \vdots \\ O \end{array} \right ] u
      \;-\; \left [ \begin{array}{c} O \\ I_d \\ O \\ \vdots \\ O
        \end{array} \right ] z_{ref}
\eqno  (10.53) $$
 The feedback dynamic regulator (10.52) is rewritten  for  the composite system (10.53) as follows
$$ u \;=\; [ K_1, \;K_2,\;\ldots, K_{\nu} ] \left [ \begin{array}{c} x \\
   \bar{q} \end{array} \right ], \qquad \bar{q}^T \;=\;
   [ \bar{q}_1^T, \bar{q}_2^T, \ldots, \bar{q}_{\eta}^T]
\eqno  (10.54) $$
 This regulator is, in fact, a linear proportional state feedback. Therefore,  Problem 3 is
reduced to  design of a proportional state feedback that stabilizes the closed-loop system
$$  \left [ \begin{array}{c} \dot{x} \\ \dot{q}_1 \\ \dot{q}_2\\ \vdots \\
  \dot{q}_{\eta} \end{array} \right ] \; = \; \left [ \begin{array}{cccccc}
     A+BK  & BK_1 & BK_2 & \cdots & BK_{\eta -1} & K_{\eta } \\
     D & O & O & \cdots & O & O \\ O  & I_d & O & \cdots & O & O \\
      \vdots& \vdots & \vdots & \cdots & \vdots & \vdots   \\
       O  & O & O & \cdots & I_d & O  \end{array} \right ] \left [
\begin{array}{c} x \\ q_1 \\ q_2\\ \vdots \\ q_{\eta} \end{array} \right ]
       \;-\; \left [ \begin{array}{c} O \\ I_d \\ O \\ \vdots \\ O
        \end{array} \right ] z_{ref}
\eqno  (10.55) $$

At first we  show  that  this problem has a solution if  the above Problem 3 has a solution. We differentiate
$\eta $ times  both sides of (10.55) with respect to $t$ and denoting  $\tilde{x}^T = [x^{(\eta)T},
\bar{q}^{(\eta)T}]$  get  the  following  linear homogeneous  differential equation in the vector $\tilde{x}$
$$  \dot{\tilde{x}}\;=\;  \left [ \begin{array}{cccccc}
              A+BK  & BK_1 & BK_2 & \cdots & BK_{\eta -1} & K_{\eta } \\
              D & O & O & \cdots & O & O \\ O  & I_d & O & \cdots & O & O \\
                 \vdots& \vdots & \vdots & \cdots & \vdots & \vdots   \\
              O  & O & O & \cdots & I_d & O  \end{array} \right ]\tilde{x}
\eqno(10.56)    $$
 If (10.55) is asymptotic stable then its the  dynamics matrix has all eigenvalues with negative
real parts and we have from (10.56): $\tilde{x} \to 0$ as $t \to \infty$ or $x^{(\eta)} \to 0, \; q^{(\eta)}_i
\to 0$, $i=1,2,\ldots,\eta$ as $t \to \infty$. In according with  (10.51)  $q^{(\eta)}_{\eta } = q^{(\eta
-1)}_{\eta -1}= \ldots = q^2_2 = \dot{q}_1 = z-z_{ref}$. Consequently, $q^{(\eta)}_{\eta } \to 0$  as  $t \to
\infty $ i.e.  $z-z_{ref} \to 0$ as $t \to \infty $ and  solvability conditions of Problem 3 are equivalent to
solvability conditions of a proportional state feedback, i.e. stabilizability  of matrices
$$\hat{A} \; = \; \left [ \begin{array}{cccccc}
        A  & O & O & \cdots & O & O \\ D & O & O & \cdots & O & O \\
        O  & I_d & O & \cdots & O & O \\ \vdots& \vdots & \vdots & \cdots &
              \vdots & \vdots   \\  O  & O & O & \cdots & I_d & O
              \end{array} \right ], \qquad  \hat{B}\;=\;
         \left [ \begin{array}{c} B \\ O \\ O \\ \vdots \\ O \end{array}
         \right ]
\eqno  (10.57)$$

     Now we show that if the matrices $A$, $B$ and $D$  (10.57)
satisfy  conditions of Assertion 10.1 then  the pair $(\hat{A},\hat{B})$ is stabilizable. Indeed, using
reasonings of the proof  of  Assertion 10.1 we should analyze a rank of the  $(n+ \eta d)\times (n +\eta d + r$)
matrix
$$  Q(\lambda) \;=\; [ \lambda I - \hat{A}, \; \hat{B}]    $$
at $\lambda = \lambda^*_i$ and  $\lambda =0 $  where $\lambda^*_i, \;i=1,2,\ldots, m \;(m\le n)$  are  unstable
eigenvalues of $A$. We consider these two cases separately.

     1. $\lambda = \lambda_i^* \neq 0$.$\;$ Using  equivalent block
        operations we get
$$ rankQ(\lambda_i^*)\;=\; rank \left [ \begin{array}{ccccccc}
            \lambda_I^* I_n-A & O & O & \cdots & O & O & B\\
          - D & \lambda_i^* I_d & O & \cdots & O & O & O \\
           O & -I_d & \lambda_i^* I_d & \cdots & O & O & O \\
           \vdots & \vdots & \vdots & \cdots & \vdots &\vdots &\vdots \\
           O & O  & O & \cdots & -I_d & \lambda_i^* I_d & O
             \end{array} \right ] \;=\;
                   d\eta + rank[\lambda_i^* I_n -A, B] $$
    Therefore, the matrix $Q(\lambda_i^*)$ has the full rank if and only if
the condition  (a) of Assertion 10.1  is held.

     2. $\lambda =0$. In this case
 $$ rankQ(0)\;=\; rank\left [ \begin{array}{ccccccc}
             -A & O & O & \cdots & O & O & B\\
           - D & O & O & \cdots & O & O & O \\
            O & -I_d & O & \cdots & O & O & O \\
            \vdots & \vdots & \vdots & \cdots & \vdots &\vdots &\vdots \\
            O & O  & O & \cdots & -I_d & O & O
         \end{array} \right ] \;=\; d(\eta -1) +
         rank \left [ \begin{array}{rr}-A & B \\ -D & O \end{array} \right ]
         \;=\; $$
$$        =\; d(\eta -1) + rank\left [ \begin{array}{rr}-A & -B \\ D & O
        \end{array} \right ]  $$
 Consequently, the rank of the  matrix   $Q(0)$ is reduced if and only if the rank of the $(n+d)\times (n+r)$ matrix $
\left [ \begin{array}{rr} -A & -B \\ D & O \end{array} \right ] $ becomes less then  $n+d$. This is not
fulfilled   if and only if $d\leq r$ and  the system $\dot{x} = Ax + Bu, \; z =Dx$  has no system zeros in
origin. This completes the proof.

We result in that   Problem 3  has  a  solution  if   conditions (a)-(c) of Assertion 10.1 are fulfilled for
matrices  $A, B, D$ of the  original system (10.48),(10.49).

\it{EXAMPLE 10.2.} \rm  For  illustration we consider the second order system with two inputs and outputs from
[P5]
$$  \left [ \begin{array}{c} \dot{x}_1 \\ \dot{x}_2 \end{array} \right ] \; =
\;  \left [ \begin{array}{rc} 0 & 1 \\ -6 & 5 \end{array} \right ] \left [ \begin{array}{c} x_1 \\ x_2
\end{array} \right ]  +  \left [
\begin{array}{cc} 1 & 1 \\ 0 & 2 \end{array} \right ] \left [ \begin{array}{c}
u_1 \\ u_2 \end{array} \right ] \eqno(10.58) $$
$$  \left [ \begin{array}{c} z_1 \\ z_2 \end{array} \right ] \; =
\;  \left [ \begin{array}{rc} 1 & 0 \\ -1 & 1 \end{array} \right ] \left [ \begin{array}{c} x_1 \\ x_2
\end{array} \right ] \eqno(10.59) $$
It is necessary to track for the following reference signal
$$ \left [ \begin{array}{c} v_1(t) \\ v_2(t) \end{array} \right ] \; =\;
 \left [ \begin{array}{c} 2t \\ t \end{array} \right ], \qquad 0 \le t < \infty
\eqno(10.60)$$
 Since in this case $\eta = 2$  then the dynamic regulator (10.51) must have the following structure
$$ \begin{array}{ccl} \dot{q}_1 & = & z - z_{ref} \\
                      \dot{q}_2 & = &q_1 \end{array}
\eqno(10.61)$$
$$             u \; = \; Kx + K_1q_1 + K_2q_2
\eqno(10.62)   $$
where $q_1,\; q_2$ are $2\times 1$ vectors, $x^T \;=\; [ x_1, x_2 ],\;\; q_i^T \;=\;[q_{i1},
q_{i2}],\;\; z_{ref}^T \;=\; [v_1, v_2],\;\; z^T\;=\; [z_1, z_2], \;\; K, K_1, K_2$ are $2\times 2$ matrices.

     Let us test  conditions of Assertion 10.1. One can  see  that
conditions  (a)  and   (b)   are  fulfilled  because   $rank[ B , AB ] = 2$  and  $r = d = 2$. Checking the
condition  (c) gives
$$ rank\left [ \begin{array}{rr} -A & -B\\ D & O  \end{array} \right ] \;=\;
       rank \left [ \begin{array}{rrrr} 0 & -1 & -1 & -1\\ 6 & -5 & -0 & -2 \\
        1 & 0 & 0 & 0 \\ -1 & 1 & 0 & 0 \end{array} \right ] \;=\; 4 $$
Therefore,  the  regulator of the structure (10.61), (10.62) may be used. Constructing  the augmented system
(10.53) with the concrete matrices $A$, $B$ , $D$  yields
$$  \left [ \begin{array}{c} \dot{x}_1 \\ \dot{x}_2 \\ \dot{q}_{11} \\
     \dot{q}_{12}\\ \dot{q}_{21}\\ \dot{q}_{22} \end{array} \right ] \; = \;
    \left [ \begin{array}{rrcccc} 0 & 1 & 0 & 0 & 0 & 0 \\
     -6 & 5 & 0 & 0 & 0 & 0 \\ 1 & 0 & 0 & 0 & 0 & 0 \\
      -1 & 1 & 0 & 0 & 0 & 0 \\ 0 & 0 & 1 & 0 & 0 & 0 \\
      0 & 0 & 0 & 1 & 0 & 0 \end{array} \right ] \left [
      \begin{array}{c} x_1 \\ x_2 \\ q_{11} \\ q_{12} \\ q_{21} \\ q_{22}
       \end{array} \right ]\;+\; \left [ \begin{array}{cc} 1 & 1 \\ 0 & 2 \\
       0 & 0 \\ 0 & 0 \\0 & 0\\ 0 & 0 \end{array} \right ] \left [
       \begin{array}{c} u_1 \\ u_2 \end{array} \right ] \; +\;
       \left [ \begin{array}{cc} 0 & 0 \\ 0 & 0 \\
       1 & 0 \\ 0 & 1 \\0 & 0 \\ 0 & 0 \end{array} \right ] \left [
       \begin{array}{c} v_1(t) \\ v_2(t) \end{array} \right ]
$$
Assigning poles of the closed-loop tracking system equal to -1 and using the modal control method we calculated
[P5]
$$  [ K, \; K_1,\; K_2 ]\;=\; \left [\begin{array}{rcrrrr}
               -9 &  269/13 & -23/26 & -23/26 & -1/6 & -1/6 \\
                0 & 1 & 3/2 & -9/2  & 1/2 & -5/2
\end{array} \right ]
\eqno(10.63) $$
 The dynamic regulator (10.61), (10.62)  with   feedback  gains  (10.63)  maintains asymptotic
tracking for the polynomial signal (10.60)
$$  \begin{array}{c} \\ lim \\ \scriptstyle t \rightarrow \infty
          \end{array} \displaystyle (v_1(t) -z_1(t)) =
         \begin{array}{c} \\ lim \\ \scriptstyle t \rightarrow \infty
         \end{array} \displaystyle (v_1(t) -x_1(t)) = 0 $$
$$  \begin{array}{c} \\ lim \\ \scriptstyle t \rightarrow \infty
          \end{array} \displaystyle (v_2(t) -z_2(t)) =
         \begin{array}{c} \\ lim \\ \scriptstyle t \rightarrow \infty
         \end{array} \displaystyle (v_2(t) +x_1(t)-x_2(t)) = 0 $$

\section[Tracking for modelled reference signal]
        {Tracking for modelled reference signal}

     We will design a dynamic feedback regulator
which maintains  asymptotic  tracking for  a reference signal of a  general  form described by several
differential equations. This regulator, known as servo-regulator (servo-compensator), has a special dynamics
matrix  with eigenvalues coincided  with    characteristic  numbers  of the reference signal. To study this
problem we  follow the approach of Davison [D3], [D5], [D6] and Ferreira [F1].

   Let's consider a completely controllable and observable system
$$ \dot{x} \;=\; Ax + Bu + Ew
\eqno (10.64) $$ with a regulated output
$$   z \; = \; Dx
\eqno (10.65)  $$ where  vectors $x \in \mbox{R}^n$ , $u \in \mbox{R}^r$, $z \in \mbox{R}^d$ have  same senses
as mentioned above, $w \in \mbox{R}^r$ is an unmeasurable disturbance, $A$, $B$, $E$, $D$ are  constant matrices
of appropriate dimensions, $rank B = r$, $rank D = d$. It is assumed that each element of  the vector $w =[w_1,
w_2, \ldots,w_p]$ satisfies similar differential equations of the order $\beta$
$$       w_i^{(\beta )} + \alpha_{\beta -1}w_i^{(\beta -1)} + \cdots +
          \alpha_ow_i \;=\;0
\eqno     (10.66)    $$
 with unknown initial conditions: $w_i(t_o), \dot{w}_i(t_o),\ldots, w_i^{(\beta
-1)}(t_o)$.
 We suppose that the  disturbance $w(t)$ is unmeasurable and  zeros of the characteristic
polynomial\footnote{Or  characteristic numbers.}  of (10.66)
$$     \tilde{\phi}(s) \;=\; s^{\beta } + \alpha_{\beta -1}s^{\beta -1} +
        \cdots + \alpha_o \;=\;0
\eqno (10.67)   $$
 have non-negative real parts. This requirement  assures that   the problem is nontrivial.
 Let
$$   z_{ref}(t) \;=\; z_{ref}\;=\; [z_{r1},\;z_{r2},\;\ldots,\; z_{rd}]
$$
is the desirable reference signal ($d$ vector),  which components satisfy  differential equations of the order
$\beta$ of the form (10.66)
$$     \psi_i (s)=z_{ri}^{(\beta )} + \alpha_{\beta -1}z_{ri}^{(\beta -1)} + \cdots +
          \alpha_oz_{ri} \;=\;0
\eqno     (10.68)    $$
with known initial conditions:
$$  z_{ri}(t_o) = z_{rio}, \;\; z^{(1)}_{ri}(t_o) = z^{(1)}_{ri0},\ldots,
z^{(\beta -1)}_{ri}(t_o) = z^{(\beta -1)}_{ri0} \eqno (10.69)   $$
 It is assumed that the characteristic
polynomial  of  (10.68) coincides with one of (10.66) and  elements of the vector $z_{ref}$ are accessible
\footnote{If characteristic polynomials of $w$ and $z_{ref}$ are different then we should find their a common
multiple and use the approach of [D7].}.

\it{PROBLEM 4.} \rm It is required to find a feedback dynamic regulator as a function of   $x$,  $z$ and
 $z_{ref}$,  $ u = u(x , z , z_{ref})$, such  that  the following asymptotic regulation
$$   z(t) \to z_{ref}, \qquad t \to \infty
\eqno   (10.70)   $$ occurs in the closed-loop system for all disturbances $w$ and for the reference signal
$z_{ref}$ with arbitrary initial conditions $x(t_0), \; z^{(j)}_{rio}, \; j=1,2,\ldots, \beta -1,\;
i=1,2,\ldots,d$.

     We use a feedback regulator of the form [D3],[F1]
$$ \dot{q}\;=\; Fq + \Gamma \epsilon
\eqno  (10.71)  $$
$$   u\;=\;K_1x + K_2 q
\eqno     (10.72) $$ having the tracking error
$$  \epsilon \;=\; z - z_{ref}\;=\; Dx - z_{ref}
\eqno (10.73)   $$
 as the input. Here $q \in \mbox{R}^{d\beta }$ is the state vector of the regulator, $K_1$ and
$K_2$ are  constant matrices of dimensions $r\times n$ and $r\times d\beta $ respectively; $F$ and $\Gamma $ are
constant quasi-diagonal $d\beta \times d\beta $ and $d\beta \times d $ matrices of the following structure
$$ F\;=\; diag(F_1, F_2,\ldots, F_d), \qquad \Gamma \;=\; diag(\gamma_1,
\gamma_2,\ldots, \gamma_d ) \eqno(10.74) $$ In (10.74) $F_i$, $i=1,2,\ldots, d$ are   $\beta \times \beta$
matrices satisfied the equality
$$   det(sI_{\beta} - F_i)\;=\; \tilde{\phi}(s)
\eqno  (10.75)  $$ and $\gamma_i, \; i=1,2,\ldots,d$ are  $\beta$-vectors  assured that pairs $(F_i,\gamma _i )$
are controllable.

Let's introduce  a new  state  vector $\tilde{x}^T \;= \;[x^T,\;q^T]$ and a new reference input $(d+p)$-vector
$\tilde{z}_r \;=\; \left [\begin{array}{c} z_{ref}\\ w \end{array} \right ]$. Uniting  equations (10.64),
(10.71), (10.72) and writing the augmented system with respect to the  vector $\tilde{x}$ we obtain
$$  \dot{\tilde{x}} \;=\; \left [ \begin{array}{cc} A & O \\ \Gamma D & F
    \end{array} \right ]\tilde{x} \;+\; \left [ \begin{array}{c} B \\ O
    \end{array} \right ] u \;+\; \left [ \begin{array}{rc} O & E \\
     -\Gamma & O \end{array} \right ] \tilde{z}_r
\eqno(10.76)  $$
Feedback control (10.73) is also written as follows
$$          u \;=\; [ K_1,\;K_2]\tilde{x}
\eqno(10.77)
 $$
 Thus Problem 4 is  reduced to a state feedback
problem involving the calculation of matrices $K_1$ and $K_2$  of  the  proportional  feedback regulator (10.77)
such that the closed-loop  system
$$  \dot{\tilde{x}} = \left [ \begin{array}{cc} A +BK_1& BK_2 \\
        \Gamma D & F \end{array} \right ]\tilde{x} \;+\;
        \left [ \begin{array}{rc} O & E \\ -\Gamma & O \end{array} \right ]
         \tilde{z}_r
\eqno (10.78)  $$
$$  \epsilon\;=\; [D,\; O]\tilde{x} \;-\;[I_d,\;O]\tilde{z}_r
\eqno   (10.79)$$
 with the  input $\tilde{z}_r $ and  the output $\epsilon$   be asymptotic stable, i.e. the dynamics matrix
of (10.78) must satisfy to the following condition
$$  Re\lambda_i \left ( \begin{array}{cc} A +BK_1& BK_2 \\ \Gamma D & F
   \end{array} \right ) < 0, \;\, i=1,2,\ldots,n+d\beta
 \eqno (10.80) $$
We need to show that if  matrices $F$ and $\Gamma $ in  (10.78)  satisfy   structural restrictions (10.74),
(10.75) then the  dynamic feedback regulator (10.72), (10.73) solves  Problem 4.

     Our reasonings  consist of two steps. At  first  we  find
conditions, which ensure  that  the  error $\epsilon =\epsilon(t)$ (10.79) vanishes as $t \to \infty$; then we
demonstrate that this condition is valid  for given matrices $F$ and $\Gamma$.

     Let $G(s)$ is the transfer function matrix  of   system
(10.78), (10.79) and  $\bar{\epsilon}(s)$ and $\bar{z}_r(s)$ are the Laplace transform of  vectors $\epsilon(s)$
and $\tilde{z}_r(s)$ respectively. We can express  the vector $\bar{\epsilon}(s)$ via $\bar{z}_r(s)$ as follows
$$    \bar{\epsilon}(s)\;=\; G(s)\bar{z}_r(s)
\eqno      (10.81)  $$

\it{ASSERTION 10.3.} \rm If the  condition (10.80) is satisfied for  system (10.78), (10.79) and all elements of
the
 $d \times (d+p)$ matrix $G(s)$ are divided into the  polynomial $\tilde{\phi}(s)$ (10.67) then $ \epsilon(t)
\to 0$   as  $t \to \infty$.

\it{PROOF.$\;$ \rm Let us apply the Laplace transform to  both sides  of  equations (10.68) and (10.66). Using
formulas (2.2) of Sect. 2.1 we calculate $\bar{z}_{ri}(s)$, $i=1,2,\ldots,d+p$,  which are the Laplace transform
of elements of the  vector
$$ \bar{z}_{ri}(s) \;=\; \frac{\psi_i(s)}{s^{\beta }+\alpha_{\beta -1}
         s^{\beta -1} + \ldots + \alpha_o}\;=\;
        \frac{\psi_i(s)}{\tilde{\phi}(s)}
\eqno(10.82) $$
     Representing  the matrix  $G(s)$  in the form\footnote{For simplicity, we
assume that  system (10.78), (10.79)  is completely controllable and observable and  poles of $G(s)$ are equal
to zeros of $det(sI -\tilde{A})$.}
$$    G(s) \;=\; \frac{1}{det(sI -\tilde{A})}\Phi(s)
\eqno (10.83)$$ where $\Phi(s)$ is an  $d \times (d+p)$  polynomial  matrix and $\tilde{A}$ is the dynamics
matrix of  system (10.78)
$$  \tilde{A}\;=\; \left [ \begin{array}{cc} A +BK_1& BK_2 \\
        \Gamma D & F \end{array} \right ]
\eqno (10.84) $$
and substituting (10.83) and (10.82) into (10.81) yields
$$    \bar{\epsilon}(s)\;=\; \frac{\Phi(s)}{det(sI -\tilde{A})}
     \frac{1}{\tilde{\phi}(s)}\left [ \begin{array}{c} \psi_1(s) \\ \vdots \\
     \psi_{d+p}(s) \end{array} \right ]
\eqno (10.85)  $$
 From the last relation it follows that if the  polynomial $\tilde{\phi}(s)$ divides all elements
of $\Phi(s)$ then the dynamical behavior of the error $\epsilon(t)$  depends on eigenvalues of $\tilde{A}$.
According to the condition (10.80) these  eigenvalues  are in  the left-hand part of the complex plan.
Therefore, $\epsilon(t) \to 0$  as $t \to \infty$. The assertion is proved.

\it{ASSERTION 10.4.} \rm  If  dynamics matrices $F$ and $\Gamma $ of  the dynamic regulator (10.72) have the
structure (10.74),  (10.75)  then  all  elements $q_{ij}(s), i=1,2,\ldots,d,\;  j=1,2,\ldots,d+p$  of $G(s)$
have the polynomial $\tilde{\phi}(s)$  as its multiple.

\it{PROOF}.$\;$ \rm Let's denote  $i$-th row of  the matrix $D$ by $d_i$ ($i=1,2, \ldots, d$),  $j$-th column of
the matrix $\Gamma $   by $[O,\; \gamma _j^T,\; O]$ ($j=1,2,\ldots,d$) and $s$-th column of the matrix
$E\;=\;[e_1,\;e_2,\;\ldots, \;e_p]$ by $e_s$. It is evident that the $ij$-th element of the matrix $G(s)$,
namely $g_{ij}(s)$ ($i=1,2,\ldots,d,\;\; j=1,2,\ldots,d+p$),  is calculated by formulas
$$ \protect \begin{array}{ccll} g_{ij}(s)\;& =\; & [d_i,\;0](sI-\tilde{A})^{-1}
       \left [ \begin{array}{c} 0 \\ \cdots  \\ 0 \\ \gamma _i \\ 0
  \end{array} \right ] +  \delta_{ij} , & \begin{array}{c} i=1,2,\ldots,d, \\
j=1,2,\ldots,d, \\ \delta_{ij} = \left \{ \begin{array}{c}1 (i=j)\\ 0(i \neq j)
            \end{array} \right. \end{array} \end{array}$$
$$ \begin{array}{ccll} g_{ij}(s) \;& = & \;[d_i,\;0](sI-\tilde{A})^{-1}\left [
\begin{array}{c} e_s \\ 0 \end{array} \right ]  , & \begin{array}{c}
    i=1,2,\ldots,d,\\j=d+1,\ldots,d+p, \\  s=j-d \end{array} \end{array}
\eqno(10.86) $$
 with  $\tilde{A}$ (10.84). By the direct calculation we reduce  formulas
(10.86) to the following ones
$$  g_{ij}(s) \;=\; \frac{1}{det(sI -\tilde{A})}det\left [
         \begin{array}{ccc} \begin{array}{c} \\sI-\tilde{A}\\ \\ \\ \end{array}
         & \vdots & \begin{array}{c} O \\ \cdots  \\ 0 \\ \gamma _i \\ 0
        \end{array} \\ \dotfill & \dotfill & \dotfill \\ d_i,\; 0 & \vdots &
         \delta {ij} \end{array} \right ], \;\;
          \begin{array}{c}i=1,2,\ldots,d, \\ j=1,2,\ldots,d, \\
         \delta_{ij} = \left \{ \begin{array}{c}1 \;\;(i=j)\\ 0\;\;(i \neq j)
          \end{array} \right. \end{array}
\eqno (10.87a)$$
$$  g_{ij}(s) \;=\; \frac{1}{det(sI -\tilde{A})}det\left [
         \begin{array}{ccc} sI-\tilde{A} & \vdots & \begin{array}{c} e_s \\
           0 \end{array} \\ \dotfill & \dotfill & \dotfill \\
          d_i,\; 0 & \vdots & O \end{array} \right ], \;\;
          \begin{array}{c}i=1,2,\ldots,d, \\ j=d+1,\ldots,d+p, \\
            s=j-d \end{array}
\eqno (10.87b)$$
 To calculate the determinant of  block matrices in right-hand sides of (10.87a,b) we substitute
blocks $F$ and $\Gamma $ (10.74) in (10.84)  and  the result  in   (10.87a), (10.87b). The appropriate matrices
become
$$       \left [ \begin{array}{ccccc}
 sI_n -A -BK_1 & \vdots & -BK_2 & \vdots & O \\ \dotfill & \vdots & \dotfill &
  \vdots & \dotfill \\ \begin{array}{c} -\gamma _1d_1 \\ \vdots \\-\gamma _dd_d
   \end{array} & \vdots & \begin{array}{ccc}(sI_{\beta} -F_1) & \cdots & O \\
     \vdots & \ddots & \vdots \\  O & \cdots & (sI_{\beta} -F_d) \end{array} &
      \vdots &\begin{array}{c} 0 \\ -\gamma _j \\ 0 \end{array} \\
      \dotfill & \vdots & \dotfill & \vdots & \dotfill \\ d_i & \vdots &
        O &\vdots & \delta_{ij} \end{array}
         \right ] ,\;\;
         \begin{array}{c}i=1,2,\ldots,d, \\ j=1,2,\ldots,d, \\
        \delta_{ij} = \left \{ \begin{array}{c}1\;\; (i=j)\\ 0\;\;(i \neq j)
         \end{array} \right. \end{array}
\eqno (10.88a)$$
$$       \left [ \begin{array}{ccccc} sI_n -A -BK_1 & \vdots & -BK_2 &
         \vdots & e_s \\ \dotfill & \dotfill & \dotfill & \dotfill & \dotfill \\
  \begin{array}{c} -\gamma _1d_1 \\ \vdots \\-\gamma _dd_d\end{array} &
  \vdots &\begin{array}{ccc}(sI_{\beta} -F_1) & \cdots & O \\
     \vdots & \ddots & \vdots \\ O & \cdots & (sI_{\beta} -F_d)\end{array} &
      \vdots & \begin{array}{c} O \\ \vdots \\  O \end{array}\\
      \dotfill & \vdots & \dotfill & \vdots & \dotfill\\ d_i & \vdots &
        O  & \vdots & O    \end{array} \right ],\;\;
         \begin{array}{c}i=1,2,\ldots,d, \\ j=d+1,\ldots,d+p, \\
          s = j - d
         \end{array}
\eqno (10.88b)$$
 Then we  premultiply the last row
by the vector $\gamma _i$ and add the  result with all  block rows of following submatrices
$$  \left [ \begin{array}{ccccc}
    \begin{array}{c} -\gamma _1d_1 \\ \vdots \\-\gamma _dd_d\end{array} &
    \begin{array}{c}\vdots \\ \vdots\\ \vdots \end{array} &
    \begin{array}{ccc}(sI_{\beta} -F_1) & \cdots & O \\
       \vdots & \ddots & \vdots \\ O & \cdots & (sI_{\beta} -F_d)\end{array} &
        \begin{array}{c}\vdots \\ \vdots\\ \vdots \end{array} &
         \begin{array}{c} O \\ -\gamma _j \\  O \end{array} \end{array}
         \right ]
\eqno(10.89a)  $$
$$  \left [ \begin{array}{ccccc}
    \begin{array}{c} -\gamma _1d_1 \\ \vdots \\-\gamma _dd_d\end{array} &
    \begin{array}{c}\vdots \\ \vdots\\ \vdots \end{array} &
      \begin{array}{ccc}(sI_{\beta} -F_1) & \cdots & O \\
       \vdots & \ddots & \vdots \\ O & \cdots & (sI_{\beta} -F_d)\end{array} &
        \begin{array}{c}\vdots \\ \vdots\\ \vdots \end{array} &
         \begin{array}{c} O \\ \vdots \\  O \end{array} \end{array}
        \right ]
\eqno(10.89b)  $$
 We obtain the  $i$-th block row of (10.89a,b) (in (10.89a) $j=i$) in the form
$$  \left [ \begin{array}{ccc} O, & (sI_{\beta} -F_i), &  O  \end{array}
      \right ]
\eqno (10.90)$$
 Thus, determinants in (10.87a,b) may be expressed  as the products
$$  (-1)^{\tau_k}det(sI_{\beta} -F_i)det\Omega _k(s)\;=\;(-1)^{\tau_k}
     \tilde{\phi}(s)det\Omega _k(s),\;\; k=1,2
\eqno   (10.91) $$
 where $\Omega _k(s) $ is a some submatrix, $\tau _k$ is a integer, $k=1,2$. Substituting
(10.91) into (10.87a,b) yields
$$   g_{ij} \;=\; \left \{ \begin{array}{ll} (-1)^{\tau_1}
       \frac{\tilde{\phi}(s)}{det(sI-\tilde{A})}det\Omega _1(s), &
              \begin{array}{c} i=1,2,\ldots,d, \\ j=1,2,\ldots,d
                      \end{array}  \\  \\ (-1)^{\tau_2}
   \frac{\tilde{\phi}(s)}{det(sI-\tilde{A})}det\Omega _2(s), &
  \begin{array}{c} i=1,2,\ldots,d,\\j=d+1,\ldots,d+p \end{array} \end{array}
  \right.
\eqno(10.92) $$
 If  zeros of the  polynomial $\tilde{\phi}(s)$ don't  coincide  with   eigenvalues of
$\tilde{A}$ (this requirement may be  always  satisfied  by appropriate choice of  feedback matrices $K_1$ and
$K_2$) then all  elements of $G(s)$  are multiple  to the   polynomial $\tilde{\phi}(s)$. The proof of the
assertion has been completed.

     It  follows from Assertions 10.3 and  10.4   that if
matrices $K_1$  and $K_2$ of the regulator  (10.77)  have  been chosen to held the  condition (10.80)  then the
tracking error tends to zero as $t \to \infty$.  Such the regulator exists if the pair of matrices
$$   \tilde{A} \;=\;\left [ \begin{array}{cc} A & O \\ \Gamma D & F \end{array}
        \right ], \qquad  \tilde{B} \;=\; \left [ \begin{array}{c} B \\ O
       \end{array} \right ]
\eqno (10.93)  $$
 is stabilizable. Let us investigate conditions on   $A$, $B$,  $D$, $\Gamma$, $F$
that assure stabilizability of the  pair ($\tilde{A}$, $\tilde{B}$).

\it{ASSERTION 10.5.} \rm  The pair $(\tilde{A},\tilde{B})$ is stabilizable if  and only if

     a) the pair $(A,B)$ is stabilizable,

     b)  $ d \le r$,

     c) transmission zeros of  the system $ \dot{x} = Ax + Bu, \;
         y =Dx$ don't coincide with  eigenvalues  of the
        matrix $F$ (or zeros of the polynomial $\tilde{\phi}(s)$ (10.67)),

     d) pairs $(F_i, \gamma_i)$, $i=1,2,\ldots,d$  are  stabilizable.

\it{PROOF}. $\;$ \rm We recall that the pair $(\tilde{A},\tilde{B})$ is stabilizable if  and only if the
$(n+d\beta )\times (n+ d\beta +r)$ matrix
$$             Q(\lambda )\;=\; (\lambda I - \tilde{A}, \tilde{B}) \;
      =\; \left [ \begin{array}{ccc} \lambda I_n-A & O & B\\
           -\Gamma D & \lambda I -F & O \end{array} \right ]
\eqno   (10.94)  $$
 has the full rank  $n+d\beta\;$  for $\lambda = \lambda_i^*$ where $\lambda_i^*$ is  an
unstable eigenvalue of  the matrix $\tilde{A}$. One can see that the set of  unstable eigenvalues of $\tilde{A}$
contains   unstable eigenvalues of $A$ and $F$. Therefore, we need to  examine two cases: $\lambda =
\lambda_i^*(A)$ and  $\lambda = \lambda_i^*(F)$ where $\lambda_i^*<.>$  denotes an unstable eigenvalue of a
matrix $<.>$. Without loss of generality we assume that $\lambda_i^*(A) \neq \lambda_i^*(F)$.

\it{CASE 1. } \rm $\lambda = \lambda_i^*(A)$. In according with the condition $\lambda_i^*(A)
\neq\lambda_i^*(F)$ the matrix  $\lambda I-F$ is nonsingular one and the inversion $(\lambda I-F)^{-1}$ exists.
Using equivalent block operations we can write
$$ rankQ(\lambda)\;=\; rank\left [ \begin{array}{ccc}
        \lambda I_n-A & O & B\\-\Gamma D & \lambda I -F & O \end{array}
      \right ] \;=\; rank\left [ \begin{array}{ccc} \lambda I_n-A & B & O\\
            -\Gamma D & O & \lambda I -F \end{array} \right ] \;=\; $$
$$      rank\left [ \begin{array}{ccc} \lambda I_n-A & B & O\\
           O & O & \lambda I -F \end{array} \right ] \;=\;
           d\beta + rank[\lambda I_n -A, B]
$$
Hence,  the matrix $Q(\lambda )$ has the  full rank $n+d\beta$ if  and  only if $rank[\lambda I_n -A, B]
\;=\;n$. This rank condition is fulfilled if and only if the pair of   matrices  $(A,B)$  is  stabilizable
because $\lambda = \lambda_i^*(A)$ is an unstable eigenvalue of $A$.

\it{CASE 2.} \rm  $\lambda = \lambda_i^*(F)$. In according with the condition $\lambda_i^*(A) \neq
\lambda_i^*(F)$ the matrix $\lambda I-A$ is a nonsingular one  and  the inversion $(\lambda I-A)^{-1}$ exists.
Using equivalent block operations  we  can  write  the series of the rank equalities
$$ rankQ(\lambda)\;=\; rank\left [ \begin{array}{ccc}
        \lambda I_n-A & O & B\\-\Gamma D & \lambda I -F & O \end{array}
      \right ] \;=\;rank\left [ \begin{array}{ccc} \lambda I_n-A & O & B\\
       -\Gamma D & \lambda I -F & \Gamma D(\lambda I-A)^{-1}B \end{array}
       \right ]  \;=\; $$
$$       =\;  n + rank[\lambda I -F, \; \Gamma D(\lambda I-A)^{-1}B]
$$
Thus  the matrix $Q(\lambda$  has the full rank  $n+d\beta $ if
$$ rank [\lambda I_{d\beta } -F,\;\Gamma D(\lambda I-A)^{-1}B] \;=\; d\beta
\eqno  (10.95) $$

Let us show that the rank equality (10.95) is true if conditions (b), (c), (d) of the assertion are carried out.
     Indeed, if $\lambda = \lambda_i^*(F)$ then  the condition (c) is
the necessary and sufficient condition for the rank fullness of the matrix $T(\lambda)\;=\; D(\lambda
I-A)^{-1}B$; $rankT(\lambda ) = min(d,r)$ because  $T(\lambda)$ is the transfer function matrix of the system $
\dot{x} = Ax + Bu, \; y = Dx$. If the condition  (b)  is  fulfilled  then  we  always can  find  a nonsingular
$r\times r$ matrix $L(\lambda)$ such that the following relationship takes place
$$ T(\lambda)L(\lambda)  \;=\; [ I_d,\;O]
\eqno (10.96) $$
 Postmultiplying the second block column of the matrix $[\lambda I_{d\beta } - F, \;\Gamma D(\lambda
I-A)^{-1}B] $ by the nonsingular matrix $L(\lambda)$ and using (10.96) we can  write  series of rank equalities
$$ rank [\lambda I_{d\beta } -F,\;\Gamma D(\lambda I-A)^{-1}B] \;=\;
   rank [\lambda I_{d\beta } -F,\;\Gamma D(\lambda I-A)^{-1}BL(\lambda )] \;=\;$$
$$   rank [\lambda I_{d\beta } -F,\;\Gamma [I_d,\;O]\;] \;=\;
   rank [\lambda I_{d\beta } -F,\;\Gamma]
\eqno (10.97)  $$
 It is evident that  rank of the matrix $[\lambda I_{d\beta } -F,\;\Gamma]$ is  equal to $d\beta$
if and only  if  the  pair  of   matrices  $(F,\Gamma)$ is stabilizable. Analysis of  matrices $F$ and $\Gamma$
shows that the  pair $(F,\Gamma)$ is stabilizable if and only if  identical pairs $(F_i,\gamma_i)$ are
stabilizable, i.e. if  the condition (d) is fulfilled. The assertion has been proved.

     From the assertion it follows

\it{THEOREM 10.2.} \rm  Sufficient conditions for  existing the servo-regulator (10.72),(10.73) that assures
asymptotic tracking  (10.70) in the system (10.64),(10.65)   for  all disturbances $w(t)$ (10.66) and all
reference signals $z_{ref}(t)$ (10.68) are conditions (a)-(d)  of Assertion 10.5.

\it{REMARK 10.1.} \rm  One can see that  above Problems 1, 3 are  particular cases of Problem 4 with $F$ is the
zero matrix.

\it{CONCLUSION}.$\;$ \rm Theorem 10.2 reveals the relationship
 between the tracking problem and the system  zeros  location,  namely,
the problem is solvable if transmission zeros  don't  coincide with characteristic numbers of the reference
signal.

               \it{EXAMPLE 10.3.} \rm

For the illustration we consider the example from [S9]. Let the completely controllable and observable system
has the single input and  output
$$  \left [ \begin{array}{c} \dot{x}_1 \\ \dot{x}_2 \end{array} \right ] \; =
\;  \left [ \begin{array}{rc} 1 & 0 \\ 2 & 1 \end{array} \right ] \left [ \begin{array}{c} x_1 \\ x_2
\end{array} \right ]  +  \left [
\begin{array}{c} 1 \\ 1  \end{array} \right ] u
\eqno(10.98) $$
$$  z \; = \;  \left [ \begin{array}{cc} 0 & 1 \end{array} \right ]
\left [ \begin{array}{c} x_1 \\ x_2 \end{array} \right ] \eqno(10.99) $$
 It is desirable to find an  output
dynamic feedback regulator, which maintains asymptotic tracking of the output $z$   for an reference signal,
which is changed  in  according  with  the increasing exponential law
$$       z_{ref}\;=\; e^{2t}z_{ro}
\eqno (10.100)   $$ where $z_{ro} $ is a nonzero real number.

     In order to employ the above results we at first  ought  to  write the
differential equation for  $z_{ref}$
$$ \dot{z}_{ref} - 2z_{ref} \;=\;0, \qquad  z_{ref}(t_o) = z_{ro}
\eqno   (10.101)$$
 Since (10.101) is the linear differential equation  of  the  first order then $\beta =1$. The
characteristic polynomial (10.67) for (10.101)
$$    \tilde{\phi}(s) = s-2   $$
has the characteristic number $s_1 =2$. Therefore,   we can choose $F= 2,\; \Gamma = 1$. As $d = 1,\; d\beta =
1$ then the feedback regulator  (10.72)  is  to  have  the  following structure
$$  \dot{q} \;=\; 2q + \epsilon\;=\; 2q + z -z_{ref}
\eqno    (10.102) $$
$$             u \; = \; k_1x_1 + k_2x_2 + k_3q
\eqno(10.103)   $$
where $q$ is the scalar variable, $k_1$,  $k_2$, $k_3$ are  constant feedback gains, which
are needed to find.

At first we analyze conditions (a)-(d) of Assertion 10.5. The fulfilment of  conditions (a), (b), (d) are
obviously. For checking  the condition  (c) we form the system matrix $P(s)$ for system (10.98), (10.99) and
calculate $det P(s)$
$$      detP(s) \;=\; det \left [ \begin{array}{ccr} s-1 & 0 & -1 \\
                      -2 & s-1 & -1 \\ 0 & 1 & 0 \end{array}
              \right ] \;=\;s+1 $$
Hence, the system has the only transmission zero  being  equal  to $-1$, which does't coincide with
characteristic number  $s_1 = 2$. Consequently, the tracking problem is solvable.

     To calculate  feedback gains $k_1$,  $k_2$, $k_3$  we unite
differential equations (10.98) and (10.102) by introducing a new vector $\tilde{x}^T\;=\;[x_1,\;x_2,\;q]$ and
representing (10.102) as follows: $\dot{q} \;=\; 2q+ Dx -z_{ref}.$ We result in
$$  \dot{\tilde{x}} = \left [ \begin{array}{ccc} 1& 0 & 0 \\ 2 & 1 & 0\\
       0 & 1 & 2 \end{array} \right ]\tilde{x} \;+\;
     \left [\begin{array}{c} 1 \\1 \\0 \end{array} \right ] u +
      \left [\begin{array}{r} 0 \\0 \\-1 \end{array} \right ]z_{ref}
\eqno (10.104)  $$
For system (10.104) we  find a proportional state feedback regulator
$$  u\;=\; [k_1,\;  k_2,\; k_3 ]\tilde{x}
\eqno (10.105)   $$
 shifting  poles  of  the  closed-loop  system  to  numbers: $ -2.148, -1.926 \pm
0.127j$. The appropriate row vector $[k_1,\;  k_2,\; k_3]$ is calculated as [S9]
$$  k\;=\; [\; -2.167, \; -7.833,\;  -21.333\;]
\eqno    (10.106) $$ Thus, dynamic servo-regulator (10.102), (10.103) with  $k_i$, $i=1,2,3$ from (10.106)
becomes
$$  \dot{q} \;=\;  2q + z -z_{ref}, \;\;\;\;\;
         u \; = \; -2.167x_1 -7.833x_2 -21.333 q \eqno(10.107)   $$

\section[Zeros and maximally
                       accuracy of optimal system]
         {Zeros and maximally 
                      accuracy of optimal system }

In the first step of control  design it is desirable to analyze properties of an  open-loop  system,  namely,
one  of  the  main question is: What can  maximally accuracy be achieved when there is no a limitation in the
power of an input action. As it has been shown by Kwakernaak and Sivan [K4], the optimal system may be
classificated into two groups:

     1. Systems having unlimited accuracy. For such systems the
performance criterion  can be reduced to zero if input amplitudes are allowed to increase indefinitely.

     2. Systems having limited accuracy. For such systems the
performance criterion  can't be reduced beyond a certain value even if  input amplitudes are allowed to increase
indefinitely.

     The problem of maximally achievable accuracy has been studied
for the optimal regulator and the optimal filtering  in [K4]. We  consider only a few questions connected  with
transmission zeros. It will be shown that the property of  maximally achievability accuracy of a linear optimal
system is  related  with  the  lack  of  right-half transmission zeros in an  open-loop system.

     Let consider the linear quadratic cost optimal regulator  problem
for the  completely  controllable  and  observable  time-invariant system
$$ \dot{x} \;=\; Ax + Bu
\eqno (10.108) $$
$$   z \; = \; Dx
\eqno (10.109)  $$ where $x \in \mbox{R}^n$,  $u \in \mbox{R}^r$, $z \in \mbox{R}^l$, $A$, $B$, $C$ are constant
matrices of
 appropriate  dimensions.  Let $J(u)$ is  the performance criterion, which is necessary to minimized
$$ J(u) \;=\; \int_{t_o}^{\infty}{(z^T\bar{N}_1z \;+\;u^TN_ou)}dt
\eqno(10. 110)    $$ where $\bar{N}_1>0$\footnote{Positive-definite and nonnegative-definite matrices are
denoted by $N>0,\; N \ge 0$ respectively.} is  an $l\times l$ and $ N_o > 0$  is an  $r\times r$ symmetric
positive-definite matrices. Substituting  (10.109)  into  (10.110) yields the following performance criterion
$$ J(u) \;=\; \int_{t_o}^{\infty}{(x^TN_1x \;+\;u^TN_ou)}dt
\eqno(10. 111)    $$
with   $N_1 = D^T\bar{N}_1 D \ge 0$ is the symmetric nonnegative-definite matrix. For this
case the Riccati equation will be
$$ -\dot{P}(t) \;=\; A^TP(t) + P(t)A - P(t)BN_o^{-1}B^TP(t) +D^T\bar{N}_1 D
\eqno(10.112)$$
 where $P(t)$ is an symmetric $n\times n$ matrix. It is known [K5] that if the pair $(A,B)$ is
controllable, the pair $(A,D)$  is observable and $N_1 \ge 0,\;\;N_o > 0$ then there exists a unique
nonnegative-definite steady-state solution of (10.112).

     Let's investigate steady-state solution properties of (10.112) when $N_o=\rho N$
as $\rho \to 0$  where $\rho$ is a real constant. Such the investigation allows to evaluate  maximally
achievable accuracy of the optimal control with unbounded input power.

We denote by $\bar{P}$  the  $n\times n$ matrix which is the steady-state solution of  the Riccati equation
(10.112). It has been shown in [K4] there exists $\begin{array}{c} \\ lim \\
\scriptstyle \rho \rightarrow 0
\end{array}\displaystyle \bar{P} = P_o $ \footnote{The exact value
of $ lim \; \bar{P}$ can be calculate by the singular optimal problem  [K5].}
 when  $N_o =\rho N,\; p \to 0$ and  for  the closed-loop optimal system
the following limit takes place:
$$ \begin{array}{c} \\ lim \\ \scriptstyle \rho \rightarrow 0 \end{array}
\displaystyle \begin{array}{c} \\ min \\ \scriptstyle u \end{array} \displaystyle J(u)\;=\; x(t_o)^TP_ox(t_o)$$
     The properties of  $P_o$  are defined by
the following theorem.

\it{THEOREM 10.3.} \rm

     a) If $l>r$  then  $P_o \neq O$,

     b) if $l \le r$   then  $P_o \neq O$  only for system (10.108),
        (10.109) having right-half transmission zeros (a non-minimum phase system).

\it{PROOF}.$\;$ \rm  Consider  the  case  (a) ($l>r$) and assume the converse: $\begin{array}{c} \\ lim \\
\scriptstyle \rho \rightarrow 0
\end{array}\displaystyle \bar{P} = O $.
For $N_o = \rho N$ we consider the appropriate algebraic Riccati equation
$$ O \;=\; D^T\bar{N}_1 D  - \frac{1}{\rho }\bar{P}BN^{-1}B^T\bar{P} +
        A^T\bar{P} + \bar{P}A
\eqno(10.113)$$ where $\bar{P}$ is  the symmetric $n\times n$ matrix. Let  $\rho \to 0$. Since the first term in
the right-hand side of (10.113) is independent of $\rho $
and  a finite one then, according to the assumption $\begin{array}{c} \\
lim \\ \scriptstyle \rho \rightarrow 0 \end{array}\displaystyle \bar{P} = O $, the last two terms approach to
zero as $ \rho \to 0$ and (10.113) becomes
$$  \begin{array}{c} \\ lim \\ \scriptstyle \rho \rightarrow 0
  \end{array} \displaystyle \frac{\bar{P}}{\sqrt{\rho }}BN^{-1}B^T
      \frac{\bar{P}}{\sqrt{\rho }} \;=\;D^T\bar{N}_1 D
\eqno(10.114)$$
Since $detN \neq 0$ then it follows from (10.114)  that the limit
$$   L\;=\; \begin{array}{c} \\ lim \\ \scriptstyle \rho \rightarrow 0
           \end{array} \displaystyle B^T\frac{\bar{P}}{\sqrt{\rho }}
\eqno(10.115)$$
must exists. Hence, the following equality takes place
$$         L^TN^{-1}L \;=\; D^T\bar{N}_1D
\eqno  (10.116) $$
 Denoting by $N^{-1/2}$ the  $r\times r$  matrix of the full rank, which satisfies the relation:
$N^{-1/2}N^{-1/2} =  N^{-1}$, we  rewrite (10.116) as follows
$$         R^TR\;=\; D^T\bar{N}_1D
\eqno  (10.117)  $$ where $R = N^{-1/2}L$ is the  $r\times n$ matrix. We now consider (10.117) as the  matrix
equation with respect to the matrix $R$.  As it has been shown in [K4] this equation has a solution for the
$n\times n$ nonnegative-definite symmetric matrix  $D^T\bar{N}D$ if and only if
$$           rankD^T\bar{N}D \le r  $$
The last inequality is equivalent to the following one
$$           rank\bar{N}_1^{1/2}D \le r
\eqno   (10.118)   $$
 where $\bar{N}_1^{1/2}\bar{N}_1^{1/2} =  N_1$. As $\bar{N}_1 >0$ is the square
positive-definite $l\times l$  matrix then  $\bar{N}_1^{1/2}$ is the square nonsingular $l\times l$  matrix.
Hence the equality (4.118) is equal to the following one: $rank D \le r$. By the assumption  of fullness rank of
$D$ we get  the  following  solvability  condition  for  the matrix equation (10.117): $l \le r$. The result
obtained is the contradiction with  assumed the condition $l >r$.  This  implies  that the assumption $P_o = 0$
was not true.

     Now we consider the case (b). Let $r=l$. Then  equation  (4.117)
has the following solution
$$     R \;=\;\bar{N}_1^{1/2}D   $$
Since  $R = N_1^{-1/2}L$ then we can present the  matrix $L$ (10.115) as follows
$$ L\;=\; N_1^{1/2}\bar{N}_1^{1/2}D
\eqno    (10.119) $$
     Let us assume the converse: $\begin{array}{c} \\ lim \\ \scriptstyle \rho
\rightarrow 0 \end{array}\displaystyle \bar{P} = O $ although system (10.108), (10.109) has left-half
transmission zeros. As it  has  been  shown above (see formula (5.11)), system zeros of a system with equal
number of inputs and outputs are defined as zeros  of the following polynomial
$$ \psi(s) \;=\; det(sI_n-A)det(D(sI_n-A)^{-1}B)
\eqno (10.120)$$
 Since we  consider the  completely  controllable   and observable system (10.108), (10.109) then
 zeros defined  from  (10.120) are  transmission zeros.

Now  we study behavior of poles of the closed-loop optimal system when  $\rho \to 0$. These poles coincide with
zeros of the following polynomial
$$ \phi(s) \;=\; det(sI_n-A +BK)
\eqno      (10.121)  $$
 where   $K\;=\; \frac{1}{\rho }N^{-1}B^T\bar{P}$ is the gain matrix of the optimal
regulator. Using Lemma 1.1 from [K5]\footnote{ Lemma 1.1: $\;$ For matrices $M$ and $N$ of  dimensions $m\times
n$ and $n\times m$ respectively the following equality $det(I_m +MN) \;=\;det(I_n +NM)$ takes place.} we can
write  the series of  equalities for $\rho \neq 0$
$$ \phi(s) \;=\; det(sI_n-A +BK)\;=\; det(sI_n-A +BK(sI_n-A)^{-1}(sI_n-A))\;=\;$$
$$     \;=\;det(sI_n-A)det(I_n +BK(sI_n-A)^{-1})\;=\;
        det(sI_n-A)det(I_r +K(sI_n-A)^{-1}B)\;=\;$$
$$     \;=\;det(sI_n-A)det(I_r+\frac{1}{\rho }N^{-1}B^T\bar{P}(sI_n-A)^{-1}B)\;=\;$$
$$      \;=\;det(sI_n-A)det(I_r+\frac{1}{\sqrt{\rho }}N^{-1}B^T\frac{\bar{P}}
      {\sqrt{\rho}}(sI_n-A)^{-1}B)\;=\; $$
$$    \;=\;(\frac{1}{\sqrt{\rho }})^rdet(sI_n-A)det(I_r\sqrt{\rho } +
       \frac{N^{-1}B^T\bar{P}}{\sqrt{\rho}}(sI_n-A)^{-1}B)
$$
If $\rho \to 0$  then using   (10.115) we get
$$ \begin{array}{c} \\ lim \\ \scriptstyle \rho \rightarrow 0
   \end{array}\displaystyle \phi(s) \;=\;(\frac{1}{\sqrt{\rho }})^rdet(sI_n-A)
   det(N^{-1}L(sI_n-A)^{-1}B) $$
Substituting  $L$  (10.119) into the right-hand side of the last expression we obtain
$$ \begin{array}{c} \\ lim \\ \scriptstyle \rho \rightarrow 0
   \end{array}\displaystyle \phi(s) \;=\;(\frac{1}{\sqrt{\rho }})^r det(sI_n-A)
   det(N^{-1/2}\bar{N}_1^{1/2}D(sI_n-A)^{-1}B) $$
 Taking  account  that  the  $r\times r$ matrix $N^{-1/2}\bar{N}_1^{1/2}$ is
nonsingular one and applying formula (10.120) we represent the  last relation as follows
$$ \begin{array}{c} \\ lim \\ \scriptstyle \rho \rightarrow 0
   \end{array}\displaystyle \phi(s) \;=\;(\frac{1}{\sqrt{\rho }})^r
          det(N^{-1/2}\bar{N}_1^{1/2})\psi(s)
\eqno (10.122) $$
 Thus, as  $\rho \to 0$  then  $r$ poles infinitely increase while remained $n-r$ poles will
asymptotically achieve locations of transmission zeros. Since the closed-loop optimal system is asymptotic
stable  then it has to have poles in the left-half of the complex plan. This restriction  is violated  as $\rho
\to 0$ if  the original open-loop system  has right-half transmission zeros. The contradiction proves the case
(b) for $r=l$.

Consider  case $l < r$. Since the set of transmission zeros of the  system $\dot{x}\;=\; Ax +Bu,\; y =Dx$ is
included in the set of transmission zeros of the squared down system $\dot{x}\;=\; Ax +Bu,\; y =TDx$ with an
$r\times l$ constant matrix $T$    then the present case is reduced to the previous one: $l=r$.

\it{CONCLUSION}.$\;$ \rm Theorem 10.3   indicates    expected possibility of the optimal regulator, namely,  it
is impossible to achieve the desirable accuracy of the regulation in a system  with  right-half transmission
zeros.

                       \it{EXAMPLE 10.4.} \rm

To illustrate Theorem  10.3 we consider the simple example  from [S13]. Let the completely controllable and
observable system of the second order with the single input/output is described as [K5]
$$   \dot{x} \;=\; \left [ \begin{array}{cr} 0 & 1 \\ 0 & -4.6 \end{array}
    \right ] x \;+\;\left [ \begin{array}{c} 0 \\ 0.787 \end{array} \right ]u
\eqno(10.123)$$

$$      z \;=\; [1 \;\;-1]x
\eqno  (10.124)  $$
 It is desirable to find an  optimal regulator, which minimized the performance criterion
(10.110) with $\bar{N}_1 = 1,\;\; N_o = 1$.

At first according Theorem 10.3 we analyze expected possibility of the system. To calculate the transmission
zero we build the system matrix $P(s)$  and determine
$$ detP(s)\;=\; det\left [ \begin{array}{ccr} s & -1 & 0 \\
          0 & s-4.6 & -0.787\\ 1 & -1 & 0 \end{array}\right ] \;=\;
           0.787(1-s)  $$
Therefore, the system has the right-half zero (1). In according with the  point b) of Theorem 10.3 the optimal
system will  have  a nonzero   maximally   achievable    error defined as
 $\begin{array}{c} \\ lim \\ \scriptstyle \rho \rightarrow 0
   \end{array}\displaystyle \int_{t_o}^{\infty}{(x^T\bar{N}_1x \;+\;u^T\rho Nu)}dt
      \;=\;=  x(t_o)^TP_ox_0(t_o) $
where $P_o$  is the solution of the algebraic Riccati equation (10.113) as $\rho \to 0$. For testing of this
fact we write the algebraic Riccati equation (10.113) for $\bar{N}_1 = 1,\; N_o =\rho N=\rho $ and
the $2\times 2$ matrix $P= \left [ \begin{array}{cc} p_{11} & p_{12} \\
 p_{21} & p_{22} \end{array} \right ]$
$$  \left [ \begin{array}{cc} 1 & 0 \\  0 & 4.6 \end{array} \right ]P \;+\;
     P\left [ \begin{array}{cc} 0 & 1 \\  0 & 4.6 \end{array} \right ]
     \;-\;P\left [ \begin{array}{c} 0 \\ 0.787 \end{array} \right ]  \rho^{-1}
        \left [ \begin{array}{cc} 0 & 0.787 \end{array} \right ] P \;+\;
         \left [ \begin{array}{r} 1 \\ -1 \end{array} \right ]
          \left [ \begin{array}{cc} 1 & -1 \end{array} \right ] \;=\; O $$
and calculate
$$p_{11}\;=\; \left (\left (\frac{4.6}{\gamma }\right )^2 \;+\;1
 \pm 2\left (\frac{1}{\gamma }\right )^{1/2}\right )^{1/2}\;+\;1,\qquad
        p_{12} \;=\;\pm\left (\frac{1}{\gamma }\right )^{1/2} $$
$$p_{22} \;=\; - \frac{4.6}{\gamma } \;\pm \; \left (\left (\frac{4.6}
       {\gamma }\right )^2 \;+\; \frac{1}{\gamma }\;\pm \;\frac{2}{\gamma }
     \left ( \frac{2}{\gamma }\right )^{1/2}\right )^{1/2} $$
where  $\gamma =(0.787)^2\rho ^{-1}$.

One can see that as  $\rho \to 0$  then $\gamma \to \infty$  and $p_{11}=2,\;p_{12}\;=\;p_{22} \to 0$. Hence,
$P_o \;=\;\left [ \begin{array}{cc} 2 & 0 \\  0 & 0 \end{array} \right ]$ and the  system  has always a nonzero
value  $x^T(t_o)P_ox(t_o)$ for $x^T(t_o) \;=\;[x_1(t_o),\; x_2(t_o)]$ with   $ x_1(t_o) \neq 0$.

     If in the above system we use  the following output
 $$           y \;=\; [1 \;\;1]x   $$
 instead (10.124) then the system zero becomes $-1$.
Calculating elements $p_{11},\;p_{12},\;p_{22} $ of the matrix $P_o$ as $\rho \to 0$  yields:
$p_{11}\;=\;p_{12}\;=\;p_{22}\;=\;0 $.

     This example confirms the connection between  maximally
achievable accuracy of an optimal system  and  locations of transmission  zeros.

\chapter*{List of symbols}
\addcontentsline{toc}{chapter}{List of symbols}
\markboth{List of symbols}{List of symbols}

\begin{list}{}{\leftmargin=5ex \labelwidth=0pt \itemsep=0pt
   \listparindent=-5ex \parsep=0.7ex
   \topsep=0pt \partopsep=0pt \parskip=0pt}
\item~ \vspace{-\baselineskip}

\par
$A,B,A_i,B_i$  -  matrices
\par
$a,b,a_i,b_i$  -  vectors
\par
$a,a_i,a_{ii}$  -  scalars
\par
$\alpha ,\beta ,\alpha _i ,\beta _i$  -  scalars or vectors
\par
$I,I^k$  -  unity matrices
\par
$I_q$   -  unity matrix of order $q$
\par
$O$     -  zero matrix
\par
$diag(a_1,\ldots,a_n)$  -  diagonal matrix with diagonal elements $a_1,
\ldots,a_n$
\par
$diag(A_1,\ldots,A_n)$  -  block diagonal matrix with diagonal blocks $A_1,
 \ldots,A_n$
\par
$A^{i_1,\ldots,i_{\eta }}$ -  matrix constructing from a matrix $A$
                        by deleting all rows expect rows $i_1,\ldots,i_{\eta }$
\par
$A_{j_1,\ldots,j_{\eta }}$ -  matrix constructing from a matrix $A$ by
                     deleting all columns expect columns $j_1,\ldots,j_{\eta }$
\par
$A^{i_1,\ldots,i_{\eta }}_{j_1,\ldots,j_{\eta }}$ -  minor constructing from a
                  matrix $A$  by deleting all rows expect rows $i_1,\ldots,
               i_{\eta }$ and all columns expect columns $j_1,\ldots,j_{\eta }$
\par
$detA$  -  determinant of  matrix $A$
\par
$rankA$  -  rank of  matrix $A$
\par
$\phi (s)$  -  characteristic polynomial of a matrix
\par
$\lambda_i, \; \lambda_i(A)$  -  eigenvalue of  matrix $A$
\par
$Y,\; Y_{AB}$ -  controllability matrix of  pair $(A,B)$
\par
$Z,\;Z_{AC}$ -  observability matrix of  pair $(A,C)$
\par
$\nu $ -  controllability index, integer
\par
$\alpha $  -  observability index, integer
\par
$ ^T$  -  symbol of transponse of a matrix
\par
$A(s),\; \Psi(s)$  -  matrices having polynomial
                     or rational functions as elements
\par
$\epsilon (s) $  -  invariant polynomials of a matrix
\par
$\mbox{N},\;\mbox{R} $   -  linear subspaces
\par
 $\emptyset $ -  empty set

\end{list}

\chapter*{References}
\addcontentsline{toc}{chapter}{References}
\markboth{References}{References}

\begin{list}{}{\leftmargin=5ex \labelwidth=0pt \itemsep=0pt
   \listparindent=-5ex \parsep=0.7ex
   \topsep=0pt \partopsep=0pt \parskip=0pt}
\item~ \vspace{-\baselineskip}

\hyphenation{
  al-ge-b-ra-ic Ana-ly-sys ana-ly-sis ap-pli-ca-tion
  ap-p-ro-ach Av-to-ma-ti-ka Auto-ma-ti-ca
  cal-cu-la-tion ca-no-ni-cal co-ef-fi-ci-ents
  com-p-lex Con-cept Con-t-rol cont-rol-la-bi-li-ty
  De-fi-ni-tion de-ter-mi-na-tion De-ter-mi-na-tion
  ge-o-met-ric Ge-o-met-ric
  li-ne-ar Li-ne-ar
  mat-rix Mat-ri-ces mul-ti mul-ti-va-ri-ab-le
  Mul-ti-va-ri-ab-le mul-ti-con-t-rol
  out-put
  Pha-se Po-les po-ly-no-mi-al Prob-lems Pro-per-ties
  Spa-ce squ-are sur-vey syn-the-sis sys-tem Sys-tem sys-tems Sys-tems
  Te-le-me-kha-ni-ka the-o-ry The-o-ry
  trans-for-ma-tion trans-mis-sion
  va-ri-ab-le
  Uni-ver-si-ty
  ze-ro ze-ros Ze-ros
  de-scrip-tion
   li-near
   dy-na-mi-cal
   track-ing
   exis-tence measu-rements
   con-di-tions ca-no-ni-cal  in-comp-lete
  }
\par
   [A1]  Andreev Yu.N. Control of multivariable linear objects.
         Moscow: Nauka, 1976 (in Russian).
\par
   [A2]  Anderson B.D.O. A note on transmission zeros of a transfer
        function matrix. IEEE Trans. Autom. Control, 1976, AC-24,
        no.4, p.589-591.
\par
   [A3]  Amosov A.A.,Kolpakov V.V. Scalar-matrix  differentiation
        and  its application to constructive problems of
        communication theory.  Problemi peredachi informatsii.
        1972. v.7, no.1, p.3-15 (in Russian).
\par
   [A4]  Asseo  S.J.  Phase-variable  canonical  transformation of
        multicontroller systems. IEEE Trans. Autom. Control, 1968,
        AC-13, no.1, p.129-131.
\par
   [A5]  Athans M. The matrix minimum principle. Information and
         Control, 1968, v.11, p.592-606.
\par
   [B1]  Barnett S. Matrices, polynomials and linear time-invariant
        systems. IEEE Trans. Autom. Control, 1973, AC-18, no.1,
        p.1-10.
\par
   [B2]  Barnett S. Matrix in control theory. London: Van Nostrand
        Reinhold, 1971.

   [B3] Braun  M. Differential equations and their applications. New York:
                Springer-Verlag, 1983.
\par
   [D1]  D'Angelo H. Linear time-invariant  systems: analysis and
        synthesis.  Boston: Allyn and Bacon, 1970.
\par
   [D2]  Desoer C.A., Vidyasagar M. Feedback systems: input-output
        properties.  New York: Academic Press, 1975.
\par
   [D3]  Davison E.J. The output control of linear time-invariant
        multivariable systems with unmeasurable arbitrary
        disturbances. IEEE Trans. Autom. Control, 1972, AC-17,
        no.5,p.621-630.
\par
   [D4]  Davison  E.J.,Wang  S.H.  Property  and  calculation   of
        transmission zeros of linear multivariable systems.
        Automatica, 1974, v.10, no.6. p.643-658.
\par
   [D5]  Davison E.J. A generalization of the output control of
        linear multivariable system  with  unmeasurable  arbitrary
        disturbances. IEEE Trans. Autom. Control, 1975, AC-20, no.6,
        p.788-791.
\par
   [D6]  Davison E.J. The robust control of a servomechanism
        problem for linear time-invariant multivariable system.
        IEEE Trans. Autom. Control, 1976, AC-21, no.1, p.25-34.
\par
   [D7]  Davison E.J.  Design  of  controllers  for  multivariable
        robust servomechanism problem using parameter optimization
        methods. IEEE Trans. Autom. Control, 1981, AC-26, no.1,
        p.93-110.
\par
   [F1]  Ferreira P.G. The servomechanism problem and method of
        the state-space in frequency domain. Int. J.Control.
        1976, v.23, no.2,p.245-255.
\par
   [G1]  Gantmacher F.R. The theory of matrices. v.1,2.  New York:
        Chelsea Publishing Co.,1990.
\par
   [G2]  Gohberg I., Lancaster P. Matrix polynomials. New York:
        Academic Press, 1982.
\par
   [H1]  Hse C.H., Chen C.T. A proof of the stability of
        multivariable  feedback  systems. Proc IEE, 1968, v.56,
        no.1, p.2061-2062.
\par
   [K1]  Kalman R.E. Mathematical description of linear dynamical
        systems. SIAM J. Control, 1963, Ser. A, v.1, no.2, p.152-192.
\par
   [K2]  Kouvaritakis B., MacFarlane A.G.J. Geometric  approach  to
        analysis  and  synthesis  of  system  zeros. Part 1. Square
        systems. Int J. Control, 1976, v.23, no.2, p.149-166.
\par
   [K3]  Kouvaritakis B., MacFarlane A.G.J. Geometric  approach  to
        analysis  and  synthesis  of  system  zeros. Part 2. Non-
        square systems. Int J. Control, 1976, v.23, no.2,p.167-181.
\par
   [K4] Kwakernaak H. Sivan R. The maximally  achievable  accuracy
       of optimal regulators and linear optimal filters. IEEE
        Trans. Autom. Control, 1972, AC-17, no.1, p.79-86.
\par
   [K5] Kwakernaak H. Sivan R. Linear optimal control systems.
        New-York: Wiley, 1972.
\par
   [L1] Lancaster P. Lambda-matrices and vibrating systems.
        London:Pergamon Press, 1966.
\par
   [L2] Lancaster P. Theory of matrices. New York: Academic Press,
        1969.
\par
   [L3] Laub A.J.,Moore B.C. Calculation  of  transmission  zeros
        using QZ techniques. Automatica, 1978, v.14, no.6, p.557-566
\par
   [M1] MacFarlane A.G.J.,Karcanias N. Poles and  zeros of linear
        miltivariable  systems: a survey of the algebraic,
        geometric and complex variable theory. Int.J.Control, 1976,
        v.24, no.1, p.33-74.
\par
   [M2] MacFarlane A.G.J. Relationships between recent
        developments in linear control theory and classical design
        techniques. Control system design by pole-zero assignment.
        London: Academic Press. 1977, p.51-122.
\par
    [M3] MacFarlane A.G.J.Complex-variable  design  methods. Modern
        approach to control system design. London: Proc. IEE. 1979, ch.7.
        p.101-141.
\par
   [M4] Maroulas J.,Barnett S. Canonical forms for time-invariant
        linear control systems: a survey with extensions. Part 1.
        Single-input case. Int.J.Syst.Sci, 1978, v.9, No.5, p.497-514.
\par
   [M5] Maroulas J.,Barnett S. Canonical forms for time-invariant
        linear control systems: a survey with extensions. Part 2.
        Multivariable case. Int.J.Syst.Sci, 1979, v.10, No.1, p.33-50.
\par
   [M6] Moler C.B., Stewart G.W. An algorithm  for  generalized
        matrix eigenvalue problem. SIAM J.Numer.Anal., 1973, v.10,
        no.2, p.241-256.
\par
   [O1] O'Reilly J. Observers for linear systems. London:
        Academic Press,1983.
\par
   [O2] Owens D.H. Feedback and multivariable systems. Stevenage:
        Peter Peregrinus, 1978.

   [P1] Barnett B.N.The symmetric eigenvalue problem. Prentice-
        Hall: Englewood Cliffs, 1980.

   [P2] Paraev Yu.I. Algebraic methods in linear control system
        theory. Tomsk: Tomsk State University, 1980 (in Russion).

   [P3] Patel  P.V.  On  transmission  zeros  and  dynamic  output
        feedback. IEEE Trans. Autom. Control, 1978, AC-23, no.4,
        p.741-749.

   [P4] Porter B.,Crossley R. Modal control. Theory and application.
        London: London Taylor and Francis, 1972.

   [P5] Porter B.,Bradshow A.B. Design of linear multivariable
        continuous-time tracking systems. Int.J.Syst.Sci, \-1974, v.5.
        no.12, p.1155-1164.

   [P6] Porter B. System zeros and invariant zeros. Int.J.Control.
        1978, v.28, no.1, p.157-159.

   [P7] Porter B. Computation of the zeros of linear multivariable
        systems. Int.J.System Sci, 1979,  v.10,  no.12,\- p.1427-1432.

   [R1] Rosenbrock H.N. State-space and multivariable theory.
        London: Nelson, 1970.

   [R2] Rosenbrock H.H. The zeros of a system. Int.J.Control.
        1973, v.18, no.2, p.297-299.

   [R3] Rosenbrock H.H. Correction to 'The zeros of a system'.
         Int.J.Control, 1974, v.20, no.3, p.525-527.

   [S1] Samash J.Computing the invariant zeros of multivariable
        systems. Electron. Lett., 1977, v.13, no.24, p.722-723.

   [S2] Schrader C.B., Sain M.K. Research on system zeros:
        a survey. Int.J.Control, 1989, v.50, no.4, p.1407-1733.

   [S3] Smagina Ye.M. Modal control in multivariable system by
        using generalized canonical representation. Ph.D, Tomsk State University,
        Tomsk, Russia, 1977 (in Russian).

   [S4] Smagina Y.M. Computing the zeros of a linear
        multi-dimensional systems. Transaction on Automation and Remote Control,
        1981,  v.42, No.4, part 1, p. 424-429 (Trans. from Russian).

   [S5] Smagina Ye.M. To the problem of squaring down of outputs
        in linear system. Moscow, 1983, Deposit in the All-Union
        Institute of the Scientific and Technical Information,
        no. 5007-83Dep., p.1-10 (in Russian).

   [S6] Smagina Ye.M. Design of multivariable system with assign
        zeros. Moscow, 1983, Deposit in the All-Union Institute of
        the Scientific and Technical Information no. 8309-84Dep.,
        p.1-15 (in Russian).

   [S7] Smagina Y.M. Zeros of multidimensional linear  systems.
        Definitions,  classification, application (Survey). Transaction on
        Automation and Remote Control, 1985, v.46, No.12, part 1, p.1493-1519
        (Trans.  from Russian).

   [S8] Smagina Ye.M. Computing and specification of zeros in
        a linear multi-dimensional systems. Avtomatika i
        Telemekhanika, 1987. no.12, p.165-173 (in Russian).

   [S9] Smagina Y.M. Problems of linear multivariable system analysis
        using the concept of system zeros. 1990, Tomsk: Tomsk  State University,
         159p (In Russian).

   [S10] Smagina Y.M. Determination of the coefficients of zero polynomial
         in terms of the output controllability matrix. Trans. on Automat. and
         Rem. Contr., 1991, v.52, p.1523-1532 (Trans. from Russian).

   [S11] Smagina Ye.M. A method of designing of observable  output
         ensuring given zeros locations.  Problems of  Control  and
         Information Theory, 1991, v.20(5), p.299-307.

   [S12] Smagina Y.M. Existence conditions of  PI-regulator for
         multivariable system with incomplete measurements.
          Izv. Acad. Nauk SSSR. Tekhn. Kibernetika, 1991, no.6, p.40-45
          (In Russian).

   [S13]  Smagina Ye.M., Sorokin A.V. The use of the concept
          of a system zero when weight matrices are selected
          in the analytic design of  optimal regulators.
          J. Comput. Systems Sci.Internat., 1994, v.32, No.3, p.98-103
          (Trans. from Russian).

   [S14]  Smagina Ye.M. Influence of system zeros on stability of the
          optimal filtering with colored noise. J. Comput. Systems Sci.
          Internat., 1995, v.33, no.3, p.21-25, 1995 (Trans.  from Russian).

   [S15]  Smagina Y.M. System zero determination in large scale system.
          Proc. Symposium IFAC "Large Scale Systems: Theory and Applications",
          11-13  July, 1995, London, UK, v.1, p.153-157.

   [S16]   Smagina Y.M. Definition, calculation and application of system zeros
            of multivariable systems. D.Sc. thesis,
            Tomsk State University, Tomsk, Russia, 1995 (In Russian).

   [S17]  Smagina Y.M. The relationship between the transmission zero
           assignment problem and modal control method. Journal of
           Computer and  System Science International, 1996, v.35, No.2,
            p.39-47 (Trans. from Russian).

   [S18]  Smagina Y.M. Tracking for a polynomial signal for a system
          with incomplete information. Journal of Computer and System
          Science  International, 1996, v.35, No.1,  p.53-37
          (Trans.  from Russian).

    [S19] Smagina Y.M. System zero definition via low order linear pencil.
          Avtomatika i Telemekhanika, 1996, no.5, p.48-57 (In Russian).

    [S20] Smagina Ye.M. New approach to transfer function matrix factorization.
           Proc. Conf. IFAC on Control of Industrial Systems,
          22-22 May, 1997, Belfort, France, v.1/3, p.413-418.

    [S21] Smagina Ye.M. Solvability conditions for the servomechanism
           problem using a controller with the state estimation.
           Engineering Simulation. 1998, v.15, p.137-147.

   [S22] Stecha J. Nuly a poly dynamickeho systemu. Automatizace,
           1981, v.24, no.7.p.172-176.

   [S23] Strejc V. State-space theory of discrete linear control.
                Prague: Academia,1981.

   [V1] Voronov A.A. Stability, controllability, observability.
        Moscow: Nauka, 1979.

   [W1] Wolowich W.A. Linear multivariable system.New York,
        Berlin: Springer-Velag, 1974.

   [W2] Wolowich W.A. On the numerators and zeros of rational
        transfer matrices.IEEE Trans. Autom. Control, 1973, AC-18,
        no.5, p.544-546.

   [W3] Wonham W.M. Linear multivariable control. A geometric
         approach. New York: Springer-Verlag, 1980.

   [Y1] Yokoyama R. General structure of linear multi-input
        multi-output systems. Technol.Report Iwata Univ, 1972,
        p.13-30.

   [Y2] Yokoyama R.,Kinnen E.Phase-variable canonical forms for
        the multi-input, multi-output systems. Int.J.Control, 1976,
        AC-17, no.6, p.1297-1312.
\end{list}

\chapter*{Notes and references}
\addcontentsline{toc}{chapter}{Notes and references}
\markboth{Notes and references}{Notes and references}

 In accordance with the purpose of this  book  some  references
are omitted in the text. The following notes will acquaint with the works used by the authors:

\begin{list}{}{\leftmargin=5ex \labelwidth=0pt \itemsep=0pt
   \listparindent=-5ex \parsep=1.4ex
   \topsep=0pt \partopsep=0pt \parskip=0pt}
\item~ \vspace{-\baselineskip}
\par
Chapter 1: [A1],$\;$[A4],$\;$[B1],$\;$[G2],$\;$[K1],$\;$[K5],$\;$[L1],$\;$[M3],
            $\;$[M4],$\;$[M5],$\;$[O2],$\;$[R1],$\;$[S1],$\;$[S3],$\;$[S17],\\
            $\;$[S23],$\;$[V1],$\;$[W1],$\;$[Y1],$\;$[Y2]
\par
Chapter 2: [A2],$\;$[D2],$\;$[G1],$\;$[K5],$\;$[L2],$\;$[M1],$\;$[O2],$\;$[P2],
            $\;$[S22],$\;$[W1],$\;$[W2],$\;$[Y1]
\par
Chapter 3: [G1],$\;$[M1],$\;$[M3],$\;$[R1],$\;$[S4]
\par
Chapter 4: [B1],$\;$[M1],$\;$[M2],$\;$[S16],$\;$[S20],$\;$[S22],$\;$[W2]
\par
Chapter 5: [M1],$\;$[M3],$\;$[P6],$\;$[R1],$\;$[R2],$\;$[R3],$\;$[S22]
\par
Chapter 6: [D4],$\;$[A1],$\;$[K2],$\;$[M1],$\;$[M3]
\par
Chapter 7: [D1],$\;$[S4],$\;$[S8],$\;$[S9],$\;$[S10],$\;$[S16],$\;$[S19]
\par
Chapter 8: [D4],$\;$[H1],$\;$[K2],$\;$[L3],$\;$[M6],$\;$[P1],$\;$[P7],$\;$[S1],
           $\;$[S8],$\;$[S9],$\;$[S15],$\;$[S16]
\par
Chapter 9: [A3],$\;$[A5],$\;$[K3],$\;$[S5],$\;$[S6],$\;$[S11],$\;$[S16],$\;$[S17]
\par
Chapter 10: [B2],$\;$[D3],$\;$[D4],$\;$[D5],$\;$[D6],$\;$[D7],$\;$[F1],$\;$[K4],
            $\;$[O1],$\;$[P4],$\;$[P5],$\;$[R3],$\;$[S12],$\;$[S13],$\;\\$[S14],$\;$[S18],$\;$[S21],$\;$[W3]
\end{list}

\end{document}